\documentclass[a4paper,11pt,reqno]{amsart}
\usepackage{amsthm}
\usepackage{amsmath}
\usepackage{amsfonts}
\usepackage{amssymb}
\usepackage{enumerate}
\usepackage{pbox}
\usepackage{xcolor}
\usepackage{placeins} % \FloatBarrier
\usepackage{tabularray}
\usepackage{tikz}
\usetikzlibrary{calc}
\usetikzlibrary{cd}
\usepackage{graphicx}
\graphicspath{{img/}}
\usepackage{array}
\usepackage{dsfont} % better looking mathbb letters
\usepackage[utf8]{inputenc}
\usepackage[english]{babel}
\usepackage[T1]{fontenc}
\usepackage{float}

\usepackage{titletoc}
\usepackage[a4paper, margin = 3cm, bottom = 3cm]{geometry}
\usepackage{marginnote}
\usepackage{yhmath} % big delimiters, wide accents
\usepackage[normalem]{ulem}
\usepackage{ifthen}
\usepackage{hyperref} % put last
\usepackage{cleveref} % put after hyperref
\hypersetup{pdfborder=0 0 0,
    colorlinks=true,
    citecolor=black,
    linkcolor=blue,
    urlcolor=red,
    pdfauthor={Guillaume Tahar}
}

\newtheorem{thm}{Theorem}[section]
\crefname{thm}{theorem}{theorems}
\newtheorem*{thm*}{Theorem}
\newtheorem{cor}[thm]{Corollary}
\crefname{cor}{corollary}{corollaries}
\newtheorem{prop}[thm]{Proposition}
\newtheorem*{prop*}{Proposition}
\crefname{prop}{proposition}{propositions}
\newtheorem{lem}[thm]{Lemma}
\newtheorem*{lem*}{Lemma}
\crefname{lem}{lemma}{lemmas}

\crefname{conj}{conjecture}{conjectures}
\newtheorem{comp}[thm]{Complement}
\crefname{comp}{complement}{complements}
\newtheorem*{pbm}{Problem}
\crefname{pbm}{problem}{problems}
\theoremstyle{definition}
\newtheorem{defn}[thm]{Definition}
\newtheorem*{defn*}{Definition}
\crefname{defn}{definition}{definitions}

\crefname{nota}{notation}{notations}
\newtheorem{ex}[thm]{Example}
\crefname{ex}{example}{examples}

\crefname{hypo}{hypothesis}{hypotheses}
\theoremstyle{remark}
\newtheorem{rem}[thm]{Remark}
\newtheorem*{rem*}{Remark}
\crefname{rem}{remark}{remarks}
\numberwithin{equation}{section}

\newcommand{\Aut}{\operatorname{Aut}}
\newcommand{\AutC}{\operatorname{Aut}(\C)}

\newcommand{\C}{\mathbb{C}}
\newcommand{\D}{\mathbb{D}}
\newcommand{\Hp}{\mathbb{H}}
\newcommand{\N}{\mathbb{N}}
\newcommand{\Q}{\mathbb{Q}}
\newcommand{\R}{\mathbb{R}}
\newcommand{\Z}{\mathbb{Z}}

\newcommand{\ie}{i.e.\ }
\newcommand{\cet}{\tilde\C^*}

\DeclareMathOperator{\res}{res} % empêche d'écrire z^\res
\newcommand{\ssm}{\smallsetminus}
\newcommand{\dlog}{\mathrm{d}\log}

\newcommand{\Hol}{\on{Hol}} % l'ensemble d'holonomie
\newcommand{\RS}{\mathbb{CP}^{1}} % la sphère de Riemann
\newcommand{\tra}{\mathcal{T}} % le groupe des translations de C
\renewcommand{\Re}{\operatorname{Re}}
\renewcommand{\Im}{\operatorname{Im}}
\newcommand{\eps}{\epsilon}
\newcommand{\dom}{\operatorname{Dom}}
\newcommand{\tends}{\longrightarrow}
\newcommand{\Id}{\operatorname{Id}}
\newcommand{\ds}{\displaystyle}
\newcommand{\wt}[1]{\widetilde{#1}}
\newcommand{\wh}[1]{\widehat{#1}}
\newcommand{\on}[1]{\operatorname{#1}}
\newcommand{\cal}[1]{\mathcal{#1}} 
\newcommand{\ov}[1]{\overline{#1}}

\newcommand{\setof}[2]{\big\{{#1}\,;\,{#2}\big\}}

\newcommand{\loctitle}[1]{\addvspace{\medskipamount}\noindent\textit{#1.}}

\newcommand{\loctitlebold}[1]{\addvspace{\medskipamount}\noindent\textbf{#1.}}

\title{The affine geometry of meromorphic connections with irregular singularities}

\author{Xavier Buff}
\address[Xavier Buff]{Institut Mathématique de Toulouse; UMR 5219; Université de Toulouse; CNRS; UPS, F-31062 Toulouse Cedex 9, France}
\email{xavier.buff@math.univ-toulouse.fr}

\author{Arnaud Chéritat}
\address[Arnaud Chéritat]{Institut Mathématique de Toulouse; UMR 5219; Université de Toulouse; CNRS; UPS, F-31062 Toulouse Cedex 9, France}
\email{arnaud.cheritat@math.univ-toulouse.fr}

\author[G.~Tahar]{Guillaume Tahar}
\address[Guillaume Tahar]{Beijing Institute of Mathematical Sciences and Applications,  No. 544, Hefangkou Village,
Huaibei Town, Huairou District, Beijing, China, 101408}
\email{guillaume.tahar@bimsa.cn}

\begin{document}

\begin{abstract}
A meromorphic connection on the tangent bundle of a Riemann surface induces a complex affine structure on the complement of the poles. Local models for Fuchsian singularities are already known. In this paper, we introduce a complete set of local invariants for a meromorphic connection and provide local models for a complex affine structure in a punctured neighborhood of an irregular singularity. Generalizing a construction attributed to Veech, we introduce the Delaunay decomposition of a compact Riemann surface endowed with a meromorphic connection with irregular singularities. In particular, we give upper bounds on the complexity of the decomposition.
\end{abstract}

\maketitle

\makeatletter
\def\l@subsection{\@tocline{2}{0pt}{2.5pc}{5pc}{}}
\makeatother
\setcounter{tocdepth}{3}
\startcontents[all]
\centerline{\textsc{Contents}}
\printcontents[all]{}{1}[2]{}

\section{Introduction}

Let $\nabla$ be a meromorphic connection on the tangent bundle of a Riemann surface $X$.
(We adopt the convention that Riemann surfaces are connected.)
Let $\mathcal{S} \subset X$ be a closed discrete subset of $X$ containing every pole of $\nabla$ and possibly some additional marked points.
Denoting by $X^{\ast}$ the complement of the \textit{singular set} $\mathcal{S}$ in $X$, $\nabla$ induces a \textit{complex affine structure} on $X^{\ast}$ i.e.\ an atlas for which the change of charts are complex affine.
If $X$ is compact then $\cal S$ is finite and we will refer to such a triple $(X,\nabla,\mathcal{S})$ as a \textit{finite type affine surface}.
Note that finite type affine surfaces are in particular connected.
When the context makes it unambiguous, we will sometimes omit the mention of \emph{finite type} and refer to $(X,\nabla,\mathcal{S})$ as $X$.
\par
Two finite type affine surfaces $(X_{1},\nabla_{1},\mathcal{S}_{1})$ and $(X_{2},\nabla_{2},\mathcal{S}_{2})$ are isomorphic if there is a biholomorphism $f:X_{1}\rightarrow X_{2}$ such that $f^{\ast} \nabla_{2}=\nabla_{1}$. This is equivalent to requiring that $f$ is an isomorphism for the complex affine structures defined by $\nabla_{1}$ and $\nabla_{2}$ on $X_{1}^{\ast}$ and $X_{2}^{\ast}$ respectively.
\par
On affine surfaces there is a notion of geodesics. On finite type affine surfaces, though geodesics naturally extend across marked points, we decide to stop them there, i.e.\ geodesics must live in $X^*$.
\par
The complex affine structure in the punctured neighborhood of a simple pole (also called a Fuchsian singularity) of a meromorphic connection $\nabla$ is already well understood. The analytic classification up to local isomorphism is given by local invariants, each corresponding to a canonical model in complex affine geometry, as described for example in \cite{CT}. In contrast, for multiple poles, also called \textit{irregular singularities}, no systematic classification was known to date.
\par
In this paper, we introduce a complete set of local isomorphism invariants for multiple poles and provide systematic local models for their neighborhoods in the language of complex affine geometry. As an application, we extend the \textit{Delaunay decomposition} to the settings of meromorphic connections on compact Riemann surfaces in full generality (see \Cref{sub:IntroDelaunay} for the mention of some special cases where the decomposition was already constructed).

\subsection{Notes on terminology}

By convention, the total number $n$ of singularities in $\mathcal{S}$ is the sum of the number of marked points and the total order of the poles of $\nabla$ counted with multiplicity. Marked points are interpreted as conical singularities of angle $2\pi$ and though they are erasable, in this article they are counted with multiplicity $1$, not $0$.

Our surfaces have two structures, as a Riemann surface and as an affine surface. Though the second one induces the first, the second misses local charts near the singularities. We call \emph{Riemann chart} any chart of any atlas of the Riemann surface, and \emph{affine chart} any chart of any atlas of the affine surface.

Concerning the residue of singularities, there are two meaningful but opposite sign conventions.
We choose the one such that the sum of residues equals $2-2g$ in a genus $g$ finite type affine surface. Note that this is the opposite of the convention adopted in \cite{CT}.

The set of complex-affine maps from $\C$ to $\C$ that are invertible is denoted $\Aut(\mathbb{C})$.

Differential forms, of expression $f(z)\,dz$ in Riemann charts, are also known under the name of abelian differentials, and of $1$-forms.

Given a collection of sets we will often abbreviate \emph{pairwise disjoint} as \emph{disjoint}.

\subsection{Irregular singularities}\label{sub:irrsing}

Let us put ourselves in a local situation and assume we are given a meromorphic connection $\nabla$ defined on a punctured neighborhood of a point $p$ in a Riemann surface and that $\nabla$ has an irregular singularity at $p$.
Consider a Riemann chart sending $p$ to the origin in $\C$ and denote $U$ the image of the chart: it is an open subset of $\C$ containing $0$.
The connection has a Christoffel symbol $\Gamma$ in this chart which is a meromorphic function on $U$.
The hypothesis that $\nabla$ has an irregular singularity at $p$ means $\Gamma$ has pole at $0$, of order at least $2$ (see \Cref{sub:minicours}).
Given two points of two Riemann surfaces locally endowed with a meromorphic connection, if there is a local isomorphism then call them \emph{equivalent}.

By \Cref{sub:minicours}, the degree and residue of the pole of $\Gamma$ at $0$ are independent of the choice of chart, and hence form an invariant for local isomorphism. For some reasons explained in \Cref{note:minus}, we define the residue of the singularity, which we denote as $\res$, to be the the residue of $-\Gamma$, i.e.\ the \emph{opposite} of the residue of $\Gamma$.

\medskip

The notion of \emph{formal equivalence} is developed in \Cref{sub:formal}.
The degree and residue are not only invariant under analytic equivalence (local isomorphism), but also under formal equivalence. Moreover:

\begin{thm}\label{thm:FormalEquivalence}
The degree and residue of the pole are the only formal invariants of irregular singularities. In other words, if two multiple poles have the same degree and residue, then their connections in any two Riemann charts near the poles are formally equivalent.
\end{thm}

The theorem above, proved in \Cref{sub:formal},
is also valid for, and hence unifies with, Fuchsian singularities (in the generalized sense, i.e.\ including erasable singularities).

\medskip

Actually, for Fuchsian singularities, the residue is almost the unique analytic invariant, i.e.\ formal equivalence implies analytic equivalence in most cases.
First, all erasable singularities (Fuchsian singularities of residue $0$) are analytically equivalent (there always is a chart in which $\Gamma$ vanishes locally).
More generally two Fuchsian singularities of residue $\res \in\C \setminus \{2,3,\ldots\}$ are analytically equivalent (locally isomorphic).
And for a fixed $\res \in \{2,3,\ldots\}$, there are only two analytic equivalence classes.

For irregular singularities of order $d=2$, we will see that each formal equivalence class has only two analytic classes.
But for $d\geq 3$, analytic equivalence classes are further away from formal equivalence classes: in the theorem below we will see that for poles of degree $d \geq 3$ and residue $\res \in \mathbb{C}$, there is a space of complex dimension $d-2$ of formally equivalent but analytically distinct singularities. These singularities are distinguished by an invariant taking values in a space $\mathcal{I}_{d,\res}$ defined below.

\newcommand{\defIsText}{
For $d\geq 2$ and $\res \in \mathbb{C}$:
\begin{itemize}
\item Let $\cal{U}_{d,\res}$ be the set of families $u \in \mathbb{C}^{\mathbb{Z}}$ such that there exists $b \in \mathbb{C}$ for which $u_{k+(d-1)} = e^{-2\pi i \res} u_k + b$ for any $k \in \mathbb{Z}$.
\item Let $\cal I_{d,\res}$ be the quotient of $\cal U_{d,\res}$ by the equivalence relation $u \sim u'$ iff there exists an affine bijection $g$ of $\C$ such that $\forall k\in\Z$, $g(u_k) = u'_k$.
\end{itemize}
We call $\mathcal{I}_{d,\res}$ the \emph{space of invariants} and denote $\Pi: \mathcal{U}_{d,\res} \to \cal \mathcal{I}_{d,\res}$ the quotient map.
We denote by $\cal U_d$ the disjoint union of $\cal U_{d,\res}$ for $\res\in \C$ and $\cal I_d$ the disjoint union of $\cal I_{d,\res}$.
}
\begin{defn}\label{def:is}
\defIsText
\end{defn}

To any multiple pole of order $d \geq 2$ and residue $\res \in \mathbb{C}$, we associate in \Cref{sub:constr:inv} an element of $ \mathcal{I}_{d,\res}$ called its \textit{asymptotic values invariant}.
It is defined by looking at the limit, along paths tending to the pole tangentially to $d-1$ half-lines we call the repelling axes, of the developing map associated to a germ of affine chart.

\newcommand{\thmMainText}[1]{
The invariant is complete and effective:
\begin{enumerate}
\item\ifthenelse{#1}{\label{item:1}}{} Consider two Riemann surfaces with a meromorphic connection, each with a pole of the same order $d \geq 2$ and the same residue $\res \in \mathbb{C}$ together with a choice of reference divergent axis.
They have the same \textit{asymptotic values invariant} if and only if they are isomorphic near the punctures by an isomorphism matching the divergent axes.
\item\ifthenelse{#1}{\label{item:2}}{} For any $d \geq 2$, any residue $\res \in \mathbb{C}$ and any element $\iota \in \cal I_{d,\res}$, there exists a meromorphic connection on a punctured disk having a pole of order $d$, residue $\res$ and invariant $\iota$. at the puncture.
\end{enumerate}
}
\begin{thm}\label{thm:main}
\thmMainText{\boolean{true}}
\end{thm}

\Cref{thm:main} is proved in \Cref{sub:model} by giving an affine geometric local model for any prescribed invariant.

\subsection{Delaunay decomposition}\label{sub:IntroDelaunay}

The geometry of affine immersions of disks and half-planes into an affine surface captures much of its qualitative structure. These immersions are organized into an ordered structure: the Delaunay category (\Cref{defn:DelaunayCategory}). We classify its maximal elements into several types: $A$, $B$, $C$, $I$, $II$, $III$, $IV$ and $V$ depending on their boundary behavior. Maximal immersions of type A define the \textit{Delaunay segments}. Cutting along the Delaunay segments yields the \textit{Delaunay decomposition} of the affine surface. This construction generalizes several previously known decompositions:
\begin{itemize}
    \item The classical Delaunay triangulation in the Euclidean plane, dual to the Voronoi tessellation and widely used in computational geometry, see \cite{Del34}.
    \item Its global generalization to translation surfaces, see for example \cite{GJ00} or more recently \cite{BSTW25} for an application to dynamics of differential operators associated to meromorphic differentials.
    \item Veech's construction for affine surfaces with conical singularities that presents specific difficulties arising from the fact that affine structures are generally not geodesically complete, and the Delaunay decomposition is not always a triangulation. See \cite{Vee93} for the original construction and \cite{DFG19} for a detailed proof in the context of dilation surfaces (affine surfaces where transition maps are constrained to dilations and translations). The construction has been crucially used to prove of existence of closed geodesics in dilation surfaces in \cite{BGT25}.
\end{itemize}
In particular, the present construction addresses the case of affine surfaces arising from meromorphic connections with irregular singularities.

\begin{defn}\label{defn:DelaunayEdge}
Given a finite type affine surface $(X,\nabla,\mathcal{S})$, a disk immersion of \textit{type A} is the affine immersion $f:\Delta \to X^{\ast}$ of an open disk $\Delta$ into $X^{\ast}$ that extends continuously to its boundary circle and such that exactly two points of the boundary circle are mapped to the singular set $\mathcal{S}$.
\par
A \textit{Delaunay segment} of an affine surface $S$ is the image of the chord drawn between the two singular boundary points under the extension of a disk immersion of type A.
\end{defn}

\begin{rem}
We recover the classical Delaunay decomposition with respect to a finite set of points in the flat plane by marking these points.
\end{rem}

We prove in \Cref{thm:ExceptionalSurfaces} that affine surfaces with no Delaunay segments are precisely those for which no finite time maximal trajectory exists. We call them exceptional affine surfaces, and their dynamical behavior is particularly simple.

\begin{thm}\label{thm:ExceptionalSurfaces}
For a finite type affine surface $X$ of genus $g$ with $n$ singular points (counted with multiplicity), the following statements are equivalent:
\begin{enumerate}[(i)]
    \item every maximal geodesic on $X^{\ast}$ is infinite in the future or in the past;
    \item $X$ has no Delaunay segment;
    \item $2g+n \leq 2$;
    \item the universal covering of $X^{\ast}$ is either the flat plane $\mathbb{C}$ or the exponential-affine plane $\mathcal{E}$ (see \Cref{sub:ExpoAffine});
    \item there exists an affine immersion $\cal E\to X^*$;
    \item the core (defined in \Cref{sub:core}) is contained in $\cal S$.
\end{enumerate}
\end{thm}

\Cref{thm:ExceptionalSurfaces} is proved in \Cref{sub:Exceptional}.
We refer to these surfaces as \emph{exceptional affine surfaces}. A complete list is given in \Cref{prop:NoDelaunayEdge}.

\medskip

In \Cref{thm:ClassificationComponents} below we prove that for non-exceptional affine surfaces, the Delaunay segments form a finite system of disjoint arcs drawn between singularities of the meromorphic connection, which by \Cref{lem:tips} are either Fuchsian conical or irregular (in which case the segment actually tends to what we call a \emph{focus} of the singularity, see \Cref{sub:helix}), and we provide the following complete classification of the connected components of the \emph{Delaunay components}, defined as the connected components of the complement of the union of the Delaunay segments.

\begin{thm}\label{thm:ClassificationComponents}
In any non-exceptional affine surface $(X,\nabla,\mathcal{S})$, Delaunay segments interiors are a finite collection of disjoint embedded geodesics. % joining points of $\mathcal{S}$.
Every conical Fuchsian singularity and every irregular singularity focus is the end of at least one Delaunay segment.
The Delaunay components are of the following types (described in details in \Cref{sub:cylindersAndSkewCones,subsub:UnboundedDomains,sub:core}):
\begin{itemize}
    \item[(i)] \textbf{Delaunay polygons}: embedded convex polygons with $p \geq 3$ sides drawn between points of $\mathcal{S}$;
    \item[(iia)] \textbf{Finite angle Reeb cylinders}: dilation cylinders of angle at least $\pi$;
    \item[(iib)] \textbf{Semi-infinite Reeb cylinders} and their associated singularity: exactly one semi-infinite Reeb cylinder for each Fuchsian singularity of residue $\res$ with $\Re(\res)=1$ and $\Im(\res) \neq 1$;
    \item[(iii)] \textbf{Semi-infinite translation cylinders} and their associated singularity: exactly one for each Fuchsian singularity of residue $1$;
    \item[(iv)] \textbf{Anti-conical domains} and their associated singularity: exactly one for each Fuchsian singularity whose residue $\res$ satisfies $\Re(\res)>1$;  
    \item[(v)] \textbf{Swath domains}: exactly $d-1$ for each irregular singularity of degree $d$.
\end{itemize}
Topologically, Delaunay polygons and swath domains are contractible. Reeb cylinders, semi-infinite translation cylinders and anti-conical domains are homeomorphic to twice-punctured spheres.
\end{thm}

\Cref{thm:ClassificationComponents} is proved in \Cref{sub:DelaunayComponents}.

\begin{rem}
In particular, we recover Veech's criterion (see Theorem~4 in~\cite{DFG19}). A dilation surface can be decomposed into Euclidean triangles by a family of saddle connections if and only if the following two conditions hold:
\begin{itemize}
    \item all of its singularities are conical singularities (Fuchsian singularities whose residue $\res$ satisfies $\Re(\res)<1$) or marked points;
    \item the surface does not contain any dilation cylinder of angle at least $\pi$.
\end{itemize}
\end{rem}

\Cref{thm:ClassificationComponents} yields an even more elementary decomposition of the surface, which is nevertheless particularly useful for studying the dynamics of the geodesic flow: for a non-exceptional surface,
\begin{itemize}
    \item \textit{Exterior domains} are the connected components of the union of the images of all the immersed open half-planes into the affine surface i.e. components of type (ii) to (v) in \Cref{thm:ClassificationComponents}, and of the non-conical Fuchsian singularities;
    \item The \textit{core} is the union of the Delaunay polygons and the Delaunay segments (its interior is the part of the affine surface that can be decomposed into Euclidean triangles drawn between the singularities).
\end{itemize}
In particular, any geodesic entering an exterior domain through its boundary never exits.

The notion of \textit{core} has been introduced in \cite{HKK17} and developed in \cite{Ta18} in the context of flat surfaces.

\begin{rem}
The dichotomy \textit{core/exterior domains} has also been used in \cite{T21} for an analogous decomposition of surfaces endowed with a cone spherical metric. In this settings, exterior domains are the connected components of the union of the images of the \textit{isometrically embedded half-spheres}. As spherical metrics and complex affine structures are special cases of complex projective structures (also called M\"{o}bius structures). It is conceivable that such decompositions could still be generalized to surfaces locally modeled on $\mathbb{CP}^{1}$, in which disks and half-planes (or half-spheres) are treated on an equal footing.
\end{rem}

The following result is in particular a bound on the topological complexity of the core.

\begin{thm}\label{thm:ComplexityBound}
Consider a finite type affine surface $(X,\nabla,\cal S)$ of genus $g$ with $n$ singular points (counted with multiplicity) that is non-exceptional, i.e.\ where $2g+n \geq 3$. We introduce the following numbers:
\begin{itemize}
    \item $\beta$ is the total number\footnote{See \Cref{sub:sides} about side counting.} of boundary sides of exterior domains (counted twice if both sides of a saddle connection are boundary edges of such domains).
    \item $t$ is the number of triangles in any triangulation\footnote{Some of the triangles may be self-folded, i.e.\ have two sides that coincide.} of the core (where the vertices are points of $\mathcal{S}$).
\end{itemize}
Then $t$ and $\beta$ satisfy $t+\beta = 4g-4+2n$.
\end{thm}

\Cref{thm:ComplexityBound} is proved in \Cref{sub:ComplexityBound}.

\subsection{About graph complements}\label{sub:sides}

\begin{figure}[htbp]
\begin{center}
\includegraphics[scale=0.9]{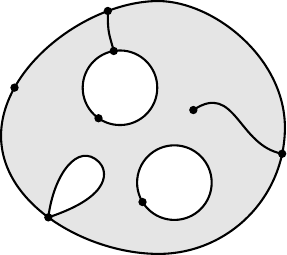}
\end{center}
\caption{In gray: a face in a graph complement. It is bounded by 10 edges, two of which bound it on both sides, so its side count is 12.}\label{fig:sides}
\end{figure}

Consider a compact orientable topological surface $X$ and an embedded finite graph $G$, with its vertices and edges. The complement $X-G$ consists in finitely may connected components that are open subset of $X$ and can be considered as generalized topological polygons which may or may not be homeomorphic to disks.
They are called \emph{faces}.
The boundary of a face consist in finitely many edges of $G$.
Each edge has two sides and an edge bounding a face $F$ may do it from one side or both, see \Cref{fig:sides}.
The \emph{side count} (or side number) of a face is the number of edge sides bounding it. It is greater or equal to the number of edges in its boundary.

If every face of a graph complement is 3-sided, then $3t=2e$ and $t=4g-4+2v$ where $t$ is the number of faces, $e$ the number of edges and $v$ the number of vertices.
This is an easy consequence of Euler's characteristic formula. Compare with \Cref{thm:ComplexityBound}.

\subsection{Relations with other works}

In recent years, complex affine structures have attracted interest from several mathematical communities.

\subsubsection{Affine interval exchange maps and dilation surfaces}

Dilation surfaces (a special case of complex affine structures in which the transition maps of the atlas are either translations or positive homotheties) form a class of geometric objects of major interest to specialists in dynamical systems, for the following reasons:
\begin{itemize}
    \item Like translation surfaces, they are equipped with directional foliations for which the first return map can be studied as a one-dimensional dynamical system. Just as directional flows on translation surfaces appear as suspensions of interval exchange transformations (IETs), directional flows on dilation surfaces appear as suspensions of affine interval exchange maps (AIETs).
    \item Still in analogy with the case of translation surfaces, moduli spaces of dilation surfaces are endowed with a natural action of $\text{GL}^{+}(2,\mathbb{R})$ playing the role of a \textit{renormalization operator} for the directional flows. The dynamics of the diagonal flow (also called the Teichm\"{u}ller flow) is closely related to the \textit{Rauzy-Veech induction} that renormalizes IETs and AIETs.
    \item Affine interval exchange maps form an intermediate class of one-dimensional dynamical systems, lying between IETs and generalized interval exchange transformations (GIETs), which are piecewise increasing homeomorphisms with finitely many discontinuities of the derivative. AIETs share with GIETs certain dynamical behaviors that do not occur in IETs, such as the existence of wandering intervals, as proved by Marmi–Moussa–Yoccoz in \cite{MMY10}.
    \item More recently, the phenomenon of \textit{affine shadowing} where the orbit of a GIET under renormalization is approximated by the orbit of an AIET appears in the proof by Ghazouani-Ulcigrai of a rigidity theorem for foliations on surfaces of genus two (see \cite{GU23}). 
\end{itemize}
Many conjectures and open problems related to dilation surfaces are discussed in \cite{G17,DFG19}. One of the most important among them concerns the existence of attracting closed geodesics on every dilation surface that has only conical singularities and is not a translation surface. To date, we only have a proof of the existence of closed geodesics (see \cite{BGT25}), without further information about the holonomy they carry.

\subsubsection{Moduli spaces of flat structures}

Complex affine structures are, in particular, complex structures, which makes it possible to construct an algebro-geometric moduli space. In~\cite{ABW}, affine structures are interpreted as meromorphic sections of a line bundle twisted by a character. The space of affine surfaces of genus~$g$ with $n$ Fuchsian singularities is constructed as an affine bundle over the moduli space~$\mathcal{M}_{g,n}$ of complex algebraic curves.
\par
In the special case of \textit{dilation surfaces}, an $\text{SL}(2,\mathbb{R})$-invariant measure, analogous to the Masur-Veech measure of the moduli space of translation surfaces, has been constructed in \cite{AS25}, drawing on a variant of period coordinates adapted to the twisting by a character.

\subsubsection{Compactification of moduli spaces of affine structures}

To date, no satisfactory compactification of the moduli space of affine structures is known. Constructing such a compactification would open up two particularly promising research directions:
\begin{itemize}
    \item It would enable the development of enumerative geometry and intersection theory on the moduli space, including the study of divisors corresponding to special loci and the computation of the corresponding cohomology classes.
    \item In the context of dilation surfaces, such a compactification would provide a natural setting for arguments based on degeneration toward the boundary of a stratum, allowing one to reduce to surfaces of lower complexity. This strategy could be key to establishing the existence of attracting cycles.
\end{itemize}

For translation surfaces (defined by meromorphic differentials), satisfactory compactifications exist from the algebro-geometric point of view. These involve adding, as degenerate objects, translation structures on the irreducible components of \textit{nodal curves}, subject to certain compatibility conditions at the nodes (see \cite{BCGGM,BCGGM2}). Such an approach is only possible because the degeneration of a translation structure always entails the degeneration of the underlying complex structure. No such relation holds for affine structures. As shown in \cite{ABW}, the moduli space is an affine bundle over $\mathcal{M}_{g,n}$ so degenerations may occur both in the base and in the fiber.
\par
In a recent preprint (see \cite{ABW1}), an algebro-geometric compactification was proposed. Degenerations in the base correspond to a notion of affine structure on a stable curve. In contrast, the fibers being affine spaces are compactified abstractly by adjoining a projective space. At present, there is no known geometric interpretation of these boundary points.
\par
Multiple poles of meromorphic connections have to play a role in the compactification of moduli spaces of affine surfaces with only simple poles. Indeed, the cone of infinite angle (also called the \textit{exponential-affine plane} in \Cref{sub:ExpoAffine}) has an affine structure induced a meromorphic connection having a double pole.
It is a credible candidate for the limit of the family of infinite cones $\mathcal{C}_{\theta}$ of angle $\theta$ as $\theta \tends + \infty$ (each corresponding to a meromorphic connection on $\mathbb{CP}^{1}$ with a pair of Fuchsian singularities of residues $1\pm\frac{\theta}{2\pi}$ that would merge into a double pole of residue $2$).
\par
In \cite{C1}, a situation is studied where an affine surface is defined via gluing polygons. It uniformizes to the Riemann sphere with five punctures. The corresponding connection is meromorphic and has simple poles at the punctures. When one of the polygons is deformed so as to degenerate, four poles merge into a pair of double poles, and the underlying affine surface is proved in \cite{C2} to converge in a specific sense.

\subsubsection{Bidimensional holomorphic dynamics}

Another remarkable connection between complex affine surfaces and holomorphic dynamics arises from the work of Abate–Tovena (see \cite{AT11}) or from the work of Guillot-Rebelo (see \cite{GR}). Given a homogeneous polynomial vector field on $\C^2$, one understands fairly well its complex-time trajectories which are either lines passing through the origin or Riemann surfaces which cover $\RS$ minus finitely many points. However, understanding the real-time trajectories of such a vector field remains a challenge. Abate-Tovena and Guillot-Rebelo proved that those trajectories project to geodesics of a meromorphic connection on $\RS$. 
\par
A related problem is the classification of holomorphic vector fields on $\C^2$ in the neighborhood of a singularity. As of today this classification is not complete, even if one restricts to a neighborhood of a separatrix passing through such a singularity. Abate and Tovena gives a classification when the separatrix corresponds to a Fuchsian singularity of the meromorphic connection. The classification in the general case requires classifying meromorphic connections in the neighborhood of irregular singularities. 
\par
Another related problem is the study of the discrete dynamics of germs tangent to the identity in the neighborhood of a fixed point in $\C^2$. It is strongly related to the study of the real-time dynamics of homogeneous polynomial vector fields in $\C^2$. On one hand, the time-$1$ map of such a vector field is tangent to the identity at the origin (when the degree of homogeneity is at least $2$). On the other hand, if $f:(\C^2,0)\to (\C^2,0)$ is of the form $f(x,y) = (x,y) + (P(x,y),Q(x,y)) + h.o.t.$ for some homogeneous polynomial map $(P,Q):\C^2\to \C^2$ of the same degree $\geq 2$, then, near $0$, the discrete dynamics of $f$ shadows the real-time trajectories of the vector field $P(x,y)\partial/\partial x + Q(x,y)\partial/\partial y$. And as mentioned previously, understanding the dynamics of $f$ near $0$ will require to understand the geodesic flow of a meromorphic connection associated to this vector field.

\subsubsection{Log-Riemann surfaces}\label{subsub:LogRiemann}

The log-Riemann surfaces have been introduced by Biswas and Perez-Marco
(see \cite{BP} and \cite{BP2}), motivated by questions in differential Galois theory. They are ``identity surfaces'' in the sense that their atlases have their changes of charts equal to the restrictions of the identity $\C\to\C$.
They are a fortiori translation surfaces hence a fortiori affine surfaces. In the aforementioned articles is studied a particular case, whose straightening is given by $\int Qe^P$ with $P$ and $Q$ polynomials. If $\deg P \geq 1$, the corresponding affine surface is an example for which the corresponding connection has a multiple pole at infinity, of order $1+\deg P$ i.e.\ $d-1=\deg P$, and residue $2+\deg Q$. Since they are identity surfaces, the affine holonomy is the identity. Their asymptotic value families are thus periodic and are the infinite ramification points mentioned in the above articles.

\subsubsection{Spaces of stability conditions}

One of the major research themes related to the homological mirror symmetry conjecture is to understand the similarities between stability condition spaces on Fukaya categories (and related variants) and moduli spaces of Abelian or quadratic differentials on Riemann surfaces. One of the major results was the explicit realization by Bridgeland and Smith in \cite{BS15} of spaces of quadratic differentials with simple zeros as the space of stability conditions on certain classes of triangulated categories defined using quivers.
\par
However, for a certain class of triangulated categories, namely partially wrapped Fukaya categories of surfaces, the quadratic differentials that appear may have exponential-type singularities, that is, locally of the form $f(z)e^{g(z)}dz^{2}$ where $f,g$ are meromorphic functions (see \cite{HKK17}). At poles of $g$, the flat structure (which is in particular an affine structure) admits irregular singularities.
\par
In this respect, the foundational work we carry out here on the local geometric models of higher-order poles of meromorphic connections will make it possible to include these exponential-type singularities within the standard, readily accessible objects of flat geometry.

\subsection{Organization of the paper}

\begin{itemize}
    \item In \Cref{sec:Background}, we provide background on the notion of complex affine surfaces and their link with holomorphic connections. We define meromorphic connections, order and residue of the poles, and finite type affine surface. We recall the classification of simple poles up to holomorphic equivalence, and also state and prove the classification of poles of any order up to \emph{formal} equivalence (\Cref{thm:FormalEquivalence}).
    Finally we present the exponential-affine plane and list all finite type affine surfaces with $2g+n\leq 2$. 
    \item In \Cref{sec:AsymptoticValues}, we introduce the asymptotic values invariant and construct local models of irregular singularities in affine geometry, realizing every possible asymptotic values invariant. We also prove in \Cref{thm:main} that this provides a complete set of invariants of multiple poles up to local biholomorphism.
    We define foci associated to irregular singularities and describe the behavior of geodesics near them.
    In \Cref{subsub:parabo} make an analogy with parabolic fixed points in holomorphic dynamics.
    \item In \Cref{sec:ConsequencesLocalModels}, we build on the previously established local models of irregular singularities to provide several technical results needed in the subsequent sections: counting disjoint simple saddle connections and generalizations, extension to sectors of geodesic ending on singuarities in finite or infinite time, extension to cylinders of closed geodesics. We also define anti-conical and swath domain and prove that any geodesic starting from their boundary and pointing inward tends to the associated singularity.
    \item In \Cref{sec:AffImm}, we provide classification results for affine immersions of open triangles, sectors, disks, half-planes and planes into a finite type affine surface. 
    In particular, affine embeddings of convex sets always extend to the boundary.
    \item In \Cref{sec:Delaunay}, we first introduce the Delaunay category, whose objects are affine immersions of disks, half-planes and planes, classify the maximal elements in terms of boundary behavior, organize them into a graph that we call the affine Delaunay spine. From this and the previous results, we provide proofs of \Cref{thm:ExceptionalSurfaces,thm:ClassificationComponents,thm:ComplexityBound}.
    \item In \Cref{sec:ExistenceProblem}, we determine the minimal number of Fuchsian singularities we need to add for a given isomorphism class of double pole to exist on the Riemann sphere.
\end{itemize}

\par\noindent{\bf Acknowledgements.}
Part of this research was conducted during a semester on Holomorphic Dynamics and Geometry of Surfaces funded by the ANR LabEx CIMI (grant ANR-11-LABX-0040) within the French State program “Investissements d’Avenir”. The research of the authors is supported by the French National Research Agency under the project TIGerS (ANR-24-CE40-3604). Research by G.T. is also supported by the Beijing Natural Science Foundation (Grant IS23005).

\section{Background on meromorphic connections and complex affine structures}\label{sec:Background}

\subsection{Affine surfaces, geodesics and meromorphic connections}\label{sub:minicours}

Consider an \emph{affine surface} $\cal A$, i.e.\ a topological surface with an atlas whose changes of charts are locally of the form $z\mapsto az+b$, $a\in\C^*$, $b\in \C$.
One of many equivalent definitions of a \emph{geodesic} on an affine surface is: a parametrized curve $\gamma:I\to \cal A$ where $I\subset\R$ is a non-empty open interval and such that in every affine chart $A$, $\gamma$ has uniform rectilinear motion, i.e.\ 
\[\gamma_A(t) = at+b\]
for some $a\in \C^*$ and $b\in \C$ where 
\[\gamma_A:=A\circ\gamma\]
is the expression of $\gamma$ in the chart $A$.
Note that near a given parameter $t_0$, this condition is independent of the affine chart in which it is inspected.
The condition can be expressed as
\[\gamma_A''(t)=0.\]

Since the change of charts of an affine surface are holomorphic, it is also a Riemann surface.
If one considers any other Riemann chart $R$ of the surface, the curve $\gamma_R:=R\circ\gamma$ usually does not anymore have constant derivative.
The equation $\gamma_A''=0$ becomes equivalent to
\[\gamma_R'' = - \frac{\phi''\circ \gamma_R}{\phi'\circ \gamma_R}  \times (\gamma_R')^2\]
where $\phi:=A\circ R^{-1}$ is the change of chart, so $\gamma_A = \phi\circ \gamma_R$.
This is exactly the geodesic equation associated to a holomorphic connection $\nabla$ on $TX$ of Christoffel symbol $\Gamma = \frac{\phi''}{\phi'}$ in the chart $R$.

Recalling what connections are is not crucial for the present article. We only need to know how the Christoffel symbol transforms under a change of Riemann chart.
Let $R_1$ and $R_2$ be two Riemann charts and $\psi = R_2\circ R_1^{-1}$ be the change of charts, so that $\gamma_{R_2} = \psi\circ \gamma_{R_1}$.
The equation
$\gamma_{R_2}'' = - \Gamma_2\circ \gamma_{R_2} \times (\gamma_{R_2}')^2$
is equivalent to
$\gamma_{R_1}'' = - \Gamma_1\circ \gamma_{R_1} \times (\gamma_{R_1}')^2$
with
\begin{equation}\label{eq:g}
\Gamma_1 = \psi'\times\Gamma_2\circ\psi + \frac{\psi''}{\psi'}.
\end{equation}

Conversely, given any system of functions $\Gamma_R$ in charts $R$, such that \cref{eq:g} holds, one gets a notion of geodesic on $X$, which are the curves which, in charts, satisfy
\[\gamma_{R}'' = - \Gamma_R\circ \gamma_{R} \times (\gamma_{R}')^2,\]
a condition which is independent of the chart.
This actually defines a connection on $TX$.
If the functions $\Gamma_R$ are meromorphic, we say the connection is \emph{meromorphic} and we call $X^*$ the complement of the poles and possibly of a discrete set of marked points.
Note that, because of \cref{eq:g}, if a point is a pole in one chart, it is a pole in every chart.
Moreover, the residues are the same, since \cref{eq:g} can be reinterpreted into: the differential forms $\Gamma_1(z)dz$ and $\psi^{*}(\Gamma_2(z)dz)$ differ by a holomorphic term $\frac{\psi''}{\psi'}dz$.
The pole orders are the same too.

\begin{ex}\label{ex:C}
Consider $X=\RS$. On its canonical chart $\C$, call the variable $z$ and consider the trivial connection, defined by $\Gamma_z=0$.
This defines on $\RS$ a meromorphic connection: indeed in the chart $w=1/z$, \cref{eq:g} gives $\Gamma_w = -2/w$.
\end{ex}

More generally a connection on $\RS$ is meromorphic if and only if its Christoffel symbol $\Gamma$ in the canonical chart $\C$ is a rational function.

\begin{defn}\label{def:res}
  The \emph{residue} of a pole of a meromorphic connection on a Riemann surface is defined as the \emph{opposite} of the residue of the corresponding pole of the Christoffel symbol $\Gamma$, in any chart:
  \[\res := -\mathrm{residue}(\Gamma) = -\frac{1}{2\pi i} \oint \Gamma(z) dz\]
\end{defn}

Hence in \Cref{ex:C}, the residue we assign to the pole is $2$ (not $-2$).

\begin{rem}\label{note:minus}
There are at least three reasons why we added the minus sign.
First, note that though the Christoffel symbol $\Gamma$ is central in the expression of the $\nabla$ operator and the study of connections, it is actually $-\Gamma$ that appears in the geodesic equation.
Second, in \cite{BR}, the collection of 1-forms $NR = -\Gamma_R\, dR$ where $R$ varies among the Riemann charts is considered instead of the collection of their Christoffel symbols $\Gamma_R$.
Third, with this convention the sum of residues of a meromorphic connection on a compact Riemann surface of genus $g$ is equal to the Euler characteristic $\chi = 2-2g$. 
In view of the analogy with the Gauss-Bonnet formula which says that the integral of the Gaussian curvature against the area of a Riemannian metric on a compact surface is equal to $2\pi \chi$, this is more coherent.
Actually, for a meromorphic connection with an invariant metric (i.e.\ one for which the monodromies are only rotations and translations), and with only conical Fuchsian singularities, the curvature is concentrated at the cone points, and the Gauss-Bonnet formula is known to take the form of an integral over a finite sum of Dirac masses.
The weight of the Dirac mass is in this case exactly the residue of the connection at this point, as we defined above.
\end{rem}

The fact that we are in (complex) dimension one has a very strong consequence: any connection on $TX$ can be locally trivialized on $X^*$. In terms of geodesics, this means that for any system of compatible functions $\Gamma_R$ as above, there exists an atlas of charts $A$ for which $\Gamma_A = 0$, i.e.\ where the geodesics are solutions of $\gamma_A''=0$.
These charts are obtained by composing Riemann charts $R$ by local solution $\phi$ of the equation
\[\phi''/\phi' = \Gamma_R
,\]
i.e.\ one takes $A = \phi\circ R$.
This atlas turns $X^*$ into an affine surface, since the only holomorphic changes of variables preserving the equation $\gamma''=0$ are the affine ones.

\medskip

Let us mention that the resolution of the O.D.E.\ $\phi''/\phi' = \Gamma$ is relatively simple:
\begin{equation}\label{eq:intODE}
  \phi = \int \exp({\textstyle \int \Gamma}),
\end{equation}
where $\int$ refers to an antiderivative (so is defined up to addition of a constant).

\medskip

Given two holomorphic/meromorphic connections $\nabla^1$ and $\nabla^2$ on the same Riemann surface $X$, one can define in Riemann charts a holomorphic/meromorphic differential form $\omega$ by letting
\[\omega_R = (\Gamma_R^2-\Gamma_R^1) dz\]
When changing the Riemann chart, the term $\psi''/\psi'$ cancels-out in \Cref{eq:g} and $\omega_R$ changes as a differential form.
This defines the \emph{difference} $\omega = \Gamma^2-\Gamma^1$.
(More generally connections in higher dimensional spaces live in an affine vector bundle directed by the bundle of $(1,2)$-tensors). 
From the definition in charts, it is immediate that
\begin{equation}\label{eq:resD}
\res(\omega,p) = \res(\Gamma_1,p)-\res(\Gamma_2,p)
\end{equation}

\begin{rem}
There is a very interesting interpretation of all this in terms of an operator called \emph{nonlinearity} in \cite{BR}.
\end{rem}

A map between affine surfaces is called \emph{affine} if its expression in charts is locally non-constant and complex-affine.
Such a map is locally injective, but is not necessarily injective.

An analogue of the identity theorem holds: two maps from a \emph{connected} affine surface $\cal A$ to an affine surface $\cal A'$ and that coincide on some open subset of $\cal A$ will coincide on the whole set $\cal A$. We call this the \emph{affine identity theorem} or principle. Of course it is an immediate consequence of the classical identity theorem but it can also easily be proved directly.

Any geodesic $\gamma:I\to \cal A$, with $I\subset \R$ an interval, extends to an affine map $U\to\cal A$ where $U$ is a connected open subset of the canonical affine plane $\C$ and $I\subset U$. This for instance a consequence of the identity theorem in complex analysis.\footnote{Indeed, geodesics are analytic maps, so $\gamma$ has an analytic extension to a connected open set $U\subset\C$ containing $I$. In a chart of $\cal A$ near $\gamma(t_0)$, $t_0\in I$, since $\gamma$ is locally affine along $I$ it affine in a neighborhood of $t_0$ in $U$; then by the identity theorem for analytic maps the set of points of $U$ where the extension is locally affine is open and closed in $U$, so equals $U$.}

A \emph{germ of affine map} at $x\in \cal A$ to another affine surface $\cal A'$ is an equivalence class among the affine map whose domains contain $x$ and that map in $\cal A'$, under the relation of being equal in a neighborhood of $x$.
A germ retains its origin: if $x\neq x'$ no germ at $x$ is equal to a germ at $x'$.

\begin{lem}\label{lem:scAffExt}
Consider an affine surface $\cal A$, and a 
universal covering $\tilde {\cal A} \overset{\pi}\longrightarrow \cal A$.
\begin{enumerate}
\item There is an affine atlas on $\tilde {\cal A}$ such that $\pi$ is affine and it is unique up to atlas equivalence.
\item On $\tilde {\cal A}$, any germ of affine chart is the germ of a unique global affine (but not necessarily injective) map $\tilde {\cal A}\to\C$ called its \emph{developing map}.
\end{enumerate}
\end{lem}
\begin{proof}
The first claim is standard.
The second claim can be found for instance in \cite{Thurston} Definition~3.4.2.
(Alternatively we can argue that since $\nabla$ is locally flat and $\tilde {\cal A}$ simply connected, there is a global section of the connection $\nabla$ on $T\tilde{\cal A}$, that it is holomorphic, and that since we are in complex dimension one, this vector field is the inverse of a unique differential form $\omega$ whose local antiderivatives are straightening coordinates of the vector field, hence affine. By simple connectivity, $\omega$ is the differential of a global function, which is affine. It is then easy to adjust the global chart to coincide with any given germ.)
\end{proof}

We end this with a subtle point of terminology: given an affine surface, we define a \emph{closed geodesic} as a geodesic that eventually comes back to its initial point with the same initial direction, but possibly a different speed.

\subsection{Parallel transport, holonomy and turning number}\label{sub:Turning}
The notions developed here will be used in \Cref{sec:ConsequencesLocalModels,sec:AffImm,sec:Delaunay} to study finite type affine surfaces and construct their Delaunay decomposition.

\loctitle{Linear holonomy}
Holonomy of a connection, such as $\nabla$, is defined through parallel transport: given a $C^1$ path with non-vanishing derivative $\gamma:[0,1]\to X$ and an initial vector $v(0)\in T_{\gamma(0)}X$, there is a unique function $t\in[0,1]\mapsto v(t)\in T_{\gamma(t)} X$ such that, in affine charts $A$ the vector $v(t)$ has constant expression.
In Riemann charts $R$, this condition takes the form
\[v_R'(t)=-\Gamma_R(\gamma_R(t))\, \gamma'_R(t)\, v_R(t),
\]
which is maybe more readable it we remove the index $R$ and variable $t$:
\[v'=(-\Gamma\circ\gamma) \times \gamma' \times v, \quad \text{($v=v_R$, $\gamma=\gamma_R$, $\Gamma=\Gamma_R$)}.
\]
It is called the \emph{parallel transport of $v(0)$ along $\gamma$}.
The map $v(0)\mapsto v(1)$ is a $\C$-linear map from $T_{\gamma(0)} X$ to $T_{\gamma(1)}X$ that we denote $\on{tra}\gamma$.
Since the connection is locally flat, the condition that $\gamma$ is $C^1$ can be relaxed to $C^0$ (this is obvious if one looks in affine charts), and the map $\on{tra}_\gamma : T_x X \to T_{x'}X$ only depends on the homotopy class of $\gamma$ as a path for $x$ to $x'$.
Concatenation of paths is compatible with parallel transport in that, for instance:
\begin{equation}\label{eq:tra}
\on{tra}_{\gamma \cdot \gamma'} = \on{tra}_{\gamma'}\circ\on{tra}_{\gamma}.
\end{equation}

Given a closed loop $\gamma$ based in $x\in X$, its \emph{holonomy} is the function from $T_{x}$ to itself defined by the parallel transport above, i.e.\ that takes $v(0)\in T_x$ and maps it to $v(1)\in T_x$ where $t\mapsto v(t)$ is the parallel transport of $v(0)$ along $\gamma$. This gives us an element of $\on{GL}(T_x X)$ associated to the homotopy class of $\gamma$.
Since we are in complex dimension one, $\on{GL}(T_x X)$ is canonically isomorphic to $(\C^*,\times)$, so in this setting holonomy can be considered as taking values in $\C^*$.
We denote this complex number
\[\on{hol}(\gamma).\]
By \cref{eq:tra}, and again using that we are in complex dimension $1$, $\on{hol}(\gamma)$ only depends on the free homotopy class of $\gamma$ in $X^*$.

\begin{rem}
Since holonomy only depends on homotopy classes in our context, it is often called a \emph{monodromy}. We will not use this term here.
\end{rem}

The transport equation in a single Riemann chart $R$ rewrites as
\[\frac{v_R'}{v_R} = -\Gamma_R(z) dz\]
i.e.
\[\frac{v_R(t_2)}{v_R(t_1)} = \exp \int_{\gamma_R|_{[t_1,t_2]}} -\Gamma_R(z) dz \]
In particular, for an isolated singularity of a connection, expressed in a Riemann chart, the holonomy around any loop circling counterclockwise around it is equal to the multiplication by
\[\exp\left(-2\pi i \res\right)\]
where $\res$ is the residue of the singularity (see \Cref{def:res}).

\loctitle{Affine holonomy and development}
Above, we defined a linear version of holonomy. There is an affine version, which can also be related to the notion of $(G,X)$-structure of Ehresmann (see \cite{Thurston}, Chapter~3).

Let $(X,\nabla,\mathcal{S})$ be an affine surface. 
Consider a universal covering $\wt X \overset{\pi} \longrightarrow X^*$ and endow $\wt X$ with the affine structure such that $\pi$ is affine.
Given a path $\gamma$ on $X^*$ and germ of affine chart $g_0$, they lift by $\pi$ into a path $\tilde \gamma$ on $\wt X$ and a germ $\tilde g_0$ of affine chart at $\tilde \gamma(0)$. 
We saw in \Cref{lem:scAffExt} that on $\wt X$, $\tilde g_0$ is the germ of a unique global affine (but not necessarily injective) map $\phi:\wt X\to\C$.
We define the \emph{development} of $\gamma$ from germ $g_0$ as the function $\phi\circ\tilde\gamma$. It is independent of the choices of lifts.

If $\gamma$ is a loop, consider the unique deck transformation $h:\wt X\to \wt X$ of $\pi$ that sends $\tilde\gamma(0)$ to $\tilde \gamma(1)$.
Then
\begin{equation}\label{eq:ah}
\phi\circ h = L\circ\phi
\end{equation}
for some $L\in\AutC$ (it holds near $\tilde \gamma(0)$ for some $L$ and hence everywhere), and $L$ is called the \emph{affine holonomy} associated to $\gamma$ and $g_0$.
It can instead be considered as associated to $\phi$ and $h$.
The affine holonomy only depends on $g_0$ and the homotopy class of the loop $\gamma$.
Changing the germ $g_0$ to another germ $g_1$ but fixing the loop has the following effect: let $f\in\AutC$ such that $g_1=f\circ g_0$; then the affine holonomy for $g_1$ and $\gamma$ is the   conjugate by $f$ of  the affine holonomy for $g_0$ and $\gamma$.
The conjugacy class in $\Aut(\C)$ of the affine holonomy of $\gamma$ is thus independent of the choice of the germ $g_0$.
Moreover, it only depends on the free homotopy class of the loop: consider a path $\delta$ between two basepoints; if one transports the germ $g_0$ along the path, one gets a germ $g_1$ and the affine holonomy for $g_1$ and $\delta^{-1}\cdot\gamma\cdot\delta$ is equal to the affine holonomy for $g_0$ and $\gamma$.

\loctitle{Turning number}
The classical definition of the turning number of a $C^1$ loop $\gamma:\R/\Z\to \C$ with non-vanishing derivative is the winding number of $\gamma'$ around $0$. It is an element of $\Z$.
This notion generalizes to affine surfaces but is now $\R$-valued: it is the total algebraic (a.k.a.\ signed) lifted angle by which any development of $\gamma$ turns, divided by $2\pi$.
We can define it too for a $C^1$ path $\gamma:[0,1]\to X^*$ with non-vanishing derivative (immersion). We denote it $\omega(\gamma)$.
For any germ of chart at $\gamma_0$, denote by $t\mapsto \hat\gamma(t)\in\C$ a development. Then $\omega(\gamma)$ is also the classical turning number associated to the path $\hat\gamma(t)$.
It follows that if $\gamma$ is a $C^1$ loop with non-vanishing derivative then
\begin{equation}\label{eq:tnHol}
\omega(\gamma) \equiv \frac{1}{2\pi}\Im(\log( \on{hol}(\gamma))) \bmod \Z
\end{equation}
(note that $\log(\on{hol}(\gamma))$ is only defined modulo $2\pi i\Z$).
We generalize this below to loops with corners.

It can also be expressed with integrals, and this time without the modulo $\Z$.
For this we introduce the \emph{nonlinearity} of a curve $\gamma :I\subset \R \to X$. 
It is the $\C$-valued differential $N(\gamma)$ on $[0,1]$ which for every affine chart $A$ takes the expression
\[ \frac{\gamma_A''(t)}{\gamma_A'(t)} dt
.\]
One can check that it does not depend on the choice of affine chart.
The turning number can then be expressed as:
\begin{equation}\label{eq:tn}
\omega(\gamma) = \frac{1}{2\pi}\Im\int_{[0,1]} N(\gamma).
\end{equation}
If we are given a Riemann chart $R$ then $N(\gamma)$ is equal to
\begin{equation}\label{eq:tn0}
\left(\frac{\gamma_R''(t)}{\gamma_R'(t)} + \Gamma_R\circ \gamma_R(t) \times \gamma_R'(t) \right) dt.
\end{equation}
So the turning number of a portion $\gamma$ contained in a chart expresses as
\begin{equation}\label{eq:tn1}
\omega(\gamma) = \omega_\C(\gamma_R) + \int_{\gamma_R} \Gamma_R(z)dz
\end{equation}
where $\omega_\C$ refers to the classical turning number of a curve (here $\gamma_R$) traced in $\C$.
It follows in particular that the turning number of a small simple close $C^1$ loop winding once anticlockwise around a pole of $\nabla$ of residue $\res$ has a turning number equal to
\begin{equation}\label{eq:tn2}
\omega(\gamma) = 1 - \Re(\res)
\end{equation}
(recall that $\res$ is defined as the residue of $-\Gamma$).

On the space of $C^1$ paths (parametrized by $[0,1]$) or loops (parametrized by $\R/\Z$) with non-vanishing derivative (immersions) we put the $C^1$ topology, for which a path converges iff it converges uniformly and its differential too.
Immersions of $\R/\Z$ are equivalent to immersions of the circle, so we will call them this way, but keep in mind that for the turning number we look at the derivative from $\R/\Z$.

\begin{lem}\label{lem:tninv}
Given an affine surface $(X,\nabla,\mathcal{S})$, for any $C^1$ immersion $\gamma$ of the circle into $X^{\ast}$, there is a neighborhood $V$ of $\gamma$ for the $C^1$ topology, such that $\forall \delta \in V$, $\omega(\delta) = \omega(\gamma)$.
\end{lem}

\begin{proof}
The conditions imply that the turning numbers are close.
On the other hand, the curves $\delta$ and $\gamma$ are in particular free homotopic, so \cref{eq:tnHol} implies 
$\omega(\gamma)\equiv \omega(\delta) \bmod\Z$.
\end{proof}

A continuous path in the space of immersions with the $C^1$ norm is called a \emph{regular homotopy}. The statement above tells us that in the space of immersions of the circle to $X^*$, the turning number is invariant under regular homotopies.

\medskip

The turning number can be generalized to piecewise $C^1$ paths with non-vanishing derivative, but some attention shoud be paid. At a corner, the left and right derivatives cannot be opposite, i.e.\ the curve is not allowed to have a cusp. We say that the path is \emph{cusp-free}.
The convention is to define the turning number of such a path/closed loop as the sum of the turning numbers of its parts and of its corners, where the turn at a corner is the representative in $(-\pi,\pi)$ of the signed angle from the left derivative to the right derivative.
\Cref{eq:tnHol} generalizes in this case.

\begin{rem}
Topological questions are even more subtle as for instance with our initial definition of regular homotopy, we cannot change the parameter $t$ at which a (non zero angle) corner occurs. This means that to be useful, the definition of the topology has to be adapted. We will not try to do this here.
\end{rem}

Recall that a loop is called \emph{simple} when it is injective.

\begin{lem}\label{lem:GB}
If a piecewise $C^1$ simple and cusp-free loop in $X^*$ is the counterclockwise oriented boundary of a topological disk of $X$ containing poles\footnote{This would also work with essential singularities.} of residues $\res_j$, then its turning number is $1-\sum_j\Re(\res_j)$.
\end{lem}

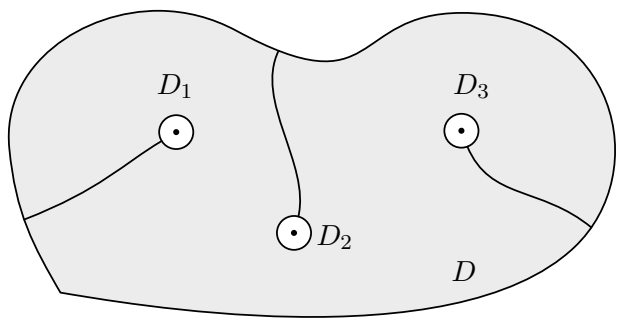
\begin{figure}[htbp]
  \begin{center}
    \begin{tikzpicture}
      \node at (0,0) {\includegraphics[width=8cm]{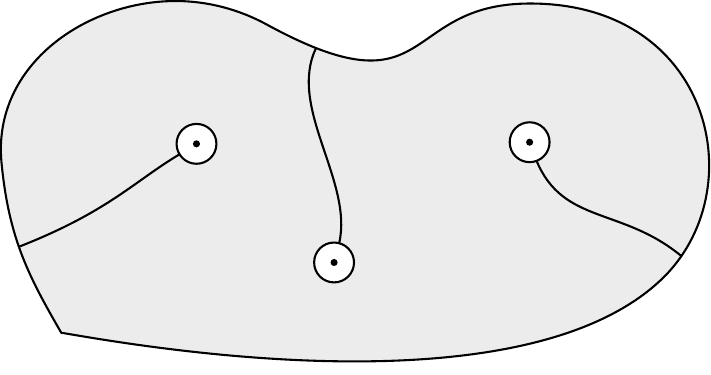}};
      \node at (2,-1.35) {$D$};
      \node at (-1.8,1.1) {$D_1$};
      \node at (0.3,-0.9) {$D_2$};
      \node at (2.1,1.1) {$D_3$};
    \end{tikzpicture}    
  \end{center}
  \caption{Illustration of the proof of \Cref{lem:GB}}
  \label{fig:GB}
\end{figure}

\begin{proof}
Call $D$ the (open) topological disk and $p_1$, \ldots, $p_m$ the poles of $\nabla$ in $D$.
Consider a chart near each pole in $D$, sending it to $0$, and a small circle in this chart.
Take the circles small enough so that the closed disks they bound correspond to disjoint subsets $D_j$ of $X$, contained in $D$ and containing only the pole $p_j$.
It is possible to link each disk to the boundary of $X$ by disjoint $C^1$ arcs $\gamma_j$ that are perpendicular to these disk and to $\partial D$, as in \Cref{fig:GB}.
The complement in $D$ of the closed disks $D_j$ and of the arcs is simply connected and hence can be given a flat conformal Riemannian metric that is constant in the affine charts.
We then apply the Gauss--Bonnet formula to this open subset: the total turning number $\tau$ of its contour (followed counterclockwise) is $2\pi$.
In the definition of $\tau$, $\gamma_j$ is followed in both direction and hence their turning number cancel-out. The endpoints of $\gamma_j$ add a total turning of $4\times 1/4=1$.
Hence
\[ 1 = \tau = \tau (\gamma) + m - \sum_{j=1}^m \tau(\partial D_j)
\]
where $\partial D_j$ is followed counterclockwise.
Now in the chosen chart around $p_j$ we apply \cref{eq:tn0,eq:tn} to get, denoting $\delta(t) r\,e^{2\pi i t}$ a parametrization of the image of $\partial D_j$ in the chart
\[
\begin{split}
\tau(\partial D_j) & = \frac{1}{2\pi}\Im\int_{0}^1 \left(\frac{\delta''}{\delta'} + \Gamma_j\circ \delta \times \delta' \right)
= 1 + \frac{1}{2\pi}\Im \int \Gamma(\delta) d\delta
\\
& = 1+\frac{1}{2\pi}\Im (2\pi i \res(\Gamma_j,0))
= 1+\Re \res(\Gamma_j,0)
= 1-\Re \res(p_j)
\end{split}
\]
as we defined the residue of a pole as the residue of $-\Gamma$ in a chart.
A simple computation then gives the result.
\end{proof}

The result above can be seen as a (known) generalization of the Gauss-Bonnet theorem. Another generalization leads to the following statement:

\begin{prop}\label{prop:sumRes}
If $(X,\nabla,\cal S)$ is a finite type affine surface, then the sum of the residues of all the singularities is equal to $2-2g$.
\end{prop}
\begin{proof}
We saw that the difference between two meromorphic connections is a meromorphic differential.
The sum of residues if a meromorphic differential is always $0$.
By \cref{eq:resD}, this implies that any two meromorphic connections on $X$ have the same sum of residues. 
It is thus enough to prove the claim for one meromorphic connection.
The compact Riemann surface $X$ admits a meromorphic vector field $V$. A version of a theorem of Poincaré says that the sum of all orders of its zeroes minus all orders of its poles is its Euler characteristic $2-2g$.
Denote $\cal S_V$ the singularities of $V$ (zeroes and poles).
Consider the translation surface structure on $X-\cal S_V$ given by the charts on which $V$ is constant of coefficient $1$.
This defines an affine surface structure on $X-\cal S_V$, which is actually meromorphic on $X$, with only simple poles, of residues equal to the order at the zeroes and minus the order at the poles.
\end{proof}

\subsection{Cylinders and skew cones}\label{sub:cylindersAndSkewCones}

We describe here a few classes of affine surfaces. They are not all of finite type according to our definition.\footnote{But they may be considered as finite type in a class of affine surfaces with geodesic boundary.}

\loctitlebold{Translation cylinders}
The \emph{full translation cylinder} is the quotient $\C/\Z$ of the canonical affine plane $\C$
by the group generated by $T:z\mapsto z+1$.
The notation $\C/\Z$ comes from $\Z$ being seen as a sub-group of $(\C,+)$.
Its subset $\mathbb{H}/\Z = \{[z]\in\C/\Z\,;\,\Im(z)>0\}$ is called the \emph{semi-infinite translation cylinder} and
$C_h:=\{[z]\in\C/\Z\,;\,0<\Im(z)<h\}$ is called the \emph{translation cylinder of height $h$}.
Note that the quotient by $\Z$ of the lower half plane bounded by $\R$ is isomorphic by $z\mapsto -z$ to the quotient of the upper half plane.

We define a boundary completion of translation cylinders with the following convention: by replacing the strict inequality by a large one in the definition of $C_h$.
Note that if an end is infinite then by we do not complete it, so the completion will not be compact in this case. 

\loctitlebold{Universal covering of the punctured affine plane} Consider $\C^*$ as a subset of the canonical affine plane $\C$, and let $\pi: \cet \to \C^*$ be a universal covering.
We endow $\cet$ with the pull-back by $\pi$ of the affine structure.
On $\cet$ the polar coordinates $(r,\theta)$ give a (non-conformal) bijection to $(0,+\infty)\times \R$.
The function $z\mapsto \log z$ has a determination $(r,\theta)\mapsto \log(r)+i\theta$ on $\cet$, which gives a conformal bijection to $\C$. In \Cref{sub:ExpoAffine} we use this to push the affine structure to a special affine structure on $\C$, which we call the \emph{exponential-affine plane} and study it further.
The group of affine automorphisms of $\cet$ consists in the lifts of linear maps: $(r,\theta)\mapsto (sr,\theta+\alpha)$, for some $s>0$ and $\alpha\in\R$, is a lift of $z\mapsto se^{i\alpha} z$.
It is simply transitive (the product of complex numbers naturally lifts to $\cet$ to give a commutative group law and its affine automorphisms coincide with  multiplication).
The support of maximal geodesics, i.e.\ lifts to $\cet$ of straight lines in $\C^*$, form two orbits under $\on{Aut}(\cet)$, whose representatives can be taken as: a half-line through $0$, a straight line not containing $0$.

\loctitlebold{Dilation cylinders}
Let $\lambda\in\R$ with $\lambda>1$.
We call \emph{full Reeb cylinder} of factor $\lambda$ the quotient of the affine surface $\cet$ above by the group $\Lambda$ generated by the affine map $(r,\theta)\mapsto (\lambda r,\theta)$.
We call \emph{semi-infinite Reeb cylinder} of factor $\lambda$ the subset of the full Reeb cylinder defined by $\theta>0$ and the one defined by $\theta<0$. They are \emph{not isomorphic}.
We call \emph{dilation cylinder of angle $\alpha$} and factor $\lambda$ the subset of the full Reeb cylinder defined by $\theta\in(0,\alpha)$.
We reserve the denomination \emph{Reeb cylinder} to the cylinders whose opening angle is $>\pi$, or infinite. Reeb cylinders will play a special role in the Delaunay decomposition of a finite type affine surface, while smaller cylinders do not appear in this decomposition.
A semi-infinite Reeb cylinder is illustrated in \Cref{fig:reeb}. 

We define a boundary completion of dilation cylinders in a similar way to translation cylinders: by replacing the strict inequality by a large one in the definition of the subset.
Note that if an end is infinite then by we do not complete it, so the completion will not be compact in this case. 

\medskip

Any affine surface isomorphic to one of the named surfaces above is called by the same name.

\medskip

Recall that a \emph{closed geodesic} is a geodesic that eventually comes back to its initial point with the same initial direction, but possibly a different speed.
\emph{Periodic} geodesics are the special case where the speed is the same.
Cylinders are foliated by the support of closed geodesics.
More precisely translation cylinders are foliated by the horizontal circles of equation $\Im(z) = \mathrm{cst}$ which are the support of periodic geodesics.
While dilation cylinders are foliated by circles that are quotients of half lines radiating from $0$: they are the support of the closed geodesics $t>0\mapsto (t,\theta)/\Lambda$.
These geodesics are infinite lived and progressively slow down as $t\to+\infty$, while they accelerate and blow up in finite time as $t\tends 0$.

\medskip

The full Reeb cylinder $R$ has an automorphism group that consists in the quotient of the automorphisms of $\cet$ and is simply transitive.
The support of maximal geodesics form two orbits under $\on{Aut}(R)$: support of closed geodesics (quotient of radial lines from $0$), support of the other ones, that all accumulate a closed geodesic in the future and another one in the past, separated by an angle of $\pi$.

\begin{figure}
\begin{tikzpicture}
\node at (0,0) {\includegraphics[scale=0.75]{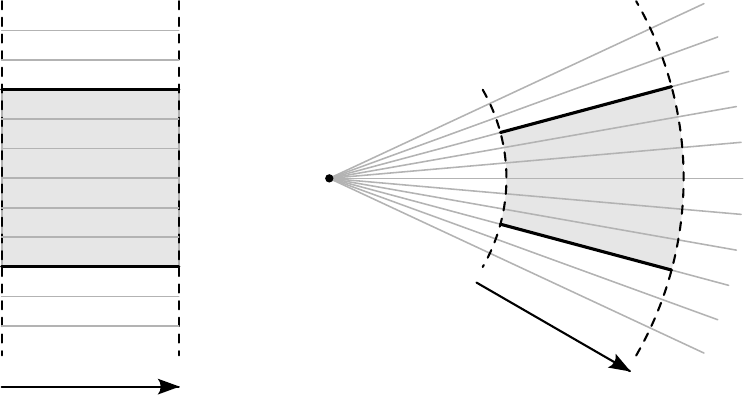}};
\end{tikzpicture}
\caption{A translation cylinder of finite height and a dilation cylinder of finite angle, obtained by quotienting the gray fundamental domain by a translation and a dilation.}
\label{fig:ueras}
\end{figure}

\loctitlebold{Skew cones}
Let $s>0$ and $\alpha\in(0,+\infty)$ and let $\Lambda$ be the group generated by the affine automorphism $\lambda:(r,\theta)\mapsto (rs,\theta+\alpha)$ of $\cet$, which is a lift of the complex multiplication $z\mapsto se^{i\alpha}z$.
A fundamental domain is given by the sector $S$ of equation $\theta\in[0,\alpha]$.
The quotient affine surface $C=\cet/\Lambda$ is called the \emph{skew cone} of factor $s$ and angle $\alpha$. 
An equivalent definition is by gluing the two sides of the fundamental domain by $\lambda$: $C=S/\lambda$.

Skew cones have self-intersecting geodesics if and only if their opening angle is $<\pi$.
Two such examples (with some part cut out) are given in \Cref{fig:cone-affine,fig:unbounded}.

The affine surface $C$ has a simply transitive group of affine automorphisms $\on{Aut}(C)$ that consist in the quotients of the automorphisms of $\cet$.
The support of maximal geodesics form two orbits under $\on{Aut}(C)$, whose representatives can be taken as: a half-line through $0$, a straight line not containing $0$.

Skew cones are isomorphic to the regular set of a meromorphic connection on $\RS$: let $a \in\C^*$ such that $a\times(\log(s)+i\alpha) =2\pi i$; there is a well-defined branch of $\log z$ on $\cet$; the map $z\mapsto \exp(a\log(z)) = \exp(a\times(\log(r)+i\theta))$ quotients to a holomorphic bijection from the skew cone $C$ to the Riemann surface $\C^*$, and is a determination of $z\mapsto z^a$.
Push the affine structure of $C$ to $\C^*\subset \C\subset\RS$ by this map.
Then affine charts on $\C^*$ are given by local branches of $z\mapsto z^{1/a}$ and \cref{eq:g} gives the meromorphic expression
\[\Gamma_z = \frac{c}{z}\]
with $c = a^{-1}-1 = \frac{\log(s)+i\alpha}{2\pi i}-1$.
It has two poles, one of residue $-c = 1-a^{-1}$ at $0$, and one at $\infty$ of residue $2+c = 1+a^{-1}$.

\begin{figure}[htbp]
  \begin{center}
    \begin{tikzpicture}
      \node at (0,0) {\includegraphics[width=8cm]{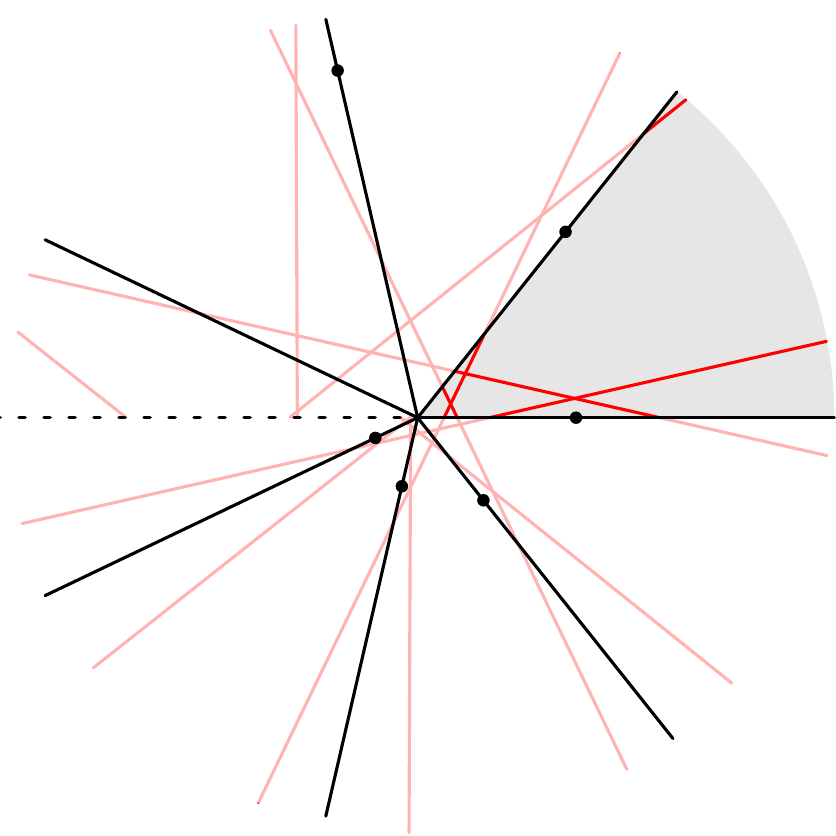}};
    \end{tikzpicture}
  \end{center}
  \caption{A skew cone of angle close to $\pi/4$ and dilation factor close to $1.4$. It is the quotient of $\cet$ by a subgroup $\Lambda$ of $\on{Aut}(\C^*)$ generated by some linear map $\lambda$.
  A portion of width $2\pi$ of $\cet$ is shown as $\C \setminus (-\infty,0]$. The orbit of a geodesic by $\Lambda$ is shown in red. Its intercept with the fundamental sector trace a single geodesic in the skew cone.
  }\label{fig:skew-cone}
\end{figure}

\subsection{Flat local models of Fuchsian singularities}\label{sub:Fuchsian}

This section will be used in \Cref{sec:ConsequencesLocalModels,sec:AffImm,sec:Delaunay}.
Its content is well-known, proofs can be for instance found in\footnote{In \cite{CT}, the convention for the residue $\res$ is the opposite complex number.} \cite{CT}.

We include Fuchsian singularities of residue zero to take into account marked points. 
Formally this means the following:

\begin{defn}
Consider a Riemann surface $X$ and a meromorphic connection $\nabla$ on $X$. Let $X^*$ be the complement in $X$ of the poles and of a (possibly empty) discrete subset of marked points.
We call \emph{Fuchsian singularity} the marked points and the poles of order $1$ of the connection.
\end{defn}

A Fuchsian singularity is a marked point if and only if its residue is $0$.
Morally, marked points are conical singularities of angle $2\pi$ and without dilation.

\medskip

A Fuchsian singularity has a punctured neighborhood isomorphic as an affine surface to one and only one of the following models, defined as quotients of subsets of the canonical affine plane $\C$. Most models are subsets of the surfaces introduced in \Cref{sub:cylindersAndSkewCones}.

\loctitlebold{Finite angle bounded sector for conical singularities}
The model is given by a neighborhood of $0$ in a the skew cone. More precisely, take the quotient, considered only near $0$, of the sector $\arg z \in [0,\alpha]$ (with $\alpha>0$ possibly $>2\pi$, in which case, work in the universal covering of $\C^*$ or in log coordinates) by identifying $x>0$ on one of its boundary lines to $xse^{i\alpha}$ on its other line, $s>0$.
We saw that
\[\res=1-\frac{\log s + i\alpha}{2\pi i} = 1-\frac\alpha{2\pi}+i\frac{\log s}{2\pi}
.\]
Topologically, the singularity is at $0$.
This kind of singularity is called \emph{conical} and is characterized by:
\[\Re(\res) \in (-\infty,1).\]

\loctitlebold{Semi-infinite translation cylinder for cylindrical singularities} It is the quotient $\mathbb{H}/\Z$ of the upper half-plane by the group of translations generated by $z\mapsto z+1$, for which a fundamental domain is $\{z\in\mathbb{H}\,;\,\Re(z) \in [0,1]\}$.
Topologically, the singularity is at the upper end.
This kind of singularity is called \emph{cylindrical} and is characterized by:
\[\res = 1.\]
Its affine holonomy is a non-zero translation.

\loctitlebold{Semi-infinite Reeb cylinder for Reeb type singularities}
There are two non-isomorphic variants.
First, the quotient by a dilation of factor $s>1$ of the sets of points in the universal covering of $\C^*$ satisfying $\arg z\in (0,+\infty)$.
A point tends to the singularity iff its argument tends to $+\infty$.
This kind of singularity is characterized by:
\[\res = 1+i\frac{\log s}{2\pi}\]
\[\Re(\res) =1,\quad \Im(\res)> 0.\]
The second variant is the quotient by a dilation of factor $s>1$ of $\arg z\in (-\infty,0)$ (the other half).
A point tends to the singularity iff its argument tends to $-\infty$.
This kind of singularity is characterized by:
\[\res = 1-i\frac{\log s}{2\pi}\]
\[\Re(\res) =1,\quad \Im(\res)< 0.\]
Both variants of singularity are called of \emph{Reeb type}.

\loctitlebold{Finite angle unbounded sector for anti-conical singularities}
This kind is characterized by:
\[\Re(\res) \in (1,+\infty)\]
and splits into two subkinds.
One kind, which we call \emph{pure anti-conical} given by a neighborhood of $\infty$ in the skew cone. More precisely take a neihborhood of $\infty$ in the sector in $\wt\C^*$ defined by $\arg (z)\in[0,\alpha]$ for some $\alpha>0$ and glue the two sides of the sector by the linear map $z\mapsto se^{i\alpha}$ for some scaling factor $s$.
We saw that
\[\res=1-\frac{\log s - i\alpha}{2\pi i} = 1+\frac\alpha{2\pi}+i\frac{\log s}{2\pi}.\]
The other kind, which we call \emph{shifted anti-conical}, only concerns $\alpha=2\pi n$,
$n\in\N^*$ and $s=1$.
The two sides of the sector are glued, near infinity, by the translation sending the point $x$ of argument $0$ to the point $x+1$ of argument $2\pi n$. For this kind,
\[\res=1-\frac{\log s - i\alpha}{2\pi i} = 1+n+i\frac{\log s}{2\pi}.\]
Its affine holonomy is a non-zero translation.
In both cases, topologically, the singularity is at infinity.

\medskip

Note that the conical singularities are the only Fuchsian singularities having a neighborhood of bounded area. Moreover:

\begin{lem}\label{lem:fuchsianTips}
  Let $p$ be a Fuchsian singularity or a marked point in an affine surface $(X,\nabla,\mathcal{S})$, then:
  \begin{enumerate}
    \item If $p$ is of Reeb type, then no geodesic can accumulate on $p$. 
    \item There exists a geodesic tending in finite time to $p$ if and only if $p$ is a conical singularity or a marked point.
    \item There exists a geodesic tending in infinite time to $p$ if and only if $p$ is anti-conical or of cylindrical type.
  \end{enumerate}
\end{lem}

Each of these models, except the two whose affine holonomy is a non-zero translation, can be obtained by taking in the exponential-affine plane $\cal E$ of \Cref{sub:ExpoAffine} the quotient by the translation $T_t$, $t=2\pi i(1-\res)$, of the $T_t$-invariant half-plane of equation $\Im (z/t)>0$.

Note that the residue characterizes the model type and that if $\res\notin\{2,3,\ldots\}$ all Fuchisian singularity with this residue are locally isomorphic, while for $\res\in\{2,3,\ldots\}$, there are two non-isomorphic possibilities.

\begin{figure}[htbp]
  \begin{center}
    \begin{tikzpicture}
      \node at (0,0) {\includegraphics[width=10cm]{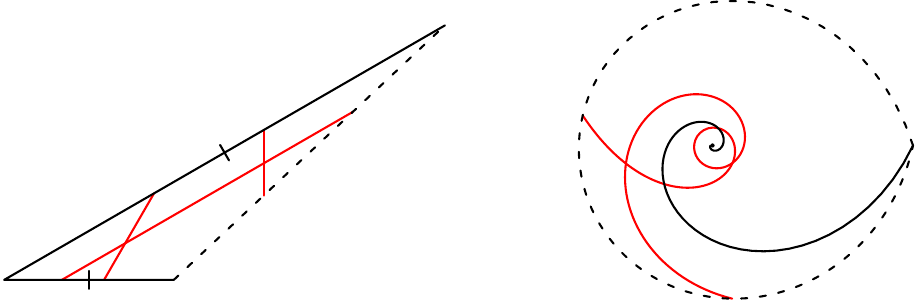}};
      \node at (0,-0.5) {$z\mapsto z^\beta$};
    \end{tikzpicture}
  \end{center}
  \caption{Finite angle bounded sector.
  We have a sector of angle 30\textdegree\ and a gluing of factor 3 between its sides. We restrict to the subset of the sector that ends at the dotted (straight) line, whose shape is not so important.\\
  Right: a Riemann chart of the flat model is given by a branch of the map $z\mapsto z^\beta$ with $\beta = 2\pi i/(\log(3) + 2\pi i\frac{30}{360})$.
  }\label{fig:cone-affine}
\end{figure}

\begin{figure}[htbp]
  \begin{center}
    \begin{tikzpicture}
      \node at (0,0) {\includegraphics[width=10cm]{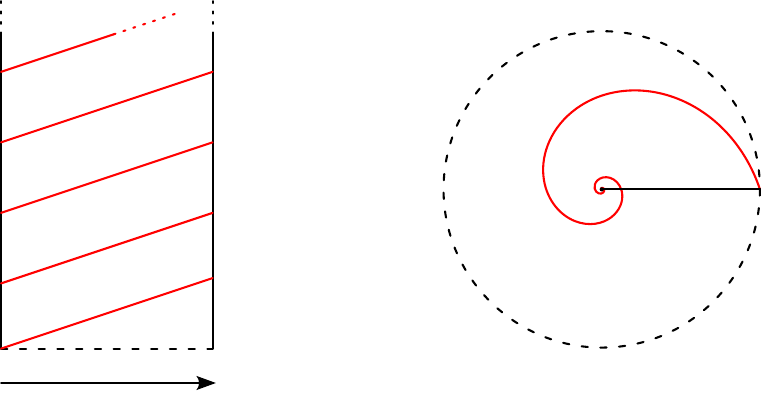}};
      \node at (-0.6,-1) {$z\mapsto \exp(2\pi i z)$};
      \node at (-3.6,-2.8) {$T_1$};
    \end{tikzpicture}
  \end{center}
  \caption{Semi-infinite translation cylinder. The quotient $\mathbb{H}/\Z$ is also the quotient of the fundamental domain ``$\Re(z)\in[0,1]$ and $\Im(z)>0$'' by identifying gluing the two sides by a horizontal translation.}\label{fig:cylindre}
\end{figure}

\begin{figure}[H]
  \begin{center}
    \begin{tikzpicture}
      \node at (0,0) {\includegraphics[width=12cm]{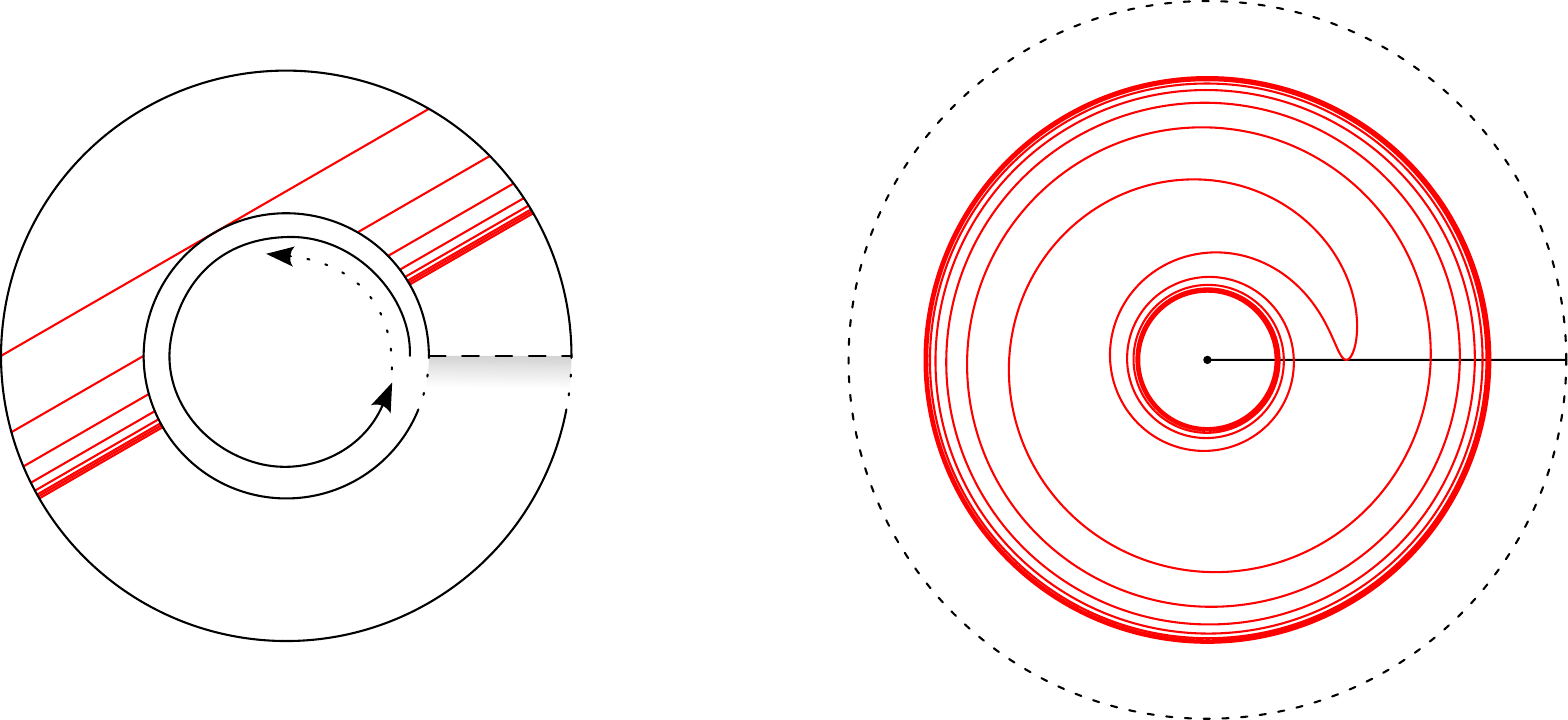}};
      \node at (-0.7,-1.8) {$z\mapsto z^{2\pi i/\log 2}$};
    \end{tikzpicture}
  \end{center}
  \caption{\small Semi-infinite Reeb cylinder. On the left, in the universal covering of $\C^*$ we took a fundamental domain of the form $|z|\in[1,2]$, $\arg z>0$. It winds infinitely may times, in one direction, over the annulus between radii $1$ and $2$. The identification is $z\sim 2z$ for $|z|=1$.
  It is mapped conformally to a pointed disk by $z=r e^{i\theta} \mapsto z'=z^{2\pi i/\log 2} = r' e^{i\theta'}$ with $r'=\exp(-a \theta)$, $a=2\pi /\log 2$ and $\theta' = 2\pi \frac{\log r}{\log 2}$.
  Since this map is too extreme (the ratio between the radii of the outer and inner limit cycles is $e^{\pi a}\approx 2.3\times 10^{12}$) we modified it on the right in a non-conformal way by using instead $r'=\exp(-a'\theta)$ with $a'=a/20$.}
  \label{fig:reeb}
\end{figure}

\begin{figure}[htpb]
  \begin{center}
    \begin{tikzpicture}
      \node at (0,0) {\includegraphics[width=13cm]{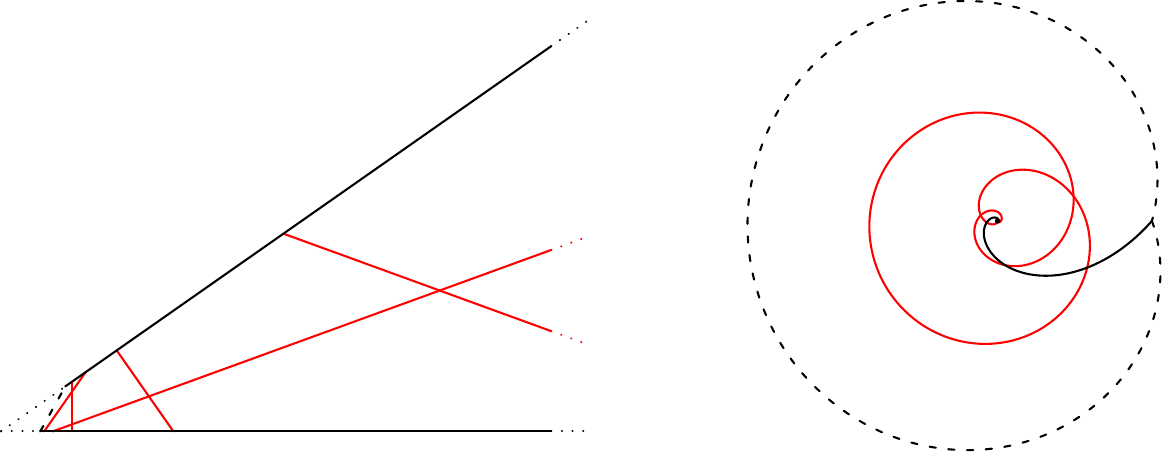}};
      \node at (0.9,-0.8) {$z\mapsto z^{\beta}$};
    \end{tikzpicture}
  \end{center}
  \caption{Finite angle unbounded sector. 
  We have a sector of angle 35\textdegree\ and a gluing of factor 2 between its sides. We restrict to the subset of the sector to the right of the dotted segment.\\
  Right: a Riemann chart of the flat model is given by a branch of the map $z\mapsto z^\beta$ with $\beta = -2\pi i/(\log(2) + 2\pi i\frac{35}{360})$.}\label{fig:unbounded}
\end{figure}

\begin{figure}[htbp]
  \begin{center}
    \begin{tikzpicture}
      \node at (-3.5,0) {\includegraphics[width=6cm]{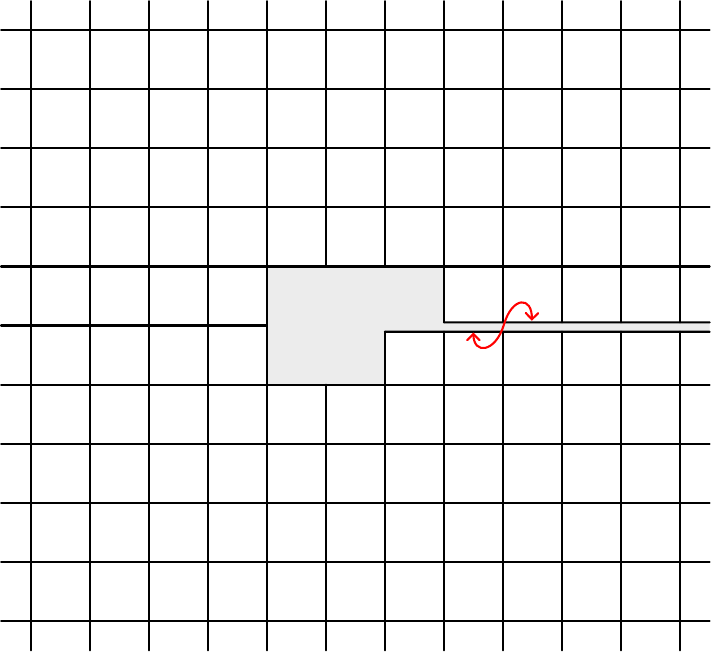}};
      \node at (3.5,0) {\includegraphics[width=6cm]{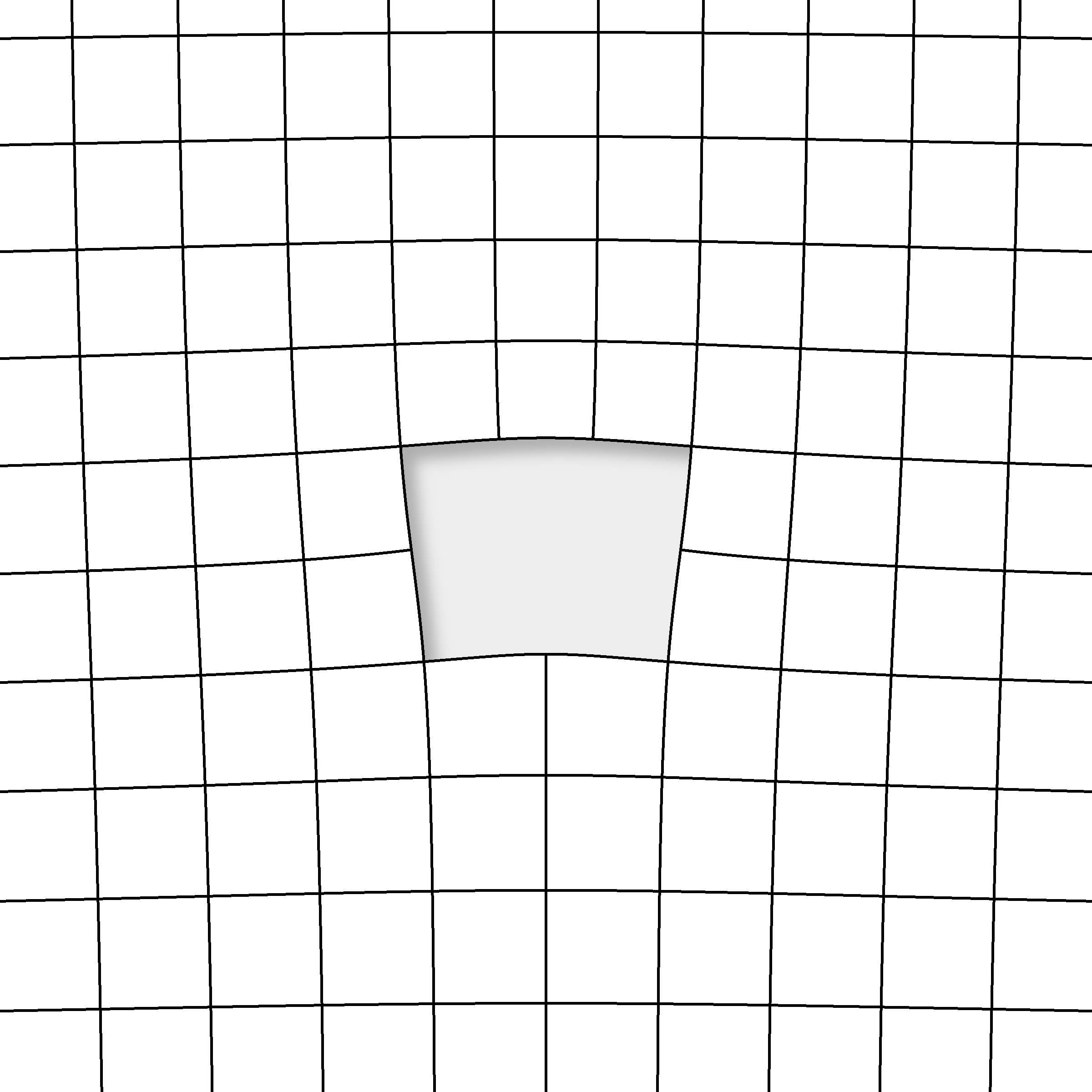}};
      \draw[->] (0.4,0) -- node[above]{$\scriptstyle f$} (-0.4,0);
    \end{tikzpicture}
  \end{center}
  \caption{Finite angle unbounded sector with a translation. Here the angle is $2\pi$. Left: affine chart, with a straight grid drawn on it. The slit is glued by a horizontal translation according to the red arrow. Right: conformal chart, but with the singularity at infinity in $\C$. For this illustration we chose $f(z) = z - \frac{1}{2\pi i}\log(z)$.}\label{fig:rec}
\end{figure}

\begin{figure}[htbp]
  \begin{center}
    \begin{tikzpicture}
      \node at (0.42,0) {\includegraphics[scale=0.863532574]{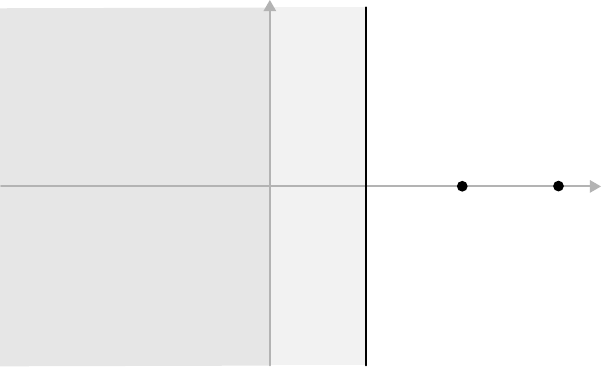}};
      \node at (3,-3.5) {unbounded};
      \node at (-0.5,-3.5) {affine cone};
      \node at (-2.1,-3.1) {$\alpha>2\pi$};
      \node at (0.7,-3.1) {\small $\alpha<2\pi$};
      \node at (1.7,1.45) {\rotatebox{90}{Reeb}};
      \node (A) at (2.4,-1.2) {cyl.};
      \draw[-] (1.38,-0.05) -- (A)[anchor = north west];
      \node (B) at (-1.4,1.2) {eras.};
      \draw[-] (0,-0.05) -- (B)[anchor = south east];
      \node at (-0.23,-0.3) {$0$};
      \node at (1.19,-0.3) {$1$};
      \node at (2.54,-0.3) {$2$};
      \node at (3.98,-0.3) {$3$};
    \end{tikzpicture}
  \end{center}
  \caption{Kind of Fuchsian singularity, given its residue $\res$, represented in the complex plane. In gray the real and imaginary axes. $\res=0$ corresponds to an erasable singularity, $\res=1$ to the semi-infinite translation cylinder. $\Re \res<1$ to a conical singularity. $\Re \res>1$ to a finite angle unbounded sector. $\res = 1+iy$ with $y\in\R^*$ to a semi-infinite Reeb cylinder. Dots are placed at $\res \in\{2,3,\ldots\}$, where there are two inequivalent isomorphism classes.
  Dilation surfaces satisfy $\Re\res\in\Z$, surfaces with an invariant metric $\res\in\R$ and translation surfaces $\res\in\Z$.}\label{fig:res}
\end{figure}
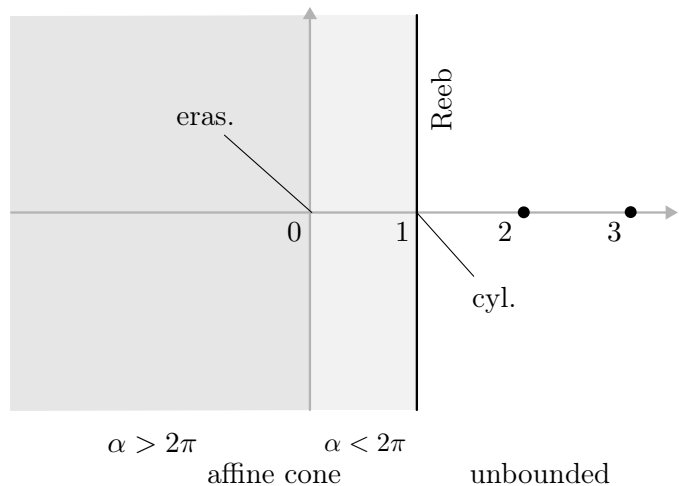

In each of \Cref{fig:cone-affine,fig:cylindre,fig:reeb,fig:unbounded,fig:rec}, we show on the left an example of local flat model for each type above, together with a geodesic (except for the last image), and on the right their image in an explicit Riemann chart.
In \Cref{fig:res}, we show a chart telling, for a given value of the residue $\res$, what kind of Fuchsian singularity we get according to the models above.

A neighborhood of a singularity that is isomorphic to one of the models above is called a \emph{standard neighborhood}.

\begin{lem}\label{lem:uniqGeodConical}
  A conical singularity $p$ has a neighborhood $V$ for which $\forall x\in V-\{p\}$, there is exactly one geodesic from $x$ to $p$ (up to reparametrization) and such that every other geodesic from $x$ eventually escapes $V$.
\end{lem}
\begin{proof}
It is enough to take a small open subset $V$ of the skew cylinder that is starlike with respect to $0$.
Then every point can be linked to $0$ by a ray.
In the skew cylinder model, any other geodesic from a point  will miss $0$ and exit $V$: this is obvious if the opening angle is $\geq 2\pi$, otherwise, it becomes obvious if one develops the fundamental sector of the skew cone, i.e.\ works in the universal covering $\cet$.
\end{proof}

We prove below that non conical Fuchsian singularities have \emph{trapping neighborhoods}.
We call \emph{trapping region} or just trap an open set $V\subset X^*$ such that every geodesic entering $V$ through its boundary never exits $V$.
Note that a union of traps is a trap and that the image of an affine immersion of a half-plane is a trap.
\begin{lem}\label{lem:traps}
  Let $p$ be a non conical Fuchsian singularity of a meromorphic connection.
  There exists trapping neighborhoods $V$ of $p$ bounded by a $C^1$ curve, and other ones bounded by a piecewise geodesic simple loop.
  If $p$ is not of Reeb type, then any geodesic entering $V$ must moreover tend to $p$.
  Otherwise any geodesic accumulating $p$ must tend to $p$.
\end{lem}
\begin{proof}
In the case of a cylindrical or Reeb type singularity, the models of neighborhoods that we described above are bounded by closed geodesics (in particular $C^1$).

In the case of a pure anti-conical singularity, a neighborhood model is a neighborhood of $\infty$ in an unbounded fundamental sector in $\cet$ for the action of a linear map $\lambda$ of non-zero argument, and we can take the following subsets: for a $C^1$ boundary take a logarithmic spiral (image by $\exp$ of a straight line) linking a pair of points on the two boundaries of the (generalized) sector; for a piecewise geodesic boundary, replace this spiral by a polygonal approximation whose ``internal'' angles are $>\pi$.
Because, in affine coordinates, the boundary of $V$ always curves/turns in the same direction, it follows that $V$ is a union of embedded or immersed (depending on whether the opening angle $\alpha$ satisfies $\alpha \geq \pi$ or $\alpha<\pi$) half-planes,\footnote{In the case $\alpha<\pi$, \Cref{fig:skew-cone} shows a relevant part of the universal covering of the skew cone.} hence a trap.
In this model, an entering geodesics develops in $\cet$ as a straight line which will only visit finitely many copies of the fundamental sector $S$, hence tends to $\infty$ if taken back to $S$ by the action of $\lambda$.

The last case is a shifted anti-conical singularity, i.e.\ with total angle $2\pi n$ and an affine holonomy which is a non-zero translation $z\mapsto z+\alpha$. The treatment is similar, with the logarithmic spiral replaced by a curve of the form $z\mapsto R e^{i t 2\pi n} - t\alpha$ for $R$ big enough, and again appropriate polygonal approximations, or simpler right angled shapes like on \Cref{fig:rec}.
Entering geodesics also tend to $\infty$.

Finally, the existence of arbitrarily small traps imply that goedesics accumulating $p$ actually converge.
\end{proof}

Note that all models above are isomorphic, as Riemann surfaces, to the punctured disk $\D^*=\D-\{0\}$ (see the next sentence) and that the corresponding meromorphic connection has a simple pole at $0$, of the indicated residue.
Indeed, in each case we provided in the figures an example of conformal isomorphism either directly to such a disk or to a punctured simply connected strict subset of $\C$; these were given in particular cases but easily generalize; for instance for the anti-conical singularity of residue $1+n$ and non-trivial affine holonomy, one can take an inverse branch of the map $z\mapsto z^{-n}+\frac{1}{2\pi i}\log z$ from a punctured neighborhood of $0$ to the flat model near $\infty$.

\begin{lem}\label{lem:lem}
An injective holomorphic map $f$ from $\D^*$ to a compact Riemann surface $X$ necessarily extends holomorphically at $0$.
\end{lem}
\begin{proof}
It is a classical result but it is hard to find it stated this way in the mathematical literature, so we provide a proof.
Let $z_0=1/2$ and $B=B(z_0,1/4)\subset \D^*$.
Since $f$ is assumed injective, the restriction $f_1$ of $f$ to $U_1:=B(0,1/4)\setminus\{0\}$ takes values in a component $Y$ of the open subset $X\setminus f(\ov{B})$ of $X$. The (non-compact) Riemann surface $Y$ is hyperbolic.
The map $f_1: U_1\to Y$ is non-expanding for the respective hyperbolic metrics.
A circular loop following the circle of center $0$ and radius $\eps$ has a hyperbolic length in the punctured disk $U_1$ that tends to $0$ as $\eps\to 0$.
For $\eps$ small enough it is contractible in $Y$.
It follows that there is a lift $g_1:U_1\to \D$ of $f_1$ under a universal covering $\D\to Y$.
Since $g_1$ is bounded, its singularity at $0$ is erasable.
\end{proof}

\begin{lem}\label{lem:embedFuchsian}
If one of the model affine surfaces above (without the singularity) is affinely embedded in a finite-type Riemann surface, then the embedding holomorphically extends by sending the singularity of the model to a singularity of $X$ of the same type.
\end{lem}
\begin{proof}
Composing the embedding with an isomorphism from the model to $\D^*$, we get a holomorphic injection from $\D^*$ to $X$, which extends at $0$ by \Cref{lem:lem}.
The extension sends $0$ to some point $y\in X$, $\D$ gives a chart in which the Christoffel symbol has a simple pole, the indicated residue and the model tell in which local affine isomorphism class we are.
\end{proof}

%\FloatBarrier

\subsection{Formal invariants of poles}\label{sub:formal}

Consider a meromorphic connection given in a chart by its coefficient $\Gamma$, a meromorphic function of the complex variable $z$, and assume that $\Gamma$ is defined near $0$ (with possibly a pole there).
Because the singularity of $\Gamma$ at $0$ is erasable or polar, but not essential, a change of variable through formal power series $\Phi \in \C[[X]]$, $\Phi = a_1 X + \sum\limits_{n\geq 2} a_n X^n$ with $a_1\neq 0$ makes sense, through the formula\footnote{There is a slight abuse of notation denoting $\Phi^* \Gamma$. We are not pulling back the function $\Gamma$ but the associated connection.}
\[ \Phi^* \Gamma := \Phi' \times \Gamma \circ \Phi + \frac{\Phi''}{\Phi'}
\]
which, in the case of analytic $\Phi$, is formula \eqref{eq:g} for analytic changes of variables.
We recall that for a formal $\Phi$ and a meromorphic $\Gamma(z) = \sum\limits_{n=-d}^{+\infty} \gamma_n z^n$ the function $\Gamma\circ \Phi$ is well-defined as limit of the series $\sum\limits_{n=-d}^{+\infty} \gamma_n \Phi^n$, which converges in the formal sense.

Actually we can also start with $\Gamma$ only formal-meromorphic, i.e.\ a formal power series with finitely many terms of negative exponents. Recall that the set of such power series is denoted as $\C((X))$ and is the fraction field associated to $\C[[X]]$ and called the field of formal Laurent series.
We say that $\Phi^*\Gamma$ and $\Gamma$ are \emph{formally equivalent}.

It is a simple computation to check that the pole of $\Phi^*\Gamma$ has the same order than the pole of $\Gamma$, call it $d\geq 0$, and that in the case $d>0$, its leading coefficient (in $X^{-d}$) is $\frac{\gamma_{-d}}{a_1^{d-1}}$.

The residue is also preserved. This can for instance be proven by noting that it only depends on finitely many coefficients of $\Phi$ and $\Gamma$, and if one replaces $\Phi$ and $\Gamma$ by a convergent series (for instance a polynomial) with the same initial coefficients, we already know that the residue is preserved (see \Cref{sub:minicours}).

Similarly the fact that, for two formal power series $\Phi_1$ and $\Phi_2$ of the form above, we have
\begin{equation}\label{eq:formalcomp}
\Phi_2^*(\Phi_1^* \Gamma) = (\Phi_1\circ\Phi_2)^* \Gamma
\end{equation}
can be deduced from the analytic case by approximations as above, or by direct computation.

The order and residue are the only formal invariants for poles of order at least $2$:
\begin{lem}
Consider two meromorphic connections given in charts by function $\Gamma_1$ and $\Gamma_2$, both defined near $0$, and assume that at both have the same order $d\geq 2$ and the same residue.
Then they are formally equivalent.
\end{lem}
\begin{proof}
It is enough to prove that if $\Gamma$ is polar at $0$ then it is conjugate to $X^{-d}+\on{res} X^{-1}$ for $d>1$ or to $\on{res} X^{-1}$ for $d=1$.
We use the classical approach: cancel-out coefficients one by one by using a sequence of changes of variables of the simple form $\Phi = X+b_n X^n$: their composition formally converges.

First, note that if $d>1$, if we use $\Phi = b X$ then the coefficient of $X^{-d}$ of $\Gamma$ becomes $\gamma_{-d}/b^{d-1}$ as already noted. Since $d-1\neq 0$ there is at least one $b\in\C^*$ such that this new coefficient is $1$.

Then assume there is some $\Phi$ such that $\Phi^* \Gamma = X^{-d}+\on{res} X^{-1}+\cal O(X^m)$ for $d>2$ or to $\on{res} X^{-1}+\cal O(X^m)$ for $d=1$, where $m\in\Z$ and $m\geq 1-d$.
If $m=-1$ we skip this step (we cannot change the residue).
Otherwise, consider a further change of variable $\Psi = X + b X^{1+d+m}$.
Then one computes that $\Psi^*(\Phi^*\Gamma) = \Phi^*\Gamma + (m+1)b X^m + \cal O(X^{m+1})$.
So one can choose $b$ so as to cancel-out the coefficient in $X^m$ of $\Phi^*\Gamma$.
\end{proof}

Consider now two poles of two connections on two affine surfaces.
Consider a Riemann chart for each pole, sending it to $0$.
The respective expressions of the meromorphic connections are given by meromorphic functions $\Gamma_1$ and $\Gamma_2$ as above.
The existence of a formal equivalence between them is independent of the choice of chart (use \cref{eq:formalcomp} with the change of charts).
Regular points of meromorphic connections are always analytically trivializable (there is a chart such that $\Gamma =0$).
For the corollary below, we define their order to be $0$ and residue to be $0$.

\begin{cor}
Consider two points in two Riemann surfaces endowed with two meromorphic connections. 
The connections are formally isomorphic at these points if and only if they have the same order and residue.
\end{cor}

\begin{rem}
It would be interesting to recover the asymptotic values invariant in terms of the Borel resummation of the formal series.
\end{rem}

\subsection{The exponential-affine plane, its geodesics and its quotients}\label{sub:ExpoAffine}

We introduce a natural example of a meromorphic connection with a double pole. It arises from a log-Riemann surface in the sense of \Cref{subsub:LogRiemann}. The geodesic flow on this affine surface, as well as on its quotients, was recently studied in \cite{BT25}.

\begin{defn}
The \textit{exponential-affine plane} $\mathcal{E}$ is the complex plane $\mathbb{C}$ endowed with an affine structure whose affine charts are the restrictions of the map $z \mapsto e^{z}$.
\par\noindent% noindent parce que sinon, visuellement, on croit qu'on est sorti de la définition.
The corresponding connection is the pull-back by the exponential map of the canonical connection of $\C$.
We call $\C$ the \emph{global log-chart} of $\cal E$.
\end{defn}

The exponential map realizes $\mathcal{E}$ as a universal covering of the punctured plane $\mathbb{C}^{\ast}$ with the canonical affine structure.
On the global log-chart $\C$ of $\cal E$ one gets a constant Christoffel symbol $\Gamma=1$.
Seeing $\cal E$ as a subset of $\RS$, we get a meromorphic connection on $\RS$ with a single singularity, at $\infty$.
The change of variable $v=1/z$ sends it to $0$ and in this coordinate the Christoffel symbol has expression
\[\Gamma_v = -\frac{1}{v^2}-\frac{2}{v}\]
so the sinuglarity is a pole of order $2$ and residue $2$ (recall that we work with the opposite of $\Gamma$ to define the residue).
We can hence view $\cal E$ as a finite type affine surface $(X,\nabla,\{\infty\})$ with $X=\RS$ and $\nabla$ with symbol $\Gamma=1$ in the canonical chart $\C$ of $\RS$.

\begin{lem}\label{lem:ExpoAffineUniqueness}
Up to (affine) isomorphism, for a meromorphic connection $\nabla$ on $\RS$ that has a double pole at infinity and no other pole, then $(\RS,\nabla)$ is isomorphic to $\cal E$.
\end{lem}

\begin{proof} In the variable $v=1/z$ the Christoffel symbol $\Gamma_v$ has expansion at $0$:
$\Gamma_v = \frac{a_{-2}}{v^2} + \cal O(\frac{1}{v})$
for some $a_{-2}\in\C^*$.
By \cref{eq:g} we have $\Gamma_z(z) = \frac{-1}{z^2} \Gamma_v(v) -\frac{2}{z}$ so $\Gamma_z(z) = a_{-2}+\cal O(\frac{1}{z})$ as $z\to \infty$.
As a holomorphic function that has a limit at $\infty$, $\Gamma_z$ is constant.
A linear change of variable $w =\lambda z$
yields $\Gamma_w(w) = \lambda^{-1} \Gamma_z(z)$ so by taking $\lambda=a_2$ we get $\Gamma_w =1 $ on $\C$.
\end{proof}

\begin{figure}
\begin{tikzpicture}
\node at (4,3.5) {\includegraphics[scale=0.43]{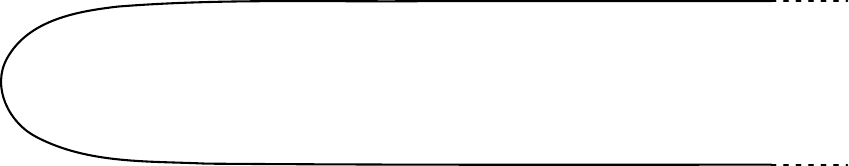}};
\node at (0,0) {\includegraphics[scale=0.43]{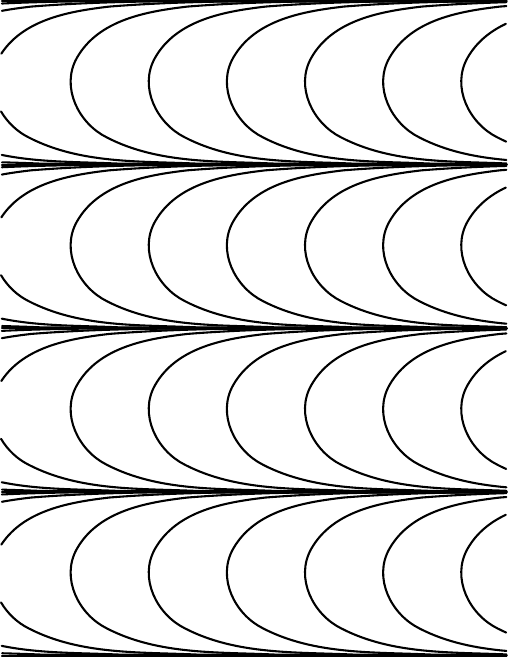}};
\node at (4,0) {\includegraphics[scale=0.43]{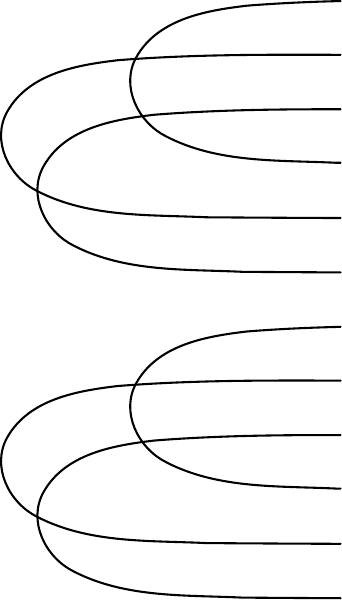}};
\node at (8,0) {\includegraphics[scale=0.43]{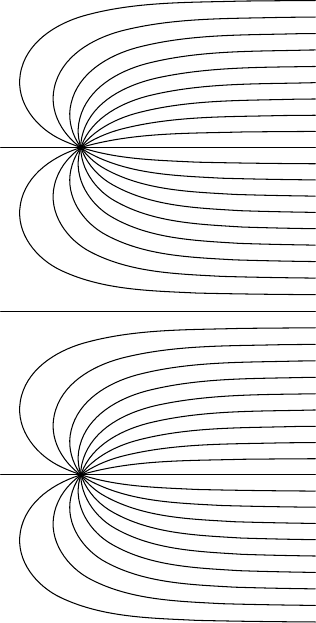}};
\end{tikzpicture}
\caption{Every geodesic of $\cal E$ either traces a translate of the curve on top or a horizontal line. Below, from left to right we have two periods (a strip of height $4\pi$) of the pull-back of: a set of parallel lines, some equilateral triangle whose sides are extended to infinity, a set of lines through a point not at the origin.}
\label{fig:cloggeod}
\end{figure}

\loctitle{Geodesics} It is interesting to figure out the geodesics in the global log-chart $\C$ of $\cal E$.
They are the pull-backs by $\exp$ of straight lines $\ell$ of $\C$.
So they are either horizontal lines if $\ell$ goes through $0$ or translates in $\C$ of the curve of equation $x=\frac{1}{2}\log(1+(\tan y)^2)$.
See \Cref{fig:cloggeod} for an illustration.
The map $y\mapsto x$ along the curve is convex.
In particular, right half-planes are \emph{trapping regions}: a geodesic entering a right half-plane never exits (this also is easily seen by passing to the exponential: a geodesic becomes a straight line and a right half-plane the complement of a disk).
The same holds for lower and upper half-planes, or half-planes whose inner normal has argument in $[-\pi/2,\pi/2]$.
These remarks will be useful in \Cref{sub:model}.

\loctitle{Special subsets} 

The following subsets of $\cal E$ will play a role in several parts of the article: the upper, lower, left and right half-planes in the global log-chart $\C$ of $\cal E$.
Let us look at their image by the exponential.
Consider first a moving horizontal line, scanning $\cal E$ from bottom to top.
Its image by $\exp$ is a radial half-line with vertex $0$ that scans $\C^*$ infinitely many times, turning around $0$ as its argument varies from $-\infty$ to $+\infty$.
Restricting $\exp$ to an upper half-plane, this half-line has argument varying from some $\theta_0\in\R$ to $+\infty$.
Restricting to a lower half-plane, the argument varies from $-\infty$ to some $\theta_0\in\R$.
A left half-plane is sent by $\exp$, as a universal covering, to the punctured disk $B(0,r)-\{0\}$ for some $r>0$.
A right half-plane is sent as a universal covering to the complement of a closed disk $\C \setminus \ov B(0,r)$.

\loctitle{Automorphisms, quotients} The group $\Aut(\mathcal{E})$ of affine automorphisms of $\mathcal{E}$ consists of the maps $z\mapsto z+\mu$ with $\mu \in \C$. 
In the universal covering of $\C^*$, they correspond to the similarity of center $0$, ratio $e^{\Re\mu}$ and angle $\Im\mu$.
Then, affine surfaces covered by $\mathcal{E}$ are quotients of $\mathcal{E}$ by discrete subgroups of the additive group of $(\mathbb{C},+)$. Depending on whether the group is of rank $0$, $1$ or $2$, we obtain $\cal E$, an affine surface of genus zero with two singularities (an affine cylinder) or an affine surface of genus one without any singularity (an affine torus).

For the next statement, recall that we consider marked points as Fuchsian singularities of residue $0$.
\begin{lem}\label{lem:AffineCylinder}
For any $\alpha \in \mathbb{C}$, the unique meromorphic connection on $\RS$ with a Fuchsian singularity of residue $1+\alpha$ at $0$, a Fuchsian singularity at $\infty$, and no other singularity, has expression
\[\Gamma_z = -\frac{1+\alpha}{z}\]
in the canonical chart $\C$ of $\RS$.
As a consequence, the unique meromorphic connection on $\RS$, up to isomorphisms, with two Fuchsian singularities of residues $1+\alpha$ and $1-\alpha$ is isomorphic on the complement of its singularities to :
\begin{itemize}
\item the quotient of $\mathcal{E}$ by $\langle z \mapsto z+ 2\rm{i}\pi \alpha\rangle$ if $\alpha\neq 0$,
\item the quotient of the flat plane by a non-zero translation\footnote{All such quotients are isomorphic.} if $\alpha=0$.
\end{itemize}
Note that if $\alpha\in\{-1,1\}$ the quotient in the first case is isomorphic to the subset $\C^*$ of the flat plane $\C$.
\end{lem}

\begin{proof}
We can make a first change of variable on $\RS$ by a homography so as to put the singularities at $0$ and $\infty$.
Again, we use the charts $z$ and $v=1/z$.
At $0$, we have $\Gamma_v = \frac{\alpha-1}{v}+\cal O(1)$.
Hence at $\infty$, $\Gamma_z = (1-\alpha)/z + \cal O(1/z^2) - \frac{2}{z} = -(1+\alpha)/z + \cal O(1/z^2)$.
The function $z\Gamma_z$ extends holomorphically at $0$ (where the pole of $\Gamma$ is simple) and has a limit at $\infty$ so is constant.
Conversely the Riemann sphere endowed with the meromorphic connection of symbol $\Gamma_z = -(1+\alpha)/z$ satisfies the asumptions.
Taking the preimage by $\exp$, i.e.\ letting $z=\exp(u)$,  logarithm, we get $\Gamma_u = -\alpha$.
If $\alpha=0$ we get $\Gamma_u=0$ so we are in the flat plane, but recall that $\log$ is only defined modulo translation by $2\pi i$.
If $\alpha\neq 0$, a further change of variable $w= -\alpha u$ gives $\Gamma_w = 1$, so we end up in $\cal E$ and the deck transformation generator $u\mapsto u+2\pi i$ becomes $w\mapsto w-2\pi i\alpha$.
\end{proof}

\subsection{Affine surfaces of low complexity}

Finite type affine surfaces with a very small genus and number of singularities can be completely classified.

\begin{defn}
Let $(X,\nabla,\cal S)$ be a finite type affine surface.
A \emph{saddle connection} is a geodesic that tends in finite time in the future and in the past to two (possibly identical) singularities.
\end{defn}

We recall that $\cal E$ denotes the exponential-affine plane of \Cref{sub:ExpoAffine}, and that its automorphisms have been described there.
Finally we recall that marked points are counted as multiplicity one singularities. Below we omit the mention of \emph{finite type} in front of \emph{affine surface}.
\begin{prop}\label{prop:LowComplexity}
If an affine surface $(X,\nabla,\mathcal{S})$ of genus $g$ with $n$ singular points (counted with multiplicity) satisfies $2g+n \leq 2$, then it is isomorphic to one of the following models:
\begin{itemize}
    \item[(a)] \textit{whole affine plane}: $X^{\ast}=\mathbb{C}$ and $\nabla$ has a unique Fuchsian singularity of residue equal to $2$ (at infinity);
    \item[(b)] \textit{translation cylinder}: $X^{\ast} = \mathbb{C}/\mathbb{Z}$ and $\nabla$ has two Fuchsian singularities with residues equal to $1$ (the two ends of the cylinder);
    \item[(c)] \textit{translation torus}: $X=\mathbb{C}/\Lambda$ where $\Lambda$ is a lattice of $(\mathbb{C},+)$ and $\nabla$ has no singularity.
    \item[(d)] \textit{infinite angle cone}: $X^{\ast}$ is affine-isomorphic to the exponential-affine plane $\cal E$ and the only singularity of $\nabla$ is a double pole;
    \item[(e)] \textit{affine cylinder}: $X^{\ast}$ is a quotient of $\mathcal{E}$ by a linear transformation; 
    \item[(f)] \textit{affine torus}: $X$ is a quotient of $\mathcal{E}$ by a multiplicative lattice of linear transformations.
\end{itemize}
In particular, these affine surfaces have no saddle connections.
\end{prop}

\begin{proof}
The relation $2g+n\leq 2$ has only four solutions in $\N\times \N$: $(g,n)=(0,0)$, $(0,1)$, $(0,2)$ and $(1,0)$.
We cannot have $(g,n)=(0,0)$ because a meromorphic connection on $\RS$ has a sum of residues equal to $2$ by \Cref{prop:sumRes}.

We split the case $(g,n)=(0,2)$ into three subcases to get the following five cases:
\begin{itemize}
    \item[(1)] $g=0$, $n=2$ and the meromorphic connection admits a unique double pole;
    \item[(2)] $g=0$, $n=2$ and the meromorphic connection admits two Fuchsian singularities, whose residues are equal to $1+\alpha$ and $1-\alpha$ for some $\alpha \in \mathbb{C}\setminus \lbrace{ 1 \rbrace}$; 
    \item[(3)] $g=0$, $n=2$ and the meromorphic connection admits a unique Fuchsian singularity whose residue is equal to $2$ and a marked point;
    \item[(4)] $g=0$, $n=1$ and the meromorphic connection admits a unique Fuchsian singularity whose residue is equal to $2$;
    \item[(5)] $g=1$ and $n=0$, i.e.\ $X$ is homeomorphic to a torus and there is no singularity.
\end{itemize}
By the Poincaré-Koebe theorem, a compact Riemann surface of genus $0$ is necessarily isomorphic to $\RS$.
\Cref{lem:ExpoAffineUniqueness} characterizes case (1) as the exponential-affine plane (d).
Case (2) is dealt with by \Cref{lem:AffineCylinder}: it corresponds to a translation cylinder (b) if $\alpha=0$, and if $\alpha \notin \lbrace{0,1\rbrace}$, $X^*$ is an affine cylinder (e).
In cases (3) and (4), the affine surface is the infinite plane with or without a marked point by \Cref{lem:AffineCylinder}. The unmarked plane is (a) while the marked case is a special case of affine cylinder (e).

Finally, in the case (5) of an affine surface of genus $g=1$, the underlying Riemann surface is isomorphic to a complex torus, which is (again by a classical consequence of the Poincaré-Koebe theorem) biholomorphic to $\C/\Lambda$ where $\Lambda$ is a co-compact discrete translation subgroup, i.e.\ a group of translations generated by $T_a$ and $T_b$ where $(a,b)$ is an $\R$-basis of $\C$.
This isomorphism sends $\nabla$ to a connection on $\C/\Lambda$, which we can pull-back on $\C$ by the quotient map $\C/\to\Lambda$ to a $\Lambda$-invariant connection.
The Christoffel symbol $\Gamma$ is then a holomorphic function (it has no poles by hypothesis) which is invariant by $\Lambda$ (given $\nabla$ on $\C$ with symbol $\Gamma$, the symbol of $T_a^*\nabla$ is $\Gamma\circ T_a$ by \cref{eq:g}) hence constant.
There are then two cases: either $\Gamma=0$, in which case we have a translation torus (c).
Or $\Gamma$ is a non-zero constant, which we can take equal to $1$ by a further linear change of variable.
Then we are in case (f).

In all of these surfaces, $X^{\ast}$ is a quotient of the canoncial affine plane $\mathbb{C}$ or the exponential-affine plane $\mathcal{E}$.
If there were a saddle connection in a quotient of one these spaces, this would lift to a geodesic of $\C$ or $\cal E$ whose lifespan is bounded both in the past and the future.
We observe that in $\cal E$ and in $\C$, there is no such geodesic.
\end{proof}

\section{Classification of irregular singularities}\label{sec:AsymptoticValues}

In this section we say that a connected affine surface $\cal A$ has a \emph{puncture} if there exists a Riemann surface $\cal R$ \footnote{Riemann surfaces are connected by definition.} containing it, and having only one more point $p$.
A neighborhood of $p$ intersected with $\cal A$ is called a neighborhood of the puncture.
If the corresponding connection is meromorphic at $p$ with a pole of order $d$, we say that we have a \emph{meromorphic puncture of order $d$}.

\subsection{Construction of the invariant}\label{sub:constr:inv}

Consider a pole $p$ of order $d\geq 2$ of a meromorphic connection defined on an open subset of a Riemann surface $X$.
Consider a Riemann chart mapping the pole to $0\in\C$.
Let $\Gamma$ be the Christoffel symbol of the connection in this chart.
By hypothesis has a power series expansion of the form
\[\Gamma(z) = \frac{a_{-d}}{z^d} + \cdots + \frac{a_{-1}}{z} + \cdots
\]
and its residue is
\[\res = -a_{-1}\]
(see \Cref{sub:minicours}).

\begin{figure}[htbp]
\includegraphics[width=7cm]{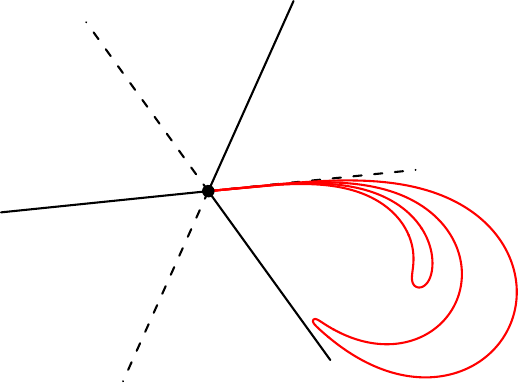}
\caption{Attracting axes as dashed lines, repelling axes as solid lines, and two typical geodesics in red, for a pole of order $d=4$, seen in a Riemann chart.}
\label{fig:axes}
\end{figure}

\begin{defn}\label{defn:axes}
We define \emph{repelling} and \emph{attracting} axes as the half-lines $[0,+\infty)e^{i\theta}$ where $\theta$ is such that $-a_{-d}e^{-i(d-1)\theta} \in\R_-$ or $\R_+$, respectively.
\end{defn}
A justification of these terms is given later in this article, related to the behavior of geodesics. 
There are $d-1$ repelling axes, alternating with the same number of attracting axes, see \Cref{fig:axes}.

\begin{rem}
Given two Riemann charts sending a given irregular pole to $0$,
the differential at $0$ of the change of chart sends repelling/attracting directions to repelling/attracting directions (see \Cref{sub:minicours}), so these directions are well defined as elements of the tangent space $T_p X$.
\end{rem}

\loctitle{Universal covering of a neighborhood}
Let $B(0,\eps)\subset \C$ be contained in the image of the chart and denote $\tilde B$ a universal covering of $B^*:=B-\{0\}$.
The argument function has a well defined $\R$-valued lift on $\tilde B$.
Endow $\tilde B$ with the lift of the connection/affine structure.\footnote{Given an affine surface $\cal A$, and a universal covering $\tilde {\cal A} \overset{\pi}\longrightarrow {\cal A}$, one can endow $\tilde {\cal A}$ with an affine atlas such that $\pi$ is affine. This can be done for instance by taking as charts injective restrictions of $u\circ \pi$ where $u$ are affine charts of $\cal A$.
Such an affine structure is unique, i.e.\ any two atlas on $\tilde {\cal A}$ such that $\pi$ is affine are compatible.
This is equivalent to taking the pull-back of the connection by the map $\pi$.}
The universal covering $\tilde B \to B^*$ has a deck transformation group isomorphic to $\Z$ and generated by the map $T:\tilde B\to \tilde B$ that sends a point of $\tilde B$ with lifted argument $\theta\in\R$ to the point with the same projection on $B^*$ but with lifted argument $\theta+2\pi$.

\loctitle{Developing map}
Given a point $x_0\in B^*$, a lift $\tilde x_0 \in \pi^{-1}(\{x_0\})$ and a germ\footnote{A \emph{germ} at $x$ is here understood as a class of maps defined in neighborhoods of $x$, where two maps are equivalent if they coincide on a possibly smaller neighborhood of $x$.} of affine chart $\phi$ at $x_0$, then the germ $\tilde \phi = \phi\circ \pi$ at $\tilde x_0$ has a unique continuation to a map $\tilde\phi:\tilde B \to \C$ that is affine.
This map is called the \emph{developing map} of $\phi$ (w.r.t.\ $\tilde x_0$). 
It can be constructed by general methods (see \cite{Thurston}, Section~3.4) but also, since it satisfies\footnote{For a map $f:\tilde B\to \C$ and any $\tilde x\in \tilde B$, the quantity $f'(\tilde x)$ is well-defined, and expresses as $g'(\pi(x))$ for $g=f\circ \pi^{-1}$ and a local inverse of $\pi$.} the O.D.E.\ $\tilde\phi''/\tilde\phi' = \tilde\Gamma$ where $\tilde\Gamma :=\Gamma\circ \pi$, it expresses as $\tilde\phi = \int \exp(\int \tilde\Gamma)$ for well-chosen additive constants in the antiderivatives $\int$.
From this expression one deduces that there exists a complex affine map
\[L:z\mapsto az+b,\]
where $b\in \C$ and
\[ a = \exp(-2\pi i \res)\in\C^*,\]
called the \emph{holonomy factor} of $\phi$ and such that
\begin{equation}\label{eq:T}
\tilde \phi \circ T = L\circ \tilde \phi
\end{equation}
holds on $\tilde B$.
It is independent\footnote{This is essentially a consequence of the commutativity of the fundamental group of $B^*$.} of the choice of the lift $\tilde x_0$ of $x_0$.
Note that we can skip the choice of an initial $x_0$, $\wt x_0$ and $\phi$ and directly start from any affine map $\wt \phi : \wt B\to \C$, which we still call a \emph{developing map} since any such map is the developing map of any germ $\wt\phi\circ\pi^{-1}$ where $\pi^{-1}$ is any inverse branch germ of $\pi$.

\begin{prop}[Asymptotic value family]\label{prop:di}
Choose any repelling direction in the universal covering $\tilde B$ and denote $\theta_0\in\R$ its lifted argument.
For every direction of argument $\theta_n = \theta_0+n\frac{2\pi}{d-1}\in\R$ in $\tilde B$, with $n\in\Z$, any path tending to $0$ in $\wt B$ with lifted argument tending to $\theta_n$ has an image by $\tilde{\phi}$ that converges to some complex number $u_n$.
\end{prop}

\begin{rem} The hypothesis on the path can be weakened, see \Cref{prop:avmg}.
\end{rem}

\Cref{prop:di} will be proved in \Cref{sub:pf:prop:di}.
The direction of lifted argument $\theta_0$ is called here a \emph{lifted reference direction} and its projection to a direction of $\C$ a \emph{reference direction}.
The $\Z$-indexed sequence $(u_n)$ is called the \emph{asymptotic value family} (see \Cref{ss:term:avf} for a comment on this choice of terminology) associated to the developing map $\tilde\phi$ and lifted reference direction.
The asymptotic value family depends on the developing map $\tilde \phi$, on the lifted reference direction and on the choice of the lift of the argument function $z\mapsto \arg z$.

\begin{rem}
The denomination ``asymptotic value family'' is chosen in reference to the notion of \emph{asymptotic value} in holomorphic dynamics, which we recall in \Cref{ss:term:avf}.
\end{rem}

\loctitle{Example}
Consider the exponential-affine plane $\cal E$ defined in \Cref{sub:ExpoAffine} and let $m\in\N^*$.
In its canonical chart $\C$ its Christoffel symbol is $\Gamma_z=1$.
A developing map is $z\mapsto \exp z$.
The formula $z=w^m$ defines an $m$-fold covering from $\C^*$ to itself.
This defines a pulled-back affine surface (equivalently, a pulled-back connection) on $\C^*$, which is a meromorphic connection $\nabla$ on $\RS$, of symbol one computes using \cref{eq:g}:
\[\Gamma_w = m w^{m-1} + \frac{m-1}{w}\]
so if $m\geq 2$ it has a simple pole of residue $1-m$ at $0$, which corresponds to the fact that $0$ a conical point of angle $2\pi m$ and without dilation.
In the coordinate $u=1/w$ one finds
\[\Gamma_u = - \frac{m}{u^{m+1}} - \frac{m+1}{u},\]
so $\nabla$ has at $\infty$ a multiple pole of degree $m+1$ and residue $m+1$.
Its asymptotic value family in this case is just the constant family $u_n=0$, $n\in\Z$, since its repelling axes directions $\theta$ in $u$-coordinates are given according to \Cref{defn:axes} by $me^{-im\theta} \in\R_-$ and a developing map is $u\mapsto \exp(u^{-m})$.

\medskip

The holonomy equation \eqref{eq:T} implies
\begin{equation}\label{eq:L}
\forall k\in\Z,\ u_{n+(d-1)} = L(u_n) = a u_n + b.
\end{equation}
\begin{rem*}
From $m$, $\res$ and the family $(u_n)$, one can of course deduce the value of $a=\exp(-2\pi i\res)$ but also of $b = u_{d-1}-a u_0$. On the other hand, one cannot deduce solely from the family $(u_n)$ which determination of $-\log(a)/2\pi i$ is equal to $\res$.
\end{rem*}

Given a developing map $\tilde \phi_1$, the other developing maps obtained for various germs are the maps of the form $\tilde{\phi}_2 = a'\tilde{\phi}_1+b'$ for some $a'\in\C^*$ and $b'\in \C$.
Their associated sequence (for the same choice of $\theta_0$ in \Cref{prop:di}) are given by 
\begin{equation}\label{eq:dd}
  \forall k\in\Z,\ u_n[\tilde{\phi}_2] = a'u_n[\tilde{\phi}_1]+b' .
\end{equation}

Changing the choice of $\theta_0$ for another reference direction $\theta'_0 = \theta_0 + s \frac{2\pi}{d-1}$ with $s\in\Z$ just shifts the sequence $u_n$ by $s$. More precisely the new sequence $u'_n$ satisfies: 
\begin{equation}\label{eq:shift}
u'_n = u_{n-s}.
\end{equation}
In the case $s$ is a multiple of $m$, meaning that we take the same axis but another lift of its direction, we get
\begin{equation}\label{eq:ol}
u'_n = L^{-k}(u_n)
\end{equation}
where $s = mk$.

We repeat now \Cref{def:is}:
\begin{defn*}
\defIsText
\end{defn*}

By \cref{eq:ol}, any two lifts of the same reference direction yield asymptotic sequences $u$ and $u'$ such that $u\sim u'$.

\begin{defn}\label{def:ia}
  Given a pole of order $d\geq 2$ of a meromorphic connection, and a choice of reference direction (i.e.\ of repelling axis), their \emph{asymptotic values invariant} is the element of $\cal I_d$ which is the equivalence class $\Pi(u)$ as in \Cref{def:is}, where  $u$ is the asymptotic value family of any developing map $\tilde\phi$ as described in \Cref{prop:di}, for any choice of lift of the chosen divergent direction. 
\end{defn}

The asymptotic values invariant, more briefly called the invariant, depends on the choice of reference axis but is independent of the choices of: the Riemann chart,\footnote{Developing maps are naturally defined on the surface, independently of the considered Riemann chart.} the lift of the reference direction, and the developing map.

\begin{rem}[Working without Riemann charts]
We used a Riemann chart near the pole but we could work directly on $X$.
The developing maps of a germ of affine chart is a well defined notion on the affine surface $X^*$, independently of Riemann charts.
As we saw earlier the divergent directions are also well defined independently of the Riemann chart, in particular it makes sense to say that a path in $U$ tends to the puncture tangentially to a given divergent direction, without referring to a particular choice of Riemann chart.
\end{rem}
We repeat here \Cref{thm:main}:
\begin{thm*}
\thmMainText{\boolean{false}}
\end{thm*}

Part (\ref{item:2}) of \Cref{thm:main} will be proved as \Cref{thm:model2} in \Cref{sub:model} and Part~(\ref{item:1}) is the object of \Cref{sub:momo}: the direct implication follows from \Cref{thm:model,prop:isocan} while the easier reciprocal was stated above the theorem.

\begin{rem}
By \cref{eq:L} and the note following it,
the map $u\in\cal U_{d,\res} \mapsto (u_0,\ldots,u_d)\in \C^{d+1}$ is a bijection.
The quotient $\cal I_d=\cal U_d/\sim$ can be identified with $\C^{d+1}/\Aut(\mathbb{C})$ where $\Aut(\mathbb{C})$ is the group of affine automorphisms of $\C$ that we let act simultaneously on every component of $\C^{d+1}$.
Another interesting identification is as follows: first send 
$u$ to $(u_0,v)$ with $v_n = (u_{n+1}-u_n)\exp(2\pi i \frac{\res}{d} n)$. Then we have a bijection between $\cal U$
and $\C\times P$ where $P$ is the set of sequences $v\in\C^\Z$ with $v_{n+d}=v_n$, i.e.\ of period dividing $d$.
Then map it to $[v]$ where $[v] = \C^* v$ denotes the class of $v$ under the action of $\C^*$ by scalar multiplication on $\C^d$. Since $P$ is in bijection with $\C^d$ via $v\mapsto (v_0,v_1,\ldots, v_{d-1})$, this gives a bijection between $\cal I_{d,\res}$ and $\C^d/\C^*$.
\end{rem}

We saw in \Cref{eq:shift} that changing the reference divergent axis shifts the index of the sequence $u_n$ by the translation $n\mapsto n+s$: $u' = \sigma^{s}u$ where $\sigma$ is the shift operator that sends $n\mapsto u_n$ to $n\mapsto u_{n-1}$.
This shift action descends to the invariant $\cal I_{d,res}$ in the sense that $u \sim v \iff \sigma u\sim \sigma v$ allows to define $\sigma\Pi(u) = \Pi(\sigma u)$.

\begin{thm}
  Consider two affine surfaces, each with a pole of order $d \geq 2$. Then they are isomorphic near the poles if and only if their invariant are identical up to the action of the shift.
\end{thm}

In \Cref{sub:model,sub:helix} we complement these statements with a model consisting in an alternance of what we call conical sectors (centered on repelling axes) and unbounded sectors (centered on attracting axes), together with a geometric description of the corresponding affine surface.

\subsubsection{Asymptotic values}\label{ss:term:avf}

Let us explain our choice of terminology in the denomination \emph{asymptotic value family} for $(u_n)_{n\in\Z}$.
First, we use the word family instead of sequence because the latter is mostly used for index set $\N$.
Second, given a holomorphic map $S\to S'$ between open Riemann surfaces, in holomorphic dynamical systems and some other fields, one calls \emph{asymptotic value} of $f$ an element $v\in S'$ such that there exists a path $\gamma :[0,1)\to S$ that leaves every compact subset of $S$ and such that $f\circ \gamma(t)\tends v$ as $t$ tends to $1$.
The $u_n$ are asymptotic values of $\tilde\phi$ for which the paths tend to $0$.
Even with this restriction, the map $\tilde{\phi}$ may have another asymptotic value if the affine holonomy map $L(z)=az+b$ satisfies $|a|\neq 1$: the fixed point of $L$, obtained with a path that winds infinitely many times around $0$ while it tends very slowly to $0$.
We give a more complete treatment in \Cref{sub:other}.
In \Cref{sub:helix} we will give another interpretation of the family $(u_n)$ in terms of \emph{foci}.

\subsubsection{Case of double poles}\label{subsub:doublepoles}

The space of invariants in the case $d=2$ is particularly simple. First note that the quotient map $\cal I_1 \to \cal I'_1$ is trivial (a bijection): surjectivity is by definition and injectivity follows from the fact that $\sigma u = L^{-1} \circ u$ in the case $d=2$.
\begin{prop}\label{prop:d1}
For every
$\res\in\C$, the set $\cal I_{2,\res} \equiv \cal I'_{2,\res}$ contains exactly two elements:
\begin{enumerate}
  \item If $\res\in\Z$, i.e.\ $e^{-2\pi i \res} = 1$, there are:
  \begin{enumerate}
    \item the class of $(r,u_n=0)$, 
    \item the class of $(r,u_n=n)$.
  \end{enumerate}
  \item If $\res\notin\Z$, i.e.\ $e^{-2\pi i \res} \neq 1$, there are:
  \begin{enumerate}
    \item the class of $(r,u_n=0)$,
    \item the class of $(r,u_n=e^{-2\pi i \res n})$.
  \end{enumerate}
\end{enumerate}
\end{prop}
\noindent The proof is elementary.  Case (1a) has a affine holonomy $L$ equal to the identity, in case (1b) $L$ is a translation. In case (2a) the asymptotic value family is constant equal to the fixed point of $L$, while in case (2b) it is disjoint from this point.

We say that a multiple pole is \emph{centered} if its asymptotic value family is constant.
The two (a) cases above are centered, and the two (b) cases are not. 

\begin{rem}
An example of affine surface with meromorphic punctures of order $2$ appeared naturally in \cite{C1} as a limit of quasiconformal deformations of the square.
In that article, the punctures have order $d=2$, and residue $0$, so we are in Case~(1) of \Cref{prop:d1}. The affine holonomy is a translation, hence we are in Case (1b).
\end{rem}

\subsection{Proof of \texorpdfstring{\Cref{prop:di}}{Proposition ...}}\label{sub:pf:prop:di}

Here we prove and complement \Cref{prop:di}.

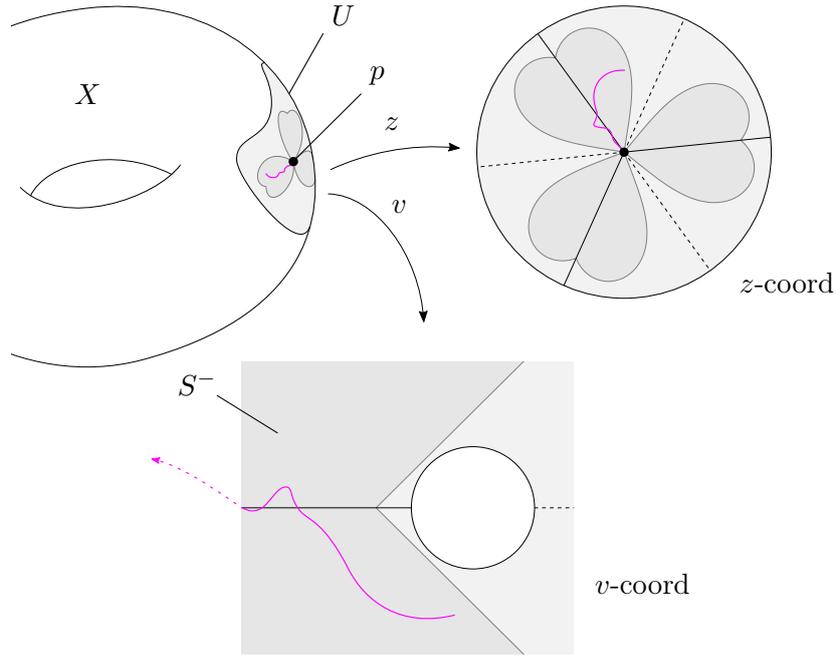
\begin{figure}
\begin{tikzpicture}
\node at (0,0) {\includegraphics[width=10cm]{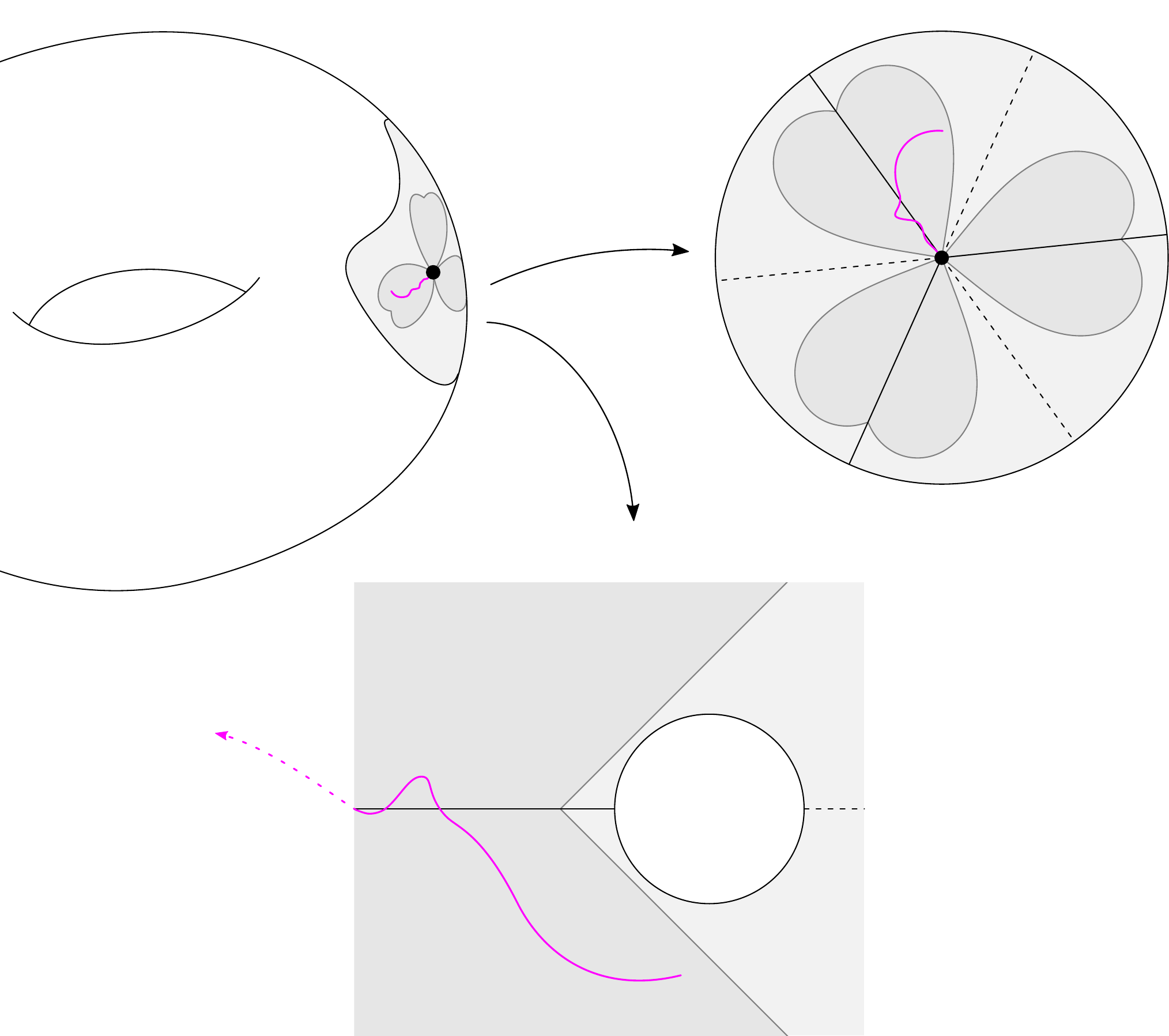}};
\node at (-4,3) {$X$};
\draw[-] (-0.4,3) node[anchor=south west, inner sep=3pt] {$p$} -- (-1.3,2.1);
\draw[-] (-0.9,3.8) node[anchor=south west, inner sep=3pt] {$U$} -- (-1.35,3);
\node at (5.2,0.5) {$z$-coord};
\node at (0,2.6) {$z$};
\node at (0.1,1.5) {$v$};
\draw[-] (-2.3,-1) node[anchor=south east, inner sep=0pt] {$S^-$} -- (-1.5,-1.5);
\node at (3.3,-3.5) {$v$-coord};
\end{tikzpicture}
\caption{Different objects involved in the proof. The magenta curve represents the path $\gamma$ on $X$ and its image in the coordinates $z$ and $v$.}
\label{fig:S1}
\end{figure}

Recall that
\[\Gamma(z) = \frac{a_{-d}}{z^d} + \cdots + \frac{a_{-1}}{z} + \cdots
\]
For convenience we denote
\[m = d-1\]
throughout this section.
It is equal to the number of repelling axes, and to the number of attracting axes, associated to the pole.

In this section we use different letters $z$, $v$, etc.\ for different coordinates and index the Christoffel symbol $\Gamma$ in these coordinates by the variable name, as in $\Gamma_z$, $\Gamma_v$, etc.
Formula~\eqref{eq:g} can be rewritten as, for any pair of coordinates $u$ and $v$ on a common open subset of $X$:
\begin{equation}\label{eq:gc}
\Gamma_v(v) dv - \Gamma_u(u) du = \dlog \frac{du}{dv}
\end{equation}
where $\dlog$ is the operator taking a function $f$ and returning the differential form $df/f$.

We first perform a non-injective change of variable
\[ v=-\frac{a_{-d}}{(d-1) z^{d-1}} = -\frac{a_{-d}}{mz^m}
\]
This change of variable is an $m$-fold covering, ramified at $0$ and maps $0$ to $\infty$.
It is motivated by the fact that since $d v = \frac{a_{-d}}{z^d} d z$ this should simplify the expression of the leading term of $\Gamma$.
Indeed by \cref{eq:gc} (equivalently, by \Cref{eq:g}), $\Gamma_v = \frac{dz}{dv}\Gamma_z + \frac{\dlog}{dv}\frac{d v}{d z} = \frac{z^d}{a_{-d}} \Gamma_z +\frac{d}{z} = 1 + \cal O(1/z)$ as $|z|\to \infty$.
In this new variable,
\[ g:=\Gamma_v - 1 = \cal O(1/|v|^{1/m})
\]
when $|v|\to\infty$. 
Note that $g$ is defined on the $m$-fold covering of a pointed disk  neighborhood $D$ of $\infty$.

Consider the sector
\[S^-= \{z\in \C\,;\,-\pi-3\pi/4 < \arg(z+R) <-\pi+ 3\pi/4\},\]
whose central axis has argument $-\pi$, apex is at $z=-R$, and opening angle is $3\pi/2$, see \Cref{fig:S1}.
We now invoke \Cref{cor:hphi2} of \Cref{sec:tech}, which we apply in the $v$-coordinate with the function $h$ and sector $S^-$ above.
(More precisely, since $z\mapsto v$ is not injective, we restrict to a subset of $X$ that is a sector in the $z$-coordinate and that is in one to one correspondence with $v\in D-[0,+\infty)$.
The central axis of the sector is actually a repelling axis of $p$ and there is one sector for each repelling axis.
On such a set, $h$ and $g$ are well-defined functions of $v$.
The standing assumption of \Cref{sub:tech} is satisfied since $h'(v) = g(v) = \cal O(1/|v|^{1/m})$ as $|v|\to\infty$ with $v\in S^-$.)
\Cref{cor:hphi2} implies that there is a developing map $\phi_0:S^-\to \C$ for the connection defined by $\Gamma_v$ and such that $\phi_0(v)\tends 0$ as $\Re(v)\to -\infty$ with $v\in S^-$.

Now for any path $\tilde \gamma : [0,1)\to \tilde U$ tending to $0$ tangentially to a repelling direction, its $z$-coordinate eventually enters and stays in the corresponding sector above and the $v$-coordinate tends to $\infty$ with an argument tending to $\pi$, hence it is eventually in $S^-$ and its real part tends to $-\infty$.
The function $\tilde\phi$ on the subset of $\tilde U$ corresponding to the sector coincides with $a\phi_0(\tilde v)+b$ for some constants $a\in\C^*$ and $b\in\C$.
It follows that $\tilde\phi\circ\tilde \gamma$ converges to $b$ and this proves \Cref{prop:di}.

\medskip

\begin{figure}
\begin{tikzpicture}
\node at (0,0) {\includegraphics[width=12cm]{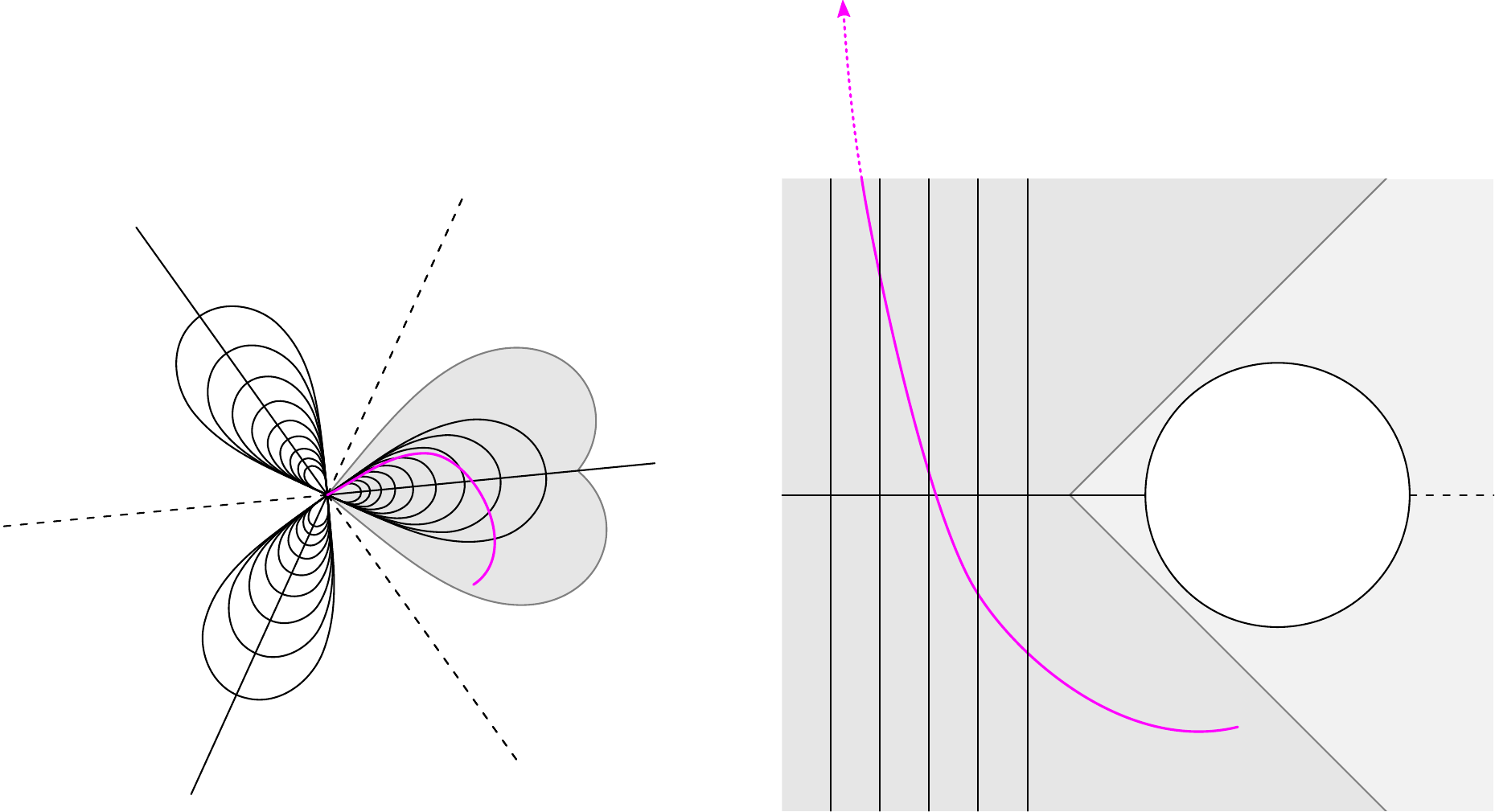}};
\node at (-3,2.7) {$z$-coord};
\node at (3,2.7) {$v$-coord};
\end{tikzpicture}
\caption{By \Cref{prop:avmg}, along a path such that $\Im v \tends -\infty$, the map $\tilde\phi$ will also converge.}
\label{fig:S2}
\end{figure}

Actually it is enough that $\Re(v\circ\gamma)$ tends to $-\infty$, so we get the slightly more general following result, illustrated on \Cref{fig:S2}:
\begin{prop}\label{prop:avmg}
Consider a path $\tilde\gamma:[0,1)\to \tilde U$ tending to the puncture in such a way that, as $t\to 1$, the $z$-coordinate of $\gamma(t)$ satisfies
\[\Re\frac{a_{-d}}{z^{d-1}} \tends +\infty.\]
Then there exists $n$ such that $\limsup|\arg z -\tilde\theta_n|\leq\pi/2d$ and $\tilde \phi$ tends along the path to the value $u_n$ associated to $\tilde\theta_n$ in \Cref{prop:di}.
\end{prop}

\subsection{Local models for irregular singularities}\label{sub:model}

%\startcontents
%\printcontents{}{2}[3]{}
%\medskip

In this section we give an isomorphic model, as an affine surface, of a neighborhood of any pole of order $d \geq 2$ of a meromorphic connection. This model will be used in to prove \Cref{thm:main} but is interesting in itself. For instance it can help to study geodesics or to describe the affine surface puncture neighborhood by cutting and pasting (infinitely many) polygons.

Throughout this section we let
\[m=d-1.\]

\subsubsection{Log-charts}\label{sub:log}

The model will be expressed in terms of gluing together charts of a type that we define and start studying here.

\begin{defn}[Log-charts]\label{def:log}
Given an affine surface we call log-charts the injective functions that are locally of the form $\log\circ \phi$ where $\phi$ is an affine chart.
\end{defn}

Log-charts are exactly the Riemann charts in which the associated Christoffel symbol take the particularly simple expression
\[\Gamma = 1.
\]
They can be considered as taking values in the exponential-affine plane $\cal E$ defined in \Cref{sub:ExpoAffine}.

\medskip

\medskip

\loctitle{Gluing log-charts} A change of variable between two log-charts necessarily locally expresses as $z\mapsto \log(a e^z+b)$, i.e.\ as $z+\log(a) + \log(1+be^{-z}/a)$ for some branches of $\log$.
From such an expression one can define an injective function on a set containing the right half-plane of equation $\Re(z)\geq \log\frac{|b|}{|a|}$ as follows.
%Ce sont également les solutions de $1 = \phi' +\frac{\phi''}{\phi'}$ c.à.d.\ $\phi'' = \phi' - {\phi'}^2$. Mais peu importe
Let $\log_p :\C \setminus (-\infty,0]\to\C$ be the principal branch of the logarithm.
For $s$ and $b$ in $\C$, we denote
\[ G_{s,b}(z) =  z + s + \log_p(1+b e^{-z-s}) ,\]
which is defined on $\C$ if $b=0$ and if $b\neq 0$ on $\C$ minus the union of the translates $T_c(A)$ of the line $A=(-\infty,0]$, where $T_c(z)=z+c$ and $c$ ranges among the different determination of $\log (-b)-s$.
The map is injective: it satisfies $\exp\circ G_{s,b}(z) = e^s \exp(z) +b$ so any two points $z$ and $z'$ with the same image differ by a multiple of $2\pi i$ but since $G_{s,b}(z+2\pi i) = G_{s,b}(z)+2\pi i$ this means they are equal.
The image of $G_{s,b}$ is a subset of $\C$ having the same form, where $c$ now ranges in the determinations of $\log(b)$.
However we will only need the restriction of $G_{s,b}$ to a right half-plane contained in its domain.

The map $G_{s,b}$ is a lift by $\exp$, in the domain and range, of the affine self-map $f : z\mapsto e^s z + b$ of $\C$ restricted to $\C \setminus [0,f^{-1}(0)]$.
The map $f$ send the latter set to $\C \setminus [0,b]$.
Similarly, the inverse of the map $G_{s,b}$ is exactly (domain included) the map $G_{-s,-e^{-s}b}$.
This is noteworthy, though maybe not essential since, again, we will only need the maps $G$ on right half-planes that are somewhat away from the boundary of their domains of definition.
We us sum this up in the following commuting diagram:
\begin{equation}\label{eq:glu}
\begin{tikzcd}[row sep=large]
\dom G_{s,b} \arrow[d,"\exp"] \arrow[r,shift left=3pt,"G_{s,b}"] & \dom G_{s',b'} \arrow[d,"\exp"] \arrow[l,shift left=3pt,"G_{s',b'}"]
\\
\C \setminus [0,b'] \arrow[r,"a\on{id}+b"'] & \C \setminus [0,b]
\end{tikzcd}
\end{equation}

We denote
\[\cal G = \setof{G_{s,b}}{s\in\C,\ b\in\C}
.\]
It is not a group for composition.\footnote{Through a notion of germ via restrictions to right half-planes, they form a group under composition. However we will not need this.}

When $\Re(z)\to +\infty$ then $G_{s,b}(z) = z + s + o(1)$.

\subsubsection{Motivation}

Consider a degree $d=m+1$ pole $p$ of a meromorphic connection on a Riemann surface, and a Riemann chart $z$ sending it to $0$.
We saw in \Cref{sub:pf:prop:di} that the change of variable $v=-\frac{a_{-{d}}}{m} z^{-m}$ changes the expression of the connection to $\Gamma_v = 1+o(1)$ as $|v|\to \infty$.
In this sense, $v$ is close to a log-chart.
Recall that $v$ sends the pole to infinity and is a $d$-fold covering from a neighborhood of $p$ to a neighborhood of infinity.
We will prove in the present section, using technical lemmas postponed to \Cref{sec:tech}, that a neighborhood of the puncture can actually be covered by log-charts mapping to simple shapes (right half-planes or unions of upper, left and lower half-planes) of the exponential-affine $\cal E$ defined in \Cref{sub:ExpoAffine}, with gluing maps thar are restrictions of maps in $\cal G$ and on which we have some control, in particular they are close to maps of the form $z\mapsto z+s$ in log-charts, $s\in\C$.

\subsubsection{Model}\label{ss:model}

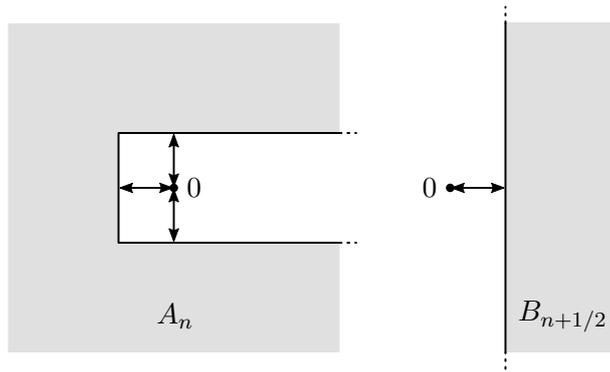
\begin{figure}
\begin{tikzpicture}
\node at (0,0) {\includegraphics[width=8cm]{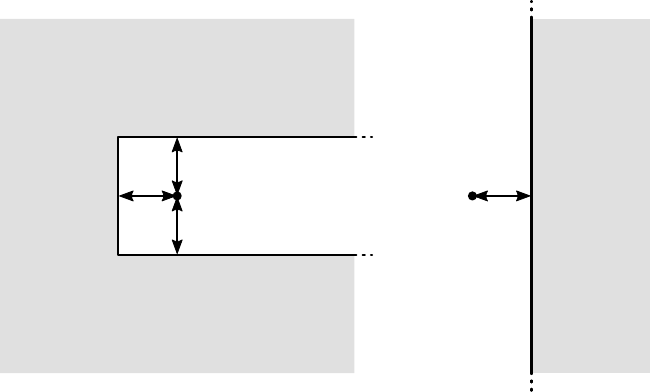}};
\node at (-1.55,0) {$0$};
\node at (1.55,0) {$0$};
\node at (-1.8,-1.7) {$A_n$};
\node at (3.3,-1.7) {$B_{n+1/2}$};
\end{tikzpicture}
\caption{The sets $A_n$ and $B_{n+1/2}$. The arrows have length $R$.}
\label{fig:AnBn}
\end{figure}

Recall that we denote
\[m=d-1.\]

We will define $2m$ open subsets of $\cal E$, denoted $A_n$ and $B_{n+1/2}$, $n\in\Z/m\Z$ where the indices live in $\Q/d\Z$.
Fix $R>0$.
For each $n\in \Z/d\Z$ we let
\[A_n = \Hp_{y>R} \cup \Hp_{y<-R} \cup \Hp_{x<-R}\]
and
\[B_{n+\frac12} = \Hp_{x>R}\]
where $\Hp$ designates a half-plane whose equation in $z=x+iy$ is written as an index.
See \Cref{fig:AnBn}.
Note that, as subsets of $\C$, the $A_n$ are identical, and so are the $B_{n+1/2}$. However, we later take a quotient of their disjoint union.

It will be convenient to denote
\[AB_n = A_n\text{ and }AB_{n+1/2}=B_{n+1/2}.\]

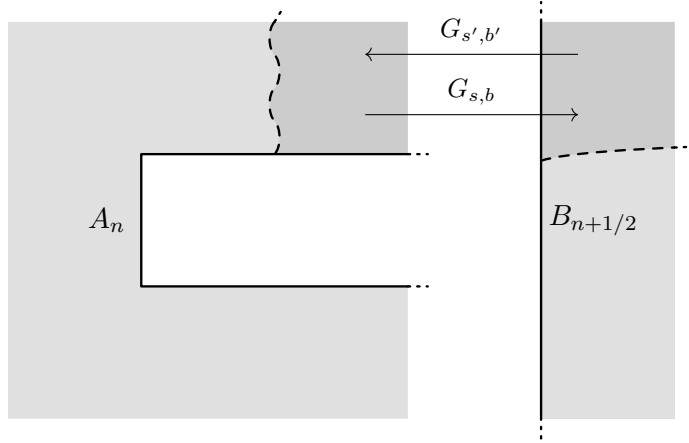
\begin{figure}
\begin{tikzpicture}
\node at (0,0) {\includegraphics[width=9cm]{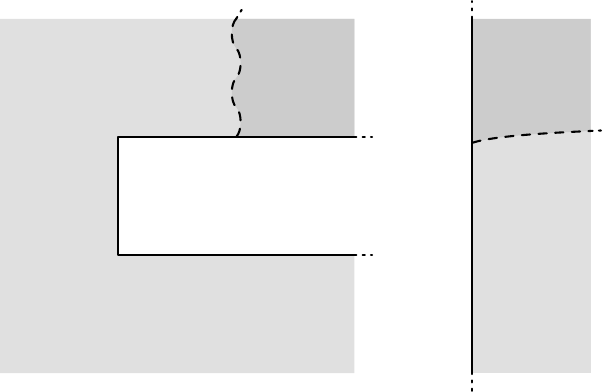}};
\draw[->] (0.2,1.4) -- node[above] {\small $G_{s,b}$} (3,1.4);
\draw[->] (3,2.2) -- node[above] {\small $G_{s',b'}$} (0.2,2.2);
\node at (-3.2,0) {$A_n$};
\node at (3.2,0) {$B_{n+1/2}$};
\end{tikzpicture}
\caption{Gluing of the top part $\Hp_{y>R}$ of $A_n$ to $B_{n+1/2}$ by an element of $\cal G$.
More precisely $z\in\Hp_{y>R}\cap A_n$ is glued to $z'\in B_{n+1/2}$ iff $z'=G_{s,b}(z)$ (iff $z=G_{s',b'}(z')$).
Here $s'=-s$ and $b'=-be^{-s}$ (i.e.\ $s=-s'$ and $b=-b'e^{-s'}$).}
\label{fig:Gglue}
\end{figure}

Each $A_n$ is glued along a subset of $\Hp_{y>R}$ to $B_{n+1/2}$ and along a subset of $\Hp_{y<-R}$ to $B_{n-1/2}$, in both case with a map $G_{s,b}\in\cal G$ for various $s$ and $b$, where $\cal G$ is defined in \Cref{sub:log}.
More precisely we choose maps $G_{n+1/4}\in\cal G$ and $G_{n-1/4}\in\cal G$ and we glue $z\in A_n\cap \Hp_{y>R}$ to $z'\in B_{n+1/2}$ iff $G_{n+1/4}(z)=z'$.
And $z\in B_{n-1/2}$ to $z'\in A_n\cap \Hp_{y<-R}$ iff $G_{n-1/4}(z)=z'$.
See \Cref{fig:Gglue,fig:ronde}.

\begin{figure}
\begin{tikzpicture}
\node at (0,0) {\includegraphics[scale=0.45]{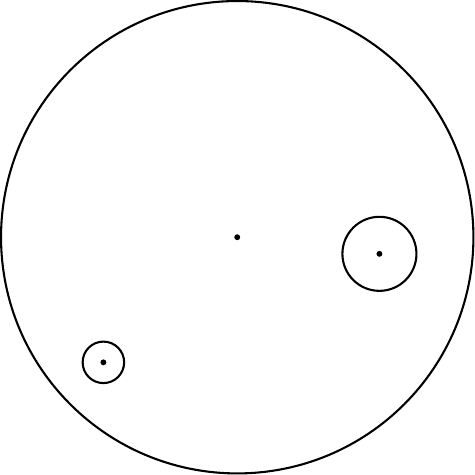}};
\end{tikzpicture}
\caption{The big circle is the image by $\exp$ of $\partial B_{n+1/2} = ``x=R"$.
The two small disks are the images by $e^s\exp +b$ of $ \Hp_{x<-R}$, part of which is in the boundary of respectively $A_n$ and $A_{n+1}$, for the  value of $(s,b)$ such that $G_{s,b} = G_{n+1/4}$ and the  value of $(s,b)$ such that $G_{s,b} = G_{n+3/4}^{-1}$.
The conditions in \cref{eq:cond} mean that the two latter disks are surrounded by the circle.}
\label{fig:gein}
\end{figure}

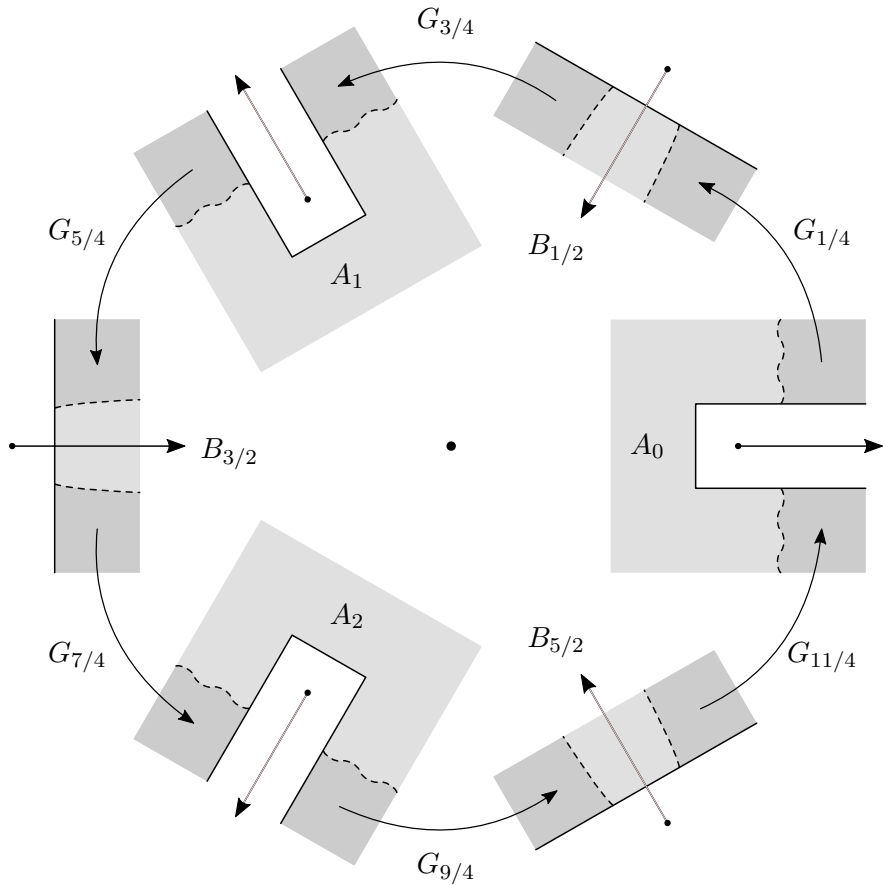
\begin{figure}
\begin{tikzpicture}
\node at (0.35,0) {\includegraphics[width=11.5cm]{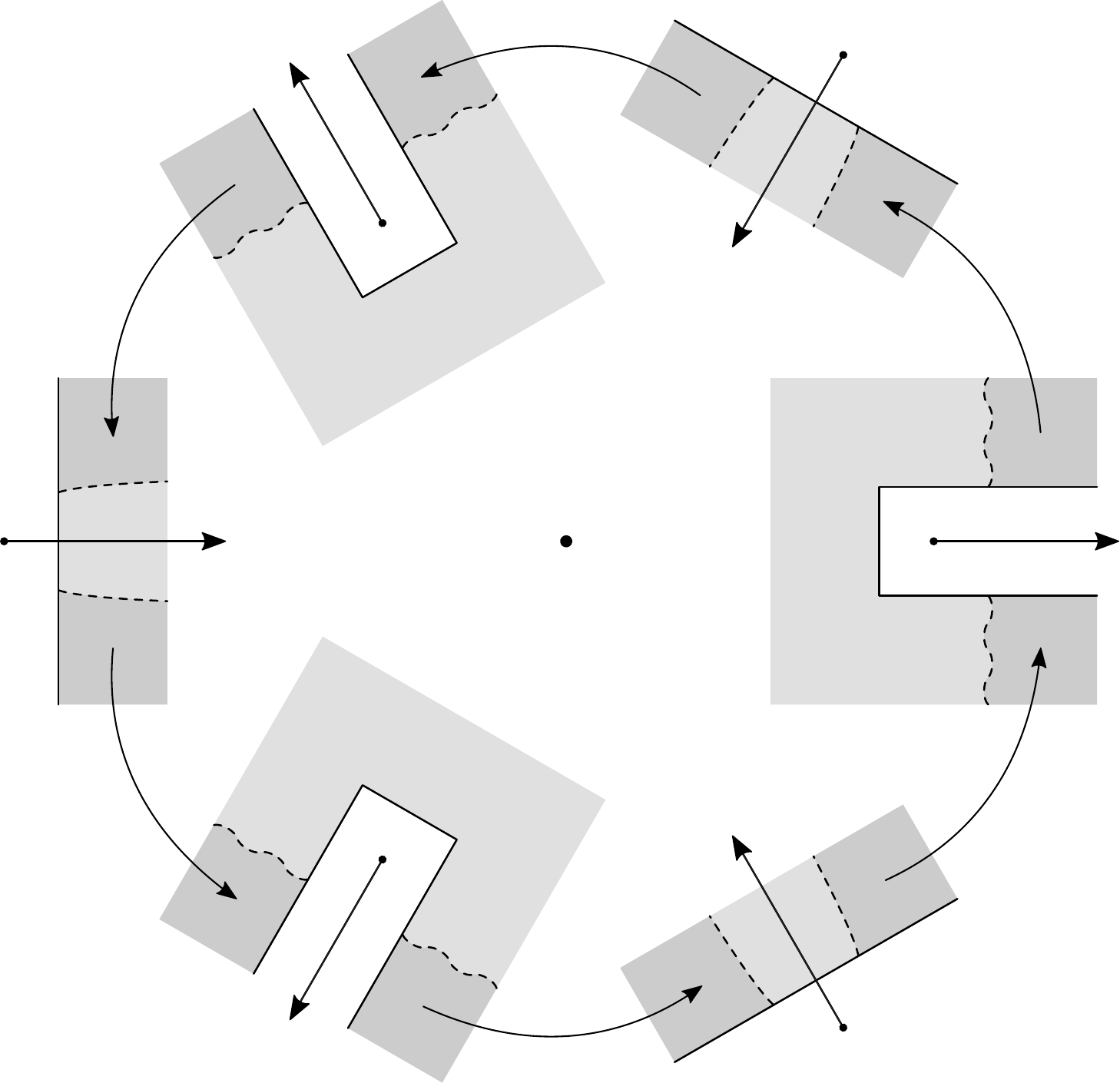}};
\node at (3,0) {$A_0$};
\node at (-0.95,2.25) {$A_1$};
\node at (-0.95,-2.25) {$A_2$};
\node at (-2.5,-0.1) {$B_{3/2}$};
\node at (1.8,2.6) {$B_{1/2}$};
\node at (1.8,-2.6) {$B_{5/2}$};
\node at (5.3,2.8) {$G_{1/4}$};
\node at (0.35,5.6) {$G_{3/4}$};
\node at (-4.5,2.8) {$G_{5/4}$};
\node at (-4.5,-2.8) {$G_{7/4}$};
\node at (0.35,-5.6) {$G_{9/4}$};
\node at (5.3,-2.8) {$G_{11/4}$};
\end{tikzpicture}
\caption{We arranged in circle the different sets $A_n$ and $B_{n+\frac{1}{2}}$. Straight arrows represent the positive real axis.}
\label{fig:ronde}
\end{figure}

For the construction to give an appropriate manifold we will need $R$ to be big enough. More precisely we assume that $G_{n+1/4}^{-1}$ and $G_{n+3/4}$ are both defined on the whole set $B_{n+1/2}$ and that they map $B_{n+1/2}$ to a set disjoint from $\Hp_{x<-R}$.
This translates for $G_{n+1/4} = G_{s,b}$ into the conditions
\begin{gather}
 e^R  > |b|, 
\\
 e^{R} > |b| + e^{\Re s} e^{-R}.\label{eq:cond}
\end{gather}
The second implies the first.
For $G_{n+3/4}$ written as an inverse : $G_{n+3/4} =G_{s,b}^{-1}$, the conditions on $(s,b)$ are the same.
A geometric interpretation is given on \Cref{fig:gein}.

\newcommand{\textModManif}{
  The topological quotient $\cal M$ of the disjoint union of the $A_n$ and $B_{n+1/2}$ by the gluings $G_{n\pm 1/4}$, where $n$ ranges in $\Z/m\Z$, is a topological manifold $\cal M$.
  Each $A_n$ and $B_n$ injects in $\cal M$ by the quotient map homeomorphically to their images.
}
\begin{lem}\label{lem:modmanif}
  \textModManif
\end{lem}
\begin{proof}
  This is standard. For sake of completeness we included a proof in \Cref{ss:pf:lem:modmanif}.
\end{proof}

We define an atlas of affine surface on the quotient $\cal M$ of the lemma above using the set $A_n$ and $B_n$ as log-charts.
More precisely consider the following restrictions of the quotient map: $\pi_n : A_n\to \cal M$ and $\pi_{n+1/2} : B_{n+1/2}\to\cal M$.
Their inverses $\pi_x^{-1}$ is considered as taking values in the log-chart $\C$ of $\cal E$.
We post-compose them with $\exp$, and restrict them to the open sets on which $\exp\circ\,\pi_x^{-1}$ is injective.
Near every point of $\cal M$ there is such a restriction.
This is our atlas: by \cref{eq:glu}, the change of charts are locally restrictions of affine maps of $\C$ so we indeed get an affine surface.

\begin{defn}\label{def:model}
  We call \emph{model} the affine manifold $\cal M$ defined above.
\end{defn}

\newcommand{\textModPunc}{
The model can be extended as a Riemann surface by addition of a point at infinity, i.e.\ we have a puncture.
}
\begin{lem}\label{lem:modpunc}
\textModPunc
\end{lem}
\begin{proof}
  See \Cref{ss:pf:lem:modpunc}.
\end{proof}

Call $\pi$ the quotient map from the disjoint union of the $A_n$ and $B_{n+1/2}$ to $\cal M$:
\[ \pi : \coprod_{x\in \frac{1}{2}\Z/m\Z} AB_x  \to \cal M
.\]
Let $w$ be a Riemann chart sending the puncture to $0$.

\newcommand{\textEstimateW}{
There exists a complex constant $\lambda\in\C^*$ and branches of $z\mapsto z^{1/m}$ on each $A_n$ and $B_{n+1/2}$ such that $w \circ \pi(z) \sim \lambda/z^{1/m}$
when $|z|\to+\infty$.
}
\begin{lem}\label{lem:estimate:w}
\textEstimateW
\end{lem}
\begin{proof}
See \Cref{ss:pf:lem:modpunc}.
\end{proof}

\begin{cor}\label{cor:mero}
\leavevmode\begin{enumerate}%\leavevmode fixes an AMS bug
\item\label{item:3} The connection is meromorphic of order $d=m+1$ at $p$.
\item\label{item:4} A path $\gamma$ in $\cal M$ tends to $p$ along a repelling direction if and only if it eventually enters some $\pi(A_n)$ and $\pi^{-1}\circ\gamma$ tends to $\infty$ with an argument tending to $\pi$.
\end{enumerate}
\end{cor}
\begin{proof}
For the first claim, we use the change of variable formula from the charts $AB_x$ to the chart $w$. For $x\in\frac{1}{2}\Z/m\Z$ denote $\pi_x$ the restriction of $\pi$ to $AB_x$.
By \Cref{lem:estimate:w}, for $z\in AB_x$, $w\circ\pi_x(z) \sim \lambda/z^{1/m}$ as $|z|\to +\infty$ for some branch of the $m$-th root. 
In particular the image of $w\circ\pi_x$ contains a sector based on $0$, of opening angle close to $2\pi/m$ and median axis a repelling axis of the pole for $A_n$, and in the case of a $B_{n+1}$ the opening angle is close to $\pi/m$ and the median axis is an attracting axis of the pole.
Let $\psi_x = (w\circ\pi_x)^{-1}$.
Then $\psi_x(z) \sim \nu z^{-m}$ for $\nu=\lambda^m$ and 
using \Cref{lem:sec2} twice:
\[ \psi_x'(z)\sim -\frac{\nu m}{z^{m+1}}
\quad\text{and}\quad
\psi_x''(z)\sim \frac{\nu m(m+1)}{z^{m+2}}
\]
both hold on slightly smaller sectors $V_x$.
We can arrange so that these smaller sectors still cover a neighborhood of the origin.
Now by \cref{eq:g} the expression of $\Gamma$ in the chart $w$ is
\[\Gamma_w = \psi_x' \times 1 + \frac{\psi_x''}{\psi_x'}\]
in the sector $V_x$.
It follows that $\Gamma_w(z) \sim \frac{-\nu m}{z^{m+1}}$ as $z\to 0$ and this implies that the function $\Gamma_w$ has a pole of order $m+1$ at the origin.

For the second claim, the repelling axes in the $w$-coordinate have direction $\theta$ defined by $\lambda^m e^{-im\theta} \in\R_-$.
The claim then follows from the already mentioned asymptotic-equivalent $\psi_x(z) \sim \lambda^m z^{-m}$. 
\end{proof}

Let
\[L_{s,b}(z) = e^s z+b.\]
and for $x=\frac{2n+1}{4}$ let $s_x,b_x$ so that
\[G_{s_x,b_x}=G_{x}\text{ and let } L_x = L_{s_x,b_x}.\]
Then
\[\exp\circ\, G_{x} = L_x \circ \exp.\]

\loctitle{Developing map} A universal covering $\wt{\cal M}$ of $\cal M$ is obtained by taking $n$ in $\Z$ instead of $\Z/m\Z$ and $x\in \frac{1}{2} \Z$ instead of its quotien by $m\Z$.
In this universal covering, if we take a germ of affine chart at a point of $\wt{\cal M}$, it extends uniquely to an affine map $\tilde\phi:\wt{\cal M}\to \C$ whose restriction to the $AB_x$ for various $x$ we denote $\tilde\phi_x$.
Then there exists $\Lambda_{x}\in\on{Aff}\C$
such that
\begin{equation}\label{eq:pLe}
\tilde\phi_{x} = \Lambda_{x}\circ \exp
\end{equation}
holds on the corresponding $A$ or $B$ and 
we have the relation
\begin{equation}\label{eq:LLL}
\Lambda_{x} = \Lambda_{x+\frac12} \circ L_{x+\frac14},
\end{equation}
which is better viewed on the commutative diagram of \Cref{fig:cd}.
Conversely any family $\Lambda_{n/2}$, $n\in\Z$ such that that \cref{eq:LLL} holds for all $x\in\frac12\Z$ defines an affine map $\tilde\phi : \wt{\cal M} \to \C$ via $\tilde \phi_{x} = \Lambda_{x}\circ\exp$.

Note that if the $L_{x+\frac14}$ are fixed, then one $\Lambda_{x}$ determines all the others.

\begin{figure}[htbp]
\begin{center}
\begin{tikzcd}
 {} \arrow{r}{G_{-1/4}} & A_0 \arrow{r}{G_{1/4}} \arrow{d}{\exp} & B_{1/2} \arrow{r}{G_{3/4}} \arrow{d}{\exp} & A_1 \arrow{r}{G_{5/4}} \arrow{d}{\exp} & B_{3/2} \arrow{r}{G_{7/4}} \arrow{d}{\exp} & A_{2} \arrow{r}{G_{9/4}} \arrow{d}{\exp} & B_{5/2} \arrow{r}{G_{11/4}} \arrow{d}{\exp} & {}
 \\
 {} \arrow{r}{L_{-1/4}}  & \C \arrow{r}{L_{1/4}} \arrow{d}{\Lambda_{0}} & \C \arrow{r}{L_{3/4}} \arrow{d}{\Lambda_{1/2}} & \C \arrow{r}{L_{5/4}} \arrow{d}{\Lambda_{1}} & \C \arrow{r}{L_{7/4}} \arrow{d}{\Lambda_{3/2}} & \C \arrow{r}{L_{9/4}} \arrow{d}{\Lambda_{2}} & \C \arrow{r}{L_{11/4}} \arrow{d}{\Lambda_{5/2}} & {}
 \\
 {} \arrow{r}{} \arrow{r}{\Id} & \C \arrow{r}{} \arrow{r}{\Id} & \C \arrow{r}{} \arrow{r}{\Id} & \C \arrow{r}{} \arrow{r}{\Id} & \C \arrow{r}{} \arrow{r}{\Id} & \C \arrow{r}{} \arrow{r}{\Id} & \C \arrow{r}{} \arrow{r}{\Id} & {}
\end{tikzcd}
\end{center}
\caption{In this diagram the maps $G_{\ldots}$ are only defined on a subset of the corresponding $A$ or $B$.
Each square of the diagram is commutative but composing two consecutive $G$ gives the empty map because of these restrictions.
The diagram is \emph{not periodic}: $\Lambda_{x+m}\neq \Lambda_x$ exept in some special situations.
Note also that the rest of the diagram is periodic, but the composition of the horizontal maps of the second row along a period is not the identity, except in particular cases.}
\label{fig:cd}
\end{figure}

\loctitle{Asymptotic value family and affine holonomy}
By \eqref{item:4} of \Cref{cor:mero} and \cref{eq:pLe}, the asymptotic value family of $\tilde\phi$ is 
\begin{equation}
u_n = \Lambda_n(0)
\end{equation}
for $n\in\Z$.
For $x\in \frac{1}{2}\Z$, the affine map $\Lambda_{x+m} \circ \Lambda_{x}^{-1} \in\on{Aff}(\C)$ is independent of $x$: indeed
$\Lambda_{x+m} \circ \Lambda_{x}^{-1} = (\Lambda_{x+1/2+m} \circ L_{x+1/4+m}) \circ (\Lambda_{x+1/2} \circ L_{x+1/4})^{-1} = \Lambda_{x+1/2+m} \circ \Lambda_{x+1/2}^{-1}$
since $L_{x+1/4+m}=L_{x+1/4}$.
So there is some $L\in\on{Aff}(\C)$ such that $\forall x\in\frac12\Z$,
\begin{equation}\label{eq:monodromyForLambda}
\Lambda_{x+m} = L\circ \Lambda_{x}.
\end{equation}
This map $L$ is the affine holonomy of $\tilde\phi$ w.r.t.\ the pole: indeed if we call $\hat\sigma$ the maps sending $x\in AB_x$ to $x\in AB_{x+d}$,  $x\in \frac{1}{2}\Z$, and $\sigma$ the induced map on $\wt{\cal M}$, such that $\pi\circ\hat\sigma = \sigma\circ \pi$, then $\sigma$ is a generator of the group of deck transformation of the covering $\wt {\cal M}\to \cal M$, which corresponds to winding once around $0$ in the positive direction in $w$-coordinates by \Cref{lem:estimate:w}, and $\tilde\phi\circ\sigma = L\circ \tilde\phi$.

Recall the notation $L_{s,b}(z) = e^s z+b$, and that 
we set $\forall x$, $L_x=L_{s_x,b_x}$ where $G_x=G_{s_x,b_x}$.
Then
\begin{equation}\label{eq:Ldecomp}
\Lambda_{x}^{-1} \circ L \circ \Lambda_{x}= L_{x+1/4}^{-1} \circ L_{x+3/4}^{-1} \circ \cdots \circ L_{x+m-1/4}^{-1}
\end{equation}
so in particular
\begin{equation}\label{eq:Lsb}
L = L_{s,b}\text{ with }s = -s_{1/4}-s_{3/4}-\cdots-s_{m-1/4}
\end{equation}
hence:
\begin{equation}\label{eq:eser}
e^{s_{1/4}}e^{s_{3/4}}\cdots e^{s_{m-1/4}} = e^{2\pi i\res}.
\end{equation}
But actually we have more precise, denoting $d=m+1$ as before:
\begin{lem}\label{lem:modres}
\[s_{1/4}+s_{3/4}+\cdots+s_{m-1/4} = 2\pi i(\res-d).\]
\end{lem}
\begin{proof}
See \Cref{ss:pf:lem:modres}
\end{proof}

This ends the description of the model.

\subsubsection{Realizability of invariants}\label{ss:realizeInv}

The following theorem states in particular that every invariant can be realized. The set $\cal U_{d,\res}$ and the set of invariants $\cal I_{d,\res}$ have been defined in \Cref{def:is}.

\begin{thm}\label{thm:model2}
For all $d\geq 2$, $\res\in\C$ and family $u_n\in\C$, $n\in\Z$ such that $u\in \cal U_{d,\res}$,
there is a model $\cal M$, as described in \Cref{ss:model}, that defines an affine surface with a meromorphic puncture of order $d$, residue $\res$ and having a developing map whose asymptotic value family is $u$. \end{thm}
\begin{proof}
Let $m=d-1$ and $a=\exp(-2\pi i\res)$.
The fact that $u\in \cal U_{d,\res}$ means that there exists some $b\in\C$ such that $\forall n\in\Z$, $u_{n+m} = au_m+b$.
Denote
\[\wideparen{L}(z)=az+b\]
so that $\forall n\in\Z$, $u_{n+m}= \wideparen{L}(u_{n})$.
In the model, we can choose the maps $G_x$, i.e.\ we can choose the complex numbers $s_x$ and $b_x$: the construction above provide some $R$ big enough so that it can be performed.
Once the $G_x$ are fixed, the choice of $\Lambda_0$ determines all the other $\Lambda_x$ (as $\Lambda_{x+1/2} = \Lambda_{x} \circ L_{x+1/4}^{-1}$ by \cref{eq:LLL}) amounts then to the choice of a developing map and we want to perform all those choices so that $\Lambda_n(0)=u_n$ for all $n\in\Z$ and so that $s_{1/4}+\cdots+s_{m-1/4} = 2\pi i(\res-d)$.
We can actually freely choose the $s_x\in\C$, as long as they satisfy the previous equation.
Then we only need the bottom two rows of \cref{fig:cd}, which we reproduce below.
\begin{center}
\begin{tikzcd}
 {} \arrow{r}{L_{-1/4}}  & \C \arrow{r}{L_{1/4}} \arrow{d}{\Lambda_{0}} & \C \arrow{r}{L_{3/4}} \arrow{d}{\Lambda_{1/2}} & \C \arrow{r}{L_{5/4}} \arrow{d}{\Lambda_{1}} & \C \arrow{r}{L_{7/4}} \arrow{d}{\Lambda_{3/2}} & \C \arrow{r}{L_{9/4}} \arrow{d}{\Lambda_{2}} & \C \arrow{r}{L_{11/4}} \arrow{d}{\Lambda_{5/2}} & {}
 \\
 {} \arrow{r}{} \arrow{r}{\Id} & \C \arrow{r}{} \arrow{r}{\Id} & \C \arrow{r}{} \arrow{r}{\Id} & \C \arrow{r}{} \arrow{r}{\Id} & \C \arrow{r}{} \arrow{r}{\Id} & \C \arrow{r}{} \arrow{r}{\Id} & \C \arrow{r}{} \arrow{r}{\Id} & {}
\end{tikzcd}
\end{center}
We claim that, then, it is possible to choose the $b_x$ and $\Lambda_0$ so that $\forall n$, $\Lambda_n(0) = u_n$.
Note that the condition only concerns half of the $\Lambda_x$: those whose indices are entire.
To simplify, denote
\[ L_{n+1/2} := L_{n+3/4}\circ L_{n+1/4}
,\]
so that the diagram becomes
\begin{center}
\begin{tikzcd}
 {} \arrow{r}{L_{-1/2}}  & \C \arrow{r}{L_{1/2}} \arrow{d}{\Lambda_{0}} & \C \arrow{r}{L_{3/2}} \arrow{d}{\Lambda_{1}} & \C \arrow{r}{L_{5/2}} \arrow{d}{\Lambda_{2}} & \C \arrow{r}{L_{7/2}} \arrow{d}{\Lambda_{3}} & {}
 \\
 {} \arrow{r}{} \arrow{r}{\Id} & \C \arrow{r}{} \arrow{r}{\Id} & \C \arrow{r}{} \arrow{r}{\Id} & \C \arrow{r}{} \arrow{r}{\Id} & \C \arrow{r}{} \arrow{r}{\Id} & {}
\end{tikzcd}
\end{center}
We are going to use the fact that we fixed the linear factor of each $L_{x+1/4}$ (to $e^{s_{n+1/4}}$), hence of each $L_{n+1/2}$ (to $e^{s_{n+1/4}}e^{s_{n+3/4}}$).
This inductively imposes the linear factor of each $\Lambda_n$ by $\Lambda_{n+1} \circ L_{n+1/2}^{-1} = \Lambda_{n}$.
Since we know the value of $\Lambda_n$ at $0$ this imposes $\Lambda_n$, and hence this imposes $L_{n+1/2}$ by $\Lambda_{n+1} \circ L_{n+1/2}^{-1} = \Lambda_{n}$ again.
There remains to check that $n\mapsto L_{n+1/2}$ is $m$-periodic.
For this we realize that $\Lambda_{n+m}\circ\Lambda_n^{-1}$: has the same dilation factor $e^{-s_{n+1/4}} e^{-s_{n+3/4}} \cdots e^{-s_{n+m-1/4}}$ as $\wideparen{L}$; sends $u_n$ to $u_{n+1}$ as $\wideparen{L}$.
Hence
\[\Lambda_{n+m}\circ\Lambda_n^{-1} = \wideparen{L}.\]
It follows that $\Lambda_{n+m}\circ\Lambda_n^{-1} \circ \Lambda_{n+1} \circ \Lambda_{n+m+1}^{-1} = \wideparen{L} \circ \wideparen{L}^{-1} = \Id_\C$, hence $\Lambda_n^{-1} \circ \Lambda_{n+1} = \Lambda_{n+m}^{-1}\circ \Lambda_{n+m}$, i.e.\ $L_{n} = L_{n+m}$.
\end{proof}

From the construction in the proof above:

\begin{prop}\label{prop:mc}
  For any $u\in\cal U_{d,\res}$, for any choice of $s_x\in\C$ such that $s_{1/4}+\cdots+s_{m-1/4} = 2\pi i(\res-d)$, and any choice of $\Lambda_0\in \Aut(\C)$ such that $\Lambda_0(0)=u_0$, there exists a unique choice of the $\Lambda_n$ and of $L_{n+1/2}:=L_{3+1/4}\circ L_{n+1/4}$ such that the corresponding model $\cal M$, defined for $R$ big enough, together with the developing map $\tilde\phi$ associated to the $\Lambda_n$, has asymptotic value family $(u_n)$.
\end{prop}

Note that knowing $L_x$ only fixes $G_x$ up to a translation by an integer multiple of $2\pi i$.

Note that in \Cref{prop:mc}, the compositions $L_{n+4/4}\circ L_{n+1/4}$ are determined but not the individual $L_{n+1/4}$ and $L_{n+3/4}$: there is still some freedom, due to the fact that the precise placement of the sets $B_{n+1}$ (corresponding to circular neighborhoods of $\infty$ after mapping by $\exp$) is not relevant, whereas the sets $A_n$ (corresponding to neighborhoods of $0$ after mapping by $\exp$) have to be positioned precisely relative to each other, so as to fit the sequence $u_n$.

If one fixes every choice in \Cref{prop:mc} but changes the linear factor $a_0$ of $\Lambda_0: z\mapsto a_0 z+u_0$ into $a_0' = \lambda a_0$, $\lambda\in\C^*$, then this multiplies the linear factor of every $\Lambda_n$ by the same number $\lambda$, and each $L_x$ is replaced by $m_\lambda^{-1} \circ L_x\circ m_\lambda $ where $m_\lambda: z\mapsto \lambda z$.

Finally, changing $u$ to $A\circ u$ for some $A\in \Aut(\C)$ can be realized without changing the model $\cal M$: just replace every $\Lambda_n$ with $A\circ\Lambda_n$.

\subsubsection{A canonical model}\label{sub:canomodel}

\begin{figure}
\begin{tikzpicture}
\node at (0,0) {\includegraphics[width=12cm]{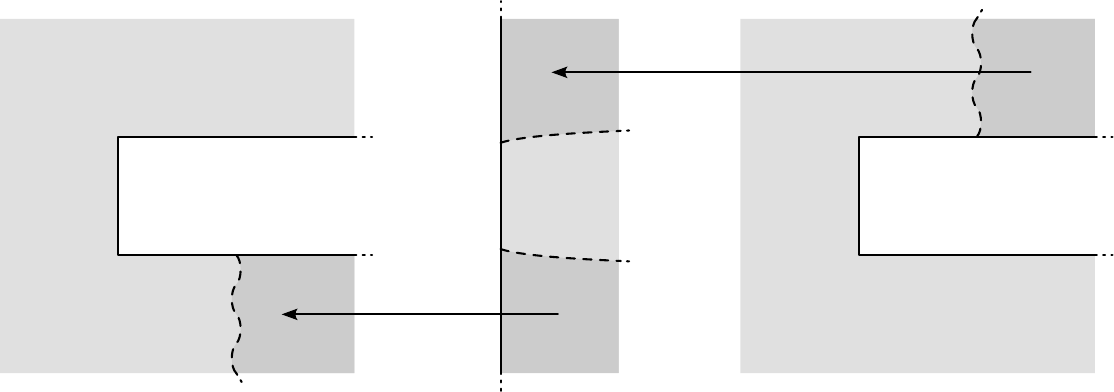}};
\node at (-4,0) {$A_n$};
\node at (-1.4,-1) {$G_{n-1/4}$};
\node at (4,0) {$A_{n-1}$};
\node at (2.2,0.9) {$G_{n-3/4}$};
\node at (0.2,0) {$B_{n-1/2}$};
\end{tikzpicture}
\caption{Illustration for \Cref{sub:canomodel}}
\label{fig:canonical}
\end{figure}

For a given invariant in $\cal I_{\res,d}$, there is not uniqueness in our construction in \Cref{ss:model} of the model $\cal M$ realizing it: first we chose a representative $u\in\cal U_{d,\res}$ of the invariant, then we freely chose the linear factor $a_0\in\C^*$ of $\Lambda_0 = z\mapsto a_0 z +u_0$, and any collection $s_x\in\C$ 
only bound by the condition $s_{1/4}+\cdots+s_{m-1/4} = 2\pi i(\res-d)$. Even with that there is still some freedom as this only fixes the compositions $L_{n+3/4}\circ L_{n+1/4}$.
In the end there is the choice of $R$ big enough, but this last choice is rather mild as far as we are concerned in constructing germs of connection multiple poles.
 
\loctitle{Canonical model}
We make here a particular choice, trying to keep things simple.
Other choices would suit as well.
We still have a choice of representative $u$.
We still leave $R$ free.
\begin{enumerate}
  \item We take each $G_{n-1/4} = G_{0,0}=\Id_\C$ (they are the  gluings of the lower parts), in particular $L_{n-1/4}=\Id_\C$.
  \item We choose the same value of $s_x=s$ for all $G_{n+1/4}$, which is necessarily
\end{enumerate}
\[s=2\pi i \frac{\res-d}{d}
.\]
We choose $\Lambda_0(z) = z+u_0$.
We saw in \Cref{prop:mc} that there is a unique choice of $L_{n+3/4}\circ L_{n+1/4}$ such that the corresponding model and developing map has asymptotic value family $(u_n)$.
Since we impose $L_{n+3/4}=\Id_\C$, this means $L_{n+1/4}$ is uniqely determined.

Let us determine everything explicitly:
From $L_{n-1/4}=\Id_\C$ the diagram implies $\Lambda_{n-1/2} = \Lambda_n$.
The map $\Lambda_n$ must map $0$ to $u_n$ so let us write it $\Lambda_n(z) = a_n z +u_n$ with $a_0=1$.
The condition $\Lambda_n = \Lambda_{n+1/2} \circ L_{n+1/4}$ reads 
$a_n z + u_n = a_{n+1} (e^s z + b_{n+1/4}) + u_{n+1}$
i.e.
\[a_{n+1} = a_n e^{-s}\quad\text{and}\quad a_{n+1}b_{n+1/4} = u_n-u_{n+1}\]
since $a_0=1$ we get
\[a_{n} = e^{-ns}\quad\text{and}\quad b_{n+1/4} = (u_n-u_{n+1})e^{-(n+1)s}
.\]

\begin{figure}[h]
\begin{center}
\begin{tikzcd}
 {} \arrow{r}{\on{Id}} & A_0 \arrow{r}{G_{1/4}} \arrow{d}{\exp} & B_{1/2} \arrow{r}{\on{Id}} \arrow{d}{\exp} & A_1 \arrow{r}{G_{5/4}} \arrow{d}{\exp} & B_{3/2} \arrow{r}{\on{Id}} \arrow{d}{\exp} & A_{2} \arrow{r}{G_{9/4}} \arrow{d}{\exp} & B_{5/2} \arrow{r}{\on{Id}} \arrow{d}{\exp} & {}
 \\
 {} \arrow{r}{\Id}  & \C \arrow{r}{L_{1/4}} \arrow{d}{\Lambda_{0}} & \C \arrow{r}{\Id} \arrow{d}{\Lambda_{1}} & \C \arrow{r}{L_{5/4}} \arrow{d}{\Lambda_{1}} & \C \arrow{r}{\Id} \arrow{d}{\Lambda_{2}} & \C \arrow{r}{L_{9/4}} \arrow{d}{\Lambda_{2}} & \C \arrow{r}{\Id} \arrow{d}{\Lambda_{3}} & {}
 \\
 {} \arrow{r}{} \arrow{r}{\Id} & \C \arrow{r}{} \arrow{r}{\Id} & \C \arrow{r}{} \arrow{r}{\Id} & \C \arrow{r}{} \arrow{r}{\Id} & \C \arrow{r}{} \arrow{r}{\Id} & \C \arrow{r}{} \arrow{r}{\Id} & \C \arrow{r}{} \arrow{r}{\Id} & {}
\end{tikzcd}
\end{center}
\caption{}
\label{fig:cd2}
\end{figure}

\begin{prop}\label{prop:isocan}
  Every model has a neighborhood of its puncture affinely isomorphic to a canonical model.
\end{prop}
\begin{proof}
Recall that $R$ must be big enough so that the models must follow the condition in \eqref{eq:cond} illustrated by \Cref{fig:gein}. Consider a model $\cal M$ and developing map. It comes with some choice of $R>0$. Let us find $R'>0$ and an isomorphism of a subset of $\cal M$ to the canonical model $\cal M'$ and developing map with the same asymptotic value family.

First we justify that we can modify and make all 
$G_{n-1/4}$ equal to (restrictions of) the identity.
Consider the maps $G_{n-3/4}$ and $G_{n-1/4}$, semi-conjugated by $\exp$ to $L_{n-3/4}$ and $L_{n-1/4}$, and whose appropriate restriction are used to glue $B_{n-1/2}$ to $A_n$ and $A_{n-1}$, see \Cref{fig:canonical}.
In the initial model, the gluing map from $B_{n-1/2}$ to $A_n$ is a restriction  of $G_{n-1/4}\in\cal G$.
Let $G=G_{n-1/4}^{-1}$, which is also in $\cal G$.
For some $R'>R$, $G(\Hp_{x>R'}) \subset B_{n-1/2}$.
Though we will not use the following fact, it is interesting to note that $R'$ is quite close to $R$ if $R$ is already big.
We replace $B_{n-1/2}=\Hp_{x>R}$ by $B'_{n-1/2} = \Hp_{x>R'}$ and $G_{n-1/4}$ by the identity restricted to $\Hp_{y<-R}\cap \Hp_{x>R'}$.
For the isomorphism between a subset of $\cal M$ and $\cal M'$, we map a point $z$ in $G(\Hp_{x>R'})\subset B_{n-1/2}$, to $G^{-1}(z) \in B'_{n-1/2} = \Hp_{x>R'}$.
We then modify the gluing map $G_{n-3/4}$ between subsets of $B_{n-1/2}$ and $A_{n-1}$ into a new gluing map from $B'_{n-1/2}$ to $A_{n-1}$ which is $G_{n-1/4}\circ G_{n-3/4}$, so that composition $L_{n-1/4}\circ L_{n-3/4}$ does not change.
Here we do not have to increase $R'$ to make it work, essentially because the first step replaced the outer circle $C$ in \Cref{fig:gein} by a circle $C'$ centered on one of the inner disks, and surrounding $C$, so it will surround both small disks.
The above procedure can be performed independently for all $n$, and we can use a common value of $R'$ for all $n$.
Then we reduce the $A_n$ to $\Hp_{y>R'} \cup \Hp_{y<-R'} \cup \Hp_{x<-R'}$.
The parts of $\cal M$ thrown away correspond to bounded subsets of the initial $A_n$ and $B_{n+1}$ so we cover a neighborhood of the puncture.
This proves our first claim.

Then to make the values of $s_{n+1/4}$ uniform, we only have to translate the pieces $AB_x$ (using the same translation for $A_n$ and $B_{n-1/2}$ to keep the previous property that their gluing is the identity) by some vector $v_n\in\C$, which replaces $G_{n-3/4}$ by $T_{v_n}\circ G_n \circ T_{v_{n-1}}^{-1}$.
We then restrict them to (unions of) half-planes with the same parameter $R''$ big enough.
The conditions of \eqref{eq:cond} are still valid because the corresponding disks in \Cref{fig:gein} decrease while the outer circle increases.

The above two procedures produce an isomorphism, in a neighborhood of the puncture, to a canonical model.
The developing map can be restricted to this neighborhood and transferred to the canonical model, and this does not change the asymptotic value family.
\end{proof}

\subsubsection{Proof of \texorpdfstring{\Cref{lem:modmanif}}{Lemma ...}}\label{ss:pf:lem:modmanif}

\begin{figure}
\begin{tikzpicture}
\node at (0,0) {\includegraphics[width=9cm]{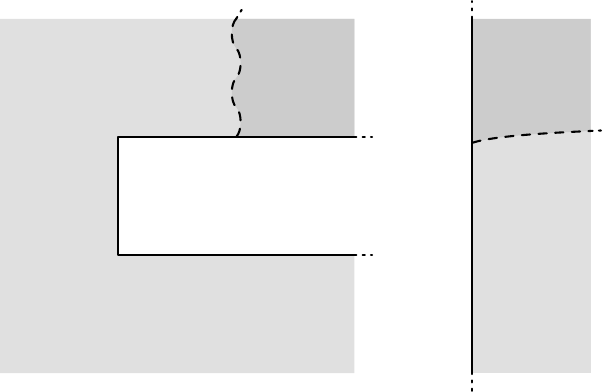}};
\draw[->] (3,2.2) -- node[above] {\small $G_{n+1/4}$} (0.2,2.2);
\node at (-3.2,0) {$A_n$};
\node at (3.2,0) {$B_{n+1/2}$};
\end{tikzpicture}
\caption{One of the cases in the proof of \Cref{lem:modmanif}}
\label{fig:separ}
\end{figure}

We repeat below \Cref{lem:modmanif}.
\begin{lem*}
\textModManif
\end{lem*}

Denote $X$ this disjoint union and $\sim$ the equivalence relation generated by the gluings.
More precisely we have the generating set $R_0 \subset X\times X$ consisting of the $(x,x')$ where ($x\in A_n\cap \Hp_{y>R}$, $x'\in B_{n +1/2}$ and $G_{n+ 1/4}(x)=x'$) and the $(x,x')$ where ($x\in A_n\cap \Hp_{y<-R}$, $x'\in B_{n +1/2}$ and $G_{n+ 1/4}(x)=x'$).

Call $A'_{n+1/4}$ the set of $z\in A_n$ such that $\Im(z)>R$ and $G_{n+1/4}(z) \in B_{n+1/2}$ and $B'_{n+1/4}$ its image by $G_{n+1/4}$.
Similarly call $A'_{n-1/4}$ the set of $z\in A_n$ such that $\Im(z)<-R$ and $G^{-1}_{n-1/4}(z) \in B_{n-1/2}$ and $B'_{n-1/4}$ its image by $G^{-1}_{n-1/4}$.
The gluings happen between the sets $A'_{n+1/4}$ and $B'_{n+1/4}$, and between $A'_{n-1/4}$ and $B'_{n-1/4}$.
Because these $4n$ sets are disjoint, it turns out that the relation $\sim$ only consists in $x\sim x'$ and $x'\sim x$ for $(x,x')\in R_0$ and of $x\sim x$ for all $x\in X$.

It follows that only points in adjacent sets $A_{n}$ and $B_{n\pm 1/2}$ can be identified by the quotient.
This first implies that $A_n$ and $B_{n+1/2}$ inject in $\cal M$ under the quotient map.
So the quotient map is locally injective.
The equivalence relation $\sim$ is generated by homeomorphisms between open subsets of $X$.
The quotient map $\pi : X\to X/\sim$ is thus necessarily open and since the quotient map is locally injective, the quotient map is a local homeomorphism.
In particular $X/\sim$ is locally homeomorphic to open subsets of $\R^2$.
Since $X$ is second countable and $\pi$ open, the quotient is second countable.

There remains to check that $X/\sim$ is separated, i.e.\ that distinct points have disjoint neighborhoods.
For points in the same piece this follows from the quotient map being injective on them, and open.
For $x\in(\frac{1}{2}\Z)/m\Z$ note that the images of $G_{x-1/4}$ and $G_{x+1/4}^{-1}$ in the piece $AB_{x}$ have disjoint closures.
So points in $X/\sim$ with representatives in non-adjacent pieces have disjoint saturated neighborhoods.
In the remaining cases, the only non-trivial one is when one point is on the boundary of the image and the other one on the range, of one of the gluing maps $G_x$, see \Cref{fig:separ}.
The map $G_x$ and its inverse extend to neighborhoods in $\C$ of the closure of these images/ranges into injective and continuous maps.
Separation then follows from this and the fact that the boundary of the domain of $G_x$, intersected with $AB_x$, is disjoint from the preimage by the extension of $G_x$ of the boundary of the range of $G_x$. And of a similar statement for $G_x^{-1}$.

\subsubsection{Proof of \texorpdfstring{\Cref{lem:modpunc,lem:estimate:w}}{Lemmas ...}}\label{ss:pf:lem:modpunc}

We repeat below \Cref{lem:modpunc}:
\begin{lem*}
\textModPunc
\end{lem*}

It is enough to prove that the Riemann surface is isomorphic near the puncture $P$ to $\D-\{0\}$ by a $K$-quasiconformal map for some $K>1$ sending $P$ to $0$.\footnote{This does not require to solve a Beltrami equation: it follows for instance from the techniques in the proof of one of the Grötzch inequalities, see \cite{Ahl06}, Section~1.B; to apply it to annuli instead of rectangles, work in the log-coordinates.}
One way is as follows: in the particular case where the $B_{n\pm 1/2}$ are glued to the $A_n$ by (a restriction of) the identity map then appropriate branches of $z\mapsto 1/z^{1/m}$ on each $AB_x$ define a conformal isomorphism to a punctured neighborhood of $0$.
In the general case, we take the same map on the $A_n$, but we first send $B_{n+1/2}$ to another set $B'_{n+1/2}$ as follow: if $z\in B_{n+1/2}$ is identified to some $z'=G_{n+1/4}^{-1}(z) \in A_n$ or to some $z' = G_{n+3/4}(z)\in A_{n+1}$ then we (have to) send $z$ to $z'\in\C$.
Otherwise we interpolate as follows: trace a vertical segment through $z$ joining the boundary of the regions.
The point $z$ is a barycenter of the ends of this segment. Send it to the point in $\C$ that is the barycenter with the same weight, of the already defined images of these endpoints.
The interpolation we defined is quasiconformal provided we perform the above construction on sub-pieces of $A_n$ and $B_{n+1}$ where $R$ is increased to a big enough value.  This induces a map $\psi$ from $\cal M$ to a quotient $\cal N$ of the disjoint union of the $A_n$ and $B'_{n+1}$, where the gluing maps are restrictions of the identity to open subsets of $\C$.
We then send $z'\in \cal N$ to a neighborhood of $0$ using appropriate branches of $1/(z')^{1/m}$. Call $s:\cal N \to \C^*$ this map.
The map $s$ is conformal so the composition $s\circ \psi$ is a quasiconformal mapping sending the puncture to $0$.

\medskip

We repeat below \Cref{lem:estimate:w}.
\begin{lem*}
\textEstimateW
\end{lem*}

Consider the map $s\circ\psi$ in the previous proof.
The map $s$ is explicit and $|s(z)|= 1/|z|^{1/m}$.
The estimate $\psi(z)\sim z$ as $|z|\to \infty$ is straightforward since $G_{r,b}(z) = z + \cal O(1)$ on the right half-plane $\cal B_{n+1}$. 
Decompose
\[w=\phi\circ s \circ \psi\]
where $\phi := w\circ \psi^{-1}\circ s^{-1}$. Then $\phi$ is quasiconformal with a quasiconformal constant that tends rapidly to $1$ near the puncture. Let us make this more precise:
First the map $\psi$ is the identity on $\cal A_n$ and equal to $G_{r,b}$ on the part (call it the top part) of $B_{n+1/2}$ glued to $\cal A_n$ and to another $G_{r,b}$ on the part (call it the bottom part) of $B_{n+1/2}$ glued to $\cal A_{n+1}$. Call intermediate part the complement in $B_{n+1/2}$ of its top and bottom parts.
The map $\psi\circ\pi$ is conformal on $A_n$ and on the top and bottom part of $B_{n+1/2}$.
We have on $B_{n+1/2}$:
\[ G_{r,b}(z) = z+r+\cal O(e^{-\Re (z)})\text{ and }G'_{r,b}(z) = 1+\cal O(e^{-\Re (z)})
.\]
The intermediate part of $B_{n+1/2}$ lays in a strip of the form $|\Im(z)| < C$ for some constant $C$ that depends on $R$ and on the values of $r$ on top and bottom.
From this and the estimates above it is standard to get that
\[\psi\circ\pi(z) = z + \cal O(e^{-\Re(z)})\text{ and }D_z\psi\circ\pi = \on{Id}+\cal O(e^{-\Re(z)})
\]
where $D_z\Psi\circ\pi$ denotes the differential of $\psi\circ\pi$ at $z$, which is $\R$-linear but usually not $\C$-linear.
Since $z$ is contained in a strip and $\Re(z)>R$, $e^{-\Re(z)}$ and $e^{-|z|}$ have bounded ratio.
The seam between the top and bottom part, where $\psi\circ\pi$ is one of the $G_{r,b}$ and the intermediate part where it is given by another smooth map, is an analytic curve, hence quasiconformally erasable.
It follows that $\psi\circ\pi$ is quasiconformal, with a constant of the form
\[K(z)=1+\cal O(e^{-|z|})\]
The image of the circular ellipse field by this map, then by $s$, is an ellipse field that is straightened by $\phi$ and hence whose ratio gives the infinitesimal quasiconformal constant of $\phi$.
Using $\psi\circ\pi(z) \sim z$ we get that $\phi$ is has, at a given $z\in\C$, a quasiconformal ratio
\[K(z) = 1+ \cal O(e^{-1/|z|^{1/m}})
\]
as $|z|\to 0$, which is indeed a rapid decrease, faster than $1+\cal O(|z|^2)$, which is already sufficient 
for $\phi$ to be $\C$-differentiable at the origin: see for instance \cite{LehtoVirtanen} Lemma~6.1 page~224.
The claim then follows.

\subsubsection{Derivative on a sector of a function close to a power map}\label{sub:tech2}

We begin with a standard application of Cauchy estimates. It is used in \Cref{cor:mero,lem:modres}.

\begin{lem}\label{lem:sec2}
Let $f$ be a holomorphic function defined on a bounded or unbounded sector of the form $|z-c|<r$ (or $|z-c|>r$) and $\arg (z-c)\in(\theta-\alpha,\theta+\alpha)$ where $\theta\in\R$, $0<\alpha<\pi$ and $c\in\C$.
Assume that there exists $\lambda\in\C^*$ and $\beta\in\C$ such that $f(z) \sim \lambda z^\beta$ for a branch of $z\mapsto z^\beta$, as $|z|\to 0$, resp.\  $|z|\to +\infty$.
Then for any $\alpha'<\alpha$, $f'(z) \sim \frac{\beta}{z}f(z)$ as $|z|\to+\infty$ within the sub-sector $\arg(z)\in (\theta-\alpha',\theta+\alpha')$.
\end{lem}
\begin{proof}
Let $S$ be the sector. We first assume it is bounded.
There is a branch of $\log f$ on $S\cap B(0,\eps)$ for $\eps$ small enough.
Let $g(z)=\log f(z)-\beta\log z - \log \lambda$.
Then $g(z) \tends 0$ as $z\to 0$, for an appropriate determination of $\log \lambda$.
We use the following Cauchy formula: $g'(z) =  \frac{1}{2\pi}\int_{0}^{2\pi} g(z+re^{i\theta})r^{-1} e^{-i\theta} d\theta$ valid whenever $\ov B(z,r)$ is contained in $S$.
If $z$ belongs to a subsector $S'$ as above, then one can take $r = a|z|$ for some $a\in(0,1)$.
This gives $|g'(z)| = |\frac{1}{2\pi} \int_{0}^{2\pi} g(z+re^{i\theta}) r^{-1} e^{-i\theta}d\theta| \leq \frac{1}{2\pi a |z|} \int_{0}^{2\pi} |g(z+re^{i\theta})| d\theta$.
It follows that $g'(z) =o(1/|z|)$ as $z\to 0$ within $S'$.
Now $g'(z) = f'(z)/f(z) - \beta/z$.
Hence $f'(z)/f(z) \sim \beta/z$ as $z\to 0$ within $S'$.

The unbounded case is actually proved exactly the same way. It can alternatively be deduced from the bounded case by applying it to $z\mapsto f(1/z)$ and using the chain rule.
\end{proof}

\subsubsection{Proof of \texorpdfstring{\Cref{lem:modres}}{Lemma ...}}\label{ss:pf:lem:modres}

We repeat below \Cref{lem:modres}.
\begin{lem*}
\[s_{1/4}+s_{3/4}+\cdots+s_{m-1/4} = 2i\pi(\res-d).\]
\end{lem*}

Here we work in $\cal M$ and not its universal covering $\wt{\cal M}$.
The index $x$ in $AB_x$, $G_x$, etc.\ will thus live in $\Q/m\Z$ and not in $\Q$.

Let $\Gamma = \Gamma_w$ be the Christoffel symbol in $w$-coordinate associated to the connection.
Then
\[ 2\pi i\res = -\oint \Gamma_w dw
\]
where the integral is taken along any path whose image in $w$-coordinate has winding number around $0$ equal to $1\in\Z$.
Recall that the $AB_x$ are log-affine charts.
Let $v=(\lambda/w)^m$: it is a coordinate that lives in the $m$-fold covering of $\C^*$ and which is close to infinity.
Then by \Cref{eq:g},
\[ \Gamma_w = \Gamma_v - \frac{m+1}{w}
.\]
In particular
\[ 2\pi i\res = -\int \Gamma_v dv + 2\pi id
\]
where the integral is taken around the image of the path in $v$-coordinates.

We now consider the following path $t\mapsto\gamma(t)\in\cal M$, that depends on a parameter $X\gg R$ that we will let tend to infinity near the end.
It starts for the lower right corner of the square $[-X,X]\times[-X,X]\subset \R^2\equiv\C$ in the chart $A_0$ and follows three sides of the square in the \emph{clockwise} sense for the chart (i.e.\ is counterclockwise in $w$-coordinate), up to the upper right corner $X+iX$.
We then pass to the chart $B_{1/2}$.
The path arrived at $z'= G_{1/4}(X+iX)$. We consider the translate of the square by the vector $G_{1/4}(X+iX)-(X+iX)$ so that its upper right corner is $z'$.
Note that this translation vector remains bounded as $X\tends \infty$ because $G_{1/4}$ is close to the identity at $X+iX$.
We then follow the left side of this square down to its lower right corner.
This point of $B_{1/2}$ corresponds to a point in $A_{1}$ to which we translate the lower right corner of the square, which we follow in $A_1$ along three sides like for $A_0$, and so on\ldots\ %
Once the path is back in Chart $A_0$, it is not back to its starting point but off by some quantity, which remains bounded as $X$ tends to infinity. We then close the path by a straight segment to its starting point $X+iX$.
The path hence consists in $m+1$ parts $\gamma_0$, $\gamma_{1/2}$, \dots, $\gamma_m$ (here the index is \emph{not} modulo $m$).

On each part of the path contained in some $AB_x$, the map $\zeta = \pi_x^{-1}: \pi(AB_x) \subset \cal M \to AB_x$ is an affine chart and
$\Gamma_v = \Gamma_\zeta + \frac{d}{dv}\log \frac{d\zeta}{dv}$.
It follows that
\[\int \Gamma_v dv = \sum_{k=0}^{m} \left( \int_{\gamma_k} \Gamma_\zeta d\zeta + \int \dlog \frac{d\zeta}{dv} \right) \]
where the sum is on every path part and $\zeta=(\pi_x)^{-1}$ for various $x\in \{0,1/2,1,\ldots,m-1/2\} \subset \Q/m\Z$.
By \Cref{lem:estimate:w}, $w\circ\pi_x(z) \sim \lambda z^{1/m}$ as $z\to\infty$ in $AB_x$, so $v\circ\zeta^{-1}(z) \sim z$.
It follows from \Cref{lem:sec2} that $(v\circ\zeta^{-1})'(z)$ tends to $1$ as $z\tends\infty$ within any sector contained in $AB_x$ with an angular margin. By margin we mean that the sector is of the form $\arg(z-c) \in (\alpha,\beta)$ with $-\pi/2<\alpha<\beta<\pi/2$ for a $B$-type piece and $0<\alpha<\beta<2\pi$ for an $A$-type piece.
As a consequence, along each of the $2m$ parts of the path, the integral $\int d\log \frac{d\zeta}{dv}$ is close to $0$ if $X$ is big, and tends to $0$ when we let $X$ tend to $+\infty$.

Now, on each piece, $\Gamma = 1$ because we are in log-charts.
So taking the sum of the integrals of $\Gamma_\zeta d\zeta$ on each part of the path, each upper side and lower side of squares will cancel-out, and each left and right sides too, and there will only remain the integral of $d\zeta$ along the last closing segment, which is equal to the difference $\Delta$ of endpoints of the segment.
This quantity is the opposite of the $2m$ translations we performed on the square $[-X,X]\times[-X,X]$, so
\[ \Delta = \sum_{x\in C_m} z_x-G_{x+1/4}(z_x)
\]
with $C_m = \{0,1/2,1,\ldots,m-1/2\}$.
Since $G_{x+1/4}(z) = z + s_{x+1/4} + o(1)$ as $\Re(z)\to +\infty$ this implies
\[ \Delta \underset{X\to\infty}= -s_{1/4} - s_{3/4} - \cdots - s_{m-1/4} + o(1)
.\]
The lemma follows.

\subsubsection{Models do model}\label{sub:momo}

Consider a meromorphic connection on a Riemann surface having a multiple pole, of order $d\geq 2$. Let $\res$ be its residue and $u_n$ the associated asymptotic value family as in \Cref{prop:di}.

\begin{thm}\label{thm:model}
A neighborhood of the puncture is isomorphic, as an affine surface, to a model $\cal M$ as defined in \Cref{ss:model}.
\end{thm}

\begin{figure}
\begin{tikzpicture}
\node at (0,0) {\includegraphics[scale=0.3]{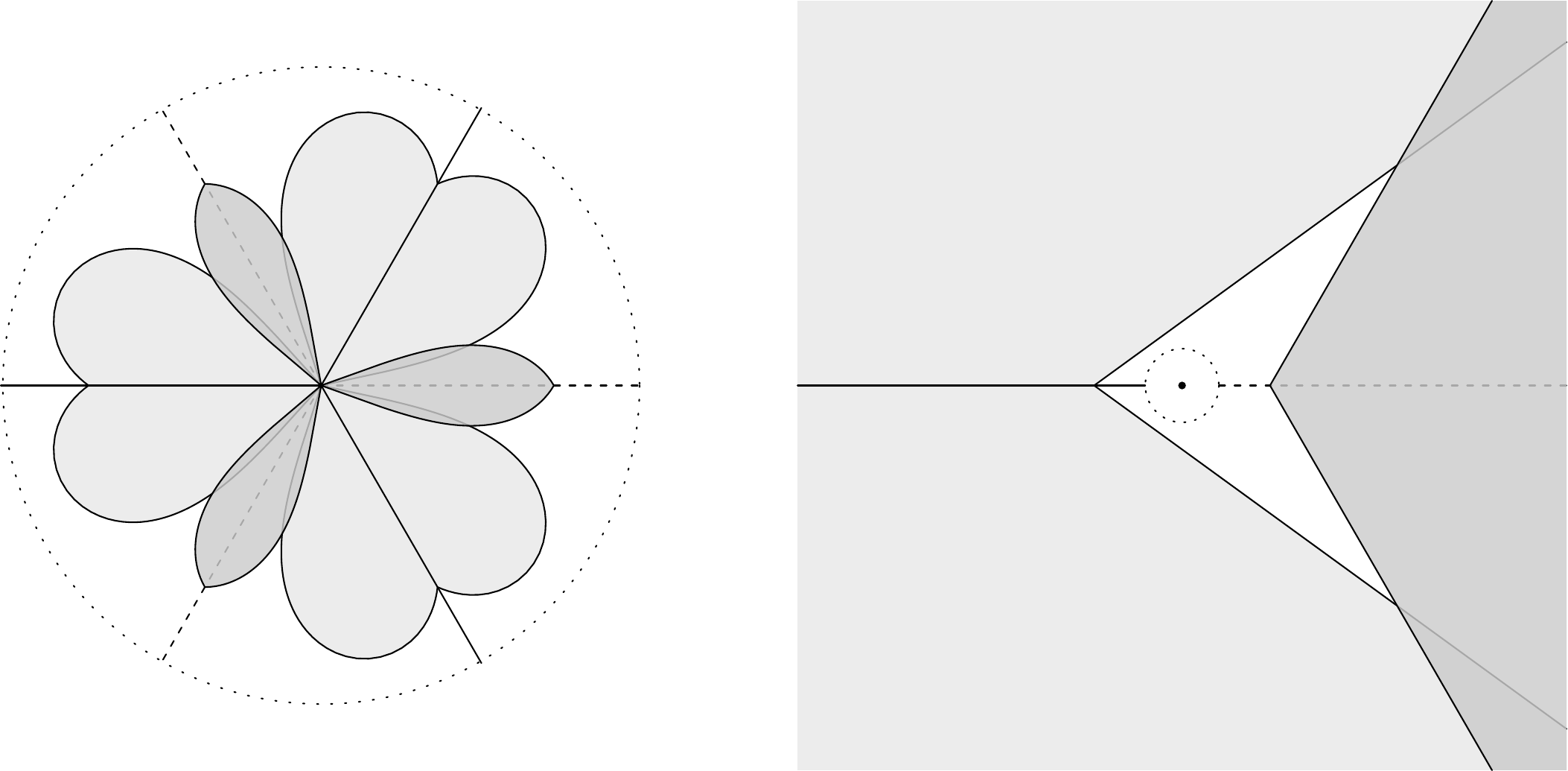}};
\node at (2.7,-2.3) {$v$};
\node at (-3.1,-2.5) {$w$};
\node at (0.9,0.7) {$S^-$};
\node at (4.9,0.7) {$S^+$};
\end{tikzpicture}
\caption{The sectors $S^-$, $S^+$ in $v$-coordinates and  the corresponding domains in $w$-coordinate.}
\label{fig:Spm}
\end{figure}
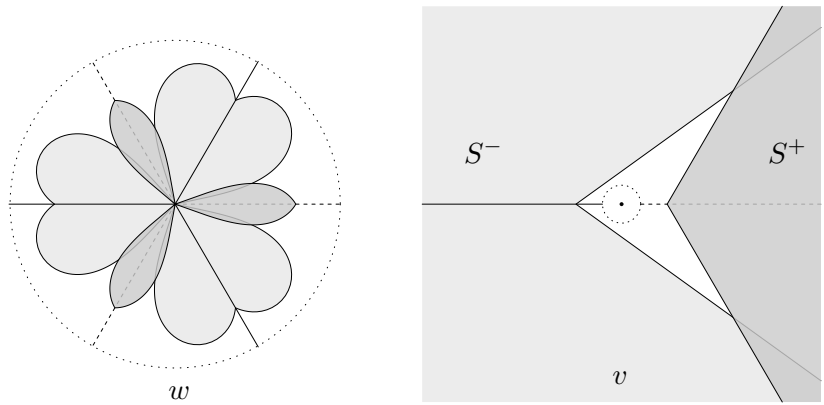

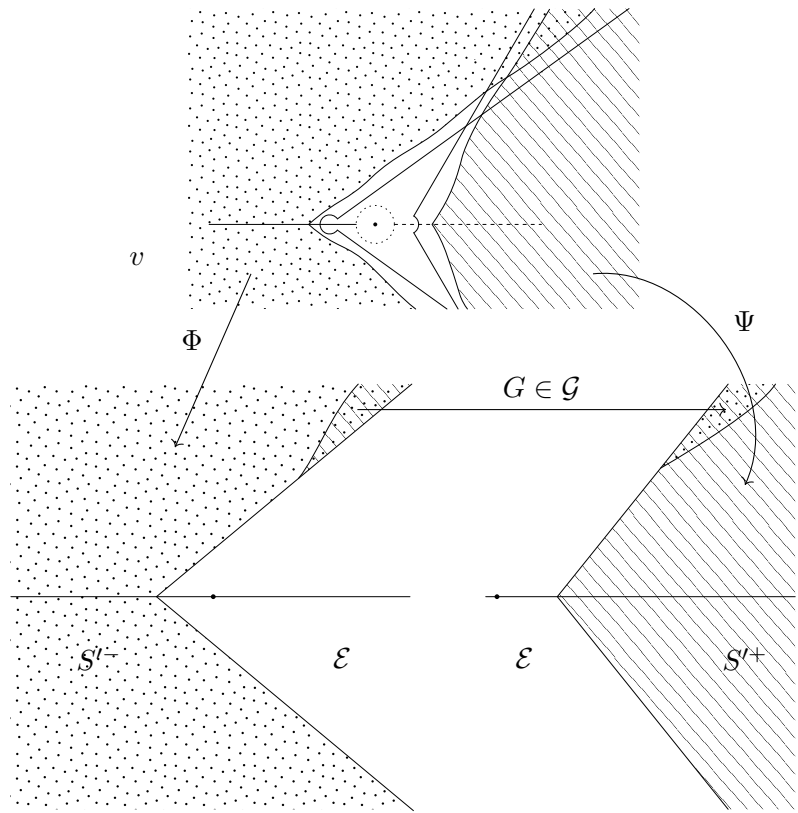
\begin{figure}
\begin{tikzpicture}
\node at (0,0) {\includegraphics[scale=0.3]{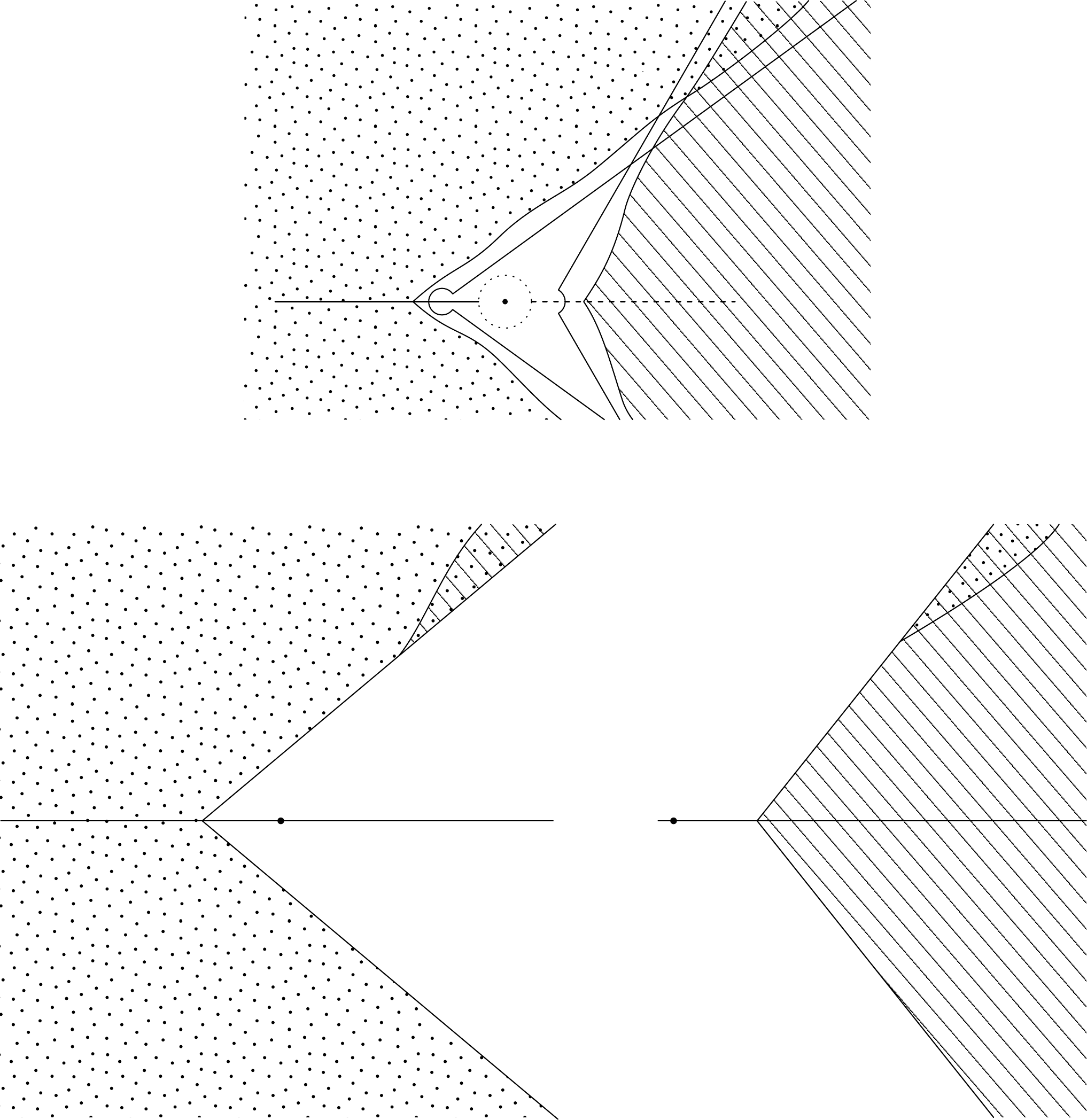}};
\node at (-3.5,2) {$v$};
\node at (-0.8,-3.3) {$\cal E$};
\node at (1.6,-3.3) {$\cal E$};
\node at (-4,-3.3) {$S'^-$};
\node at (4.5,-3.3) {$S'^+$};
\draw[->] (-0.6,0) -- node[above]{$G\in\cal G$} (4.25,0);
\draw[->] (-2,1.8) -- node[above left]{$\Phi$} (-3,-0.5);
\path[->] (2.5,1.8) edge [bend left=60] node[above right]{$\Psi$} (4.5,-1);
\end{tikzpicture}
\caption{The sectors $S'^-$, $S'^+$ in $\cal E$ are images by the affine surface isomorphisms $\Phi$ and $\Psi$ of appropriate subsets of $S^-$ and $S^+$. (These subsets are sketched in a non-realistic way so as to make the picture more readable.) The overlap between these subsets corresponds in the log-charts $S^-$ and $S^+$ to open subsets that are glued by a restriction of an element of $\cal G$ (the particular collection of log-affine maps in $\cal E$ defined in \Cref{sub:log}). The black dots mark the origin.}
\label{fig:Spm2}
\end{figure}

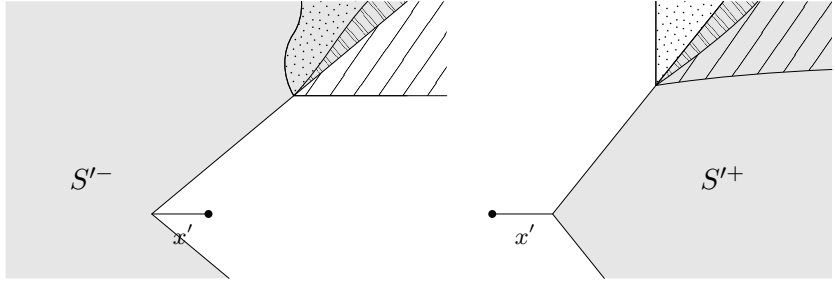
\begin{figure}
\begin{tikzpicture}
\node at (0,0) {\includegraphics[scale=0.3]{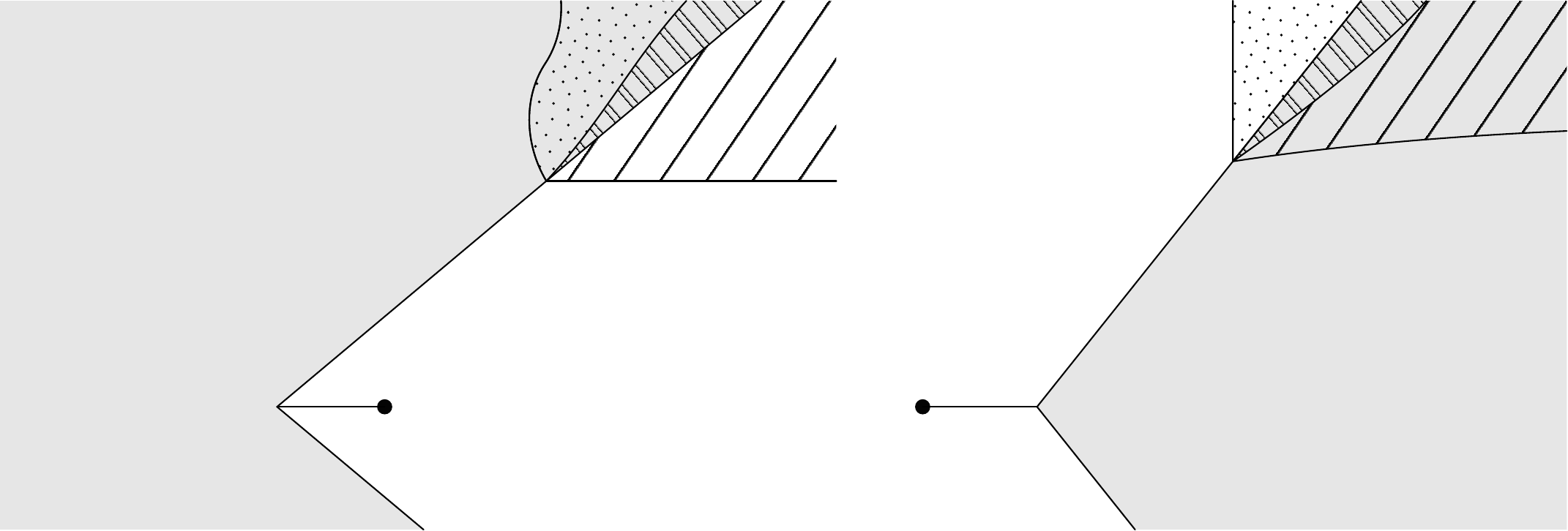}};
\node at (-4.3,-0.5) {$S'^-$};
\node at (4,-0.5) {$S'^+$};
\node at (-3.1,-1.25) {\footnotesize $x'$};
\node at (1.4,-1.25) {\footnotesize $x'$};
\end{tikzpicture}
\caption{Extending the log-charts $S'^-$ and $S'^+$ from sectors to domains containing upper half-planes for the first and a right half-plane for the second. The double hatched areas are already glued by $G\in\cal G$. Possibly using a bigger values of $x'$, the same map $G$ will map correspondingly the the dotted and single hatched areas: the single hatched is drawn on the left and is an upper half-plane minus $S'^-$. The dotted area is on the right  and is a right half-plane minus $S'^+$.}
\label{fig:Spm3}
\end{figure}

\begin{figure}
\begin{tikzpicture}
\node at (0,0) {\includegraphics[scale=0.4]{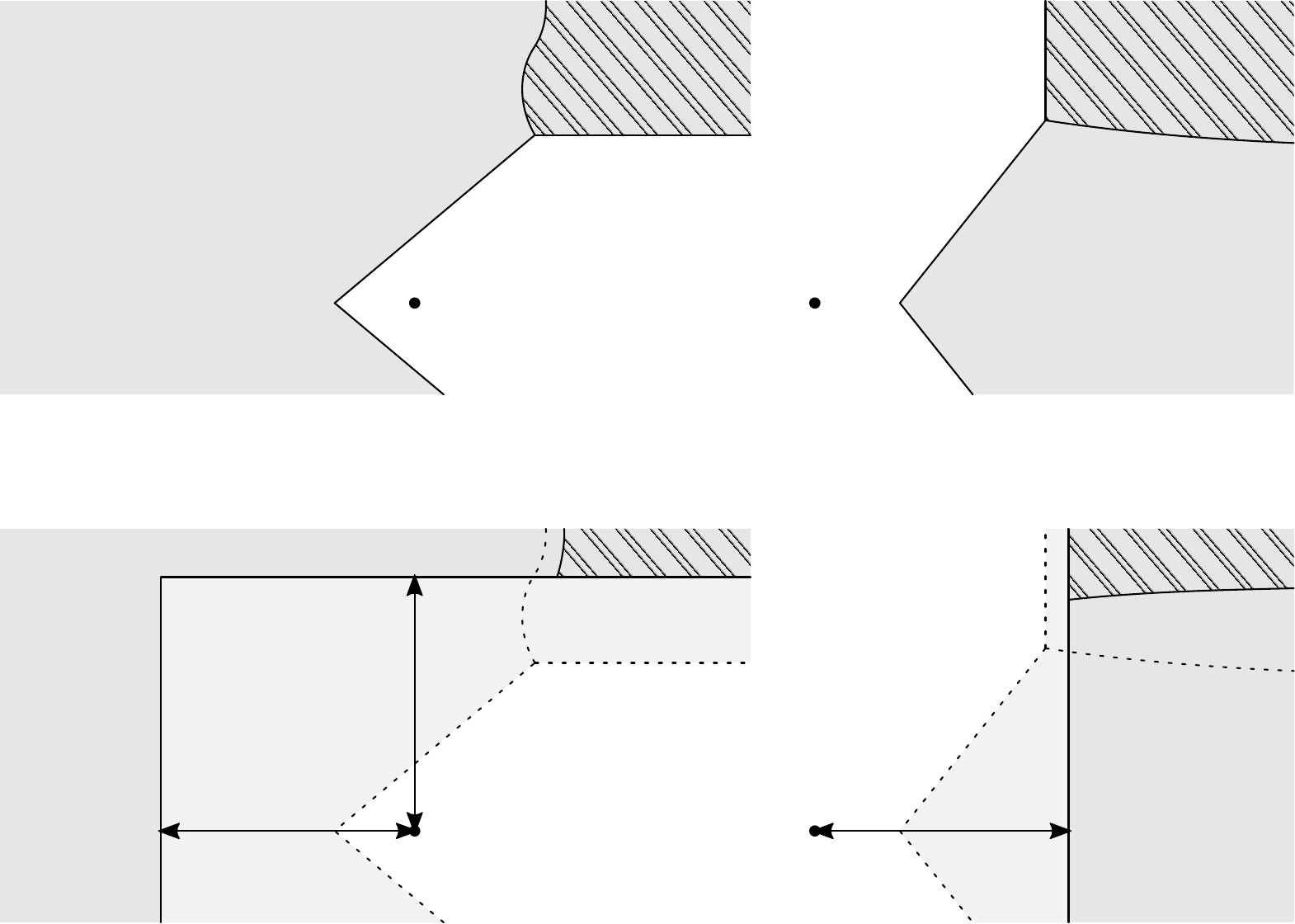}};
\end{tikzpicture}
\caption{The procedure of \Cref{fig:Spm3} gives the pieces on top, with the hatched part being glued by $G$.
By taking appropriate subsets, we end up with a model $\cal M$ of the type described in \Cref{ss:model}.
The arrows have the same length.
}
\label{fig:Spm4}
\end{figure}

Let us prove the theorem.
Let $w$ be a Riemann coordinate sending the pole to $0\in\C$.
Denote
\[m = d-1\]
and $\Gamma_w = \frac{a_{-d}}{w^d} + \cdots + \frac{a_{-1}}{w} + \cdots$.
Consider the change of variable $v=-\frac{a_{-d}}{mw^{m}}$ similar to the one in \Cref{sub:pf:prop:di}.
Recall that $w\mapsto v$ is an $m$-fold covering sending $0$ to $\infty$.
We remind that in this variable we have
\[\Gamma_v = 1+g\quad\text{with}\quad g=\cal O(1/|v|^{1/m})\]
as $|v|\to \infty$.
Let $R>0$,
\begin{gather*}
S^-= \{z\in \C~;~\frac{\pi}5<{\rm arg}(z+R) <2\pi-\frac{\pi}5\}
\\
S^+= \{z\in \C~;~-\frac{\pi}{3}<{\rm arg}(z-R) < \frac{\pi}3\}.
\end{gather*}
In $w$-coordinate, there are $m$ domains corresponding to $v\in S^-$ and $m$ domains corresponding to $v\in S^+$, and they are naturally indexed by the repelling/attracting axis that bisects each.
See \Cref{fig:Spm}.

We can apply \Cref{cor:sm2,cor:sp2}: there exists $R',x'>0$ and log-charts (affine surface isomorphisms to subsets of $\cal E$) $\Phi$ on some $S^-_{R'} = S^- - B(-R,R')$ and $\Psi$ on some $S^+_{R'} = S^+  - B(R,R')$,
whose images contain respectively a sector of the form (see \Cref{fig:Spm2})
\begin{gather*}
S'^-= \{z\in \C~;~\frac{\pi}{4+1/2}<{\rm arg}(z+x') <2\pi-\frac{\pi}{4+1/2}\}
\\
S'^+= \{z\in \C~;~-\frac{\pi}{4-1/2}<{\rm arg}(z-x') < \frac{\pi}{4-1/2}\}.
\end{gather*}
and such that $\Phi^{-1}(S'^-)$ and $\Psi^{-1}(S'^+)$ contain themselves respectively a sector of the form
$S''^-= \{z\in \C~;~\frac{\pi}{4+1/3}<{\rm arg}(z+x'') <2\pi-\frac{\pi}{4+1/3}\}$ and $S''^+= \{z\in \C~;~-\frac{\pi}{4-1/3}<{\rm arg}(z-x'') < \frac{\pi}{4-1/3}\}$, so they will overlap and cover a neighborhood of $0$ in $w$-coordinates.

By the procedure above, we covered a punctured neighborhood of $0$ in the $w$-coordinate, by $2m$ petal-like open set, half of which map bijectively to $S'^-$ by $\Phi_k$, where $k$ ranges in the set of repelling axes, and the other half to $S'^+$ by $\Psi_k$ where now $k$ ranges in the set of attracting axes.
Though not stricly necessary at this step, we arrange for simplicity so that the values of $x'$ are the same for all these sets.
These $2m$ open sets correspond to $2m$ open subsets of the original surface $X$ and their union we will call $V\subset X$.
The $2m$ functions $\Phi_k\circ v$ and $\Psi_k\circ v$ are locally affine functions from these open subsets of $X$ to the exponential-affine plane $\cal E$ whose images are respectively $S'^-$ and $S'^+$.
The affine manifold $V$ is isomorphic to the quotient of the disjoint union of $m$ copies of $S'^-$ and $m$ copies of $S'^+$, glued together along open subsets of $\cal E$ by maps that are affine isomorphisms (for the affine structure of $\cal E$).
Affine isomorphisms in $\cal E$ are locally of the form $\log(a\exp z+b)$ and pushing $x'$ to a bigger value if needed we can arrange so that the sector subsets on which the gluings occur are contained in $\Re(z)>x_1$ for arbitrarily high $x_1>0$ so we may assume that our $2m$ gluings are all restrictions of maps in $\cal G$.

The $2m$ sets $S'^-$ and $S'^+$ together with the $2m$ gluings form a log-atlas for a neighborhood of multiple pole we started from.
We now explain how to modify this atlas to obtain an atlas of a possibly smaller neighborhood and of the type $\cal M$ described in \Cref{ss:model}, which defined what we called models.
First we can always increase $x'$ a finite number of times.
Possibly using this, one can replace $S'^-$ by its union with an upper and a lower half-planes, and $S'^+$ by its union with a right half-plane, while still having a log-atlas whose changes of charts are the restriction of the same $2m$ elements of $\cal G$ but to bigger domains, as in \Cref{fig:Spm3}. \Cref{fig:Spm4} shows how to take subsets so as to get a model $\cal M$ of the desired form.

%\stopcontents

\subsection{Other asymptotic values}\label{sub:other}

Consider a multiple pole $p$ of a meromorphic connection, of degree $d$ and residue $\res$.
Consider a punctured simply connected neighborhood $U$ of $p$ and denote $\pi:\tilde U\to U$ a universal covering.
Let $\tilde \phi: \tilde U\to \C$ be a developing map (see \Cref{sub:constr:inv}, where we worked in a chart).
Let $(u_n)_{n\in\Z}$ be the asymptotic value family associated to $\tilde \phi$ in \Cref{prop:di}: these are in particular asymptotic values (as defined in \Cref{ss:term:avf}) of $\tilde\phi$.
In this section we elucidate the other possible asymptotic values of $\tilde\phi$, for paths that tend to the pole $p$.

\begin{defn}
We call \emph{$p$-asymptotic} value of $\tilde \phi$ the complex numbers $z\in\C$ such that there exists a path $\tilde\delta : [0,+\infty) \to \tilde U$ such that $\delta(t):=\pi\circ\tilde\delta(t)\tends p$ and $\tilde\phi\circ\tilde\delta(t) \tends z$, as $t\to+\infty$.
\end{defn}

The first remark is that, in the Riemann sphere $\wh\C = \C\cup\{\infty\}$, infinity\footnote{This point at infinity may be considered as artificial in the setting of affine surfaces.} is a $p$-asymptotic value of $\tilde\phi$: indeed one can take for $\tilde \delta$ any path in a attracting petal $\pi(B_{n+1/2})$ whose image in $B_{n+1/2}$ has real part tending to infinity.

Let $\gamma$ be a closed loop in $U$ winding once around $0$ in the anticlokwise direction.
Recall that $\tilde \phi \circ T = L \circ \tilde \phi$, where $T$ is the generator, associated to $\gamma$, of the deck transformation group of $\tilde U\to U$ and $L(z)=az+b$ is such that $a = e^{2\pi i \res}$.
Note that $\Im \res \neq 0 \iff |a|\neq 1$.

\begin{thm}\label{thm:av}
The $p$-asymptotic values of $\tilde\phi$ are:
\begin{enumerate}
    \item the terms $u_k$ of its asymptotic sequence,
    \item $\infty$,
    \item if $\Im \res \neq 0$, the fixed point of its affine holonomy map $L$.
\end{enumerate} 
\end{thm}

The statement means that are no other $p$-asymptotic values of $\tilde\phi$.

\medskip

In the rest of this section, we prove the theorem above.
We start by proving that if $\Im \res \neq 0$, the fixed point of its affine holonomy map $L$ is a $p$-asymptotic value.
Let $\delta:(-\infty,+\infty)$ be the periodized generator path $\gamma$ of $\pi^1(U,\gamma(0))$, i.e.\ $\delta (t) = \gamma(t -\lfloor t\rfloor)$.
Consider a lift $\tilde\delta:(-\infty,+\infty)\to\tilde U$ of $\delta$.
Then $\tilde\delta(t+n) = T^n \circ \tilde\delta(t)$ so $\tilde\phi\circ \tilde\delta(t+n) = L^n \circ \tilde\phi\circ \tilde \delta(t)$ so $\tilde\phi\circ \tilde\delta(t)$ tends to the fixed point of $L$ when $t\tends +\infty$ if $|a|<1$, or when $t\tends -\infty$ if $|a|>1$.
However $\delta$ does not tend to $p$ so we need to modify it.
Let us treat the case $|a|<1$. By composing $\wt\phi$ with a translation we can assume that the fixed point of $L$ is $0$, i.e.\ that
\[ L(z)=az
.\]
Let us work in a chart where $p$ is at $0$.
Assume the image of the chart contains $B(0,r_0)$.
We change $\gamma$ to the curve which, in the chart, is a circle of radius $r_1$.
Let $\eps_n>0$ be decreasing sequence tending to $0$, with $\eps_1=r_1$.
For all $n\geq 1$, the image by $\tilde\phi$ of the compact set $K_n$ which in polar coordinates $(r,\theta)$ is defined by $r\in[\eps_n,\eps_{n+1}]$ and $\theta\in[0,2\pi]$ is contained in some ball $B(0,R_n)$.
We choose inductively $n_1$, $n_2$, etc.\ in $\N$ so that 
\[ |a|^{n_1+\cdots+n_k} R_k\leq 1/k
\]
for all $k\geq 1$.
We start defining $\tilde \delta$ by starting at $(r,\theta)=(\eps_0,0)$ and then let $\theta$ increase to $2\pi n_1$ while leaving $r$ constant. This is part $0$ of the path.
For part $1$ we let $\theta$ constant and decrease $r$ to $\eps_2$, 
then we leave $r$ constant and let $\theta$ go from $n_1$ to $2\pi (n_1+n_2)$.
We define parts $2$, \ldots\ similarly.
We cannot control the distance to $0$ of part $0$ but for every other part $k$, we know that it is contained within distance $1/k$ from $0$.
The case $|a|>1$ is treated similarly, by winding in the other direction.

\medskip

We must now prove that all $p$-asymptotic values are of the type enumerated in \Cref{thm:av}.
Consider a path $\tilde \delta:[0,+\infty)\to\tilde U$ such that $\pi \circ \tilde \delta$ tends to the pole $p$ and such that $\tilde\phi\circ\tilde \delta$ converges.
By \Cref{thm:model}, some neighborhood $V$ of $p$ is affine isomorphic to a model, obtained by gluing finitely many log-charts $A_n$ and $B_{n+1/2}$ with $n\in\Z/mZ$.
By taking a smaller neighborhood if necessary, we can assume that tending to $p$ is equivalent to tending to infinity in the log-charts.
A model for $\tilde V$ is obtained by gluing the same sets but with index $n\in\Z$.
To distinguish $\tilde V$ from $V$, we denote them here by $\tilde A_n$ and $\tilde B_{n+1/2}$ and similarly as in \Cref{sub:model} we let $\tilde{AB}_x$ denote the first or the second according to whether $x=n$ or $x=n+1/2$.
To avoid too much notation, we identify the sets $\tilde A_n$ and $\tilde B_{n+1/2}$ with their projection in $\tilde V$, the context making it clear which one is meant.
Given $\tilde v\in \tilde V$, it there is a complex number $z\in\C$ that is projected to $\tilde v$ by $\tilde{AB}_x\mapsto \tilde V$, we call $z$ \emph{the affix of $\tilde v$ in $\tilde{AB}_x$}.
We have $B_{n+1/2} = \Hp_{x>R} \subset\C$.
 
Assume by way of contradiction that $\exists n_0\in \Z$  and $t_k\tends+\infty$ such that, denoting $x_0={n_0+1/2}$, $\tilde \delta(t_k)\in\tilde B_{x_0}$.
Let $z_k$ be the affix of $\tilde \delta(t_k)$ in $B_{x_0}$. 
From $\pi\circ\tilde\delta \tends p$ we get $|z_k| \tends \infty$.
On the other hand $\tilde\phi\circ\tilde\delta(t_k) = a_{x_0} \exp(z_k) +b_{x_0}$ for some $a_{x_0}\in\C^*$ and $b_{x_0}\in\C$.
By hypothesis $\tilde\phi\circ \tilde\delta(t)$ converges to some limit $\alpha\in\C$.
This implies that $\exp(z_k)$ converges. Its limit $\ell$ cannot be $0$ because $\Re z_k$ is bounded from below (it is in $B_{x_0}$).
Since $\tilde\phi\circ \tilde\delta(t) \tends a_{x_0}\ell+b_{x_0}$, we know $\exists T$ such that $\forall t\geq T$, $\tilde\phi\circ\tilde\delta(t) \in a_{x_0} B(\ell,r)+b_{x_0}$ where $r$ is chosen small enough so that $\exp^{-1}(B(\ell,r))\subset B_{x_0}$ and $0\notin B(\ell,r)$, which implies that connected components of $\exp^{-1}(B(\ell,r))$ are bounded.
Let $k_0$ such that $t_{k_0}\geq T$ and $W$ the component of $\exp^{-1}(B(\ell,r))$ that contains $z_{k_0}$.
Then the continuous path $\tilde \delta$ cannot leave $W$ for $t\geq t_{k_0}$.
In particular $\forall k\geq k_0$, $z_k \in W$.
This contradicts the fact that $|z_k|\tends +\infty$.
So the premisses were false, which means that for all $n_0\in\Z$, there exists a time after which the path $\tilde \delta$ cannot visit $B_{x_0}$.

So either $\tilde\delta$ enters an $A_n$ that it never leaves, or $\tilde\delta$ must wind more and more in one direction, the latter meaning that ($\forall x\in\frac{1}{2}\Z$) ($\exists T$) ($\forall t\geq T$) $\tilde\delta(t) \notin \bigcup_{x'\leq x} \tilde{AB}_{x'}$ or ($\forall x\in\frac{1}{2}\Z$) ($\exists T$) ($\forall t\geq T$) $\tilde\delta(t) \notin \bigcup_{x'\geq x} \tilde{AB}_{x'}$.

In the case $\tilde \delta_n$ enters an $A_n$ that it never leaves, recall that it must also eventually avoid $B_{n+1/2}$ and $B_{n-1/2}$, so the real part of the affix of $\tilde \delta_n$ in $A_n$ cannot tend to infinity.
Since $\tilde\phi = a_n \exp+b_n$ on $A_n$, 
$\tilde\phi\circ\tilde\delta$ can only tend to the term $u_n$ of the asymptotic sequence of $\tilde\phi$.

In the other case, we first note that $\Im \res \neq 0$. If it were, then each times the path $\tilde\delta$ has a part that traverses a $B_{n+1/2}$, the image by $\tilde\phi$ of this part would have a diameter that is bounded from below independently of $n$, preventing $\tilde\phi\circ\tilde\delta$ to converge.
So $\Im(\res)\neq 0$, i.e.\ $|a|\neq 1$, and the same argument shows that $\tilde \delta$ cannot wind in the direction that makes $\tilde\phi$ grow.
So $n\tends+\infty$ if $|a|<1$ and $n\tends-\infty$ if $|a|>1$.
Now each times the path traverses an $A_n$, it has to pass close to the corresponding $u_n$, more precisely for $z\in\Hp_{x<-R}\subset A_n$ we have $|\tilde\phi\circ\pi_n(z)-u_n|\leq |a_0| |a|^n e^-R $ where $\tilde\phi\circ\pi_0 = a_0\exp+b_0$.
Since $u_n$ tends to the fixed point of $L$ as $n\tends +\infty$ or $-\infty$ (according to whether $|a|<1$ or $|a|>1$), it implies that the latter is an accumulation point of $\tilde\phi\circ\tilde\delta$. Q.E.D.

\subsection{Perturbation of the exponential-affine plane on sectors}\label{sec:tech}

This section develop intermediate tools, grouped in \Cref{ss:csq} and used in \Cref{sub:pf:prop:di,sub:momo}.
As the title tells, we start from the global log-chart $\C$ of the exponential-affine plane, with constant Christoffel symbol $\Gamma =1$.
We consider only a sector in $\C$ and perturb $\Gamma(z)$ on this sector into a holomorphic function that is close to $1$ when $z$ is big.
In this situation, we compare the affine charts of the perturbed connection with the original ones (which are of the form $z\mapsto a\exp z + b$).

We start with technical lemmas.

\subsubsection{Estimates on some integrals}\label{sub:tech}
The inequalities in \Cref{lem:sm,lem:sp} are used in \Cref{ss:csq} to control maps sending near log-charts to log-charts.

\medskip
\noindent\textbf{Standing assumption:} In this section, $\eps\in (0,1)$, $C\in \R_+$ and $h:U\to \C$ is holomorphic on $U\subset \C\ssm\{0\}$ with 
\[\forall z\in U,\quad |h'(z)|\leq \frac{C}{|z|^\eps} .\]

\begin{lem}\label{lem:h}
    If $[z_1,z_2]\subset U$ and $\bigl|\arg(z_2)-\arg(z_1)\bigr|<\frac{5\pi}6$, then 
    \[\bigl|h(z_1)-h(z_2)\bigr| \leq \frac{C'}{|z_1|^\eps} |z_1-z_2|\quad \text{with}\quad C':=\frac{2^\eps C}{(1-\eps)}.\]
\end{lem}

\begin{proof}
We start from 
\[h(z_1) - h(z_2) = (z_1-z_2)\int_0^1 h'\bigl(tz_1 + (1-t)z_2\bigr){\rm d}t.\]
In addition, 
\[\left|\int_0^1 h'\bigl(tz_1 + (1-t)z_2\bigr){\rm d}t\right|\leq 
\int_0^1 \frac{C}{|tz_1+(1-t)z_2|^\eps}{\rm d}t = \frac{C}{|z_1|^\eps} \int_0^1 \frac{1}{t^\eps\left|1+\frac{(1-t)z_2}{tz_1}\right|^\eps}{\rm d}t.\]
Since $\ds\left|{\rm arg}\left(\frac{(1-t)z_2}{tz_1}\right)\right|= \bigl|{\rm arg}(z_2)-{{\rm arg}(z_1)}\bigr|<\frac{5\pi}{6}$, we have that 
\[\left|1+\frac{(1-t)z_2}{tz_1}\right|\geq \frac{1}{2}.\]
Thus, 
\[\int_0^1 \frac{1}{t^\eps\left|1+\frac{(1-t)w}{tz}\right|^\eps}{\rm d}t\leq 2^\eps\int_0^1 \frac{{\rm d}t}{t^\eps} = \frac{2^\eps}{(1-\eps)}.\qedhere\]
\end{proof}

For simplicity we chose a bound of $5\pi/6$ in the lemma above, but of course we may replace it by any bound less than $\pi$: this would only change the constants.

We now assume that $R>0$, $\theta\in (0,\frac{5\pi}{6})$ and $U$ is the sector bisected by $(-\infty,R)$ and of half-opening $\theta$, i.e.:
\[U = S^-= \{z\in \C~;~\pi-\theta<{\rm arg}(z+R) <\pi+ \theta\}.\]
From \Cref{lem:h} we get: 
\begin{cor}\label{cor:hphi}
Under the standing assumption on $h'$, the map $h$ satisfies
\[h(z) = {\mathcal O}\bigl(|z|^{1-\eps}\bigr)\quad \text{as } z\to \infty \text{ in } S^-.\]
In particular, the function $z\mapsto \exp\bigl(z+h(z)\bigr)$ admits an antiderivative which tends to $0$ as $\Re(z)$ tends to $-\infty$ while $z\in S^-$.%je préfère préciser car le fait que $h$ n'est pas défini en dehors de $S^-$ peut avoir échappé au lecteur, surtout qu'on va appliquer ça à des restrictions d'une fonction qu'on va appeler $h$
\end{cor}

Let us denote by $\phi:S^-\to \C$ this antiderivative. 

\begin{lem}\label{lem:sm}
For the particular choice of $\phi=\int e^{z+h(z)}dz$ above, we have:
\[\left|\frac{\phi(z)}{\phi'(z)} -1\right| =  {\cal O}\left( \frac{1}{|z|^\eps}\right)\quad \text{as }z\to \infty\text{ in }S^-.\]
\end{lem}

\begin{proof}
Note that
\[\phi(z) = \int_0^{+\infty} {\rm e}^{h(z-t)} {\rm e}^{z-t} {\rm d}t.\]
Thus,
\[\left|\frac{\phi(z)}{\phi'(z)} - 1\right| =  \left|\int_0^{+\infty} {\rm e}^{h(z-t)-h(z)}{\rm e}^{-t}{\rm d}t - 1\right|
 =\left|\int_0^{+\infty} \left({\rm e}^{h(z+t)-h(z)}-1\right){\rm e}^{-t}{\rm d}t\right|.\]
Let $C'$ be the constant provided by \Cref{lem:h} and set $c_z:=\frac{C'}{|z|^\eps}$. Then, for all $t\in \R_+$, we have that $\bigl|h(z-t)-h(z)\bigr|\leq c_z t$. If $|z|$ is sufficiently large so that $c_z<\frac{1}{2}$, then 
\[\left|\frac{\phi(z)}{\phi'(z)} - 1\right| \leq \int_0^{+\infty} ({\rm e}^{c_zt}-1)e^{-t}{\rm d}t = \frac{c_z}{1-c_z} < \frac{c_z}{2} = \frac{C'}{2|z|^\eps}.\]
\end{proof}

Let us now assume that $x_0>0$, $\theta\in(0,\frac{\pi}{2})$ and $U$ is the sector 
\[S^+= \{z\in \C~;~-\theta<{\rm arg}(z-x_0) < \theta\}.\]
Note that the condition on $\theta$ is more restrictive for $S^+$ than for $S^-$.
Choose a point $z_0\in (x_0,+\infty)\subset S^+$ and let $\psi:S^+\to \C$ be the antiderivative of the function $z\mapsto \exp\bigl(z+h(z)\bigr)$ that vanishes at $z_0$.

\begin{lem}\label{lem:sp}
Under the standing assumption on $h'$ and the assumptions above:
\[\left|\frac{\psi(z)}{\psi'(z)} -1\right| =  {\cal O}\left( \frac{1}{|z|^\eps}\right)\quad \text{as } z\to \infty\text{ in }S^+.\]
\end{lem}

\begin{proof}
Note that 
\[\psi(z) = \int_{z_0}^{z} {\rm e}^{h(w)} {\rm e}^{w} {\rm d}w,\]
where the integral is taken along the segment going from $z_0$ to $z$. 
Thus,
\begin{align*}\left|\frac{\psi(z)}{\psi'(z)} - 1\right| &=  \left|(z_0-z)\int_{0}^{1} {\rm e}^{h\bigl(z+t(z_0-z)\bigr)-h(z)}{\rm e}^{t(z_0-z)}{\rm d}t - 1\right|\\
 &=\left|(z_0-z)\int_{0}^1 \left({\rm e}^{h\bigl(z+t(z_0-z)\bigr)-h(z)}-1\right){\rm e}^{t(z_0-z)}{\rm d}t - {\rm e}^{z_0-z}\right|.
 \end{align*}
 Since $\theta<\frac{\pi}{2}$, there exists a constant $c'\in(0,1)$ such that if $z\in S^+$ with $|z|$ is large enough, then $\Re(z-z_0)\geq c'|z-z_0|$. 
Let $C'$ be the constant provided by \Cref{lem:h} and set $c_z:=\frac{C'}{|z|^\eps}$. Choosing $|z|$ larger if necessary, we may assume that $c_z<c'/2$. In that case, 
\begin{align*}
\left|\frac{\psi(z)}{\psi'(z)} - 1\right| &\leq |z_0-z|\int_0^1 \left({\rm e}^{c_z|z_0-z|t}-1\right){\rm e}^{-c'|z_0-z|t} {\rm d}t + {\cal O}\left({\rm e}^{-c'|z_0-z|}\right)\\
& \leq \frac{c_z}{c'-c_z} + {\cal O}\left({\rm e}^{-c'|z_0-z|}\right)\\
 &\leq \frac{2C'}{c'|z|^\eps} + {\cal O}\left({\rm e}^{-c'|z_0-z|}\right) = {\cal O}\left( \frac{1}{|z|^\eps}\right).\qedhere
\end{align*}
\end{proof}

\subsubsection{Injectivity lemma}

\begin{lem}\label{lem:inj}
Let $U$ be a connected open subset of $\C$ and $f:U\to \C$ a holomorphic map.
Assume that:
\begin{itemize}
\item $\exists \tau>0$, $\forall a,b\in U$, there exists a piecewise $C^1$ path from $a$ to $b$ of length $\leq \tau |a-b|$,
\item $\exists K>0$, $|f'-1|\leq K$,
\item $\tau K <1$.
\end{itemize}
Then $f$ is injective.
\end{lem}
\begin{proof}
Let $\gamma$ be a path as in the statement and assume $a\neq b$.
Then $f(b)-f(a) = b-a + \Delta$ with $\Delta = \int_\gamma (f'(z)-1)dz$ so $|\Delta|\leq \tau|a-b| K < |a-b|$.
\end{proof}

\subsubsection{Consequences}\label{ss:csq}

Let $x_0>0$, $\theta\in (0,\frac{5\pi}{6})$ and $S^-$ be the sector bisected by $(-\infty,x_0)$ and of half-opening $\theta$, i.e.:
\[S^-= \{z\in \C~;~\pi-\theta<{\rm arg}(z+x_0) <\pi+ \theta\}.\]

\begin{cor}\label{cor:hphi2}
Assume a holomorphic connection is given on $S^-$, with Christoffel symbol $\Gamma$ satisfying for some $\eps\in(0,1)$:
\[\Gamma(z) = 1 + \cal O(1/|z|^\eps)\]
as $|z|\to\infty$.
Then any developing map $\phi$ satisfies: $\phi(z)$ converges as $z\in S^-$ and $\Re(z)\to -\infty$.
\end{cor}
\begin{proof}
  Apply \Cref{cor:hphi} to the antiderivative $h=\int (\Gamma-1)$ on the simply connected set $S^-$.
  Then $\phi$ is an antiderivative of $z\mapsto \exp(z+h(z))$ (see \cref{eq:intODE}).
\end{proof}

Recall $S^-$ has apex $-x_0$.

\begin{cor}\label{cor:sm2}
Then for every $\eps'>0$ there is some $R_0>0$ such that, on the subset $S^-_{R_0} = S^- - B(-x_0,R_0)$ there exists a global log-chart of the connection defined by $\Gamma$, i.e.\ an isomorphism of affine surfaces: $\Phi$ from $S^-_{R_0}$ to an open subset $U$ of $\cal E$, and satisfying $|\Phi'-1|<\eps'$.
In particular, for all $\theta',\theta''$ such that $0<\theta''<\theta'<\theta$, $U$ contains a sector of apex $-x'$ and of the form
\[ S'^-= \{z\in \C~;\pi-\theta'<{\rm arg}(z+x') <\pi+\theta'\}
\]
for some $x'>0$, whose preimage $\Phi^{-1}(S'^-)$ contains a sector of the form
\[ S''^- = \{z\in \C~;\pi-\theta''<{\rm arg}(z+x'') <\pi+\theta''\} \]
for some $x''>0$.
\end{cor}
\begin{proof}
There exists developing maps $\phi$ on $S^-$ since it is simply connected so, using the same notations as in the proof of \Cref{cor:hphi2}, antiderivatives $\phi = \int \exp(z+h(z))$ are well-defined.
By \Cref{cor:sm2}, we can choose $\phi$ so that it tends to $0$
as $z\in S^-$ and $\Re(z)\to -\infty$.
By \Cref{lem:sm}, $\left|\frac{\phi(z)}{\phi'(z)}-1\right| = \cal O(1/|z|^\eps)$ as $|z|\to \infty \in S^{-}$ (note that the condition on $z$ is broader).
In particular, for some $R_0$ big enough, $\phi$ it does not vanish on the simply connected subset $S^-_{R_0} = S^- - B(0,R_0)$ for some $R_0>0$ big enough.
Let $\Phi$ be a determination of $\log \phi$ on this subset.
Then
\[ |\Phi'-1|=\cal O(1/|z|^\eps)
\]
as $|z|\to \infty$ within $S^-_{R_0}$.
This estimate implies that $\Phi$ is injective for $R_0$ big enough (by \Cref{lem:inj}, on a subset like $S^-_{R_0}$, it is enough that $|\Phi'-1|<1/2\pi$), so it is a log-chart, and that $\Phi(S^-_{R_0})$ contains sectors $S'^-$ as in the statement.
\end{proof}

Let $x_0>0$, $\theta\in(0,\frac{\pi}{2})$ and
\[S^+= \{z\in \C~;~-\theta<{\rm arg}(z-x_0) < \theta\}\]
a sector of apex $x_0$.

\begin{cor}\label{cor:sp2}
Assume a holomorphic connection is given on $S^+$, with Christoffel symbol $\Gamma$ satisfying for some $\eps\in(0,1)$:
\[\Gamma(z) = 1 + \cal O(1/|z|^\eps)\]
as $|z|\to\infty$.
Then for every $\eps'>0$, there exists $R_0>0$ such that, on the subset $S^+_{R_0} = S^+ - B(x_0,R_0)$ there exists a global log-chart of the connection defined by $\Gamma$, i.e.\ an isomorphism of affine surfaces: $\Psi$ from $S^+_{R_0}$ to an open subset $U$ of $\cal E$, and for all $\theta',\theta''$ such that $0<\theta''<\theta'<\theta$, the set $U$ contains a sector of the form
\[ S'^+= \{z\in \C~;-\theta'<{\rm arg}(z-x') <\theta'\}
\]
for some $x'>0$, and $\Psi^{-1}(S'^+)$ contains a sector of the form
\[ S''^+= \{z\in \C~;-\theta''<{\rm arg}(z-x'') <\theta''\}
\]
for some $x''>0$.
\end{cor}
The proof is similar as for \Cref{cor:sm2}, but simpler, using \Cref{lem:sp} in place of \Cref{lem:sm}.

\subsection{Similarities with parabolic points in dynamics}\label{subsub:parabo}

For the reader that knows the classification of parabolic fixed points in one complex dimensional holomorphic dynamical systems, it is interesting to compare with our description of multiple poles of meromorphic connections.

In both cases we have petals and axes.
The analogue of the Fatou coordinates are the log-affine charts.
The repelling axes of a multiple pole $p$ of a connection are similar to repelling axes of a parabolic point $q$, while the attracting axes are similar to attracting axes.
We restrict here to the case where the attracting and repelling petals of $q$ are sectors in Fatou coordinates, of opening angle between $0$ and $2\pi$.
Define \emph{repelling petals} of $p$ as the subsets of a neighborhood of $p$ corresponding to sectors of $A_n$, of opening angle between $0$ and $2\pi$, like $S'^-$ in \Cref{sub:momo}.
Define \emph{attracting petals} of $p$ the subsets of a neighborhood of $p$ corresponding to $B_{n+1/2}$.

\smallskip

\SetTblrInner{hspan=minimal}
\centerline{
\begin{tblr}{width=0.9\textwidth, colspec={|X|X|}}
\hline
Parabolic point $q$ & Irregular singularity $p$
\rule[-1ex]{0pt}{1ex} \\ \hline
any geodesic entering an attracting petal of $p$ remains in it
& any orbit entering an attracting petal of $q$ remains in it
\\ \hline
no orbit starting in a repelling petal of $q$ remains in it in the future
& few geodesics starting from a point in a repelling petal of $p$, or entering it, remain completely in it in the future
\\ \hline
\SetCell[c=2]{l} orbits/geodesics quitting a repelling petal of $q$ either quits the union of the petals or enters an attracting petal (in the second case it never reenters the repelling petal)
\\ \hline
\SetCell[c=2]{l} formal classification strictly coarser than the analytic classification
\\ \hline
\end{tblr}}

\medskip

There are however some differences. 

\smallskip

\centerline{
\begin{tblr}{width=0.9\textwidth, colspec={|X|X|}}
\hline
Parabolic point $q$  & Irregular singularity $p$
\\ \hline
repelling petals of $f$ are attracting petals of $f^{-1}$ and vice-versa
& time reversal does not permute the roles of attracting and repelling petals
\\ \hline
attracting and repelling petals have similar properties as long as one studies local dynamics
& attracting and repelling petals have different properties concerning the local dynamics of their geodesics
\\ \hline
attracting and repelling petals can be extended to ``wide'' petals whose images in Fatou coordinates are unions of an upper, lower and left (or right) half-plane & attracting petals can but repelling ones only can if the two asymptotic values of the two adjacent attracting axes differ
\\ \hline
the analytic invariant lives in an infinite dimensional complex space
& the analytic invariant lives in a finite dimensional complex space
\\ \hline
\end{tblr}}

\medskip

\subsection{The model in affine charts: helixes, foci}\label{sub:helix}

The exponential-affine plane $\cal E$ (where $\Gamma=1$ is constant over $\C$) is the universal covering of the subset $\C^*$ of the plane $\C$ with its canonical affine structure ($\Gamma=0$) via the map 
\[\exp : \cal E \to \C^*
.\]
The imaginary part provides a lift of the argument function on $\C^*$.

Another model of this universal covering is by taking copies $U_n$ of the slit plane $\C^*-(0,+\infty)$, indexed by $n\in\Z$ and gluing consecutive ones along copies of $(0,+\infty)$. Or to stick to a definition where the gluing happens along open subsets, one can start from a finite covering of $\C^*$ with overlapping sectors and take infinitely many copies of them and glue them appropriately.

Another way to look at it is through a horizontal line scanning $\cal E$ from bottom to top.
Its image by $\exp$ is a radial half-line with vertex $0$ that scans the second model by turning infinitely many times around $0$ as its argument varies from $-\infty$ to $\infty$.

We call this the \emph{helix model} of the universal covering.

\begin{figure}
    \centering
    \includegraphics[width=11cm]{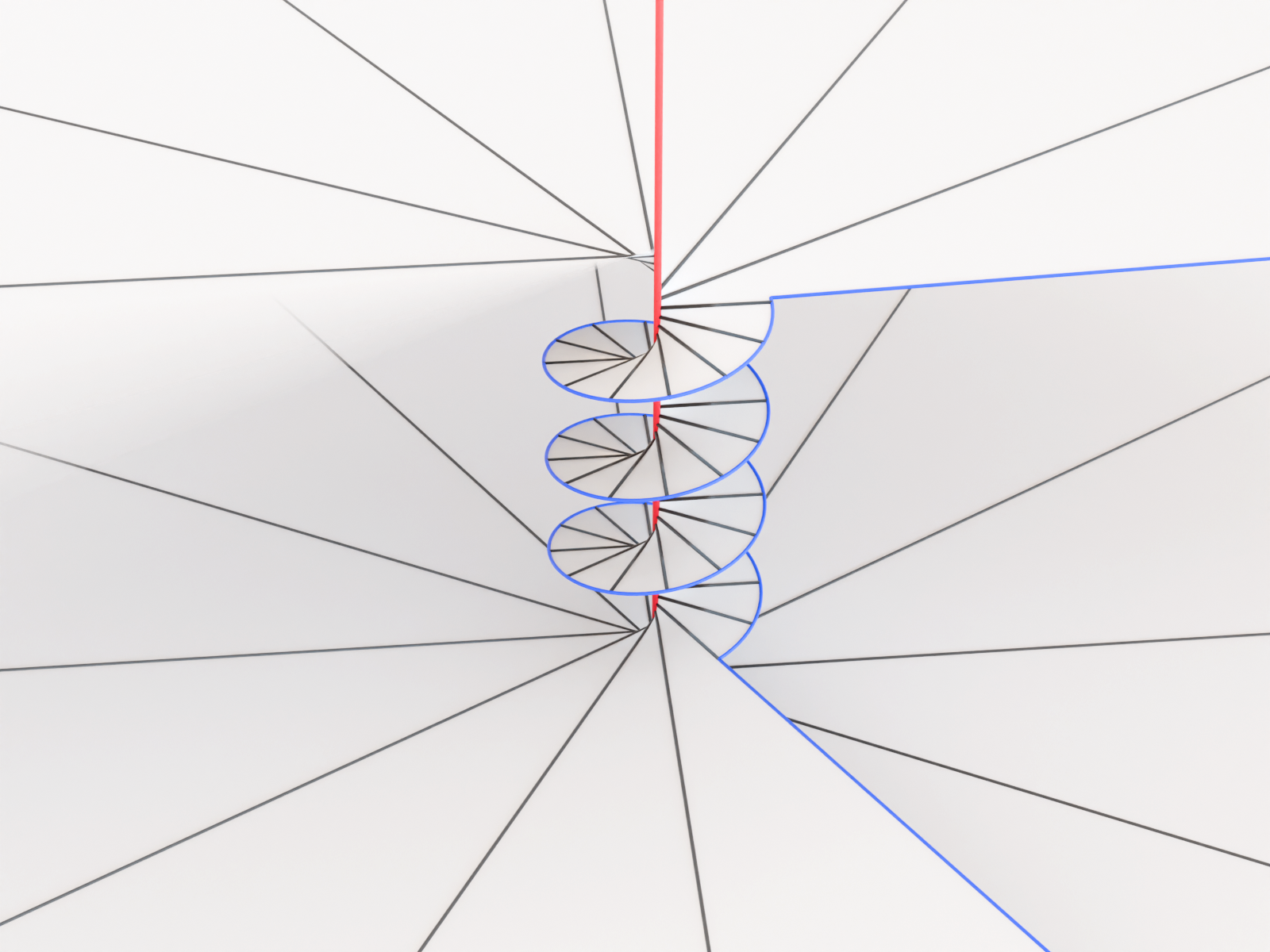}
    \caption{A 3D representation of an $A_n$ helix. The affine structure is the pull-back by the projection to a horizontal plane endowed with the canonical one. We see only a small portion of the part on which the rays  are unbounded: these rays actually make an infinite number of turns around the red vertical axis. The blue line corresponds to the boundary of $A_n$. The vertical axis can be thought of as the focus.}
    \label{fig:enter-label}
\end{figure}

The subsets $A_n$ and $B_{n+1/2}$ of this universal covering, used in the model $\cal M$, can also be described at the level of the helix model as follows:
\begin{itemize}
  \item $B_{n+1/2}$ corresponds to the preimage by $\exp$ of the complement of a closed disk: $\C \setminus \ov B(0,e^R)$. It is spanned by the sub half-line of the scanning half-line that starts from distance $e^R$ from the origin.
  \item The left part $\Hp_{x<-R}$ of $A_{n}$ corresponds to the preimage by $\exp$ of the pointed disk $B(0,e^{-R})-\{0\}$. It is spanned by the subset of the half-line between $0$ and distance $e^{-R}$.
  \item The upper part $\Hp_{y>R}$ corresponds to points of lifted argument $>R$ and the lower part $\Hp_{y<-R}$ of lifted argument $<-R$. They are spanned by the turning half-line where the argument ranges in $(R,+\infty)$ and $(-\infty,-R)$, respectively.
\end{itemize}

\loctitle{Foci}
This point of view motivates the association to $A_n$, $n\in\Z/m\Z$, of a \emph{focus}.\footnote{Alternative names for a focus may include: endpoint, repelling end, tip, apex, cone point, conical point, virtual point, \ldots}
This also holds for the three half-planes $\Hp_{y>R}$, $\Hp_{x<-R}$ and $\Hp_{y<-R}$ in $\cal E$ whose union forms $A_n$. The focus $F_n$ is a virtual point to which the affine atlas does not extend, yet for any germ of affine chart at a point of $A_n$ this germ extends to all $A_n$ into the restriction $\phi$ to $A_n$ of $a\exp+b$ for some $a\in\C^*$ and $b\in\C$ and we extend $\phi$ to $F_n$ by declaring $\phi(F_n) = b$. In particular
\[\exp(F_n)=0.\]
Any horizontal line in $A_n$ is the support of a geodesic going to the left, and is mapped by $\phi$ to a straight line tending to $\phi(F_n)$.

The foci can also be defined at the level of the universal covering, we denote them $\tilde F_n$ where now $n\in\Z$.
The asymptotic sequence is then the image of all the foci $\tilde F_n$, $n\in\Z$ by the developing map $\tilde\phi$, where now the index $n$ needs to be in $\Z$ because the developing map lives on $\wt{\cal M}$, not $\cal M$.

\begin{rem}
It is tempting to associate a focus to $B_{n+1/2}$, mapped to $0$ by $\exp$ as above, but we will not for two reasons: first it lives somewhat away from $B_{n+1/2}$; second $B_{n+1/2}$ contains many subsets affinely isomorphic to $B_{n+1/2}$ but for which the focus is offset by a translation (this is very easy to see with the helix model, and corresponds in $\cal E$ to applying some map in $G_{s,b}$). The fact that the asymptotic sequence is enough to characterise irregular singularities up to isomorphism is also an indication that introducing foci for the $B_{n+1/2}$ is not so interesting.
\end{rem}

\loctitle{Extending the topology at foci and swathes}
Given a meromorphic connection $\nabla$ on a Riemann surface $X$ with an irregular singularity $p$, 
let us explain how, near $p$, we extend the topology on the union of the affine surface $X^*$ and of the the $d-1$ foci.
A basis of open neighborhoods of $F_n$, $n\in\Z/m\Z$ is given by the sets formed by the union of $F_n$ with the subset of the affine surface corresponding to a left half $\mathbb{H}_{x<x_0}\subset A_n$, where $x_0$ varies in $(-\infty,-R)$.
(This topology is the same as the one induced by the metric completion on $A_n$ endowed with the pull-back of the canonical Riemannian metric of $\C$ by $\exp: A_n\to\C$.)

This extended topological space is not anymore a topological surface but it is still path-connected and semi-locally simply connected (it is actually locally contractible), so there is still a notion of universal covering.
Then for $n\in\Z$, $\tilde F_n$ is mapped to $F_{n\bmod m}$.

We will also complete the topological space $X^*$ at irregular points of (polar) degree $d$ by adding $d-1$ other points, that we call \emph{swathes}, one for each attracting petal. In a very similar way as for foci, a neighborhood basis at the swath $s$ associated to $\pi(B_x)$ consists in the union of $\{s\}$ with $\pi(H)$ where $H$ varies in the set of right half planes contained in $B_x$.

The topological space $X$ with all irregular singularities removed and replaced by the union of foci and swathes remains separated in the sense of Hausdorff.

\subsection{Geodesics}\label{sub:geodesics}

\begin{figure}
\begin{tikzpicture}
\node at (0,0) {\includegraphics[scale=0.75]{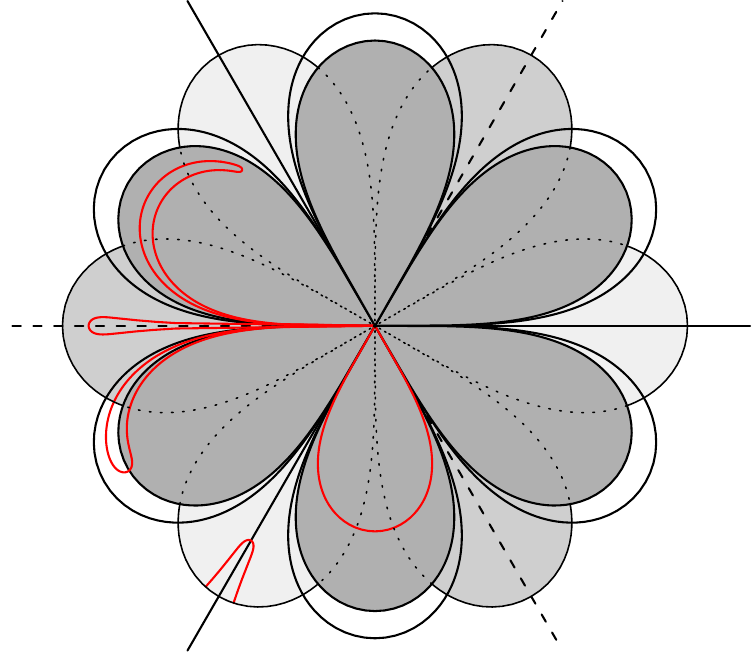}};
\node at (4.5,0.4) {$0$};
\node at (4.5,2) {$S^+_0$} edge (3.15,1.75);
\node at (4.5,2.8) {$S'^+_0$} edge (3.3,2.2);
\node at (4.5,-2) {$S^-_0$} edge (3.15,-1.75);
\node at (4.5,-2.8) {$S'^-_0$} edge (3.3,-2.2);
\node at (3,4) {$1/2$};
\node at (-2.8,4) {$1$};
\node at (-4.5,0.4) {$3/2$};
\node at (-2.8,-4) {$2$};
\node at (3,-4) {$5/2$};
\node(a) at (4.5,-0.5) {$A_0$};
\draw (a) -- (3.5,-0.3);
\draw (a) -- (2.4,-1.45);
\draw (a) -- (2.4,1.45);
\node(b) at (-4.5,-2) {$B_{3/2}$};
\draw (b) -- (-3.6,-0.8);
\end{tikzpicture}
\caption{In red, a few geodesics near an irregular singularity, as described in \Cref{sub:geodesics}. The neighborhood $V$ is the union of all gray sets. Sepals $S^\pm_n$ are in dark gray.
The repelling axes are the solid straight lines while the attracting axes are dashed.
The intersepal spaces of $A_n$ are in lighy gray while the ones of $B_n$ are in medium gray. Part of their outlines is countinued as a dotted curves, to indicate the projection by $\pi$ of $\Hp_{x<-R}$ or $\Hp_{x>R}$.}
\label{fig:sepals}
\end{figure}
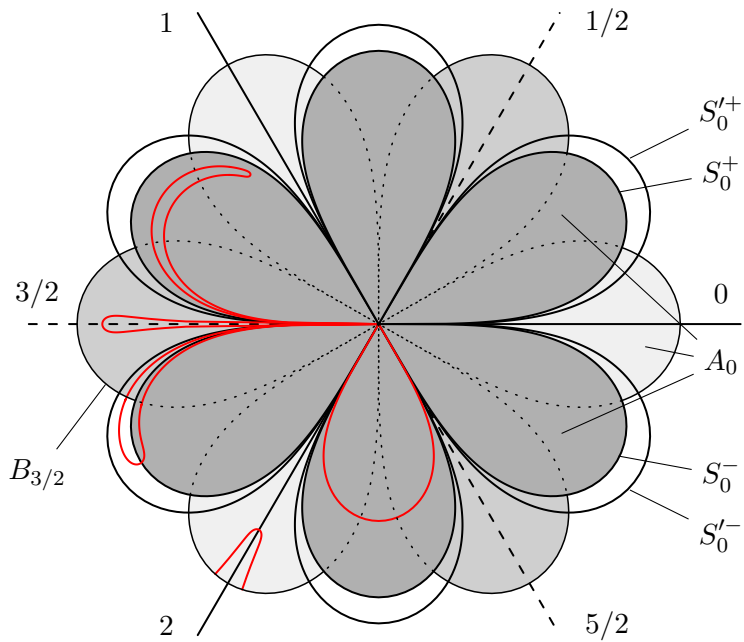

Consider a neighborhood $V$ isomorphic to a model $\cal M$ given by \Cref{thm:model}.
Possibly taking a smaller $V$, we assume that $R>2\pi$. %Cette condition peut n'être pas suffisante à cause de la partie translation des recollements~: à regarder de près..
Consider the projection
\[ \pi : \coprod_{x\in \frac{1}{2}\Z/m\Z} AB_x  \to V
\]
where, as usual, $AB_x$ denotes the set $A_x$ when $x\in\Z/m\Z$ and the set $B_x$ when $x-1/2\in\Z/m\Z$.
Let us first define \emph{sepals} $S^\pm_{n}$ as follows:
they are the $2m$ sets corresponding in $A_n$ to $\Hp_{y>R+\pi}$ for $S^+_n$ and $\Hp_{y<-R-\pi}$ for $S^-_n$. The sepals are pairwise disjoint and have a tangency of order $m+1=d$ with a pair of adjacent repelling and attracting axes (use \Cref{lem:estimate:w,lem:sec2}), there is exactly one sepal for each such pair.
We also consider the $2m$ slightly bigger and still disjoint sepals $S'^\pm_n$, corresponding to $\Hp_{y>R}$ and $\Hp_{y<-R}$ in $A_n$.
See \Cref{fig:sepals}.

The model implies the following properties:
\begin{itemize}
\item Every geodesic starting from $V$ has one of the following behaviors in the future: either it tends to $P$ in finite or infinite time, or it exits $V$.
\item Any geodesic entering the attracting petal $\pi(B_{n+1/2})$ remains in it in the future.
\item From any point $Q$ in the repelling petal $\pi(A_n)$, there is exactly one direction in $T_Q V$ for which the corresponding geodesic will hit $P$ in finite time.
\item Every geodesic starting in a sepal $S=S^\pm_n$ remains in the $S'^\pm_n$ that contains $S$.
\item Every geodesic entering a sepal remains in it forever.
\item No geodesic can visit two different $\pi(A_n)$ or two different $\pi(B_{n+1/2})$ or two different sepals.
\end{itemize}

The $2m$ spaces in $V$ between the $S^\pm_{n}$ are thin, because they are bounded by two curves with a tangency of order $m+1$.
Call them \emph{intersepal} areas.
They consist in images by $\pi$ of parts of $A_n$ (of equation $x<-R$ and $|y|\leq R+1$) and parts of $B_{n+1/2}$ (bounded by geodesics, which are not necessarily straight in the log-coordinate), and both parts have bounded height.
A geodesic starting in this part of $A_n$ is likely to escape $V$: actually if the starting point is at distance $\geq 2\pi$ of the boundary of the part, then unless it is horizontal and points to the left, it will escape $V$.
A geodesic starting in the intersepal area of $B_{n+1/2}$ can escape $V$ but the set of angles that allow this decrease as the real part of the initial point in $B_{n+1/2}$ gets bigger.

One of the conclusions to draw is that a lot of geodesics that start from $V$ never exit.

\medskip

Of the observations above, we single out the following ones, for easier later reference.
\begin{lem}\label{lem:irregularTips}
\ \begin{enumerate}
  \item Any geodesic tending in finite time to an irregular singularity is eventually contained in one of the $d-1$ repelling petals $P=\pi(A_n)$ for some $n\in\Z/m\Z$ and tends to the associated focus $F_n$.
  \item Any geodesic tending in infinite time to an irregular singularity is eventually contained in one of the $d-1$ attracting petals $P=\pi(B_{x})$ with $x-1/2\in\Z/m\Z$ and it tends to the associated swath.
\end{enumerate}
\end{lem}

\section{Geometric constructions from local models}\label{sec:ConsequencesLocalModels}

We begin by drawing upon the classical local models of Fuchsian singularities and our models from \Cref{sub:model} for irregular singularities, to derive some necessary technical results before turning to the main constructions.
First, \Cref{lem:fuchsianTips,lem:irregularTips} imply:

\begin{lem}\label{lem:tips}
Let $(X,\nabla,\cal S)$ be a  finite type affine surface.
Consider a geodesic that converges in finite time to a singularity. Then this singularity is either
\begin{itemize}
\item Fuchsian conical or 
\item irregular, and the geodesic converges to one of its foci through a repelling petal.
\end{itemize}
Consider a geodesic that converges in infinite time to a singularity. Then this singularity is either
\begin{itemize}
\item Fuchsian anti-conical,
\item Fuchsian cylindrical or 
\item irregular, and the geodesic converges to one of its swathes through an attracting petals.
\end{itemize}
\end{lem}

Given the first part of the statement above, it makes sense to group the two possibilities under a common terminology:

\begin{defn}\label{def:apex}
  We call \emph{apexes} the conical Fuchsian singularities of $\nabla$ (they form a subset of $X-\cal S$) and the foci of irregular singularities.
  We denote $\hat X$ the union of $X-\cal S$ and the apexes, endowed with the topology defined in \Cref{sub:helix}.
\end{defn}

This term apex comes from the name of the tip of cones in classical geometry. Note that marked regular points of $\nabla$ are in the set of apexes.

\medskip

Recall that a \emph{saddle connection} is a geodesic in $X-\cal S$ that tends to singularities in finite time in the future and in the past. By \Cref{lem:tips} it necessarily connects two apexes or one apex to itself.

\subsection{System of arcs and saddle connections}\label{sub:apicalSystems}

We first provide a purely topological lemma concerning systems of arcs connecting apexes on finite type affine surfaces. This result is one of the main ingredients in the complexity bound for Delaunay decompositions that we will establish in \Cref{thm:ComplexityBound}.

To motivate it, let us first mention an easy and well-known consequence of Euler's characteristic identity: given a triangulation of a genus $g$ compact orientable surface, the number $e$ of edges and the number $v$ of vertices are related by $e=6g-6+v$. To prove it, use Euler's identity $2-2g = v-e+t$ where $t$ is the number of triangles and note that $2e=3t$.
If instead we have a polygonal decomposition and every polygon has at least 3 sides then $e\leq 6g-6-v$.
See \Cref{sub:sides} for our convention about side counting (a single boundary edge is counted once or twice depending on the situation).

In this section, we consider a finite type affine surface $(X,\nabla,\mathcal{S})$.
We recall that swathes, together with a neighborhood basis, have been defined in \Cref{sub:helix}.
Recall that an embedded open arc in a topological space is a continuous injection of $(0,1)$.

\begin{defn}\label{def:aas}
  We call \emph{apical arc system} a collection $\cal A$ of embedded open arcs in $X-\cal S$ with disjoint interiors such that every arc tends on both ends to an apex of $\hat X$.
  The \emph{support} of $\cal A$ is the union of the closure of all arcs in $X$; it is a subset of $X$.
  We call \emph{faces} the connected components of the complement in $X$ of the support. They may or may not contain singularities.
  A face is called
  \begin{itemize}
  \item \emph{ordinary} if, as an open set, it is a topological disk and does not contain any interior singularity nor any swath neighborhood;
  \item \emph{special} if, as an open set, it is a topological disk and contains a neighborhood of a unique swath or simple pole that is not a conical singularity;
  \item \emph{defective} if it is neither ordinary nor special.
  \end{itemize}
  The side count of a face is the number of edge sides bounding it: a repeated arc is thus counted twice, see \Cref{fig:sides}.
\end{defn}

We recall that conical Fuchsian singularities include marked points.

\begin{lem}\label{lem:SystemTopologicalArcs}
Denote by $g$ be the genus of $X$ and $n$ the number of singular points in $\mathcal{S}$ \emph{counted with multiplicity}.
Consider an apical arc system $\mathcal{A}$ as in \Cref{def:aas} satisfying the following supplementary conditions:
\begin{itemize}
  \item each apex is the end of at least one arc of $\mathcal{A}$,
  \item there is no defective face.
\end{itemize}
Since all faces are ordinary or special, they are in particular topological disks.
Denoting $o$ the number of ordinary faces, the number $|\mathcal{A}|$ of arcs in $\mathcal{A}$ satisfies:
\begin{equation}\label{eq:first}
|\mathcal{A}| = 2g - 2 + n + o.
\end{equation}
If moreover every ordinary face is 3-sided (with possibly a repeated side), we can express this identity in terms of the the total number $\sigma$ of boundary sides of special faces, instead of $o$:
\begin{equation}\label{eq:second}
|\mathcal{A}| = 6g - 6 + 3n - \sigma.
\end{equation}
\end{lem}

\begin{proof}
We assume that $\mathcal{S}$ contains:
\begin{itemize}
    \item $n_{1}$ conical Fuchsian singularities;
    \item $n_{2}$ other Fuchsian singularities;
    \item $k$ multiple poles of $\nabla$ of orders $d_{1},\dots,d_{k} \geq 2$.
\end{itemize}
We have $n=n_{1}+n_{2}+\sum\limits_{i=1}^{k} d_{i}$.
\par
The arcs of $\mathcal{A}$ define a decomposition of the topological surface $X$ into:
\begin{itemize}
    \item $n_{1}+k$ vertices;
    \item $o$ ordinary faces;
    \item $n_{2}+\sum\limits_{i=1}^{k} (d_{i}-1)$ special faces;
    \item $|\mathcal{A}|$ edges.
\end{itemize}
Computing the Euler characteristic of $X$ with respect to this decomposition, we obtain:
$$
2-2g
=
(n_{1}+k)+o+(n_{2}+\sum\limits_{i=1}^{k} (d_{i}-1))-|\mathcal{A}|.
$$
In other words, $|\mathcal{A}| = 2g - 2 + n + o$, which is \cref{eq:first}.
\par
If each ordinary face is $3$-sided we get $2|\mathcal{A}| = \sigma+3t$, so $o=\frac{2}{3}|\mathcal{A}|-\frac{1}{3}\sigma$. Together with \cref{eq:first}, this implies \cref{eq:third}.
\end{proof}

In the next result we relax the conditions on $\cal A$, and get an inequality instead of an identity: 

\begin{lem}\label{lem:SystemTopologicalArcs2}
Denote by $g$ be the genus of $X$ and $n$ the number of singular points in $\mathcal{S}$ counted with multiplicity.
Consider an apical arc system $\mathcal{A}$ as in \Cref{def:aas} with the supplementary conditions:
\begin{itemize}
\item $\cal A$ is non-empty;
\item every ordinary face has at least three sides.
\end{itemize}
Then number $|\mathcal{A}|$ of arcs in $\mathcal{A}$ satisfies:
\begin{equation}\label{eq:third}
|\mathcal{A}| \leq 6g - 6 + 3n.
\end{equation}
\end{lem}

\begin{proof}
We assume that $\mathcal{S}$ contains:
\begin{itemize}
    \item $n_{1}$ conical Fuchsian singularities;
    \item $n_{2}$ other Fuchsian singularities;
    \item $k$ multiple poles of $\nabla$ of orders $d_{1},\dots,d_{k} \geq 2$.
\end{itemize}
We have $n=n_{1}+n_{2}+\sum\limits_{i=1}^{k} d_{i}$.
Write $n_1 = a_0'+n_1'$ where $a_0'$ is the number of conical Fuchsian singularities that are also endpoints of  arcs in $\cal A$.
The arcs of $\mathcal{A}$ define a decomposition of the topological surface $X$ into:
\begin{itemize}
    \item $a_{0}=a_0'+a_0''$ vertices, where $a_{0}''$ is the number of irregular singularities that are endpoints of arcs of $\cal A$;
    \item $a_1=|\mathcal{A}|$ edges.
    \item $a_{2} = a_2^o+a_2^n+a_2^r$ faces, where $a_2^o$ counts the ordinary faces where, $a_{2}^n$ the faces that are topological disks but not ordinary, and $a_{2}^r$ the faces that are not topological disks;
\end{itemize}
Faces counted by $a_2^n$ contain a Fuchsian singularity or a swath neighborhood (if it contains an irregular singularity then it contains a swath neighborhood) and a multiple pole of order $d_i$ has $d_i-1$ swathes, in particular:
$$
a_2^n\leq n'_1+n_2+\sum_{i=1}^k (d_i-1)
.
$$
The Euler characteristic of $X$ satisfies
$$
\chi(X) = 2-2g
= a_0 - a_1+ \sum_F \chi(F)
$$
where the sum is over all faces $F$, and $\chi(F)$ is the Euler characteristic of $F$ (as an open set).
So
\[ |\cal A| = 2g-2+a_0+\sum_F \chi(F)
.\]
Since $\cal A$ is non-empty, no face is a sphere, hence the faces homeomorphic to disks satisfy $\chi(F)= 1$ while the others $\chi(F)\leq 0$.
We thus get
\[|\mathcal{A}| \leq 2g - 2 + a_0 + a_2^o+a_2^n.\] 
Since each ordinary face is at least $3$-sided we get $2|\mathcal{A}|\geq 3a_2^o$, so $a_2^o\leq\frac{2}{3}|\mathcal{A}|$, whence $|\mathcal{A}| \leq 2g - 2 + a_0 + \frac{2}{3}|\cal A|+a_2^n$ i.e.
\[ |\mathcal{A}| \leq 6g - 6 + 3(a_0  + a_2^n).
\]
Finally,
$a_0+a_2^n = a_0''+a_0'+a_2^n \leq a_0''+a_0'+n_1'+n_2+\sum_{i=1}^k (d_i-1) = (a_0''-k)+n_1+n_2+\sum_{i=1}^k d_i = (a_0''-k) +n \leq n$.
\end{proof}

\begin{rem}
Observe in particular that \Cref{eq:third} implies $2g+n > 2$ ($\cal A$ is non-empty by hypothesis).
\end{rem}

We deduce an upper bound on the maximal number of non-intersecting saddle connections.
We recall that a saddle connection is a geodesic connecting two (possibly identical) singular points in finite time.
We call a saddle connection \emph{simple} if it is injective, except that the two ends may or may not coincide.

\begin{cor}\label{cor:UpperBoundSaddleConnections}
Consider a finite type affine surface $(X,\nabla,\mathcal{S})$ where $X$ is a compact Riemann surface of genus $g$ while $\mathcal{S}$ has $n$ singular points (counted with multiplicity).
\begin{itemize}
\item If $2g+n \geq 3$, any family $\mathcal{F}$ of \emph{simple} saddle connections with disjoint interiors has at most $6g-6+3n$ elements.
\item If $2g+n \leq 2$, then $(X,\nabla,\mathcal{S})$ does not admit any saddle connection, simple or not.
\end{itemize}
\end{cor}

\begin{proof}
The case $2g+n \leq 2$ is handled by \Cref{prop:LowComplexity}, so we assume that $2g+n \geq 3$.
Then $\cal F$ is an apical arc system.
Without loss of generality we may also assume that $\cal F$ is finite.
Then $\cal F$ is an apical arc system.
If $\cal F$ is  then $|\cal F| = 0 $ while $6g-6+3n \geq 3$ so the inequality holds. In the rest of the proof, we assume $\cal F$ non-empty.

We claim that the system $\cal F$ then satisfies the assumptions of \Cref{lem:SystemTopologicalArcs2}.
Otherwise, a developing map of the affine structure would map any ordinary face with less than 3 sides to a connected open subset $U$ of $\C$ of finite area and whose boundary is contained in the union of two straight segments. There is no such domain $U$.

Then \Cref{lem:SystemTopologicalArcs2} gives the result.
\end{proof}

\begin{lem}
If $2g+n\geq 3$ then $(X,\nabla,\mathcal{S})$ has at least one apex.
\end{lem}
\begin{proof}
If $\nabla$ has a multiple pole then it correspond to an irregular singularity, which has at least one focus, so $X$ has an apex.
If every pole of $\nabla$ is simple.
Recall (see \Cref{sub:Fuchsian}) that a simple pole is a conical singularity if and only if its residue $\res$ satisfies $\Re\res<1$.
From $2g+n \geq 3$ we get three cases:
\begin{itemize}
    \item (a): $g \geq 2$ and at least one pole is a conical singularity;
    \item (b): $g=1$ and $n\geq 1$ so at least one singular point (possibly a marked point) is a conical singularity;
    \item (c): $g=0$, $n \geq 3$, the total residue is $2$ so one of them has a real part of at most $\frac{2}{3}$ and determines a conical singularity.
\end{itemize}
This proves the claim.
\end{proof}

\subsection{Neighborhoods of geodesics}

Before entering in the subject, we prove a preliminary lemma that will be useful here and elsewhere.
In this section we consider a Riemann surface $S$ with a meromorphic connection $\nabla$. By an affine immersion we mean a locally affine map from an open subset of an affine surface and taking values in the complement of the poles of $S$.

\begin{lem}\label{lem:SlitDisk0}
Let $S$ be a Riemann surface with a meromorphic connection $\nabla$.
Let $p\in S$ be (1)~a regular point, (2)~a conical singularity or (3)~a focus of an irregular singularity.
For any $R>0$ and any neighborhood $V'$ of $p$, there exists a punctured neighborhood $V\subset V'$ of $p$ in cases~(1) and~(2), a repelling petal neighborhood of the focus in case~(3), such that $\forall x\in V$ in cases~(1) and~(2) and $\forall x\in V\cap P$ in case~(3), there is an affine immersion (possibly injective) $f$ from the slit disk
\[ B(0,R)-(-R,0] = \dom(f)
\]
to $V'$, sending $1$ to $x$ and such that $f(z)$ tends to $p$ as $z$ tends to $0$ within $\dom(f)$ 
\end{lem}

\begin{proof}
For a regular singularity the claim is very easy and we do not need the slit. Below we assume $p$ is not regular.

Recall that the exponential-affine plane has two models, one as a universal covering ${\wt\C}^*$ of $\C^*$ endowed with the restriction of the canonical affine structure of $\C$, in which we use name the coordinate $u$. 
One as $\cal E=\C$ where we name the coordinate $w$, endowed with the connection with $\Gamma_w=1$, or equivalently for which the map $w\mapsto u=\exp(w)$ is affine.

For a conical singularity, a neighborhood $V_0$ of $p$ is isomorphic, for some $\tau\in\C$ with $\Im \tau>0$ to the quotient (see \Cref{sub:Fuchsian,sub:ExpoAffine}) of the subset $H$ of $\cal E$ of equation $\Im(w/\tau)>0$ by the group of translations generated by $T_\tau:w\mapsto w+\tau$.
Being close to $p$ means $\Im(w/\tau)$ is high.
By shifting the boundary of $H$ to $\Im(w/\tau)=h_0>0$ we may assume that $V_0\subset V'$.
For any $R>0$, there exists $h>0$ such that for all $w_0\in \C$ with $\Im(w_0/\tau)>h+h_0$, then $H$ contains the half-strip $S$ of equation $\Re (w-w_0)< + \log R$ and $|\Im (w-w_0)|<\pi$.
The map $w\mapsto \exp(w-w_0)$ sends bijectively $S$ to the the slit disk $B(0,R)-(-R,0]$ and its inverse, followed by $\langle T_\tau \rangle$ and the isomorphism to a neighborhood of $p$, is the sought for affine immersion.

For an irregular singularity, the repelling petal in which the geodesic enters is isomorphic to a left half-plane in $\cal E$ and the proof is the same as above (without quotienting by a translation).
\end{proof}

In an affine surface, a geodesic can be viewed as an $\R$-affine immersion of an interval $I\subset\R$. In the next lemmas, when these trajectories encounter singularities, we draw upon the local models to extend the immersions to a suitable open subset of $\C$ containing $I$.

\begin{lem}\label{lem:SlitDisk}
Given a meromorphic connection $\nabla$ on a Riemann surface $S$, consider a geodesic arc $\gamma$ defined on $(0,1)$ and such that $\gamma(t)$ converges to a point $p\in S$ as $t \tends 0$. Then, $\gamma$ extends to the affine immersion of a slit disk
\[C := \lbrace{ z \in \mathbb{C}\,;\,0<|z|<\epsilon,\ -\pi<\arg(z)<\pi\rbrace}\]
for some $\epsilon >0$.
This immersion converges to $p$ at $z\tends 0$.
\end{lem}
\begin{proof}
By \Cref{lem:tips}, $p$ is either regular, conical Fuchsian or irregular, in which case $\gamma$ tends to one of the foci of $p$ and let $p$ denote this focus instead.
Apply \Cref{lem:SlitDisk0} to this $p$ and $R=2$: we get a neighborhood $V$ with some property.
There exists $t_0\in(0,1)$ such that $\forall t\geq t_0$, $\gamma(t)\in V$.
The property gives an affine immersion $f:\dom(f) = B(0,2)-(-2,0] \to S$ sending $1$ to $\gamma(t_0)$ and such that $f(z)$ tends to $p$ as $z$ tends to $0$ within $\dom(f)$.
On the straight segment $[t_0,0)$, the map $f_0:z\mapsto f(z/t_0)$ defines a geodesic from $\gamma(t_0)$ to $p$. By \Cref{lem:uniqGeodConical,sub:geodesics}, it must coincide with $\gamma$ on $[t_0,0)$.
The map $f_0$ provides the sought for immersion.
\end{proof}

Note that the extension is not necessarily injective in the Fuchsian conical case.

\begin{lem}\label{lem:DeformCone}
Given a meromorphic connection on a Riemann surface $S$, consider an affine immersion $\gamma:\mathbb{R}_{>0} \to S$ such that $\gamma(t)$ converges to some $p\in S$ as $t \to +\infty$.
Then $p$ is a pole of $\nabla$ and $f$ extends to the affine immersion of an unbounded sector \[C:=\lbrace{ z \in \mathbb{C}\,;\,|z|>R,\ -\epsilon<\arg(z)<\epsilon\rbrace}\]
for some $R>0$ and $\epsilon>0$.
If $p$ is not a cylindrical Fuchsian singularity, we can take $\eps=\pi$.
Moreover, all the geodesics $t \mapsto \gamma(t e^{i\theta})$ for $-\epsilon < \theta< \epsilon$ converge to $p$ as $t \tends + \infty$ (to the same swath if $p$ is is a multiple pole).
\end{lem}

\begin{proof}
By \Cref{lem:fuchsianTips}, the limit $p$ of $\gamma(t)$ can only be an irregular singularity, or a Fuchsian one, either cylindrical or anti-conical.

In the irregular case, \Cref{lem:irregularTips} tells us the geodesic ends up, say for all $t>t_0$, in an attracting petal, which is isomorphic to the universal covering $\tilde U$ of $U$ where $U$ is complement in the canonical affine plane $\C$ of a closed disk of center $0$.
In $U$, the geodesic takes the form for $t>t_0$ of an affinely parametrized straight line $at+b$, and for $\eps$ small enough and $R$ big enough, the set $C$ of the statement, with $\eps=\pi$, is mapped injectively by $z\mapsto az+b$ to a single sheet of the universal covering $\tilde U\to U$.

In the cylindrical Fuchsian case, we end up in a model $\mathbb{H}/\Z$ with an affine line $t>t_0 \mapsto at+b$ whose imaginary part tends to $\infty$ and the argument is similar, but we need the image of $C$ by $z\mapsto az+b$ to remain contained in $\mathbb{H}$, so $\eps$ may need to be small.

In the anti-conical Fuchsian case, we distinguish two cases: in the pure case, we end up as in the proof of \Cref{lem:SlitDisk} in a model $H/\langle T_t\rangle$ in the exponential-affine plane $\cal E$ with this time $\Im t<0$. In the coordinate $u=\exp w$, the geodesic take the form $t\mapsto at+b$.
For $R$ big enough, the map $z\mapsto az+b$ lifts on the set $C$ with $\eps=\pi$ to a map $C\to H$ which projects under the quotient by $T_t$ to an affine extension of the geodesic.

In the shifted anti-conical case, given any geodesic $t\mapsto at+b$ in the model, we can embedded the image of $C$ with $\eps=\pi$ by the extension $z\mapsto az+b$ provided $R$ big enough.
\end{proof}

Actually in the non-cylindrical case, we can even take $\eps>\pi$ (a sector of opening angle bigger than a full turn, so that does not embed in $\C$).

\subsection{Cylinders and closed geodesics}\label{subsub:CylindersAroundGeodesics}

Given an affine surface $\cal A$, we define a \emph{closed geodesic} as a geodesic that eventually comes back to its initial point with the same initial direction, but possibly a different speed.

Consider a closed geodesic $\gamma:I\to \cal A$, with $[a,b]\subset I$, $I\subset\R$ an \emph{open} interval, such that, at $a$ and $b$, $\gamma$ has the same value and direction.
Since the vectors $\gamma'(a)$ and $\gamma'(b)$ in the tangent space to $X^*$ at $\gamma(a)=\gamma(b)$ are related by $\gamma'(b)=s\gamma'(a)$ for some real $s>0$. 

Note: the holonomy of the parallel transport along $\gamma$ from $a$ to $b$ is the multiplication by $s$.

We saw in \Cref{sub:cylindersAndSkewCones} that the supports of all closed geodesics of translation/dilation cylinders define a horizontal/radial foliation.

\begin{lem}\label{lem:closedGeodInCylinder}
There exists an immersed translation cylinder if $s=1$, an immersed dilation cylinder if $s\neq 1$, such that $\gamma$ is a closed geodesic of this cylinder.
More precisely $\gamma = \psi\circ\tilde\gamma$ where $\psi:C\to \cal A$ is an affine immersion, $C$ is the cylinder, $\tilde \gamma$ is a closed geodesic of $C$.
\end{lem}
\begin{proof}
By uniqueness of solution of regular O.D.E.\ (Cauchy-Lipschitz), 
$\gamma\circ \phi=\gamma$ holds on $I\cap \phi^{-1}(I)$, which contains $b$, where $\phi(t)= s(t-b) + a$, $\phi\in\AutC$.
If $s\neq 1$ the map $\phi$ has a unique fixed point $z_0$, it is real and cannot belong to $I$ (geodesics are locally injective).
In particular, the sets $\phi^n([a,b])$ for $n\in\Z$ cover an open set containing $I$ and this is the case too if $s=1$ as then they cover $\R$.

Consider an affine extension $\hat\gamma: U\to X^*$ where $U$ is an open subset of $\C$ containing $I$.
If $s=1$ let $S$ be an open rectangle
\[ \setof{z\in\C}{\Re(z)\in (a-\eps,b+\eps) ,\ \Im z\in(-h,h)}
\]
contained in $U$.
If $s\neq 1$ let $z_0$ be the fixed point of $\phi$ and $S_0$ a set of the form
\[ \setof{z\in\C}{\Re(z)\in [a-\eps,b+\eps],\ \arg(z-z_0)\in (-\alpha,\alpha)}
\]
if $z_0>b$ or 
\[ \setof{z\in\C}{\Re(z)\in [a-\eps,b+\eps],\ \arg(z-z_0)\in (\pi-\alpha,\pi+\alpha)}
\]
if $z_0<a$.
Let $C$ be the affine surface defined as the quotient of $S$ by the relation $z\sim \phi(z)$ and denote $\pi : S\to C$ the quotient map.
The relation $\hat\gamma\circ\phi = \hat \gamma$ holds on $S\cap \phi^{-1}(S)$ by the affine identity theorem.
We can thus factor $\hat\gamma = \psi\circ\pi$.
We define $\tilde \gamma(t)$, first for $t\in [a-\eps,b+\eps]$ as $\pi(t)$, and extend it to $I$ 
is a geodesic of $C$ using $\phi$. 
\end{proof}

So $\gamma$ embeds into a one-parameter family of closed geodesics (up to reparametization) with the same holonomy factor $s$.

The closed geodesics of a cylinder can be parametrized by an interval in a way that is meaningful in affine geometry. In the case of a translation cylinder, this corresponds to the parametrization of a geodesic transverse to all the geodesics in the cylinder. In the case of a dilation cylinder, it corresponds to the angular parametrization of one of the concentric circular arcs of the annulus. In both cases, the interval can be finite, semi-infinite or infinite. If the one-parameter family of closed geodesics cannot be extended further, we say that the cylinder is \textit{maximal}.
\par
An infinite end of a cylinder corresponds to a Fuchsian singularity of residue $\res$ satisfying $\Re(\res)=1$, see \Cref{sub:Fuchsian}. We have $\Im(\res)=0$ or not depending on whether the cylinder is a translation or dilation cylinder. Finite ends of cylinders will be described in \Cref{sub:diskAndCylBoundary}.

Examples of self-intersecting attracting closed geodesics forming an immersed dilation cylinder are given in Section~3.3 and Figure~2 of \cite{NST23}.

\subsection{Anti-conical and swath domains}\label{subsub:UnboundedDomains}

Consider a finite type affine surface $(X,\nabla,\cal S)$.
We recall that saddle connections are finite length geodesics between points of $\cal S$.

\begin{defn}\label{defn:Unbounded}
An \textit{antic-conical domain} is a topological open disk $F$ in $X$, punctured at a unique point, which is an anti-conical Fuchsian singularity (i.e.\ of residue $\res$ satisfying $\Re(\res)>1$) and such that $F$ is bounded by simple saddle connections that do not cross each other, is $m$-sided for some $m\geq 1$ and such that its $m$ corners have angles 
$\theta_{1},\dots,\theta_{m} \geq \pi$.
\par
A \textit{swath domain} is a topological open disk $F$ containing no singularity, bounded by simple saddle connections that do not cross each others, such that $F$ is $m$-sided for some $m\geq 1$, with $m-1$ corners having angles $\theta_{1},\dots,\theta_{m-1} \geq \pi$ while the remaining corner is an irregular singularity, of which a unique swath has a neighborhood included in $F$ (it makes sense to declare that this corner has infinite angle).
\end{defn}

%The reason why these domains are qualified Delaunay will be clarified in \Cref{lem:uddadc}.

\begin{rem}
In the case of anti-conical domains, we have
\[\sum\limits_{k=1}^{m} (\theta_{k}-\pi) = 2\pi(\Re(\res)-1).\]
as follows from \Cref{lem:GB}.
\end{rem}

\begin{lem}\label{lem:Unbounded}
In an anti-conical or swath domain, any geodesic starting at a point of the boundary and pointing towards the interior converges in infinite time to the associated anti-conical Fuchsian singularity or swath.
\end{lem}

We have in particular a form of concavity: no inner geodesic can link two points of the boundary.

\begin{proof}
Orient the sides of the face $F$ counterclockwise (two sides on the same edge give the latter two opposite orientations).

\begin{figure}[htbp]
\begin{center}
\begin{tikzpicture}
  \node at (0,0) {\includegraphics[scale=0.66]{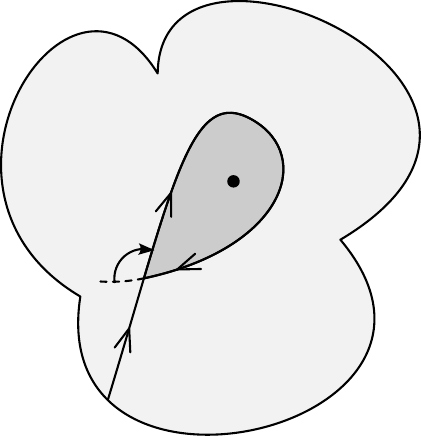}};
  \node at (-1,-0.2) {$\alpha$};
  \node(A) at (2.5,-0.5) {$U$};
  \draw (A) -- (0.3,-0.3);
  \node at (0.5,-2) {$F$};
  \node at (0,-3.5) {\begin{minipage}{5.5cm}\small On this example, the turning angle $\alpha$ and the winding number $s$ are both negative.\end{minipage}};
  \node at (7,0) {\includegraphics[scale=0.66]{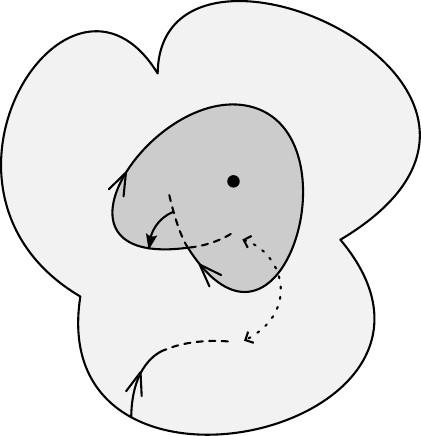}};
  \node at (7,-3.5) {\begin{minipage}{5.5cm}\small If the signs of $\alpha$ and $s$ do not match, then following the geodesic back in time up to $\partial F$ we must cross $\partial U$ again.\end{minipage}};
\end{tikzpicture}
\end{center}
\caption{}\label{fig:alpha}
\end{figure}

Consider a geodesic $\gamma:I\to X$ defined an interval $I$ such that $0\in I\subset [0,+\infty)$, $\gamma(\on{Int})\subset F$ and $\gamma(0)\in\partial F$. We first prove that $\gamma$ cannot self-intersect.
Otherwise, consider the minimal $a$ such that there is another $b>a$ for which $\gamma(a)=\gamma(b)$. We also take $b$ minimal.
The simple closed loop $\gamma|_{[a,b]}$ cuts a topological disk $U$ out of the face.
Let $s = 1$ if the restricted geodesic follows $\partial U$ it the positive (counterclockwise) orientation and $s=-1$ if it is  the negative orientation: $s$ is a also the winding number of the restricted geodesic in a topological chart sending $F$ to a disk in the plane).
Let $\alpha$ be the turning angle of the geodesic at its origin. Note that $|\alpha|\in(0,\pi)$. 
The fact that the geodesic starts on $\partial F$ implies that $\alpha$ has the same sign as $s$. Indeed, if the signs did not match, following the geodesic in the past we would enter $U$ and would not be able to exit: see \Cref{fig:alpha}).
If $U$ does not contain a singularity, \Cref{lem:GB} implies
$\alpha = s 2\pi$, which contradicts $|\alpha|<\pi$.
If $U$ contains the anti-conical point, \Cref{lem:GB} implies $\alpha= s 2\pi (1-\Re(\res))$ but the signs do not match since $1-\Re(\res)<0$.

We then prove that the geodesic cannot exit.
This geodesic would cut the face into two topological disks, one of them, call it $F'$, does not contain the singularity or some swath neighborhood.
Denote $\alpha,\beta\in(0,\pi)$ the angle of $F'$ at the ingress and egress points of the geodesic.
\Cref{lem:GB} implies
\[
2\pi = \pi-\alpha+\pi-\beta + \sum_{i\in I} \pi-\theta_i
\]
where $I\subset\{1,\ldots,m\}$ is the set of corners of $F$ still in $F'$.
However, $\theta_i\geq \pi$ so the right hand side is $<2\pi$, leading to a contradiction.

Similarly: (*) \emph{the geodesic cannot hit in finite time $\partial F$ on an apex (conical point or focus of an irregular singularity)}.
Indeed the same contradiction would hold with the only difference that $\beta$ may or may not be $\geq \pi$. 

An anti-conical singularity and a swath have a trapping neighborhood (see \Cref{lem:traps,sub:geodesics}), such that every geodesic entering it necessarily stays, and tends to the singularity or swath.
By way of contradiction we assume that geodesic $\gamma$ never enters it.
Consider the maximal interval of the form $(0, t_{\max})$ on which it is contained in $F$, with $t_{\max} \in (0,+\infty]$.
Consider the accumulation set $K$ of $\gamma$ at $t_{\max}$.
It is a compact and connected subset of $\ov F$.
Hence either
\begin{enumerate}
\item $K$ contains a non-singular point of $X$;
\item $K$ is reduced to a singular point.
\end{enumerate}

In the first case, using injectivity of $\gamma$ on $[0,t_{\max})$, we get again a contradiction from \Cref{lem:GB}, as follows.
Consider a local affine chart of the accumulation point whose image is a disk centered on the point.
The traces of the geodesic in this disk are a sequence of straight chords that accumulate $0$. The direction, taken modulo $\pi$, of such chords must converge. For any $\eps>0$ we can find two consecutive chords on which the respective geodesic speed vectors make an angle $<\eps$ modulo $\pi$.
By taking a short transversal segment between them we obtain a closed loop made of two pieces of geodesic and total turning number close to $0$ or to $s/2$, where $s$ is as a few paragraphs above. It cannot be $-s/2$ for the same reason as previously.
This contradicts \Cref{lem:GB} both in the case it contains an anti-conical singularity (the turning number must be $(1-\Re(\res)s$), and in the case it does not (the turning number must be $s$).

In the second case the geodesic tends to the singular point $p$, recall that: either $F$ contains an anti-conical singularity, no other singularity, and $\partial F$ contains only regular points and conical singularities; or $F$ contains a neighborhood of a unique swath, no singularity, and $\partial F$ contains only regular points, conical singularities and the irregular singularity of the swath.
If $p$ is anti-conical we are done.
If $p$ is a conical singularity, the local model of \Cref{sub:Fuchsian} implies that any geodesic tending to $p$ must hit it in finite time, contradicting (*).
If $p$ is an irregular singularity, then by \Cref{lem:irregularTips} either the geodesic tends to the swath and we are done, or it tends to a focus in finite time, again contradicting (*).
\end{proof}

\section{Geometry of affine immersions}\label{sec:AffImm}

The main goal of this section is to describe affine immersions of planes, open half-planes and open disks into a finite type affine surface, proving all the geometric results underlying the Delaunay decomposition constructed in \Cref{sec:Delaunay}. We more generally study immersions of convex subsets of $\C$.

\medskip

We recall that the word affine is to be understood in this article as complex-affine.
An affine map between affine surfaces is defined as a map that is locally affine in charts. It is necessarily an immersion and injective ones are necessarily embeddings. We choose here to use the terminology affine immersion to insist on the absence of injectivity hypothese, and the terminology affine embedding by coherence.

We recall that marked points of a finite type affine surface $(X,\nabla,\cal S)$ are erasable singularities, and stress that that we decide to stop geodesics at marked points even though we may naturally extend them beyond.

\subsection{Affine identity principle and relation between affine immersions}\label{sub:acparbi}

We recall that the affine indentity principle is the statement that if two affine maps are defined on the same \emph{connected} subset of an affine surface $S$ and take values in some affine surface $S'$, and if they coincide in the neighborhood of some point, then they are equal.
An immediate consequence, bearing the same name, concerns two affine maps defined on two open subsets $U_1$, $U_2$ of $S$ with non-empty intersection and taking values in $S'$: if they coincide in the neighborhood of some point $s$, then they are equal on the connected component of $U_1\cap U_2$ containing $s$.
We will use several times the following consequence this principle.

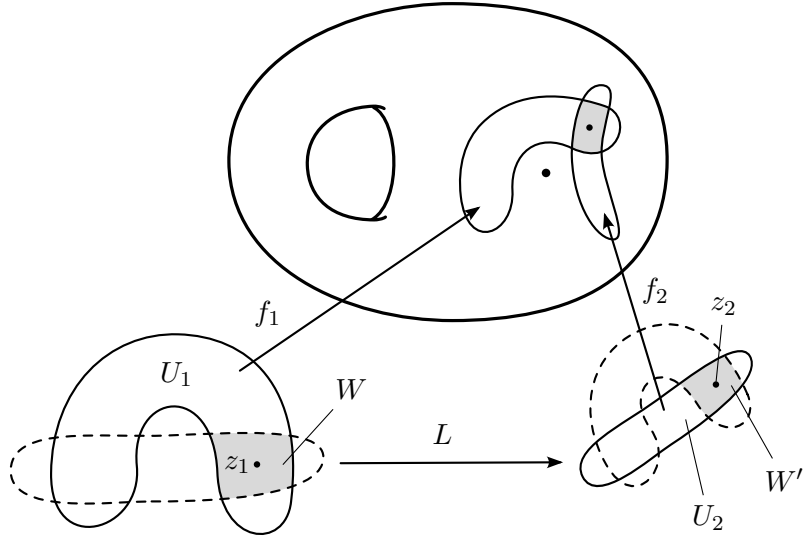
\begin{figure}[htbp]
\begin{center}
\begin{tikzpicture}
\node at (0,0) {\includegraphics[scale=0.8]{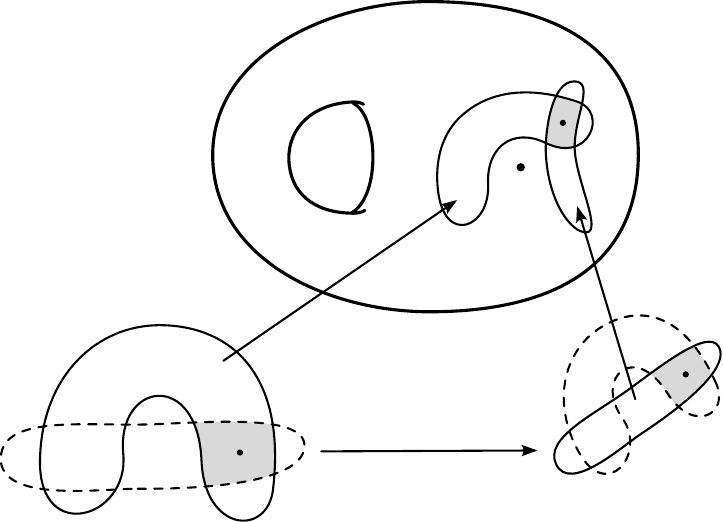}};
\node at (-1.5,-0.6) {$f_1$};
\node at (3.6,-0.3) {$f_2$};
\node at (0.8,-2.2) {$L$};
\node(A) at (-0.4,-1.6) {$W$};
\draw (A) -- (-1.3,-2.6);
\node at (-2.7,-1.4) {$U_1$};
\node(C) at (4.3,-3.3) {$U_2$};
\draw (C) -- (4,-2);
\node at (-1.9,-2.6) {$z_1$};
\node(B) at (4.5,-0.5) {$z_2$};
\draw (B) -- (4.4,-1.5);
\node(D) at (5.3,-2.8) {$W'$};
\draw (D) -- (4.6,-1.55);
\end{tikzpicture}
\end{center}
\caption{Illustration of \Cref{lem:irel} in a case where the immersions are embeddings.}
\label{fig:olu}
\end{figure}

\begin{lem}\label{lem:irel}
  Let $S$ be an affine surface (not necessarily of finite type).
  Let $f_1: U_1\to S$ and $f_2:U_2\to S$ be two affine immersions with $U_1$, $U_2$ open subsets of the canonical affine plane $\C$.
  Assume that there exists $z_1\in U_1$ and $z_2\in U_2$ such that $f_1(z_1)=f_2(z_2)$.
  Consider the unique element $L\in \AutC$ sending $z_1$ to $z_2$ and such that $f_1=f_2\circ L$ holds near $z_1$.
  Denote by $W$ be the connected component of $U_1\cap L^{-1}(U_2)$ that contains $z_1$
  and $W'=L(W)$ the connected component of $U_2\cap L(U_1)$ that contains $z_2$.
  Then 
  \[ f_1=f_2\circ L\text{ holds on }W
  .\]
  In particular, $f_2(U_2)\cap f_1(U_1)$ contains $f_1(W)=f_2(W')$. See \Cref{fig:olu}.
\end{lem}

For instance, in the particular case where $L^{-1}(U_2) \subset U_1$, we get $f_2(U_2) = f_1(L^{-1}(U_2))\subset f_1(U_1)$.

We will apply this lemma several times.
A typical situation is the following (we call strip any image by an element of $\AutC$ of the set $\setof{z\in\C}{\Re(z)\in(0,1)}$):

\begin{lem}\label{lem:onde}
If a geodesic $t\geq 0 \mapsto\gamma(t)$ tends at $+\infty$ to an anti-conical Fuchsian singularity $p$ of a meromorphic connection, for which we have a punctured neighborhood $V$ isomorphic to a standard neighborhood in a model of $p$ as in \Cref{sub:Fuchsian}, and $\gamma$ extends to an affine immersion $f$ of a convex open set $\Delta\subset\C$ containing $[0,+\infty)$ and not equal to a strip between two parallel lines nor to the whole plane, then there is a neighborhood $W'$ of $\infty$ in $\Delta$ that is mapped to $V$.
\end{lem}

The skew cylinder model of $V$ allows to control the immersion $f$ and the geodesics, in this part $W'$ of $\Delta$.

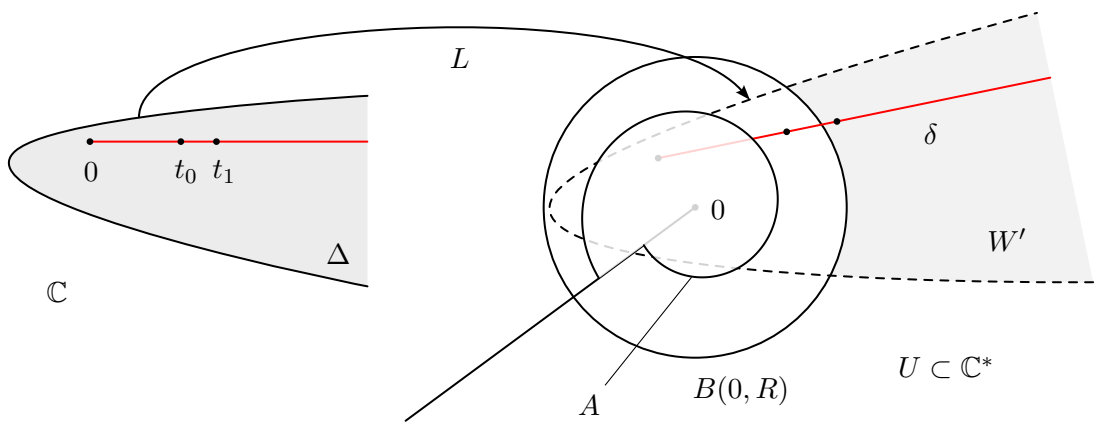
\begin{figure}[htbp]
\begin{center}
\begin{tikzpicture}
\node at (0,0) {\includegraphics[scale=0.75]{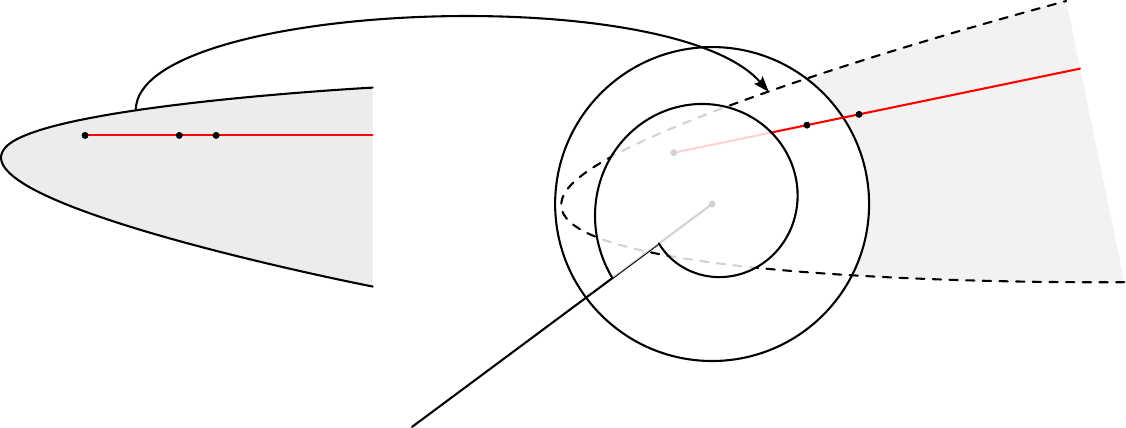}};
\node at (-6.05,0.6) {$0$};
\node at (-4.8,0.6) {$t_0$};
\node at (-4.3,0.6) {$t_1$};
\node at (-2.8,-0.5) {$\Delta$};
\node at (-6.5,-1) {$\C$};
\node at (-1.2,2.1) {$L$};
\node at (5,1.1) {$\delta$};
\node at (2.2,0.1) {$0$};
\node(A) at (0.5,-2.5) {$A$};
\draw (A) -- (1.85,-0.8);
\node at (2.5,-2.3) {$B(0,R)$};
\node at (6,-0.3) {$W'$};
\node at (5.2,-2) {$U\subset\C^*$};
\end{tikzpicture}
\end{center}
\caption{Illustration of the proof of \Cref{lem:onde}. The part of the skew cylinder that maps in $V$ omits the semi-transparent area labeled $A$.}
\label{fig:onde}
\end{figure}

\begin{proof}
There is a punctured neighborhood $V$ of $p$ isomorphic: in the pure anti-conical case to a neighborhood of $\infty$ in a skew cylinder model $\cet\to\cet/\langle \lambda\rangle$;
in the shifted anti-conical case, to a neighborhood of $\infty$ in the $n$-fold covering of the punctured canonical affine plane $\C^*$, modified by a translation shift along a radial line.
Let $t_0$ be such that $\gamma([t_0,+\infty))\subset V$.
In the model, this part of $\gamma$ lifts to a straight path $\tilde\delta:[t_0,+\infty)\to \cet$.
This path $\tilde\delta$ projects to $\C^*$ to a geodesic $\delta$ of the canonical affine plane $\C$. Let $L\in\AutC$ be such that $\delta(t)=L(t)$.
From the model, there is an affine immersion $f_2:U\to V$ for a set of the form $U=\C\setminus(B(0,R)\cup D')$, where $D'$ is a half-line from $0$ and $R$ big enough so that $L(\Delta)\cap U$ has one and only one connected component $W'$ (here we use that $\Delta$ is not equal to a strip between two parallel lines).
It is necessarily unbounded and contains $\delta([t_1,+\infty))$ for some $t_1\geq t_0$.
See \Cref{fig:onde}.
It follows that $W:=L^{-1}(W')$ is a connected neighborhood of $\infty$ in $\Delta$.
We then conclude with \Cref{lem:irel} applied to $f_1=f$, $z_1=t_1$, $U_1=\Delta$, $f_2$, $z_2=\delta(t_1)$ and $U_2=U$ that $f=f_2\circ L$ holds on $W=L^{-1}(W')$, so $f(W)\subset V$.
\end{proof}

As in the case of analytic maps, there are stronger version of the affine identity theorem: for instance, two affine maps from a connected affine surface $S$ to an affine surface $S'$ to be equal, it is enough that they coincide on a set that has an accumulation point in $S$.
We exploit this in the following statement.

\begin{lem}\label{lem:ContinExt}
Let $f:\Delta\to S$ be an affine immersion of a convex subset $\Delta$ of $\C$. Let $z_0\in \partial \Delta$ and $z_1 \in\Delta$.
Assume that $f(z)$ converges in $X$ as $z$ tends to $z_0$ along the segment $(z_0,z_1]$, in other words, that the geodesic $t\mapsto f((1-t)z_0+tz_1)$ converges at $t=0$. Denote $y$ this limit.
Then either $y$ is a regular point (i.e.\ $y\in X^*$) and $f$ has an affine extension to a disk $B(z_0,r)$, or $y$ is a conical Fuchsian singularity or an irregular singularity and 
$f$ has an affine extension to a slit disk $B(z_0,r)\setminus I$ where the slit $I$ is a radius of the disk in the direction opposed to $z_1$.
Moreover, $f$ extends continuously to $\{z_0\}$.
\end{lem}

\begin{proof}
If $y$ is regular, by post-composing a chart near $y$ with an element of $\AutC$ we can find a chart $\psi : V\to B(z_0,r)$ for some $r>0$ and such that $\psi\circ f$ is the identity for every $z\in B(z_0,r)\cap(z_0,z_1]$. The domains of the maps $f$ and $\psi^{-1}$ both contained the convex, hence connected, set $B(z_0,r)\cap \Delta$, and coincide on $B(z_0,r)\cap(z_0,z_1]$, so are equal on $B(z_0,r)\cap \Delta$ by the version of the identity principle stated just before.
We can thus extend $f$ with $\psi^{-1}$.

If $y$ is not regular then \Cref{lem:tips} implies it is a conical Fuchsian singularity, or an irregular singularity.
Then \Cref{lem:SlitDisk} implies there is some $g:U:=B(z_0,r)-I \to $ for some $r>0$ and with $I$ a radius in the direction opposed to $z_1$.
By the same argument as in the regular case, $f$ and $g$ must coincide on $\Delta\cap U= \Delta\cap B(z_0,r)$, so $g$ allows to extend $f$ as wished.
\Cref{lem:tips} also states that $g$ converges to $y$ at $z_0$.
\end{proof}

Independently, we state and prove below a useful result.

\begin{lem}\label{lem:invcc}
  Let $S$ be an affine surface (not necessarily of finite type) and $U$ an open subset of the canonical affine plane $\C$.
  Assume $f:U \to S$ is an affine immersion, $V$ is a connected open subset of $S$, and $W$ a connected open subset of $U$ on which $f$ is a bijection to $V$.
  Then $W$ is a connected component of $f^{-1}(V)$.
\end{lem}

\begin{proof}
  Let $W'$ be the connected component of $f^{-1}(V)$ that contains $W$. Then $W'=W$ for otherwise there would be a point $z$ in the boundary of $W$ relative to $W'$.
  This point would be mapped in $V$ by $f$, so there would be another point $z'\in W$ mapped to the same point by $f$.
  But then points in $W$ close to $z$ would have at least two preimages in $V$, one near $z$ and one near $z'$, contradicting injectivity of $f$ in $W$.
\end{proof}

\subsection{Immersions of planes}\label{sub:AffImmPlanes}

Few affine surfaces admit an affine immersion of the whole affine plane. All of them are translation surfaces. 

\begin{prop}\label{prop:PlaneImmersion}
Consider a finite type affine surface $(X,\nabla,\mathcal{S})$ such that there exists an affine immersion $f:\mathbb{C} \rightarrow X^{\ast}$. Then, $f$ is a universal covering and $(X^*,\nabla)$ is isomorphic to one of the following surfaces:
\begin{enumerate}[(i)]
    \item \textit{complex flat plane}: $X^{\ast}\cong \mathbb{C}$, $X$ has genus $0$ and $\nabla$ has on $X$ a unique Fuchsian singularity of residue equal to $2$;
    \item \textit{full translation cylinder}: $X^{\ast} \cong \mathbb{C}/\mathbb{Z}$, $X$ has genus $0$ and $\nabla$ has on $X$ two Fuchsian singularities with residues equal to $1$ (the two ends of the cylinder);
    \item \textit{translation torus}: $X^*=X\cong\mathbb{C}/\Lambda$ where $\Lambda$ is a lattice of $(\mathbb{C},+)$ and $\nabla$ has no singularity.
\end{enumerate}
In each case above there is an isomorphism $\phi: \C/\Lambda \to X^*$ and $L\in \AutC$ such that $f = \phi\circ \pi \circ L$, where $\pi:\C\to \C/\Lambda$ and $\Lambda < (\C,+)$ is respectively the trivial group, $\Z$, or the involved lattice.
\end{prop}

\begin{proof}
Consider a point $x\in X^*$. There is an affine chart $\phi: V\to U$ where $V$ is an open subset of $X^*$ containing $x$, $U\subset \C$ and where $V$ and $U$ are connected.
We will prove that $V$ is well covered.
If $V$ and $f(\C)$ are not disjoint then let $W$ be a connected component of $f^{-1}(V)$.
We claim that $f$ is an affine bijection from $W$ to $V$. We provide two proofs:

Argument~1: consider the map $g:= \phi\circ f: W\to U$.
For any point in $f(\C)$, maximal geodesics containing this point are bi-infinite (i.e.\ they extend to $\R$).
Any two points of $U$ can be connected by a finite chain of straight segments.
It follows that $g$ is surjective.
Since $g$ is locally affine and $W$ connected, $g$ is the restriction to $W$ of a single element of $\AutC$, hence it is injective.
So $g$ is a bijection from $W$ to $U$.

Argument~2: apply \Cref{lem:irel} to $f_1 = f$, $U_1=\C$, $z_1\in W$, $f_2 = \phi^{-1}$, $U_2 = U$ and $z_2=\phi(f(z_1))$. We obtain that there exists $L\in\AutC$ such that, denoting $W^0=L^{-1}(U)\subset \C$, we have: $V=f(W^0)$, $f = \phi^{-1}\circ L$ holds on $W^0$ and $z_1\in W^0$. By \Cref{lem:invcc}, $W^0$ is a connected component of $f^{-1}(V)$, so $W^0=W$.

Since $X^*$ is connected and $\C$ non-empty, the covering $f$ is surjective.
So $X^{\ast}$ is a quotient of $\mathbb{C}$ by a discrete subgroup $\Lambda$ of $\text{Aff}(\mathbb{C})$ that acts freely discontinuously.
Such a group can only contain translations.
So $\Lambda$ is a discrete subgroup of $(\mathbb{C},+)$, and there are three cases:
\begin{itemize}
    \item[(i)] $\Lambda$ is trivial;
    \item[(ii)] $\Lambda$ is generated by a translation $z \mapsto z + \alpha$ for some $\alpha \in \mathbb{C}^{\ast}$;
    \item[(iii)] $\Lambda$ is a lattice of $(\mathbb{C},+)$.
\end{itemize}
These three classes of subgroups correspond to the three cases of the statement. 
In cases (i) and (ii), there are explicit embeddings of the affine surfaces $\C$ and $\C/\Z$ in $\RS$ for which the affine structure on $X^{\ast}=\mathbb{C}/\Lambda$ extends to its end(s) as Fuchsian singularities of the meromorphic connection on $\RS$ with the residues indicated in the statement: see \Cref{ex:C,lem:AffineCylinder}.
\end{proof}

The affine surfaces classified in \Cref{prop:PlaneImmersion} have no apexes (conical singularity or focus of an irregular singularity), so no saddle connections. They form a subclass of the exceptional affine surfaces that we will classify in \Cref{sub:Exceptional}.

\medskip

The following result uses techniques analogue to \Cref{prop:PlaneImmersion} and will be needed in the proof of \Cref{lem:de0,thm:ExceptionalSurfaces}.

\begin{prop}\label{prop:EImmersion}
Consider a finite type affine surface $(X,\nabla,\mathcal{S})$ such that there exists an affine immersion $f:\cal E \rightarrow X^{\ast}$. Then, either $X^*$ admits a whole plane immersion, or $(X^*,\nabla)$ has the exponential-affine plane $\cal E$ as a universal covering (hence is isomorphic to $\cal E$, a full Reeb cylinder, a skew cone or an affine torus).
\end{prop}
\begin{proof}
As a Riemann surface $\cal E$ can be identified with $\C$.
If $f$ is injective, then we can use \Cref{lem:lem} (an injective holomorphic map from $\D^*$ to a compact Riemann surface $X$ necessarily extends holomorphically at $0$) to the map $z\in \D^*\mapsto f(1/z)$: we get that $f$ extends to a holomorphic map $\hat f:\RS\to X$.
The point $\hat f(\infty)$ is an irregular singularity (a pole of order $2$ of the connection). It follows that $\hat f$ is injective.
Holomorphic maps between compact Riemann surfaces are surjective.
Hence $X^*$ is isomorphic to $\cal E$.

In the rest of the proof we assume that $f$ is non-injective.
The affine surface $\cal E$ can be identified with the universal covering $\cet$ of the punctured canonical affine plane $\C^*$.
Let $(r_1,\theta_1)$ and $(r_2,\theta_2)$ be two points, in polar coordinates, mapped by $f$ to a same point $x_0\in X^*$.
Let $z_1 = r_1 e^{i\theta_1}$, $z_2 = r_2 e^{i\theta_2}$ and $\iota_1$, $\iota_2$ two local sections of $\pi:\cet\to \C$ sending $z_i$ to $(r_i,\theta_i)$.
There exists $\phi\in \AutC$, uniquely determined by $(r_1,\theta_1)$ and $(r_2,\theta_2)$, such that $\phi(z_1)=z_2$ and $f\circ \iota_1 = f\circ\iota_2\circ\phi $ holds near $z_1$.

Case 1: $\phi$ does not fix $0$. 
Then either $0\in\phi(B(z_1,|z_1|))$ or $0\in\phi^{-1}(B(z_2,|z_2))$. Up to permuting the roles of the $z_i$, we assume this is the first option.
Let $\iota_1$ be the section on $B(z_1,|z_1|)$ of $\pi$ sending $z_1$ to $(r_1,\theta_1)$.
By affine continuation, we get that $f = f\circ\iota_1\circ\phi$ holds on $\pi^{-1}(\phi(B(z_1,|z_1|)))$. In particular, there is a point of $\cet$ near which $f \circ \zeta = f$ where $\zeta:(r,\theta)\mapsto(r,\theta+2\pi)$ is a generator of the deck transformation group of $\pi$.
By affine extension, the relation $f \circ \zeta = f$ holds on all $\cet$.
So $f$ factors through $\pi$ as $f=g\circ\pi$ for some affine immersion $g:\C^*\to X^*$ and with $\phi$ we can extend $g$ into some affine immersion $\hat g:\C\to X^*$.

Case 2: for any pair $(r_1,\theta_1)$ and $(r_2,\theta_2)$, the map $\phi$ fixes $0$. 
Such a $\phi$ lifts to a global automorphism $\tilde \phi:\cet\to\cet$ sending $(r_1,\theta_1)$ to $(r_2,\theta_2)$ and by analytic continuation, $f\circ\tilde \phi = f$ holds on all $\cet$.
The set of all such $\tilde \phi$, together with the identity, then forms a subgroup $\Lambda$ of the automorphism group of $\cet$ and this group must be discrete for otherwise it would contradict local injectivity of $f$.
Hence the map $f$ factors to an embedding of the quotient manifold $\cet/\Lambda$ in $X^*$ and
$\Lambda$ is a lattice or a cyclic subgroup of the translation group of $\cal E$, so the quotient is a) an affine torus, b) a skew cylinder or c) a full Reeb cylinder.
For a) the image $f(\cal E)$ is compact, hence equal to $X$.
In cases b) and c) we can proceed as in the injective case, where we used \Cref{lem:lem}: for b) $X-f(\cal E)$ consists in a conical and an anti-conical singularity and $X$ is a sphere, for c) $X-f(\cal E)$ consists in two Reeb type singularities and $X$ is a sphere too.
\end{proof}

\subsection{Boundary behavior}

We prove here that affine embeddings of open convex subsets of $\mathbb{C}$ in $X^*$ always extend to their boundary.
Default of continuous extension of affine immersions are hence always related to their non-injectivity.
\par

We start with the following lemma that concerns immersions (not only embeddings) and that will be useful later too.

\begin{lem}\label{lem:Avoid2}
Consider an affine immersion $f$ of an open convex set $\Delta \subset \mathbb{C}$ into $X^*$.
Let $z_\infty\in \partial \Delta$ and $z_0\in \Delta$.
Consider the semi-open segment $I=[z_0,z_\infty)$.
Its image by $f$ is the support of a semi-open arc of geodesic.
Let $y$ be an accumulation point of this geodesic, i.e.\
$y=\lim f(z_n)$ with  $z_n\in I$ and $z_n\tends z_\infty$.
Then $y$ cannot be a non conical Fuchsian singularity.
If moreover $y \notin f(\Delta)$ then
either:
\begin{enumerate}[(i)]
\item\label{item:a2yreg} $y$ is regular and $f$ has an affine extension to a neighborhood of $z_\infty$;
\item\label{item:a2ysing} $y$ is singular and $f$ has a continuous extension to $\Delta\cup\{z_\infty\}$.
\end{enumerate}
In both cases, the extension maps $z_\infty$ to $y$.
\end{lem}

\begin{proof}
Without loss of generality we may assume that $z_\infty =0$ and $z_0=1$, so $I=[1,0)$ and $z_n$ is a positive real.
Let $r_0 \in (0,1)$ such that $B(1,r_0)\subset \Delta$.
Then for all $n>0$, $B(z_n,r_0z_n)\subset \Delta$.

Note that $f$ restricted to $I$ is a geodesic.
Assume that $y$ is a non conical Fuchsian singularity.
By \Cref{lem:traps}, $y$ has arbitrarily small trap neighborhoods, such that any geodesic entering the neighborhood through its boundary never leaves, but we get a contradiction because either it is impossible to reach the singularity in the case of a Reeb type one, or it takes an infinite amount of time to accumulate it.
Thus $y$ can only be a regular point of $X$, a conical Fuchsian singularity or an irregular singularity.

If $y$ is an irregular singularity, there are also arbitrarily small traps, described in \Cref{sub:geodesics}: they are the union of sepals and attracting petals. After crossing the boundary of this trap, a geodesic takes infinite time to accumulate the singularity.
The complement of the traps consists in interpetal spaces that converge to the apexes (foci) of $y$.
Hence a geodesic accumulating $y$ must accumulate the set of apexes. Up to taking a subsequence we can assume that $f(z_n)$ tends to one of the apexes. In this case we denote $P$ a repelling petal associated to this apex.

For any $R>0$ by \Cref{lem:SlitDisk0} (applied to $R+1$ composing $\phi$ by $1-z$ and taking a restriction) there exists a punctured neighborhood $V$ of $y$ in case~(1), a repelling petal neighborhood of the apex in case~(2), such that $\forall x\in V$ and $\forall x\in V\cap P$, there is an affine immersion (possibly injective) $\phi$ from the slit disk
\[ B(0,R) \setminus [1,R) = \dom(\phi)
\]
to $X^*$ sending $0$ to $x$ and such that $\phi(z)$ tends to $y$ as $z$ tends to $1$ within $\dom(\phi)$.
We choose $R=2/r_0$ (note that $R>2$) and consider such a $V$.
There exists $n$ such that $f(z_n)\in V$.
Consider the map $\phi$ associated to $x=f(z_n)$ as above.
Let us apply \Cref{lem:irel} to 
\[ \begin{split}
& f_1 = \phi,\ z_1 = 0,\ U_1=\dom(\phi), \\
& f_2=f,\ z_2=z_n \text{ and } U_2 = \Delta\cap B(z_n,2|z_n|)
.\end{split} \]
It involves some $L\in\AutC$ such that $L(0)=z_n$ and that the identity $f_1=f_2\circ L$ holds on the connected component $W$ of $U_1\cap L^{-1}(U_2)$. See \Cref{fig:aro}.

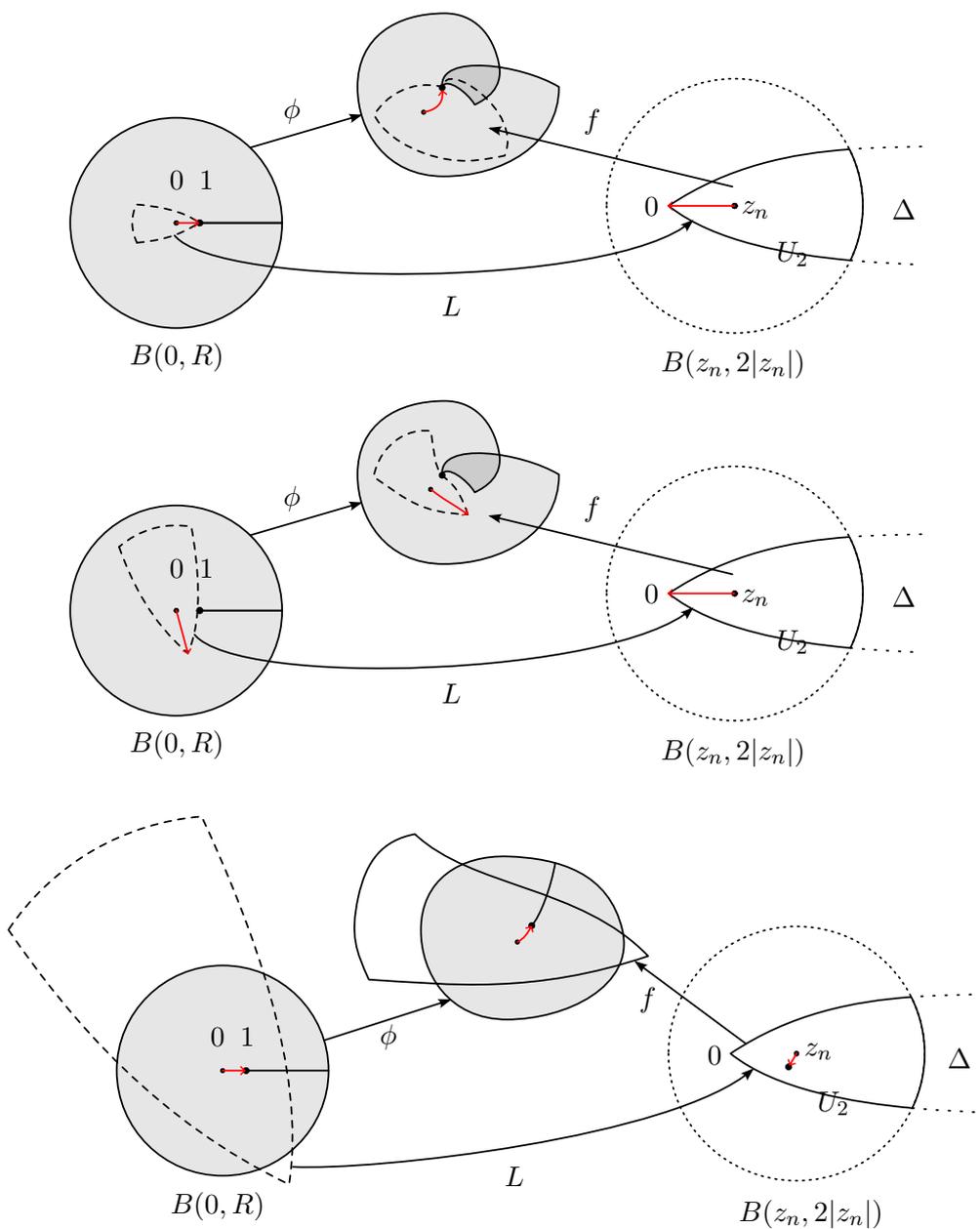
\begin{figure}[htbp]
\begin{center}
\begin{tikzpicture}
\node at (-4.5,-1.3) {or};
\node at (0,0) {\includegraphics[scale=0.66]{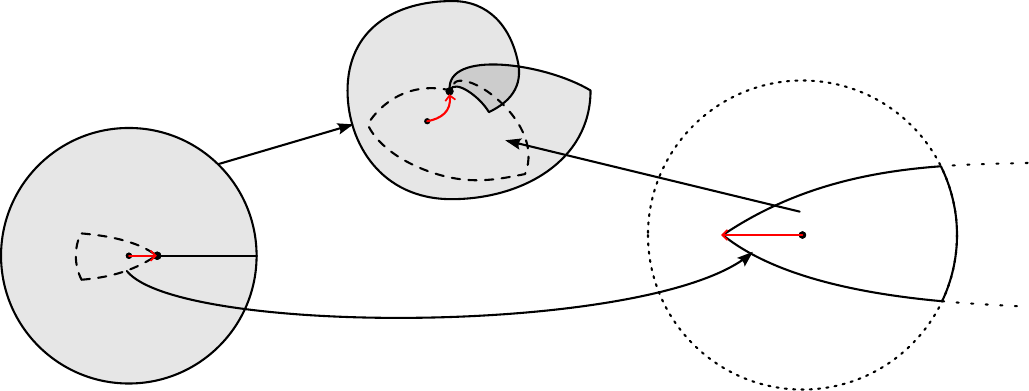}};
\node at (-2.75,0.85) {$\phi$};
\node at (1.3,0.7) {$f$};
\node at (-4.3,-0.1) {$0$};
\node at (-3.9,-0.1) {$1$};
\node at (2.1,-0.45) {$0$};
\node at (3.5,-0.5) {$z_n$};
\node at (-4.3,-2.5) {$B(0,R)$};
\node at (-0.6,-1.8) {$L$};
\node at (4,-1.1) {$U_2$};
\node at (5.5,-0.5) {$\Delta$};
\node at (3.2,-2.6) {$B(z_n,2|z_n|)$};
\end{tikzpicture}
\begin{tikzpicture}
\node at (-4.5,-1.3) {or};
\node at (0,0) {\includegraphics[scale=0.66]{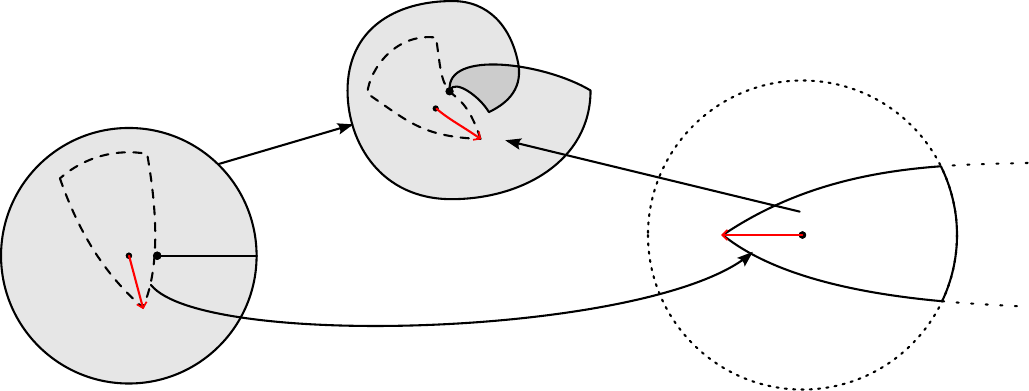}};
\node at (-2.75,0.85) {$\phi$};
\node at (1.3,0.7) {$f$};
\node at (-4.3,-0.1) {$0$};
\node at (-3.9,-0.1) {$1$};
\node at (2.1,-0.45) {$0$};
\node at (3.5,-0.5) {$z_n$};
\node at (-4.3,-2.5) {$B(0,R)$};
\node at (-0.6,-1.8) {$L$};
\node at (4,-1.1) {$U_2$};
\node at (5.5,-0.5) {$\Delta$};
\node at (3.2,-2.6) {$B(z_n,2|z_n|)$};
\end{tikzpicture}
\vskip1em
\begin{tikzpicture}
\node at (-4.5,-1.3) {or};
\node at (-0.55,0.3) {\includegraphics[scale=0.66]{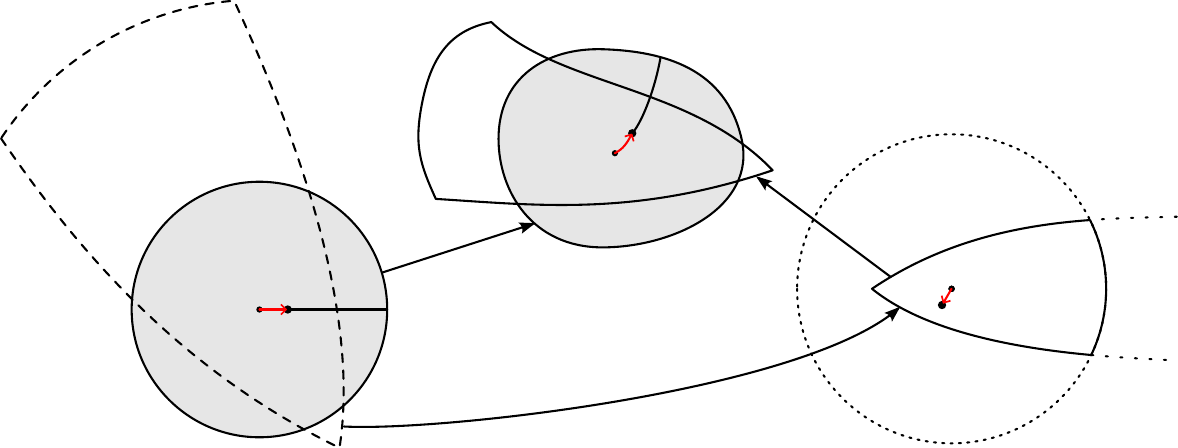}};
\node at (-2,-0.2) {$\phi$};
\node at (1.5,0.3) {$f$};
\node at (-4.3,-0.2) {$0$};
\node at (-3.9,-0.2) {$1$};
\node at (2.4,-0.45) {$0$};
\node at (3.8,-0.4) {$z_n$};
\node at (-4.3,-2.5) {$B(0,R)$};
\node at (-0.3,-2.1) {$L$};
\node at (4,-1.1) {$U_2$};
\node at (5.7,-0.5) {$\Delta$};
\node at (3.7,-2.6) {$B(z_n,2|z_n|)$};
\end{tikzpicture}
\end{center}
\caption{Illustration of the proof of \Cref{lem:Avoid2}.}
\label{fig:aro}
\end{figure}

If $1\in L^{-1}(U_2)$ then note that since $U_2$ is convex, the whole set $[0,1)$ is contained in $L^{-1}(U_2)\cap \dom\phi$ hence in $W$.
Then $f=\phi\circ L^{-1}$ on the ray $L([0,1))$, so it tends to $y$ along this ray, which is a semi-open segment from $z_n$ ending at some inner point of $U_2$.
By continuity, $f$ takes the value $y$.

If $1\notin L^{-1}(U_2)$, then since $R\geq 2/r_0$, the convex set $L^{-1}(U_2) = L^{-1}(\Delta\cap B(z_n,2|z_n|))$ is contained in $B(0,R)$.
It contains $0$ and is disjoint from the slit $[1,R)$, for otherwise, by convexity, it would contain $1$, contradicting our assumption.
So $L^{-1}(U_2) \subset U_1$ and the identity $f = \phi \circ L^{-1}$ holds on $U_2$, which is the intersection of $\Delta$ with a neighborhood of $0$ in $\C$.
This provides a continuous extension to $f$ at $0$, and actually an affine extension to a neighhorhood of $0$ in $\C$ if either $L^{-1}(0)\neq 1$ or $y$ is regular (for, then, we can take $\phi$ that extends to the whole disk $B(0,R)$).
\end{proof}

In Case~\eqref{item:a2ysing}, the proof actually provides an extension to a slit neighborhood of $z_\infty$, but note that \Cref{lem:ContinExt} also gives it.

\medskip

\Cref{lem:Avoid2} above allows us to draw an interesting consequence for embeddings of convex sets.

\begin{cor}\label{cor:Embedding}
Any affine embedding $f$ of an open convex set $\Delta \subset \mathbb{C}$  into $X^{\ast}$ extends continuously to the boundary of $\Delta$ in $\C$.
\end{cor}

\begin{proof}
Consider as in \Cref{lem:Avoid2} a point $z_\infty\in \partial \Delta$, any $z_0\in \Delta$ and the semi-open segment $I=[z_0,z_\infty)$. By compactness of $X$ there is  an accumulation point $y$ of $f$ at $z_\infty$ along $I$, i.e.\ $y=\lim(z_n)$, $z_n\in I$, $z_n\tends z_\infty$.
\par
By \Cref{lem:Avoid2} it is enough to prove that $y$ cannot belong to the image of $f$. Indeed, if there were $z \in \Delta$ such that $y=f(z)$, then consider a neighborhood $V$ of $z$ such that $z_{n} \notin V$ provided $n$ is large enough. Since $f$ is a local diffeomorphism in a neighborhood of $z$, we can find an element of $V$ that has same image as $z_{n}$ under $f$. This contradicts the injectivity of $f$ and proves the claim.
\end{proof}

The continuous extension actually holds in a stronger sense: one can replace the range $X$ of the extension by the union of $X^*$ and all the apexes (conical Fuchsian singularities and foci of irregular singularities), with the appropriate topology near foci described in \Cref{sub:helix}.

We end with a general result concerning continuous extension of immersions of convex sets.

\begin{lem}\label{lem:b}
Consider an affine immersion $f$ of an open convex set $\Delta \subset \mathbb{C}$ into $X^*$.
Let $C_0$ be the set of points in $\partial \Delta$ where $f$ has a continuous extension.
Then $C_0$ is an open subset of $\partial \Delta$, the extension maps in the set of regular points, conical Fuchsian and irregular singularities,
the points mapped to non-regular points form a closed and discrete subset of $C_0$.
\end{lem}
\begin{proof}
Denote $\hat f$ the extension.
Consider $z_0\in C_0$ and choose any $z_1\in \Delta$: the conditions of \Cref{lem:ContinExt} are in particular met, and the lemma tells us that $y:=\hat f(z_0)$ cannot be a non-conical Fuchsian singularity.
Moreover \Cref{lem:ContinExt} tells us that there is an affine extension of $f$ to a neighborhood of $z_0$ or to a slit neighborhood of $z_0$, with the slit in the direction opposite to $z_1$. The extension region contains a neighborhood of $z_0$ in $\partial \Delta$, so $C_0$ is open and points mapped to non-regular points are isolated.
Closedness follows from continuity.
\end{proof}

\subsection{Immersed sectors around geodesics}\label{sub:immersedConesLemmas}

In this section, $(X,\nabla,\mathcal{S})$ is a finite type affine surface. For the affine immersion of an infinite convex open sector $\Delta$ of vertex $z_{0}$, we refer to a geodesic contained in $\Delta$ and starting at $z_{0}$ as a \textit{radial geodesic}.

\subsubsection{Behavior at infinity}

After studying limits on boundary points we turn in the next two statements to limits  at infinity.
An infinite open sector 
\[\Delta=\setof{z\in\C\setminus\{z_0\}}{\arg(z-z_0) \in (\alpha,\beta)}\]
is convex iff $|\beta-\alpha|\leq \pi$.
We call $(\alpha,\beta)\subset\R/2\pi\Z$ the \emph{angular span} of $\Delta$.
The first statement is about limits along geodesics.

\begin{lem}\label{lem:preAvoid}
Consider an affine immersion $f$ of a convex open infinite sector $\Delta$ in the canonical affine plane $\C$ into $X^*$.
Consider a geodesic $\gamma:t\in(0,+\infty) \to a_\gamma t+b_\gamma\in\Delta$ not parallel to $\partial \Delta$ and $y\in X$ an accumulation point as $t\to+\infty$ of the geodesic $t\in(0,+\infty)\mapsto f\circ \gamma(t)$ of $X$.
Then either:
\begin{enumerate}[(i)]
\item\label{item:pa_reg} $y$ is regular and $y\in f(\Delta)$;
\item\label{item:pa_ncf} $y$ is Fuchsian non conical and $f\circ \gamma$ tends to $y$ at infinity;
\item\label{item:pa_irr} $y$ is irregular and $f\circ \gamma$ tends to a swath of $y$ at infinity.
\end{enumerate}
\end{lem}

\begin{proof}
By performing an affine change of variable on $\Delta$, we may assume that $\gamma(t)=t$.
The sector $\Delta$ is based on some $z_0\in\C$.
Since the the two boundary lines of $\Delta$ are non horizontal, there exists $r_0>0$ such that $\forall t>0$, $B(t,r_0t)\subset \Delta$.
We recall that no geodesic can accumulate a Reeb type singularity, so $y$ is not such a point.
\par
Assume that $y$ is a regular point.
Assume by contradiction that $y\notin f(\Delta)$.
Consider a local chart $\phi$ defined on an open neighborhood $V$ of $y$, with $\phi(V)=B(\phi(y),r)$.
Consider $t\in (0,+\infty)$ such that $f(t)\in \dom(\phi)$ and denote $z_1 = \phi(f(t))$.
There is a unique $L\in \AutC$ such that $L(z_1)=t$ and $\phi^{-1} = f\circ L$ holds near $z_1$. This map $L$ may depend on $t$.
We will apply \Cref{lem:irel} to
\[ \begin{split}
& f_1 = \phi^{-1}, z_1,\ U_1 = B(\phi(y),r)
\\
& f_2 = f,\ z_2=t \text{ and } U_2 = \Delta
.\end{split}\]
Denoting $W$ the connected component of $B(\phi(y),r)\cap L^{-1}(\Delta)$, we get that $\phi^{-1}=f\circ L$ holds on $W$ and particular that $f(\Delta)$ contains $\phi^{-1}(W)$.
The convex set $L^{-1}(\Delta)$ cannot contain $\phi(y)$ for otherwise it would contain the segment $[\phi(f(t)),\phi(y)]$, so $W$ would containt it too, hence $f(\Delta)$ would contain $y$.
From $\phi(y)\notin L^{-1}(\Delta)$ and $B(t,r_0 t)\subset \Delta$, we get that the scaling factor of $L^{-1}$ is at most $|\phi(f(t))-\phi(y)|/r_0$.
Hence each time $|\phi(f(t))-\phi(y)|\times(1+2/r_0) < r$, then $L^{-1}(\Delta\cap B(t,2t)) \subset B(\phi(y),r)$.
Since this happens for arbitrarily high values of $t$, this implies that the whole set $L^{-1}(\Delta)$ is contained in $B(\phi(y),r)$.
But $L^{-1}(\Delta)$ is unbounded, leading to a contradiction. So $y$ cannot be regular.

If $y$ is a conical Fuchsian singularity, then $y\notin f(\Delta)$ since the latter is contained in $X^*$.
We get contradiction by a similar approach to the Case~\eqref{item:pa_reg}, by proving that for any neighborhood $V'$ of $y$, $f(\Delta)\subset V'$.\footnote{This is absurd on two levels: $f$ would have to map to the intersection of all $V'$, which is $\{y\}$, but $f$ maps to $X^*$; actually, mapping $\Delta$ affinely in a small enough neighborhood of $y$ is already impossible since $y$ has neighborhoods that contain no geodesic defined for infinitely long time.}
We imitate the proof for a regular point: \Cref{lem:SlitDisk0} gives us some neighborhood $V$ of $y$, which, for $t>0$ such that $f(t)\in V$, allows us to replace $\phi^{-1}$ in our argument by an affine immersion $f_1:B(0,R) \setminus (-R,0]\to V' \setminus \{y\}$ such that $f_1(1)=f(t)$.
As in the case of regular points, the choice of an $R>1+2/r_0$ ensures that $\Delta\subset V'$. We again get a contradiction.

IF $y$ is a anti-conical Fuchsian or cylindrical Fuchsian then~\eqref{item:pa_ncf} follows from \Cref{lem:traps}.

If $y$ is an irregular singularity, then consider a trapping neighborhood of its swathes consisting of the union of the attracting petals and of the sepals.
If the geodesic enters the trap, then it tends to a swath;
If the geodesic never enters a trap, then it accumulates a focus of $y$. We then get as in the conical case that the whole set $\Delta$ would be contained in a repelling petal $\pi(A_n)$, contradicting that such petal contains no geodesic defined for an infinite time.
\end{proof}

\begin{rem}\label{rem:nonInjCase}
In Case~\eqref{item:pa_reg} of \Cref{lem:preAvoid}, $f$ is non-injective: we can for instance use the same argument as in the proof of \Cref{lem:Avoid2}.
Alternatively, we can apply the lemma to the restriction of $f$ to $\Delta_t = T_{\gamma(t)}(\Delta)$ where $T_v(z)=z+v$.
This implies that $y\in f(\Delta_t)$ for all $t>0$.
By applying this to cones of smaller opening we can even ensure that there is a sequence $z_n$ where $\phi(z_n)=y$ such that $\arg(z_n)$ tends to $\arg(a_\gamma)$ (the direction of $\gamma$).
\end{rem}

In the next statement, we extend \Cref{lem:preAvoid} to limits along more sequences.

\begin{lem}\label{lem:Avoid}
We assume that $X^*$ is not a full translation cylinder ($\C/\Z$).
Consider an affine immersion $f$ of a convex infinite open sector $\Delta \subset \mathbb{C}$ into $X^*$. We denote $\on{acc}$ for ``accumulation set''.
Then exactly one of the following holds:
\begin{itemize}
\item[(i)] for any geodesic $\gamma:(0,+\infty) \to \Delta$ not parallel to $\partial \Delta$, $\underset{t\to+\infty}{\on{acc}} f\circ \gamma(t) \subset f(\Delta)$;
\item[(iia)] $f(z)$ tends to a cylindrical singularity $y\in \cal S$ as $d(z,\partial \Delta)\to+\infty$;
\item[(iib)] $f(z)$ tends to an anti-conical singularity $y\in \cal S$ as $|z|\to+\infty$;
\item[(iic)] $f(z)$ tends to a swath of an irregular singularity $y\in \cal S$ as $|z|\to+\infty$;
\end{itemize}
Note that, since $X$ is compact, the accumulation set of a geodesic is never empty.
In Case~(i), every accumulation point belongs to $X^*$, i.e.\ is a regular point.
The cases are thus mutually exclusive.
\end{lem}

\begin{rem}
In the case where $X^*$ is a full translation cylinder, if
a radial geodesic $\gamma$ of the sector $\Delta$ is mapped to a closed geodesic of the cylinder, then radial geodesics that belong to distinct connected component of $\Delta \setminus \gamma$ converge to distinct ends of the cylinder.
\end{rem}

\begin{proof}[Proof of \Cref{lem:Avoid}]
We assume that we are not in Case~(i), i.e.\ that 
there exists a geodesic $\gamma:(0,+\infty) \to \Delta$, not parallel to $\partial \Delta$, and $y\in X$ an accumulation point of $f\circ\gamma(t)$ as $t\to+\infty$, such that $y\notin f(\Delta)$.
By \Cref{lem:preAvoid} either $y$ is an irregular singularity and $f\circ\gamma$ tends to one of its swathes, or $y$ is Fuchsian non conical and $f\circ\gamma$ tends to $y$.

Assume that $y$ is anti-conical Fuchsian.
By \Cref{lem:onde}, there is a neighborhood of $\infty$ in $\Delta$ that immerses in a neighborhood of $y$ that is isomorphic to the infinite end of a skew cylinder, possibly shifted.
Any geodesic $\gamma$ in $\Delta$ ends up staying in this neighborhood and $f\circ\gamma$ is mapped to a straight line in the skew cylinder model $\C^*/\langle\lambda\rangle$ or in the shifted model.
In both cases, the line tends to the infinite point of the model, so $\gamma(t)\tends y$.
%In the non-shifted model, on this line the modulus tends to $\infty$ while the lifted argument stays bounded, so it tends to the infinite end of the skew cylinder, so $\gamma(t)\tends y$.
%In the shifted model the 

In the case $y$ is a cylindrical Fuchsian singularity, consider a punctured neighborhood $V$ isomorphic via $
\phi: \mathbb{H}/\Z \to V$ to a semi-infinite translation cylinder. 
We start to apply \Cref{lem:irel} as in the anti-conical case, but here there is less angular room in the universal covering $\mathbb{H}$ of the model as there was for an anti-conical singularity, and if the sector $\Delta':=L(\Delta)$ contains a direction of argument in $(-\pi,0)$, then it will not be contained in $\mathbb{H}$, even if we remove a bounded subset from $\Delta$.
But in this case, denote $T(z)=z+1$. Since $\forall z\in \Delta' \cap T^{-1}(\Delta')$, the map
$f\circ L$ takes the same value at $z$ and $T(z)$, we can extend $f\circ L$ to an affine immersion $g:\C \to X^*$ such that $g\circ T=g$.
By \Cref{prop:PlaneImmersion} and our assumption that $y$ is cylindrical, $X^*$ is isomorphic to $\C/\Z$.
In the favorable case where $L(\Delta)$ has angular span contained in $(-\pi,0)$, then we can conclude as in the anti-conical case, except that for the convergence of $f(z)$ to $y$ we need the Euclidean distance from $z$ to $\partial \Delta$ to tend to $\infty$ to treat the case where one of the two boundary lines of $\Delta$ would be horizontal in the cylinder $\mathbb{H}/\Z$.

The last case is when $y$ is an irregular singularity and $f\circ\gamma$ tends to a swath.
In particular, it eventually enters, using the notations of \Cref{sub:geodesics}, an attracting petal $\pi(B_{n+1/2})$ and its real part in $B_{n+1/2}$ tends to $\infty$.
The set $B_{n+1/2}$ is a right half plane in the exponential-affine plane $\cal E$, and $\exp$ is an affine map on this set, realizing it as the universal covering of the complement of a disk in the canonical affine plane $\C$.
We then conclude, with an argument similar to the anti-conical case, that $f(z)$ tends to the corresponding swath as $z\tends\infty$ in $\Delta$.
\end{proof}

From the proof of~\Cref{lem:Avoid}, we get the following precisions: 
\begin{comp}\label{comp:iia}
In Case~(iia) of \Cref{lem:Avoid} above, there exists $L\in\AutC$ such that the angular span of $L(\Delta)$ is contained in $(0,\pi)$ and such that $f = \psi\circ \pi \circ L$ holds on $L^{-1}(\mathbb{H})\cap \Delta$, where $\pi:\mathbb{H}\to \mathbb{H}/\Z$ is the canonical projection and $\psi$ is an affine isomorphism from $\mathbb{H}/\Z$ to a neighborhood of $y$.
\end{comp}

\subsection{Holonomy set of an affine immersion}\label{sub:hol}

In the following, $f:\Delta \to X^{\ast}$ is the affine immersion of an open convex set $\Delta \subset \mathbb{C}$ to a finite type affine surface $(X,\nabla,\mathcal{S})$. 
If $f$ is not injective, for every pair of distinct points $x,y \in \Delta$ such that $f(x)=f(y)$, there is a complex affine automorphism that conjugates a neighborhood of $x$ with a neighborhood of $y$ in $\Delta$ while being compatible with $f$: in other words, there is $\phi \in \Aut(\mathbb{C})$ such that $y=\phi(x)$ and $f \circ \phi (z)=f(z)$ provided $z$ is close enough to $x$.
Since $\Delta$ is assumed convex, the set $\Delta\cap\phi^{-1}(\Delta)$ is connected, and by affine identity principle, the relation $f \circ \phi (z)=f(z)$ holds on all $\Delta\cap\phi^{-1}(\Delta)$.

\begin{defn}[Holonomy set]\label{def:holonomySet}
An automorphism $\phi \in \Aut(\mathbb{C})$ belongs to the \textit{holonomy set} $\Hol (f)$ of $f:\Delta \to X^{\ast}$ if it satisfies the following properties:
\begin{itemize}
    \item $\phi^{-1}(\Delta) \cap \Delta \neq \emptyset$;
    \item for any $z \in \phi^{-1}(\Delta) \cap \Delta$, $f \circ \phi (z) = f(z)$.
\end{itemize}
\end{defn}

\begin{rem}\label{rem:HolInjDefect}
The holonomy set encapsulates the defect of injectivity of $f$: the sets of pair $(x,y)$ such that $f(x)=y$ is retrieved as the set of $(x,\phi(x))$ for $\phi\in \Hol(f)$ and $x\in \Delta\cap \phi^{-1}(\Delta)$. 
In particular $f$ is injective if and only if its holonomy set is reduced to the identity.
\end{rem}

\begin{rem}\label{rem:compo}
Though $\on{Id}_\C\in \Hol(f)$ and $\forall \phi\in \Hol(f)$, $\phi^{-1}\in \Hol(f)$, the subset $\Hol(f)$ of $\Aut(\mathbb{C})$ does not have to be a subgroup.
Given $\phi,\psi\in\Hol(f)$, a sufficient condition for $\psi\circ\phi\in\Hol(f)$ is that there exists $z\in\Delta$ such that $\phi(z)\in \Delta$ and $\psi(\phi(z))\in \Delta$.
This automatically holds if $\Delta=\C$.
\end{rem}

We denote $\tra \cong (\mathbb{C},+)$ the translation subgroup of $\AutC$. 
We denote $\cal D$ the subgroup of translations and dilations of $\AutC$, i.e.\ of $\phi:z\mapsto az+b$ with $a>0$ and $b\in\C$.

\begin{lem}\label{lem:irg}
If $\Delta$ has infinite inner radius (i.e.\ contains disk of arbitrarily big radii) then $\tra \cap \Hol(f)$ is a discrete subgroup of $\tra$.

If $\Delta$ contains a sector, then $\cal D\cap\Hol(f)$ is a discrete subgroup of $\tra$.
\end{lem}
\begin{proof}
Given two translations $T_u$ and $T_v$ in $\Hol(f)$, the sufficient condition in \Cref{rem:compo} $T_{v}\circ T_{u}$ to be in $\Hol(f)$ applies to $z$ at the center of a ball contained in $\Delta$ and of radius $>|u|+|v|$.
We thus have a group.
If it were not discrete, considering a sequence $T_{u_n}$ in it, with $u_n\tends 0$ and $u_n\neq 0$, we would contradict that $f$ is locally injective.

In the case of $\cal D$ the proof is similar.
\end{proof}

This applies for instance to sectors.

\begin{rem}\label{rem:ige}
We can extend $f$ to an affine immersion $\hat f:\Delta \cup \phi^{-1}(\Delta)\to X^*$ by setting $\hat f(z)=f(z)$ when $z\in \Delta$ and $\hat f(z)=f(\phi(z))$ if $z\in \phi^{-1}(\Delta)$. The definition of $\Hol(f)$ ensures that these two values coincide on $\Delta\cap \phi^{-1}(\Delta)$.
Note however that trying to extend to $\Delta\cup \phi^{-1}(\Delta)\cup \phi^{-2}(\Delta)$ may fail, because $(\Delta\cup \phi^{-1}(\Delta))\cap\phi^{-1}(\Delta\cup \phi^{-1}(\Delta))$ may fail to be connected.
Consider for instance $\Delta$ to be a open half-plane not containing $0$ and $f$ to be a branch of $z\mapsto z^{5/2}$. 
\end{rem}

Given $\phi \in \Aut(\mathbb{C})$, we introduce the notations $a_{\phi}$ and $b_{\phi}$ where $\phi:z \mapsto a_{\phi} z + b_{\phi}$.

\begin{lem}\label{lem:HolFixNeg}
For any $\phi \in \Hol (f)$, the following two properties hold:
\begin{itemize}
    \item if $\phi$ has a fixed point in $\Delta$ then $\phi=\Id_\C$;
    \item the coefficient $a_{\phi}$ cannot be a negative real number.
\end{itemize}
\end{lem}

\begin{proof}
No element of the holonomy set different from the identity can fix a point $z$ of $\Delta$ because then we would find pairs of elements mapped by $f$ to the same point of $X^{\ast}$ in arbitrarily small neighborhoods of $z$, contradicting the fact that $f$ is an immersion.
\par
We deduce that no element $\phi$ of $\Hol (f)$ can be such that $a_{\phi} \in \mathbb{R}_{<0}$. Indeed, such an affine map would fix a unique point $z_{0} \in \mathbb{C}$. The intersection $\Delta \cap \phi^{-1}(\Delta)$ is a non-empty convex set any for any $z \in \Delta \cap \phi^{-1}(\Delta)$, the fixed point $z_{0}$ would lie in the interior of segment $[z,\phi(z)]$ and therefore would be a point of $\Delta$, contradicting the previous statement. 
\end{proof}

\subsection{Half-plane immersions}

We will completely classify half-plane immersions.

\medskip

We start by showing that the existence of a translation in the holonomy set of a half-plane immersion drastically constrains the geometry of the immersion.

\begin{lem}\label{lem:HPTranslation}
Consider an affine immersion $f:H \to X^{\ast}$ of an open half-plane $H$ into $X^*$.
If $\Hol (f)$ contains a non-zero translation, then at least one of the following statements holds:
\begin{itemize}
    \item[(i)] $f$ is the restriction of an affine immersion of the whole plane;
    \item[(ii)] $\Hol(f)$ is the group generated by a non-zero translation $\phi$ satisfying $\phi(H)=H$ and $f$  quotients-out to an affine embedding of the semi-infinite translation cylinder $H/\langle\phi\rangle$ as a punctured neighborhood of a cylindrical Fuchsian singularity of $X$.
\end{itemize}
\end{lem}

\begin{proof}
Note that a translation $\psi\in \AutC$ has vector parallel to $\partial H$ if and only if $\psi(H)=H$.
For any translation $\psi\in\Hol(f)$, if $b_{\psi}$ is not parallel to the boundary of $H$, then $\psi$ can be used to extend $f$ to an immersion $\hat f$ of the whole plane $\mathbb{C}$ satisfying $\hat f = \hat f \circ \psi$. Claim~(i) holds in this case.
Recall we denote $\tra$ the translation subgroup of $\AutC$ and that, according to \Cref{lem:irg}, $\tra\cap \Hol(f)$ is a (discrete) subgroup of $\AutC$.
Let $\tra(H)\subset\tra$ be the subgroup of translations preserving $H$
and $\Lambda:=\Hol(f)\cap \tra(H)$, which is also a discrete subgroup of $\AutC$.
\par
In the rest of this proof we assume that Claim~(i) does not hold: in particular, every translation in $\Hol(f)$ preserves $H$, i.e.\ $\Lambda \subset \tra(H)$.
By hypothesis, $\Hol (f)$ contains a non-zero translation, so $\Lambda$, as a non-trivial discrete subgroup of $\tra(H)$, is cyclic. 
Without loss of generality, we will assume that $\phi$ is a generator of $\Lambda$.
\par
Still under the hypothesis that Claim~(i) does not hold, we now prove that every element $\psi \in \Hol (f)$ is a translation.
We saw in \Cref{lem:HolFixNeg} that $a_\psi$ cannot be a negative number. It follows that the intersection of $H$ and $\psi^{-1}(H)$ is a sector, of angle $\leq\pi$.
Deep enough in the sector, there exists $z\in H$ such that $\phi(z)$ is also in the sector, so both $\psi(z)$ and $\psi(\phi(z))$ belong to $H$.
Since the map $\psi \circ \phi \circ \psi^{-1}\in\AutC$, sends $\psi(z)\in H$ to $\psi(\phi(z))\in H$ for all $z$ in an open set, it belongs to $\Hol(f)$.
Note that $\psi \circ \phi \circ \psi^{-1}$ is the translation $z \mapsto z+a_{\psi}b_{\phi}$.
Since $\phi:z\mapsto z+b_\phi$ is the generator of $\Lambda$,
we get that $a_{\psi} \in \mathbb{Z}^{\ast}=\Z-\{0\}$.
As the same reasoning can be carried out with $\psi^{-1}$, we deduce that both $a_{\psi}$ and $\frac{1}{a_{\psi}}$ belong to $\mathbb{Z}^{\ast}$.
By \Cref{lem:HolFixNeg}, $a_\psi$ cannot be a negative number, so $a_{\psi}=1$.
Hence $\text{Hol}(f) = \Lambda = \langle\phi\rangle$.

By \Cref{rem:HolInjDefect}, $f(z)=f(z')$ if and only if $z'\in \Lambda z$, so $f$ quotients-out to an injective map through $H\mapsto H/\Lambda$.
The affine surface $H/\Lambda$ is isomorphic to $\mathbb{H}/\Z$.
As a Riemann surface, $\mathbb{H}/\Z$ is isomorphic to the punctured unit disk $\D^*$ via the exponential map.
By \Cref{lem:embedFuchsian} the factor map is an affine isomorphism from $\mathbb{H}/\Z$ to a punctured neighborhood of a cylindrical singularity of $X$.
\end{proof}

Since the commutators of affine functions are translations, we now deduce from \Cref{lem:HPTranslation} that the elements of the holonomy set of a half-plane immersion commute.

\begin{cor}\label{cor:HPCommutation}
Consider an affine immersion $f:H \to X^{\ast}$ of an open half-plane $H$ into $X^*$. Then the elements of $\Hol(f)$ pairwise commute.
\end{cor}

\begin{proof}
If $\Hol (f)$ contains a non-zero translation, each of the two alternatives stated in \Cref{lem:HPTranslation} implies that the elements of $\Hol(f)$ pairwise commute (see \Cref{sub:AffImmPlanes} for details on affine surfaces admitting the immersion of a plane). In the rest of the proof, we will assume that $\Hol(f)$ does not contain any translation.
\par
We consider a pair $\phi,\psi$ of elements in $\Hol(f)$. We already know from \Cref{lem:HolFixNeg} that $a_{\phi}$ and $a_{\psi}$ cannot be negative real numbers. Up to replacing $\phi$ or $\psi$ by its inverse, we can assume that $\arg(a_{\phi}) \in [0,\pi)$ while $\arg(a_{\psi}) \in (-\pi,0]$.
\par
We recall that, given $\alpha,\beta\in\Hol(f)$, a sufficient condition for $\beta\circ\alpha\in\Hol(f)$ is that there exists $z\in H$ such that $\alpha(z)\in H$ and $\beta(\alpha(z))\in H$.
\par
By the restriction we put on $\arg(a_\phi)$ and $\arg(a_{\psi})$, the maps $\phi$, $\psi$ satisfy the condition above so $\psi\circ\phi\in \Hol(f)$.
Similarly, $\psi^{-1}\circ\phi^{-1}\in \Hol(f)$.
Moreover, the argument of the linear factors of $\psi\circ\phi$ belongs to $(-\pi,\pi)$ and the argument for $\psi^{-1}\circ\phi^{-1}\in \Hol(f)$ is the opposite, so again we can apply the criterion and deduce that the commutator $[\phi,\psi]$ belongs to $\Hol(f)$.
The argument of the linear factor of $[\phi,\psi]$ is $1$, i.e.\ it is a translation.
By hypothesis, the translation $[\phi,\psi]$ has to be trivial so $\phi$ commutes with $\psi$.
\end{proof}
%La preuve du lemme ci-dessus fonctionne probablement encore avec un cône convexe infini dans l'hypothèse 2 (pas de translation). On pourrait peut-être le scinder en 2.

Every set of pairwise commuting elements of $\AutC$ is contained in a maximal commutative subgroup of $\AutC$, which are the group of translations or the fixator of some point of $\C$.

\medskip

We get a complete classification of half-plane immersions in finite-type affine surfaces.

\begin{prop}\label{prop:halfPlaneImmersion}
Consider an affine immersion $f:H \to X^{\ast}$ of an open half-plane $H$ into $X^*$. Then either:
\begin{enumerate}[(i)]
\item\label{item:embed} $f$ is an embedding;
\item\label{item:exc} $X^*$ is the full translation cylinder $\C/\Z$, a translation torus $\C/\Lambda$ or an affine torus $\cet/\Lambda$;
\item\label{item:cylEmbed} $f$ factors to an embedding through a cylindrical Fuchsian singularity standard neighborhood;
\item\label{item:antiConeEmbed} $f$ factors to an embedding through an anti-conical Fuchsian singularity standard neighborhood, of opening angle $<\pi$;
\item\label{item:ReebEmbed} $f$ factors to an embedding through a Reeb cylinder of angle $\pi$.
\end{enumerate}
More precisely, denote by $\tra$ the translation subgroup of $\AutC$.
In each corresponding case:
\begin{enumerate}[(i)]
\item 
  \begin{itemize}
  \item $\Hol(f)$ is reduced to the identity;
  \item as $z$ tends to $\infty$, $f(z)$ tends to an anti-conical singularity of opening angle $\geq \pi$ or to a swath of an irregular singularity;
  \end{itemize}
\item
  \begin{itemize}
  \item $\Hol(f)$ is respectively equal to a non-trivial cyclic subgroup of $\tra$, a lattice of $\tra$, the following a subset of the projection to $\AutC$ of a lattice of the translation group acting on $\cal E$: $\setof{u^p\circ v^q}{p,q\in\Z,\ p\theta_u+q\theta_v\in(-\pi,\pi)}$, $u,v$ is a pair of translations of $\cal E$ of vectors of imaginary part $\theta_u$ and $\theta_v$ such that the affine torus $\cet/\Lambda$ is isomorphic via $\exp$ to $\cal E/\langle u,v\rangle$;
  \item $f = \pi \circ \iota$ where $\pi$ is the quotient map $\C\to\C/\Z$, $\C\to\C/\Lambda$ or $\cet\to\cet/\Lambda$ and $\iota$ an affine bijection from $H$ to a half-plane in $\C$ or $\cet$;
  \end{itemize}
\item
  \begin{itemize}
  \item $\Hol(f)$ is a non-trivial cyclic subgroup of $\tra$; 
  \item $f = g \circ \pi$ where $\pi : H \to H/\langle T \rangle$ is the quotient map and $T$ is a translation preserving $H$ (so $H/\langle T \rangle \cong \mathbb{H}/\Z$ is a semi-infinite cylinder) and $g$ is an affine embedding from $H/\langle T \rangle$ onto a neighborhood of a cylindrical Fuchsian singularity $p\in\cal S$;
  \end{itemize}
\item
  \begin{itemize}
  \item $\Hol(f)$ is of the form $\setof{\phi^n}{-k<n<k}$ for some $\phi:z\mapsto az+b\in\AutC$ with $\arg(a)\in(0,\pi)$ whose fixed point is not in $H$ and where $k=\lceil\pi/\arg(a) \rceil$;
  \item $f = g \circ \pi \circ \iota$ where $\iota$ is an affine embedding of $H$ in $\cet$,
  $\pi : \cet \to \cet/\langle \tilde\psi \rangle$ is the quotient map where $\tilde\psi$ sends $(r,\theta)$ to $(|a|r,\theta+\alpha)$ where $\alpha$ is the representative in $(0,\pi)$ of $\arg a$ (so $\cet/\langle \tilde\phi \rangle$ is a skew cone) and $g$ is an affine embedding from $\pi\circ\iota(H)$ onto a neighborhood of an anti-conical Fuchsian singularity $p\in\cal S$;
  \end{itemize}
\item
  \begin{itemize}
  \item $\Hol(f)=\langle \phi\rangle$ for some dilation $\phi:z\mapsto az+b\in\AutC$ with $a\in(0,1)$ whose fixed point does not belong to $H$.
  \item $f = g \circ \pi|_H$ where $\pi : H' \to H'/\langle \phi \rangle$ is the quotient map, $H'$ is the half-plane containing $H$ and whose boundary passes through the fixed point of $\phi$ (so $H'/\langle \phi \rangle$ is a Reeb cylinder of angle $\pi$) and $g$ is an affine embedding from $H'/\langle \phi \rangle$ into $X^*$;
  \end{itemize}
\end{enumerate}
Case~\eqref{item:exc} has non-empty intersection with Cases~\eqref{item:cylEmbed} and~\eqref{item:ReebEmbed}.
\par
Finally, $f$ extends continuously at every point of $\partial H$, except, possibly, one point in Case~\eqref{item:ReebEmbed} or in Case~\eqref{item:exc} if $X^*$ is an affine torus.
\end{prop}

\begin{proof}
By \Cref{rem:HolInjDefect}, $f$ is an embedding if and only if $\Hol(f)$ is reduced to the identity.
In this case, a half-line in $H$ is mapped to the support of a geodesic $\gamma$ of $X$. Consider an accumulation point $y\in X$ of this geodesic as in \Cref{lem:Avoid}.
If $y$ were regular, \Cref{rem:nonInjCase} would give us that $f$ is non-injective.
So we are in one of the cases~(iia), (iib) or~(iic) of \Cref{lem:Avoid}.
Case~(iia) where $y$ is a cylindrical singularity is incompatible with $f$ being injective: indeed by \Cref{comp:iia} applied to $\Delta=H$, there exists $L\in\AutC$ such that $L(H)$ is an upper half-plane and such that $f\circ L^{-1}=\psi\circ\pi$ holds on $L^{-1}(H)\cap\mathbb{H}$ where $\pi:\C\to\C/\Z$ is the canonical projection; but then $f$ cannot be injective, leading to a contradiction.
In Case~(iib) if the anti-conical singularity had an opening angle $\alpha$ strictly less than $\pi$, we would obtain a similar contradiction to injectivity: indeed by the proof of \Cref{lem:onde}, there is a neighborhood $W$ of $\infty$ in $H$ such that on $W$, $f$ factors as $\phi\circ \pi_\lambda \circ \iota \circ L$ where $L\in\AutC$, $\iota:L(W)\to\cet$ is a section of $\pi_0:\cet\to \C^*$,  $\pi: \cet \to \cet / \langle\lambda\rangle$ and $\phi$ is an affine isomorphism from a neighborhood of the infinite end of $\cet / \langle\lambda\rangle$ to an punctured neighborhood of $y$.
However, $L(W)$ is a half-plane minus a bounded set, and contains a neighborhood of $\infty$ in a sector of opening angle bigger than $\alpha$.
The map $\pi_\lambda\circ\iota$ is not injective on such a set.

In the case where $X^*$ is the full translation cylinder, a translation torus or an affine torus, the fact that $\pi$ is a universal covering and $H$ simply connected implies that $f$ factors as $\pi\circ\tilde f$. The map $\tilde f$ being locally affine, defined on a half-plane, and taking values in $\C$ or $\cet$, it must be injective.
The determination of $\Hol(f)$ in each case is then easy: in the first two cases it must coincide with the deck transformation group of $\pi$.
In the case of an affine torus, let $\tilde H=\tilde f(H)$; the transformation $\mu$ of $\cet$ associated to $u^p\circ v^q$ satisfies $\mu(\tilde H)\cap \tilde H\neq \emptyset$ if and only if the rotation part of $u^p\circ v^q$ lies in $(-\pi,\pi)$.

In the rest of this proof, except in the last paragraph, we assume that we are not in Cases~\eqref{item:embed} or~\eqref{item:exc}.
The set $\Hol(f)$ in particular contains an element different from the identity.
We saw in \Cref{cor:HPCommutation} and the paragraph after it that $\Hol(f)$ is either contained in the translation subgroup of $\AutC$ or in the fixator of a point $z_0$ in $\C$ that does not belong to $H$.

In the first case, \Cref{lem:HPTranslation} tells us that $f$ factors to an embeddingt through a cylindrical Fuchsian singularity neighborhood and that $\Hol(f)$ is a non-trivial cyclic subgroup of $\tra$.

In the second case,
we first assume that $\Hol(f)$ only contains dilations (fixing $z_0$).
By \Cref{lem:irg}, $\Hol(f)$ is a discrete subgroup of $\AutC$.
Let $\phi\in\AutC$ be a generator and denote $\Lambda=\langle\phi\rangle=\Hol(f)$.
Then by \Cref{rem:HolInjDefect}, $f(z)=f(z')$ $\iff$ $z\in\Lambda z'$, So $f(H)$ is isomorphic to the image of $H$ in $\C \setminus \{z_0\}/\Lambda$, which is equal to $S/\Lambda$ where $S$ is the half-plane containing $H$ and whose boundary passes through $z_0$.
We get a Reeb cylinder of opening $\pi$.

Finally we assume that $\Hol(f)$ contains a similarity $s_0$ fixing $z_0$ but which is not a dilation.
We pass to the exponential-affine plane $\cal E$ using the change of variable $z=z_0+\exp(w)$, pulling back $H$ to a subset $U$ of $\cal E$ that is either a horizontal strip of height $\pi$ if $z_0\in\partial H$ or a shape of equations $y_1-\pi/2<y<y_1+\pi/2$ and $x-x_1 > \ln(1/\cos (y-y_1))$, denoting $w=x+iy\in\cal E$.
Let
\[g=f\circ(z_0+\exp).\]
The set $\Hol(f)$ is in bijective correspondence with a set, which we denote $\Hol(g)$, of translations $T_u : w\mapsto w+u$ of $\cal E$ with the property that $T_u\in\Hol(g)$ $\iff$ ($T_u^{-1}(U)\cap U\neq\emptyset$ and $g\circ T_u =g$ holds on this intersection (which is connected)), and we know that $g(w)=g(w')$ iff $\exists T_u\in\Hol(g)$ such that $T_u(w)=w'$.

If $\Hol(g)$ contains two elements $T_u$ and $T_v$ with $\R u+\R v=\C$, then $g$ can be extended to an affine immersion $\hat g:\cal E\to X^*$.
It passes to the quotient to an affine immersion of an affine torus $Y$ to $X^*$: $\tilde g: Y\to X^*$.
The image being both open and compact, we get that $X^*=X=\tilde g(Y)$.
Since $Y$ and $X$ are compact and $\tilde g$ is a local diffeomorphism, it is a covering.
Hence $\hat g$, which is defined on the simply connected set $\cal E$, is a universal covering of $X$.
Its deck transformation group is a subgroup of the automorphism group of the affine surface $\cal E$, i.e.\ the group of translations, and contains $T_u$ and $T_v$.
It is thus a lattice and $X$ is an affine torus.

Otherwise, $\Hol(g) \subset  T_{\R u_0}$ for some $u_0\in\C^*$ with $\Im(u_0)\neq 0$ since $s_0\in\Hol(f)$.
We can use any non-trivial $T_u\in\Hol(g)$ to extend $g$ to an affine immersion $\hat g$ defined on $\hat U:= \bigcup_{n\in\Z}T_u^n(U)$, which contains a half-plane whose boundary is parallel to $\R u_0$.
Then by \Cref{rem:compo}, $\Hol(\hat g)$ is a subgroup of $\setof{T_{s u_0}}{s\in\R}$.
If it were a dense subgroup we would contradict local injectivity of $g$, so we have a cyclic subgroup. Let $v\in\C$ with $\Im(v)>0$ such that $T_{v}$ generates $\Hol(\hat g)$.
Denote $\pi:\cal E \to \cal E/\langle T_v \rangle$ the canonical projection. Note that the image is a skew cylinder and that $\pi(\hat U)=\pi(U)$ is a neighborhood of its infinite end.
Then $\hat g = \psi\circ\pi$ for some $\psi$ which is injective by property of $\Hol(\hat g)$ and is an embedding to a neighborhood of a pure anti-conical singularity by \Cref{lem:embedFuchsian}.
Finally, $f = \psi\circ\pi\circ {(z_0+\exp)|_U}^{-1}$.
Concerning $\Hol(f)$, it is the set of $z\mapsto z_0+(z-z_0)e^{n v}$ where $n\in\Z$ is such that $T_v^n\in \Hol(g)$, which is the case iff $-\pi<n\times \Im (v)<\pi$.

In Case~\eqref{item:embed}, \Cref{cor:Embedding} ensures that $f$ extends continuously to $\partial H$.
In the other cases, given a point $z_0\in \partial H$ and $\eps>0$, the set $\Delta' = H\cap B(z_0,\eps)$ is convex. If $f$ is injective on $\Delta'$, we can deduce by applying \Cref{cor:Embedding} to the $f|_{\Delta'}$ that $f$ extends continuously at $z_0$.
This applies at every point of $\partial H$ for Cases~\eqref{item:cylEmbed} and \eqref{item:antiConeEmbed}.
In Case~\eqref{item:ReebEmbed}, this also applies to every point of $\partial H\setminus\{z_0\}$.
Last, in Case~\eqref{item:exc}, the immersion is explicit and always extends for the full translation cylinder and the translation tori, and in the case of an affine torus $\cet/\Lambda$, the extension is possible on $\partial H\setminus\{0\}$.
\end{proof}

\subsection{Sector immersions}

We start by a study of the behavior at infinity.

\begin{lem}\label{lem:coneTurnProlong}
Let $\Delta\subset \C$ be an infinite convex open sector of opening angle $\theta\in(0,\pi)$ and $f:\Delta\to X^*$ an affine immersion. Assume that $\Hol(f)$ contains an element $\phi$ whose scaling factor $a_\phi$ satisfies $0<|\arg (a_\phi)|<\theta$.
Then exactly one of the following occurs:
\begin{itemize}
  \item[(i)] $X^*=X$ is an affine torus;
  \item[(ii)] there is a neighborhood $V'$ of $\infty$ in $\Delta$ on which $f$ quotients-out to the embedding of a standard neighborhood of an anti-conical singularity $p$ of opening angle dividing $|\arg a_\phi|$; $f(z)\tends p$ as $|z|\to +\infty$.
\end{itemize}
More precisely in Case~(ii), $f|_{V'} = g \circ \pi \circ \iota|_{V'}$ where $\iota$ is an affine embedding of $\Delta$ in $\cet$, $\pi : \cet \to \cet/\langle \tilde\lambda \rangle$ and $g$ is an affine embedding from $\pi\circ\phi(\Delta)$ onto a neighborhood of $p$.
\end{lem}

\begin{proof}
Let $\alpha=|\arg(a_\phi)|\in(0,\pi)$.
The proof of Case~(\ref{item:antiConeEmbed}) of
\Cref{prop:halfPlaneImmersion} adapts to this case, with the following precaution: we first restrict $f$ to a set $V'$ of the form $\Delta'-B(0,R)$ with $\Delta'$ a sector based on the fixed point of $\phi$ and opening angle $\theta'\in (\alpha,2\alpha)$ and such that $V'\subset\Delta$.
We get a factorization for the restriction, and we extend the factorization to one $f$ that holds in a neighborhood of $\infty$ in $\Delta$, by using the affine identity principle.
(Alternatively we can use $V'$ to extend the restriction to a half-plane containing a neighborhood of $\infty$ in $\Delta$ and then apply \Cref{prop:halfPlaneImmersion} so as to conclude with the identity theorem.)
\end{proof}

The lemma below can be deduced from its version for planes (\Cref{lem:HPTranslation}).

\begin{lem}\label{lem:HolConeTranslation}
Consider an affine immersion $f:\Delta \to X^{\ast}$ of an infinite open convex sector $\Delta$ into $X^*$.
Denote $z_0$ the apex of the sector and write $\Delta = T_{z_0}(\vec\Delta)$.

\noindent 1. The set $\tra \cap \Hol(f)$, of elements of $\Hol(f)$ that are translations, is invariant by composition, thus forms a subgroup of $\AutC$. It is discrete.

\noindent 2. If $\Hol (f)$ contains a non-zero translation, then at least one of the following statements holds:
\begin{enumerate}[(i)]
  \item $f$ is the restriction of an affine immersion of the whole plane;
  \item $\tra\cap\Hol(f)=\langle\phi\rangle$ with $\phi=T_v$, $v\notin\vec\Delta$ and $f$ quotients-out by $\langle\phi\rangle$ to an affine immersion $\psi$, i.e.\ $f=\psi \circ \pi$ holds on $\Delta$ where $\pi:\C\to\C/\langle\phi\rangle$ and $\psi$ is defined on $\pi(\Delta)$.
  Moreover let $H$ be the maximal half-plane contained in $\bigcup_{n\in\Z} \phi^n(\Delta)$;
  the immersion $\psi$ is injective in the semi-infinite cylinder $H/\langle\phi\rangle$ and sends it to a standard neighborhood of a cylindrical Fuchsian singularity of $X$.
\end{enumerate}
\end{lem}
\begin{proof}
Recall that we denote $\tra$ the translation subgroup of $\AutC$ and 
that by \Cref{lem:irg},
\[\Lambda:=\tra\cap\Hol(f)
\]
is a discrete subgroup of $\tra$.

If the vector $v$ of a translation in $\Lambda$ belongs to $\vec\Delta$, then we can use $f\circ T_{nv}$ with successive values of $n\in\N$ to extend $f$ to a whole plane immersion.

On the opposite, if no non-trivial element of $\Lambda$ has its vector $v$ that belongs to $\vec\Delta$, then $\Lambda$ is cyclic ($\vec \Delta$ intersects any lattice).
Let $\phi_0 = T_{v_0}$ be a generator:
\[ \Lambda = \langle \phi_0\rangle
.\]
We can still define an extension $\hat f$ of $f$ using $f\circ T_{nv_0}$ for all $n\in\Z$, but its domain now is not the whole plane.
It contains however a maximal half-plane $H$ with $\partial H$ parallel to $v_0$.
See \Cref{fig:cep}.
Note that $\hat f|_H$ and $f$ coincide on the non-empty convex set $H\cap \Delta$.
We can apply \Cref{lem:HPTranslation} to $\hat f|_H$:

\begin{figure}[htbp]
  \begin{center}
    \begin{tikzpicture}
      \node at (0,0) {\includegraphics[scale=0.8]{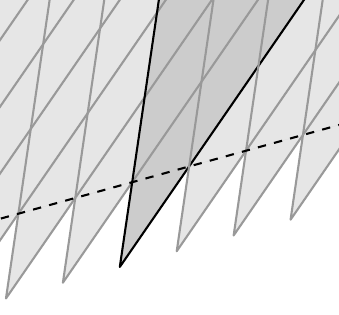}};
    \end{tikzpicture}    
  \end{center}
  \caption{The initial sector $\Delta$ on which $f$ is defined is in dark gray. We use the translation $T_{v_0}\in\Hol(f)$ to define an extension $\hat f$ that is defined on the light and dark gray set. It contains a maximal half-plane $H$ bounded by the dashed line.}
  \label{fig:cep}
\end{figure}

Either $\hat f|_H$ is the restriction of an affine immersion $\tilde f$ of the whole plane. Since $f$ and $\tilde f$ are both defined on $\Delta$ and coincide on the non-empty open subset $\Delta\cap H$, they coincide on $\Delta$, so $\tilde f$ is an extension of $f$.

Or $\hat\Lambda:=\Hol(\hat f|_H)$ is the group generated by a non-zero translation $\phi$ satisfying $\phi(H)=H$, i.e.\ of vector $u$ parallel to $v_0$, and $\hat f|_H$ quotients-out to an affine embedding of the semi-infinite translation cylinder $H/\langle\phi\rangle$ as a punctured neighborhood of a cylindrical Fuchsian singularity of $X$.
Note that $\phi$ is necessarily in $\Hol(f)$, i.e.\ in $\Lambda$, because the open set $(\Delta\cap H)\cap \phi^{-1}(\Delta\cap H)$ is non-empty.
So $u\in v_0\Z$.
Conversely, $\phi_0\in\hat\Lambda$ because $(\Delta\cap H)\cap \phi_0^{-1}(\Delta\cap H)$ is non-empty.
So $v_0\in u\Z$.
It follows that
\[ \hat\Lambda=\langle \phi_0\rangle=\Lambda
.\]
Note that $\Delta\setminus H$ has no two points related by $\langle \phi_0\rangle$; it follows that $f$ quotients-out by $\hat\Lambda$ on the whole set $\Delta$.
\end{proof}

We have an analogue of \Cref{prop:halfPlaneImmersion}, but weaker in that some of the statements only hold near $\infty$.

\begin{prop}\label{prop:coneImmersion}
Consider an affine immersion $f:\Delta \to X^{\ast}$ of an infinite open convex sector $\Delta$ of opening angle $\theta\in(0,\pi]$ into $X^*$. Then either:
\begin{enumerate}[(i)]
\item\label{item:cone:embed} $f$ is an embedding in a neighborhood of $\infty$ and $f(z)$ tends, as $|z|\to +\infty$ to an anti-conical Fuchsian singularity of opening angle $\geq \theta$ or to a swath of an irregular singularity;
\item\label{item:cone:exc} $X^*$ is the full translation cylinder $\C/\Z$, a translation torus $\C/\Lambda$ or an affine torus $\cet/\Lambda$; $f$ factors as $f=\pi\circ\phi$ where $\phi$ is an affine embedding of $\Delta$ in $\C$ or $\cet$;
\item\label{item:cone:cylEmbed} $f$ quotients-out near $\infty$ by some $\phi_*\in\AutC$ to an embedding through a cylindrical Fuchsian singularity standard neighborhood;
\item\label{item:cone:antiConeEmbed} $f$ quotients-out near $\infty$ by some $\phi_*\in\AutC$ to an embedding through a standard neighborhood of an anti-conical Fuchsian singularity $p$ of opening angle $\leq\theta$; $f(z)\tends p$ as $|z|\to+\infty$;
\item\label{item:cone:ReebImmerse} $f$ quotients-out near $\infty$ by some $\phi_*\in\AutC$ to an immersion through a dilation cylinder of opening angle $\geq \theta$.
\end{enumerate}
\end{prop}

In~\eqref{item:cone:ReebImmerse} the cylinder could very well be non-embedded (only immersed): see \Cref{subsub:CylindersAroundGeodesics}.

\begin{proof}
Case~\eqref{item:cone:embed} is treated completely similarly as in the \Cref{prop:halfPlaneImmersion}, as it involves \Cref{lem:Avoid,rem:nonInjCase,comp:iia} which are all valid on infinite open convex sectors and \Cref{lem:onde} which is valid on unbounded open convex sets.

Case~\eqref{item:cone:exc} is treated as Case~\eqref{item:exc} of \Cref{prop:halfPlaneImmersion}. In the rest of this proof we assume that we are not in Case~\eqref{item:cone:exc}.

If there is some $\phi \in \text{Hol}(f)$ such that $\arg(a_{\phi}) \in (-\theta,0)\cup(0,\theta)$, then by \Cref{lem:coneTurnProlong} we are in particular in Case~\eqref{item:cone:antiConeEmbed}.
In the rest of this proof we assume we are not in this case either, so in particular that for all $\phi\in\Hol(f)$, either $\arg a_\phi=0$ or $|\arg(a_\phi)|\geq \theta$.

In \Cref{prop:halfPlaneImmersion}, $\Delta=H$ was a plane and we could invoke \Cref{cor:HPCommutation}, but here $\Delta$ is any sector and we proceed differently.
One annoying hypothetical situation is illustrated on \Cref{fig:ci}. However we will see that there is not enough room in a finite type affine surface for this kind of immersions to exist.
Consider any radial geodesic $\gamma$ of $\Delta$, say the one whose direction bisects the sector, and apply \Cref{lem:preAvoid}: either 1.\ the accumulation set of $f\circ \gamma(t)$ as $t\to+\infty$ is contained in $f(\Delta)$, or 2.\ $f\circ\gamma(t)$ tends to a singularity, which is 2.a Fuchsian non conical or 2.b irregular (in the latter case it tends to a swath). 

\begin{figure}[htbp]
\begin{center}
\begin{tikzpicture}
\node at (0,0) {\includegraphics[scale=0.5]{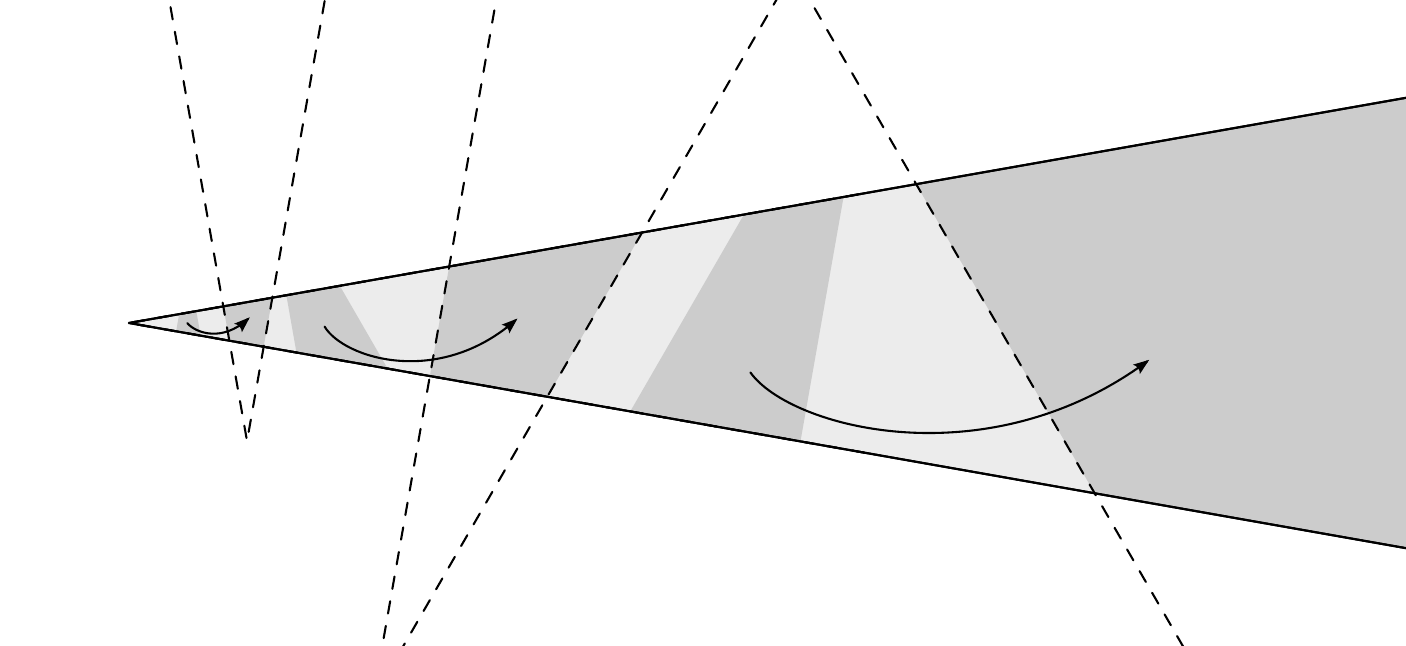}};
\node at (-4.7,-0.4) {$\phi_1$};
\node at (-3.2,-0.8) {$\phi_2$};
\node at (1.5,-1.6) {$\phi_3$};
\end{tikzpicture}
\end{center}
\caption{A hypothetical sector in $\C$, whose immersion $f$ in an affine surface would have infinitely many non-identity elements $\phi\in \Hol(f)$, all with disjoint ranges and domains, and no other ones. This situation in fact cannot happen. One could effectively define an affine surface $S$ by quotienting $\Delta$ by $\phi_1$, $\ldots$, $\phi_n$ but it is not of finite type, and actually cannot embed in a finite-type affine surface as will follow from \Cref{prop:coneImmersion}.}
\label{fig:ci}
\end{figure}
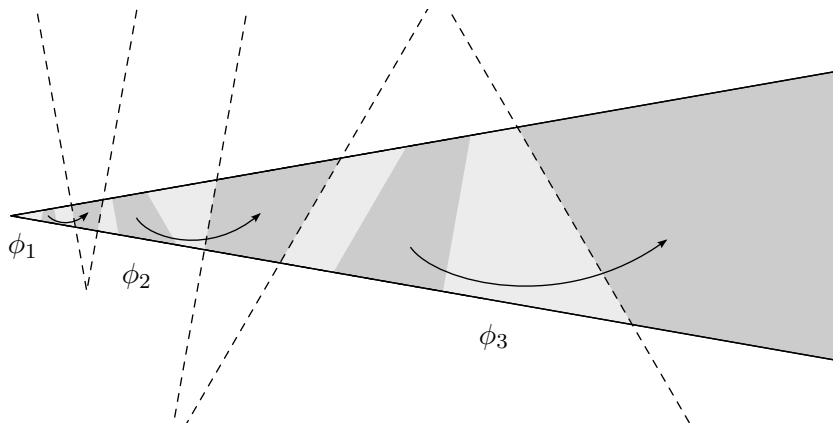

If 2.\ above holds, then we proceed as in \Cref{lem:Avoid}:
there is a neighborhood of $\infty$ in $\Delta$ that immerses into a neighborhood of $y$ that is isomorphic to the infinite end of a skew cylinder in case~2.a by \Cref{lem:onde}, or into a swath neighborhood isomorphic to a right half-plane in $\cal E$ in case~2.b.
In case~2.b, $f$ is thus actually injective in a neighborhood of $\infty$ hence we are in Case~\eqref{item:cone:embed}.
This is also the case in case~2.a if the opening angle $\theta$ of the skew cone is $< |\arg(a_\phi)|$.
By hypothesis, we excluded the case where $\theta > |\arg(a_\phi)|$.
If $\theta = |\arg(a_\phi)|$ we may be either in Case~\eqref{item:cone:antiConeEmbed} or in Case~\eqref{item:cone:embed}.

If 1.\ above holds, then consider an accumulation point $y$ of $f\circ\gamma$ and a small neighborhood $V$ around $y$ and an affine isomorphism  $\psi:V\to D=B(0,1)$, with $\phi(y)=0$.
Branches of the geodesic in $V$ accumulate on $y$, so we can extract a sequence of branches whose direction converges to some limit direction.
For each such branch there exists (apply \Cref{lem:irel}) a subset of $\Delta$ mapped by $\psi\circ f$ to the intersection of $D$ with a sector of opening $\theta$, bisected by the branch, and based on one of the endpoints of the branch, see \Cref{fig:ba}.
As soon as the branch passes close enough to $y$, we know that the intersection above will contain $0$.
Hence $y$ has arbitrarily large preimages $z_n$ by $f$.
Moreover since the branches have a convergent direction, 
if we take $m$ and $n$ big enough, the map $\phi\in \Hol(f)$ such that $\phi(z_m)=z_n$ satisfies that $\arg(a_\phi)$ is close to $0$.
From our standing assumption, we deduce that $\arg(a_{\phi})=0$,
so $\phi$ is a translation or a dilation fixing some point $z_{1}$ outside of $\Delta$ (see \Cref{lem:HolFixNeg}).
We choose such $m$ and $n$ and denote $\phi_*$ the element $\phi$ constructed above.

\begin{figure}[htbp]
  \begin{center}
    \begin{tikzpicture}
      \node at (0,0) {\includegraphics[scale=0.75,angle=90]{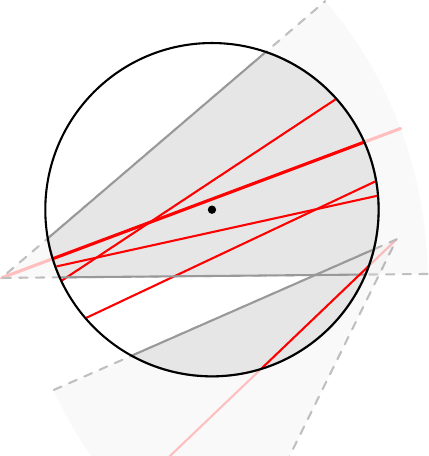}};
    \end{tikzpicture}    
  \end{center}
  \caption{The disk $D$ in the proof of \Cref{prop:coneImmersion}; in red, branches of the geodesic; in gray, images of the sector $\Delta$ around two of these branches, by the affine maps $L\in\AutC$ such that locally $\psi\circ f\circ\gamma=L$. When the branch is close to the center $0$, then $0\in L(\Delta)$.}
  \label{fig:ba}
\end{figure}

The segment $[z_m,z_n]$ is entirely contained in $\Delta$  and is mapped by $f$ to a closed geodesic, because $\phi_*'\in\R_{>0}$, that belongs to some immersed translation or dilation cylinder (see \Cref{sub:cylindersAndSkewCones}).
If $\phi_*$ is a translation, then we are in the situation of \Cref{lem:HolConeTranslation} and Claim~\eqref{item:cone:cylEmbed} holds (we cannot be in the situation~(i) of \Cref{lem:HolConeTranslation} because then by \Cref{prop:PlaneImmersion} $X^*$ is a full translation cylinder or a translation torus, which we exclude by hypothesis, or a translation plane, in which case $f$ is injective, and we already covered this case).

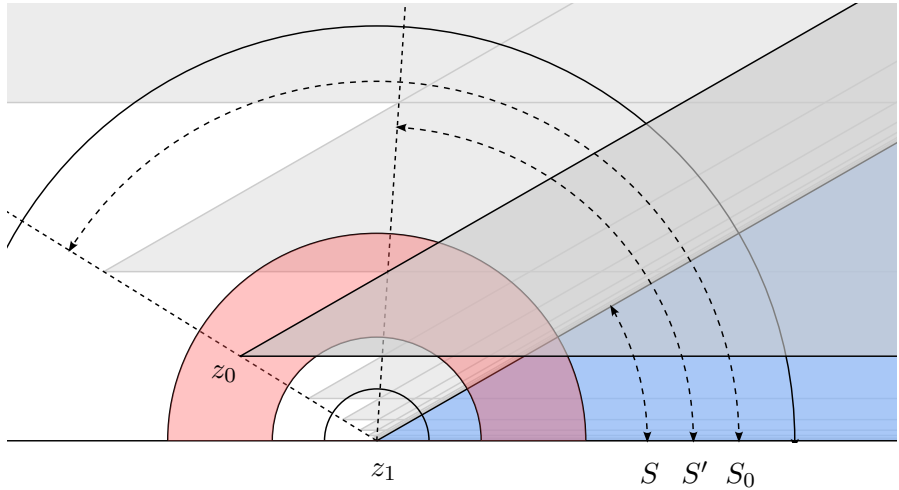
\begin{figure}[htbp]
  \begin{center}
    \begin{tikzpicture}
      \node at (0,0) {\includegraphics[scale=0.55]{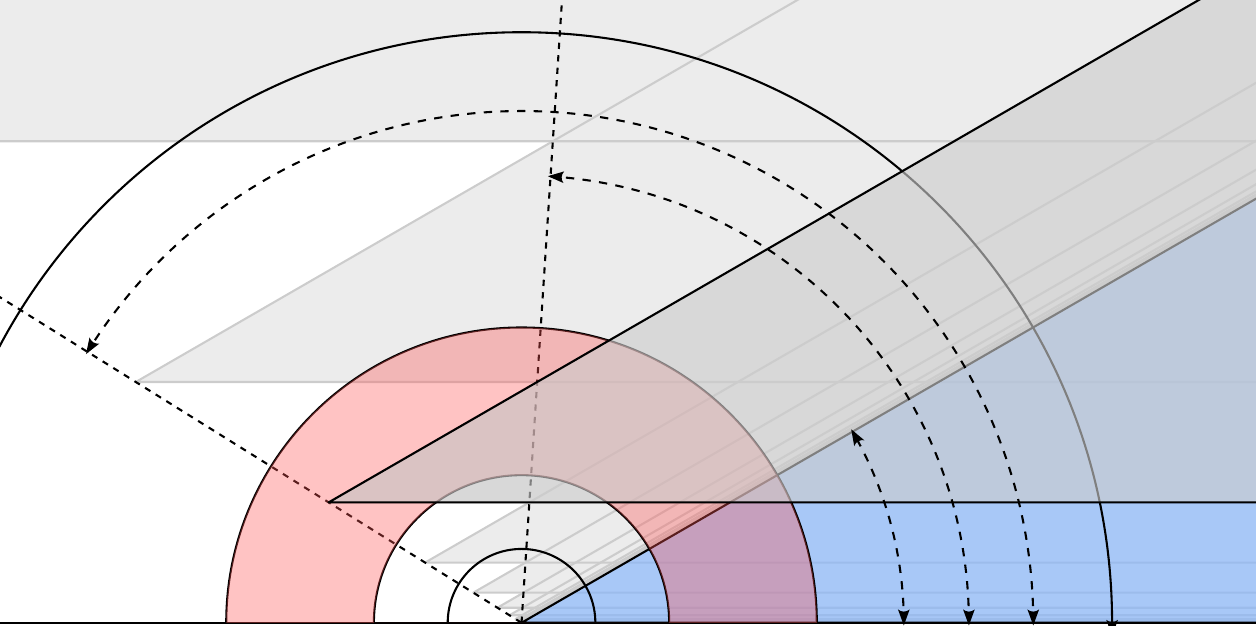}};
      \node at (-3,-2) {$z_0$};
      \node at (-0.9,-3.3) {$z_1$};
      \node at (2.6,-3.31) {$S$};
      \node at (3.2,-3.3) {$S'$};
      \node at (3.8,-3.34) {$S_0$};
    \end{tikzpicture}
  \end{center}
  \caption{Case $\phi_*$ is a dilation of factor $2$ in the proof of \Cref{prop:coneImmersion}; $\Delta$ is in dark gray; $S=T_{z_1}(\vec\Delta)$ in blue; $\bigcup_{n\in\Z}\phi_*^n(\Delta)$ in light gray, dark gray and blue; a fundamental domain of a Reeb cylinder of angle $\pi$ in red.}
  \label{fig:dc}
\end{figure}

If $\phi_*$ is a dilation, recall that its fixed point $z_1$ does no belong to $\Delta$ and
note that given $z\in\Delta$, the set of $k\in\Z$ such that $\phi_*^k(z)\in\Delta$ is a integer interval, i.e.\ of the form $\Z/\cap [a,b]$ (this follows from $\Delta$ being convex and the orbit being on a line, ordered like the exponent $k$).
From this, $f$ factors as $f=\psi\circ\pi$ where $\pi: \C \setminus \{z_1\} \to Y=(\C \setminus \{z_1\}) / \langle \phi_* \rangle$ is the quotient to a dilation torus and $\psi:\pi(\Delta)\subset Y\to X^*$ is an affine immersion.
The set $\Delta$ is actually contained in a sector $S_0$ of angle $\leq\pi$ based on $z_1$, so $\pi(\Delta)$ is contained in a dilation cylinder of angle $\leq\pi$. 
Moreover, denote $\Delta = T_{z^*}( \vec\Delta )$ where $z^*$ is the tip of the sector $\Delta$ and let $S=T_{z_1}( \vec\Delta)$.
Then $\pi(\Delta)$ contains $\pi(S)$.
Finally if we call $S'$ the biggest sector based on $z_1$ such that $\pi(S')\subset \pi(\Delta)$, then $\pi(S')$ is a dilation cylinder contained in the domain of $\psi$ and $\Delta\cap S'$ is a neighborhood of $\infty$ in $\Delta$.
So Claim~\eqref{item:cone:ReebImmerse} holds.
See \Cref{fig:dc}.
\end{proof}

\begin{comp}\label{comp:coneImmersion}
Moreover, in Cases~\eqref{item:cone:cylEmbed} and~\eqref{item:cone:antiConeEmbed}, if an element $\phi\in\Hol(f)$ satisfies that $|\arg(a_\phi)|<\theta$, then $\phi$ is a power of $\phi_*$.
In Case~\eqref{item:cone:ReebImmerse} we can choose $\phi_*$ so that this is the case too.
In Case~\eqref{item:cone:antiConeEmbed} we have moreover that $k|\arg(a_{\phi_*})|\leq\theta$.
\end{comp}

\begin{proof}
The convex set $\phi(\Delta)\cap \Delta$ contains a sector, thus is unbounded.
So $\forall R>0$, $\exists z\in \Delta$ such that $\phi^{-1}(z)\in \Delta$. By definition of $\Hol(f)$, $f(\phi^{-1}(z))=f(z)$ and this holds for nearby $z$ too.

In Case~\eqref{item:cone:cylEmbed} consider the decomposition $f=\psi\circ\pi$ of \Cref{lem:HolConeTranslation}, where $\pi:\C\to\C/\langle T_v\rangle$ and $\psi$ is defined in $\pi(\Delta)$ and injective in a half-plane $H$.
This injectivity means that $z$ and $\phi^{-1}(z)$ are related by a power of $T_v$, and this holds locally hence  $\phi\in\langle T_v \rangle$.

In Case~\eqref{item:cone:antiConeEmbed} we use the decomposition $f|_{V'} = g \circ \pi \circ \iota|_{V'}$ of \Cref{lem:coneTurnProlong} where $g$ is injective near infinity. This implies similarly that $\iota(z)$ and $\iota(\phi^{-1}(z))$ are related by a power of $\lambda$, and since this holds locally $\iota\circ\phi\circ\iota^{-1}$ is the restriction of a power of $\lambda$.
This power $k\in\Z$ is independent of $z$ and by letting $z$ tends to infinity we get that $k|\arg(a_{\phi_*})|\leq\theta$.
Note that the only dilation in $\Hol(f)$ is the identity.

Last, in Case~\eqref{item:cone:ReebImmerse}, the set of dilations that fix $z_1$ and belong to $\Hol(f)$ is a discrete subgroup of $\AutC$ (it is a subgroup of $\cal D\cap\Hol(f)$ using the notations of \Cref{lem:irg}).
We let $\phi_*$ be a generator.
The hypothesis that $\phi(\Delta)\cap \Delta$ is unbounded means in particular that $|\arg(a_\phi)|\leq \theta$.
If $\arg(a_\phi)\neq 0$ then we know from the proof of \Cref{prop:coneImmersion} that we are in Case~\eqref{item:cone:ReebImmerse} of the proposition, for which we already determined that no dilation can be in $\Hol(f)\setminus\{\on{Id}\}$.
If $\arg(a_\phi)=0$, recall that the dilations in $\Hol(f)$ form a discrete subgroup of $\AutC$.
If the fixed point of $\phi$ were different from that of $\phi_*$, by conjugating $\phi$ with powers of $\phi_*$ we would get a sequence of dilations in $\Hol(f)$ not equal to $\on{Id}$ but tending to $\on{Id}$, contradicting discreteness.
Hence the fixed points are the same, and $\phi$ is a power of the generator $\phi_*$.
\end{proof}

We deduce a description of the behavior at infinity of radial geodesics of open convex sector immersions.
The following result extends Lemma~3.3 of \cite{Ta23}, originally established for dilation structures, to the broader setting of complex affine structures.

\begin{cor}\label{cor:ConePetal}
Assume that $X^*$ is not a full translation cylinder, a translation torus, or an affine torus.
Let $\Delta$ be an infinite open convex sector $\Delta$ and denote $\cal R$ the set of radial\footnote{I.e.\ starting from the tip.} geodesics $\gamma:(0,+\infty)\to\Delta$.
Given an affine immersion $f :\Delta \to X^*$, one of the following statements holds:
\begin{enumerate}[(i)]
    \item there is an immersed dilation cylinder $C$ such that $\forall \gamma\in\cal R$, the geodesic $f\circ \gamma$ accumulates, as time goes to infinity, on a (possibly self-intersecting) attracting closed geodesic of $C$;
    \item there exists a cylindrical Fuchsian singularity $y \in \mathcal{S}$ such that $\forall \gamma\in\cal R$, the geodesic $f\circ\gamma$ converges to $y$ as times goes to infinity;
    \item there exists an anti-conical singularity $y \in \mathcal{S}$ such that $\forall \gamma\in\cal R$, the geodesic $f\circ\gamma$ converges to $y$ as times goes to infinity;
    \item there exists an irregular singularity swath $s$ such that $\forall \gamma\in\cal R$, the geodesic $f\circ\gamma$ converges to $s$ as times goes to infinity.
\end{enumerate}
\end{cor}

\subsection{Immersions that do not extend continuously to a vertex}\label{sub:TriangleImmersion}

The main difficulty in attempting to generalize the Delaunay construction to affine structures lies in the fact that they are not necessarily geodesically complete. In other words, a geodesic may blow up in finite time without hitting a singularity. This phenomenon is particularly apparent in a dilation cylinder, whose radial geodesic foliation has this property, see \Cref{sub:cylindersAndSkewCones}. In this section, we prove a partial converse.

\begin{prop}\label{prop:TriangleImmersion}
Consider an affine immersion $f$ of an open Euclidean triangle $\Delta$ to $X^*$ that does not extend continuously to one of its vertices $z_{0}$.
Then $f$ extends to the immersion of the smallest infinite open convex sector $S$ of vertex $z_{0}$ containing $\Delta$.
Moreover, one of the following statements holds:
\begin{enumerate}[(i)]
    \item $X$ is an affine torus;
    \item there exists a dilation $\phi_*$ fixing $z_0$ such that $f$ quotients to an immersion of the dilation cylinder $S/\langle\phi_*\rangle$, i.e.\ $f=\psi\circ\pi$ where $\pi: S\to S/\langle\phi_*\rangle$ is the quotient map; in particular, every radial geodesic of $S$ is mapped by $f$ to a (possibly self-intersecting) attracting closed geodesic of the cylinder.
\end{enumerate}
\end{prop}

\medskip

\noindent \textit{Proof of \Cref{prop:TriangleImmersion}:}

\loctitle{Construction of a special sequence $\phi_n\in \Hol(f)$}
We start by proving that there exists a sequence $(\phi_{n})_{n \in \mathbb{N}}$ of elements of $\Hol (f)$ satisfying the following conditions:
\begin{enumerate}
    \item\label{item:ti:moda} $\lim\limits_{n\to+\infty} a_{\phi_{n}} = 0$;
    \item\label{item:ti:arga} $\lim\limits_{n\to+\infty} \arg(a_{\phi_{n}}) = 0$;
    \item\label{item:ti:xn} there exists sequences $x_n,x'_n\in\Delta$ such that $x_n\tends 0$, $x'_n\tends 0$ and $\phi_{n}(x'_{n})=x_n$.
\end{enumerate}

Consider a radial geodesic segment $\gamma:(0,1] \to \C$ from $z_0$ to a point of the triangle $\Delta$, bisecting the angle of the triangle at $z_0$. Denote $\alpha=\arg \gamma'$ the direction of this bisector. 
Following \Cref{lem:Avoid2}, there exists $y \in f(\Delta)\subset X^{\ast}$ such that $f\circ \gamma (t)$ accumulates on $y$ as $t$ tends to $0$.

The argument is then an elaboration on an argument in the proof of \Cref{lem:Avoid2}.
Consider a neighborhood $V$ of $y$ affine isomorphic to the disk $D$ of center $0$ and radius $1$ via $\psi:V\to D$.
Below we identify $V$ with $D$.
The geodesic $f\circ \gamma$ cannot tend to $y$ nor to a point of $V$, for otherwise it would be possible to extend $f$ continuously at $z_0$.
In particular the geodesic will enter and exit $V$ infinitely many times.
Each time the geodesic $f\circ\gamma$ meets $V$, it does so along a chord of $D$, except possibly for one part which begins in $D$ corresponding to $t=1$.
We get a sequence of chords $(A_k,B_k) = \psi\circ f(I_k)$, $I_k=\gamma((a_k,b_k))$, $(a_k,b_k)\subset(0,1]$ disjoint intervals.
The chords can be oriented by the parametrization by $t$.
We know they accumulate on $0$ (i.e.\ $y$).
We can select a subsequence such that the oriented direction of the chord converges too.
For each such chord $(A_k,B_k)$, there is $L_k\in\AutC$ such that $\psi\circ f=L_k$ holds on $I_k$.
The chord length tends to the diameter of $D$ while the length of the interval $(a_k,b_k)$ tends to $0$. As a consequence, the constant $|L_k'|$ tends to $+\infty$ as $k\to+\infty$.
The image $L_k(z_0)$ of the triangle vertex cannot belong to $D$ since $0\notin (a_k,b_k)$.
This implies that for all $k$ big enough, $L(\Delta)$ contains the center $0$ of $D$.
Let
\[ x_k = L_k^{-1}(0) .\]
We just saw that $x_k\in\Delta$ for all $k$ big enough.
Note that
\[x_k\tends z_0:\]
indeed $d(x_k,[\gamma(a_k),\gamma(b_k)]) = d(0,[A_k,B_k]) / |L_k'| \leq 1 / |L_k'| \tends 0$ and $\gamma(a_k)$ and $\gamma(b_k)$ tend to $z_0$.
Since $x_k\neq z_0$ ($z_0$ is in $\Delta$ but not $x_k$) we can extract an injective subsequence, and assume $x_k\in\Delta$ for all $k$. 
Moreover, $f(x_k)=y$, so for any two terms $x_{k}$, $x_n$ of the sequence, there exists $\phi = \phi_{n,k}=L_n^{-1}\circ L_k\in\Hol(f)$ such that $\phi(x_k)=x_n$.
By selecting $n$ and $k$ big, we can ensure that $\arg (a_\phi)$ is close to $0$.
By selecting $k$ first and then $n$ big enough, we can ensure that $a_\phi$ is close to $0$.
We thus get a sequence $\phi_n$ satisfying the three required conditions by letting $x'_n=x_k$.

\loctitle{Extension of $f$ to $S$}
Denote $S_0 = T_{-z_0}(S)$: it is an infinite sector whose tip is $0$.
For $n$ big enough, the triangle $\Delta$ contains the bounded sector $T_{x_n}(C')$ where $C'=S_0\cap B(0,r)$, for some uniform $r>0$.
The triangle $\phi_{n}^{-1}(\Delta)$ thus contains the bounded sector $\setof{x'_n+a_{\phi_n}^{-1} z}{z\in C'}$ whose tip is close to $z_0$ and whose radius is $r/|a_{\phi_n}|$, hence large.
It follows that for every compact subset $K$ of $S$, there exists $N>0$ such that $\forall n\geq N$, $K\subset \phi_{n}^{-1}(\Delta)$.
Let $\hat f_n$ be the restriction to the convex open set $S\cap \phi_n^{-1}(\Delta)$ of $f\circ\phi_n$.
Note that $\Delta\cap \phi_n^{-1}(\Delta)$ is not empty since it contains $x'_n$.
If follows that for any $m$ and $n$, $f_n$ and $f_m$ coincide near $x'_n$ with $f$, so coincide near $x'_n$, and by affine continuation coincide on $\dom(f_n)\cap\dom(f_n)$, the latter set being convex, hence connected.
It follows that, given $z\in S$, the value of $\hat f_n(z)$ for the various $n$ such that $z\in\dom(\hat f_n)$ is independent of $n$.
We deduce that $f$ has an extension $\hat f: S\to X^*$ that is an affine immersion (and actually takes values in $f(\Delta)$).
The domain of $\hat f$ is the convex set $S$ and near $x'_n$, $\hat f\circ\phi_n = f \circ \phi_n = f =\hat f$, i.e.\ %
\[\phi_{n} \in\Hol(\hat f).\]

For $n$ big enough, $\phi_n(S)\cap S$ contains $\setof{x_n+z}{z\in S_0\text{ and }a_{\phi_n} z\in S_0}$, hence is unbounded.
Up to reindexing the sequence, we assume this is the case for all $n$.
We now apply \Cref{prop:coneImmersion} to $\hat f$.
Because $\phi_n\in \Hol(\hat f)$ and $\phi_n(S)\cap S$ is unbounded, it follows that $\hat f$ is not injective on any neighborhood of $\infty$, so we cannot be in Case~\eqref{item:cone:embed} of \Cref{prop:coneImmersion}.
If $X^*$ were a full translation cylinder or a translation torus, $\Hol(\hat f)$ could only contain translations, so these cases are excluded too.
Case~\eqref{item:cone:cylEmbed}, of $\hat f$ factoring near $\infty$ to an semi-cylinder embedding, is incompatible too with $\phi_n\in\Hol(\hat f)$.
Hence either 
\begin{enumerate}
\item\label{item:pc:at} $X=X^*$ is an affine torus,
\item\label{item:pc:acs} $\hat f$ factors near $\infty$ to an embedding through a anti-conical singularity of opening angle $\leq \theta$,
\item\label{item:pc:dc} $\hat f$ factors near $\infty$ to an immersion through a dilation cylinder.
\end{enumerate}

More precisely denote $\theta$ the opening angle of $\cal C$.
In cases~\eqref{item:pc:acs} and~\eqref{item:pc:dc}, \Cref{prop:coneImmersion} tells us that $\hat f=g\circ\pi\circ\iota\circ L$ holds on a neighborhood $V$ of $\infty$ in $\cal C$, where:
\begin{itemize}
\item $L\in\AutC$ is such that $0\notin L(S_0)$,
\item $\iota$ is an embedding of $L(S_0)$ in $\cet$,\item $\pi: \cet\to \cet/\Lambda$, 
\item $\Lambda$ is the group generated by the lift to $\cet$, and with rotation angle in $(-\pi,\pi)$, of $z\mapsto a_* z$
\item $a_*\in\C^*$, $a_*\neq 1$ and $\arg(a_*)\in(-\theta,\theta)$,
\item $g$ is an embedding from $\pi\circ\iota\circ L(V)$ to $X^*$.
\end{itemize}
We denote $\phi_*= L^{-1}\circ (z\mapsto a_\phi z) \circ L\in \AutC$, which fixes $z_1:=L^{-1}(0)$ and has dilation factor $a_*$.
\Cref{prop:coneImmersion} tells us that the only elements $\phi\in \Hol(\hat f)$ such that $|\arg(a_\phi)|<\theta$, must be powers of $\phi_*$.
In particular, they all fix $z_1$.
Applying this to $\phi_n$, we get $\phi_n(z_1)=z_1$.
Then $z_1-x_n=\phi_n(z_1)-\phi_n(x'_n)= a_{\phi_n} \times (z_1-x'_n) \tends 0$ so $x_n\tends z_1$. So $z_1=z_0$, i.e.\ every $\phi_n$ fixes $z_0$ and $\phi_*$ fixes $z_0$.

So in case~\eqref{item:pc:acs}, the relation $\hat f\circ\phi_*=\hat f$ holds on the sector based on the tip $z_0$ of $S_0$ and the form $S_0\cap \phi_*^{-1}(S_0)$ and $\hat f$ factors as $\hat f=g\circ\pi$ where $\pi:S_0\to S_0/(z\sim\phi_*(z))$ and $\phi_*$ fixes $z_0$ and $\arg(a_{\phi_*})\in(-\theta,0)\cup(0,\theta)$ where $\theta$ is the opening angle of $S_0$.
The map $g$ is supposed to be injective near $\infty$. However, for all $n$, the relation $f(\phi_n(z))=f(z)$ holds on an unbounded sector contained in $S_0$. For $n$, then $z$, very big, $|\arg(\phi_n(z)/z)|<\arg(a_{\phi_*})$, so $\phi_n(z)$ and $z$ project to distinct points of the skew cylinders, that are close to $\infty$, but mapped to the same point by $g$, leading to a contradiction.

Hence we are in case~\eqref{item:pc:dc} and since $z_1=z_0$, $\phi_*(S_0)=S_0$ and this proves \Cref{prop:TriangleImmersion}.

\medskip

From our results on triangle immersions and on plane immersions, we deduce the analogous result for immersions of convex sets at flat points of their boundaries, which is crucial in the case of disks for the construction of the Delaunay decomposition in Section~\ref{sec:Delaunay}.

\begin{cor}\label{cor:HPDiskImmersion}
Consider an affine immersion $f$ of an open convex subset $\Delta$ of $\C$ to $X^*$.
We assume that $f$ does not extend continuously at $z_{0}\in \partial \Delta$ and $\Delta$ has a tangent at $z_0$.
Then $f$ extends to an immersion $\hat f$ of the infinite half-plane $H$ containing $\Delta$ and delimited by the tangent to $\Delta$ at $z_0$.
\par
Moreover, one of the following statements holds:
\begin{enumerate}[(i)]
    \item\label{item:hpdi:at} $X$ is an affine torus;
    \item\label{item:hpdi:r} $\hat f$, hence $f$, factors to an embedding 
    of a dilation cylinder of angle $\pi$: $\hat f=g\circ \pi$ where $\pi:H\to H/\langle \phi\rangle$, $g$ is an affine embedding and $\phi$ is a dilation fixing $z_0$.
\end{enumerate} 
\end{cor}
In Case~\eqref{item:hpdi:r} every radial geodesic is mapped to an attracting closed geodesic of that cylinder.

\begin{proof}
\Cref{prop:TriangleImmersion} proves that for any triangle $T$ of vertex $z_{0}$ contained in $\Delta$, $f|_T$ extends to an infinite sector $\cal C$ of vertex $z_{0}$ containing $T$, and this extension actually coincides with $f$ on $\cal C\cap \Delta$ by the affine identity principle.
Using a sequence of triangles whose angle a $z_0$ is increasing, bisected by the normal to $\partial \Delta$ at $z_0$, and tends to $\pi$, the extensions extend each other by the affine identity principle, and their union defines an extension $\hat f$ of $f$ to the half-plane $H$ containing $\Delta$ whose boundary line is tangent to $\partial \Delta$ at $z_{0}$.
The extension does not extend continuously at $z_0$, for $f$ would, and the end of the statement of \Cref{prop:halfPlaneImmersion} proves that either $X$ is an affine torus, or $\hat f$, hence $f$, factors through the embedding into $X^{\ast}$ of a dilation cylinder of angle $\pi$.
The fixed point of $\phi$ must be $z_0$, for otherwise there would be continuous extension at $z_0$.
\end{proof}

\subsection{Boundary of immersed disks and cylinders}\label{sub:diskAndCylBoundary}

We draw the following summary of possibilities in the specific case of a disk, in view of the previous results.

\begin{lem}\label{lem:diskImmersion}
Given an affine immersion of a disk $\Delta$ to $X^*$, either
\begin{enumerate}[(i)]
\item\label{item:di:discont} there is a point $z_0\in\partial \Delta$ where $f$ does not extend continuously; then $f$ has a continuous extension on $\partial \Delta\setminus\{z_0\}$, which maps in $X^*$; we are in the situation of \Cref{cor:HPDiskImmersion}, i.e.\ either $X$ is an affine torus or the image of $f$ is an embedded Reeb cylinder of angle $\pi$;
\item\label{item:di:cont} $f$ has a continuous extension to $\partial \Delta$, taking values in $X$; only finitely many points map to $\cal S$.
\end{enumerate}
\end{lem}
\begin{proof}
If there is a point of the boundary where $f$ does not extend continuously then \Cref{cor:HPDiskImmersion} shows we are in Case~\eqref{item:di:discont}.
Otherwise, $f$ has a continuous extension $\bar f$ to $\partial \Delta$.
Since $\bar f^{-1}(\cal S)$ is closed and discrete by \Cref{lem:b} in the compact set $\partial \Delta$, it follows that it is finite.
\end{proof}

This also applies when $\Delta$ is any bounded convex subset of $\C$ with a tangent at every point of its boundary, i.e.\ without angular points.

\medskip

It has been proved that every (possibly self-intersecting) closed geodesic embeds in a continuous family of closed geodesics forming an immersed cylinder. Now, we can describe the boundary of a cylinder.

\begin{cor}\label{cor:CylinderBoundary}
Given an affine immersion $f$ of a translation or dilation cylinder $C$ to $xX^*$,
$f$ extends continuously to the boundary of the cylinder. The image of a boundary component is either a (possibly self-intersecting) closed geodesic or the union of finitely many (possibly intersecting and self-intersecting) saddle connections joining apexes of $\mathcal{S}$.
\par
Moreover, if $f$ is an embedding, then the closed geodesics or individual saddle connections in the image of the boundary component have no self-intersection (but two saddle connections may have non-empty intersection, as in the example below).
\end{cor}

\begin{ex}
As an example, consider the cushion, a surface formed from two copies of a square (one with inverted orientation) glued along their facing boundaries.
It haw $4$ singularities, which are Fuchsian conical.
An embedded cylinder is obtained by removing two opposite edges of the original square, which are particular saddle connections. Each boundary circle is mapped in a $2:1$ way to one of these removed edges.
\end{ex}

\begin{proof}[Proof of \Cref{cor:CylinderBoundary}]
Let $z_0\in\partial C$.
Assume by contradiction that $f$ does not have a continuous extension at $z_0$.
We defined $\partial C$ as a closed subset of some affine surface $Y$ containing $C$.
Consider an embedded (injective) disk $D$ centered on $z_0$ in $Y$ 
The set $D\cap C$ can be identified with a half-disk, which we call $\Delta$, in $\C$.
In our hypothetical situation, the restriction of $f$ to $\Delta$ is not continuous.
\Cref{cor:HPDiskImmersion} then implies that $f|_\Delta$ extends to an immersion $\hat f$ of the half-plane $H$ delimited by the line through $z_0$ supporting the straight part of $\partial\Delta$ and satisfying $\hat f \circ\phi=\hat f$ where $\phi\in\AutC$ is a dilation of ratio $s\neq 1$ and fixing $z_0$.
Following \Cref{cor:HPCommutation}, all the elements of $\Hol (\hat{f}|_{H})$ commute and therefore fix $z_{0}$.
But $H$ contains the universal covering $\tilde C$ of $C$, on which $f$ lifts to a map $\tilde f$ which has an element $\phi\in\AutC$ in its holonomy set that does not fix $z_0$.
By affine continuation, $\tilde f$ must coincide with $\hat f$ on $\tilde C$, and $\Hol(\hat f)$ contains $\phi$ which thus should fix $z_0$, leading to a contradiction.

\Cref{lem:ContinExt} implies many things:
\begin{itemize}
\item That the points of the component $B$ of $\partial C$ containing $z_0$ and mapped to a singularity form a discrete subset $A$.
It is also a closed subset, being the preimage of the closed subset $\cal S$ of $X$ by the continuous extension.
Hence $A$ is finite.
\item That the singular points must be conical Fuchsian singularities or irregular singularities.
\item That near every other point, there is an affine extension. In particular the boundary between two points of $A$ map by $f$ to the support of a geodesic of finite time span starting and ending at singularities, i.e.\ a saddle connection.
If $A$ is empty then $B$ is mapped to the support of a closed geodesic.
\end{itemize}

Finally, assume that if $f$ is not injective on a component of $B\setminus A$, meaning that two boundary points $z_{0},z_{1}$ mapped to $X^{\ast}$ have same image under $f$.
Since this component correspond to a closed geodesic or a geodesic arc, $z_1$ and $z_0$ are mapped to a self-intersection point of the geodesic, and at such a point the directions of the geodesic cannot be opposite.
It follows that we can find points in the interior of $C$ that have same image under $f$.
Thus, if $f$ is an embedding, its boundary geodesics and individual saddle connections have no self-intersection.
\end{proof}

\section{Delaunay decomposition}\label{sec:Delaunay}

In this section, we organize the affine immersions of planes, half-planes, and disks into a finite type affine surface into the \textit{Delaunay category}, introduced in \Cref{sub:DelaunayCategory}, whose structure will allow us to describe the qualitative geometry of the affine surface. We define maximal elements of the Delaunay category and prove that each element of the Delaunay category is the restriction of such an immersion. 
\par
In \Cref{sub:maximal}, we characterize the maximal elements of the Delaunay category of a given affine surface as affine immersions which exhibit a certain behavior on the boundary.
\par
In \Cref{sub:AffineSpine}, we define a deformation process of maximal elements, which we call pivoting.
This allows to give a graph structure to the set of maximal elements, which we call the Delaunay spine, and that can be naturally embedded in the affine surface, blown up at some points. %The structure of this graph and the types of the maximal elements forming its vertices and edges will play a key role in the description of the Delaunay decomposition.
\par
In \Cref{sub:delSeg} we define Delaunay segments using one of the maximal immersion types: type~A (disk, flexible).
We prove that they form saddle connections that are embedded (not self-crossing) and that they do not cross each other.
We also prove, using the counting results in \Cref{sub:apicalSystems}, that there are only finitely many Delaunay segments.
\par
In particular, in \Cref{sub:Exceptional}, we classify the affine surfaces without Delaunay segments, referred to as the exceptional affine surfaces (proving thus \Cref{thm:ExceptionalSurfaces}).
\par
In \Cref{sub:core}, we define open Delaunay polygons, as image by rigid maximal disk immersions of the interior convex hull of points mapped to singularities and prove that they are embedded and contain no point on a Delaunay segment. We also define the exterior of $X$ as the union of images of all half-plane immersions (plus the non-conical Fuchsian singularities), and the core as its complement.
We prove that the core is the union of the open Delaunay polygon interiors, the Delaunay segment interiors, the conical Fuchsian singularities and the irregular singularities.
\par
In \Cref{sub:DelaunayComponents}, we describe the connected components in the Delaunay decomposition of non-exceptional surface (proving \Cref{thm:ClassificationComponents,thm:ComplexityBound}).

\subsection{Delaunay category}\label{sub:DelaunayCategory}

Affine immersions of standard open domains of the complex plane into a finite type affine surface assemble into an ordered structure that we present below.
\par
We recall that affine \emph{immersions} are just affine maps, i.e.\ maps that are locally affine in charts. We emphasize that they do not need to be injective (they only have to be locally injective).
We call affine \textit{embedding} the injective immersions.
For brevity, we will often abbreviate affine immersion as immersion, and affine embedding as embedding. 

\begin{defn}\label{defn:DelaunayCategory}
For a finite type affine surface $(X,\nabla,\mathcal{S})$, the \textit{Delaunay category} $\mathcal{D}el$ is defined as follows:
\begin{itemize}
    \item objects of $\mathcal{D}el$ are pairs $(\Delta,f)$ where $\Delta$ is an open disk in the canonical affine plane $\C$, an open half-plane in $\C$ or the whole plane, and $f$ is an affine immersion from $\Delta$ to $X^{\ast}$;
    \item an arrow between two pairs $(\Delta_{0},f_{0}) $ and $(\Delta_{1},f_{1})$ is an affine embedding $\phi:\Delta_{0} \rightarrow \Delta_{1}$ such that $f_{0}=f_{1} \circ \phi$.
\end{itemize}
\end{defn}

The existence of an arrow between two objects of the Delaunay category is a preorder relation: it is reflexive and transitive.
It is not antisymmetric for two reasons.
First, for any automorphism $\phi$ of $\C$, the element $(\phi^{-1}(\Delta),f\circ\phi)$ is isomorphic to $(\Delta,f)$ in this category.
Second, in the case $X^*$ is a full cylinder $\C/\Z$, for any two half-planes $H$ and $H'$ with parallel non-horizontal boundary lines, there is an arrow from $\pi|_H$ to $\pi|_{H'}$ and an arrow from $\pi|_{H'}$ to $\pi|_{H}$: assuming $H$ contains a neighborhood of $-\infty$ in $\R$, and that $H'=T_x(H)$ with $x\in\R$, then for any $n\in\Z$ such that $x+n\leq 0$, we have $\pi_H' = \pi_H \circ T_n|_{H'}$, while for any $n\in\Z$ such that $n\leq x$ then $\pi_H = \pi_H' \circ T_n|_{H}$.

The notion of \emph{maximal elements} (up to isomorphism) for this relation support important geometric information on the complex affine structure.

\begin{defn}\label{defn:DelaunayMaximal}
For a finite type affine surface $(X,\nabla,\mathcal{S})$, a \textit{maximal element} of the Delaunay category is a pair $(\Delta,f)$ such that
for any $(\Delta_1,f_1)$ for which there is an arrow $\phi:(\Delta,f)\rightarrow (\Delta_1,f_1)$ then $\phi^{-1}$ is an arrow $(\Delta_1,f_1)\to (\Delta,f)$ (i.e.\ $\phi(\Delta)=\Delta_1$).
\end{defn}

\begin{rem}
In the more general setting of complex projective structures (also called M\"{o}bius geometry), we do not distinguish between disks and half-planes. The space of immersed disks in a complex projective surface is locally modeled on the space of oriented circles in $\mathbb{CP}^{1}$, which identifies with a three-dimensional de Sitter space. The relation defined by the existence of an arrow between two immersions then corresponds to the causal structure induced by the Lorentzian metric, see \cite{Sc99} for details.
\end{rem}

\begin{lem}\label{lem:charMax1}
An element $(\Delta,f)$ of the Delaunay category is maximal if and only if there is no affine extension of $f$ to a disk, half-plane or whole plane, containing $\Delta$.
\end{lem}
\begin{proof}
If $(\Delta,f)$ has an extension $(\Delta_1,f_1)$ with $\Delta\subsetneq \Delta_1$, then the injection $\phi(z)=z$ from $\Delta$ to $\Delta_1$ is an arrow from the first to the second but $\phi(\Delta)\neq\Delta_1$, so $(\Delta,f)$ is not maximal.

Conversely if $(\Delta,f)$ has no such extension and if there is some arrow $\phi:(\Delta,f)\to(\Delta_1,f_1)$ then $\phi$ is the restriction of an element $\hat\phi\in\AutC$ and $f_1\circ \hat\phi$ provides an extension of $f$, so $\hat\phi^{-1}(\Delta_1) = \Delta$ and $\phi^{-1}$ is an arrow from $(\Delta_1,f_1)$ to $(\Delta,f)$.
\end{proof}

We will characterize maximal elements in \Cref{sub:maximal}.
We first state and prove below the existence of maximal elements above a given one. 

\begin{lem}\label{lem:charMax2}
Consider a finite-type affine surface $(X,\nabla,\mathcal{S})$.
Then every element $(\Delta,f)$ of the Delaunay category is the restriction of a maximal element.
\end{lem}
\begin{proof}
In view of \Cref{lem:charMax1}, it is enough to prove that the set of open disks/half-planes/whole planes in $\C$ containing $\Delta$ and on which $f$ has an affine extension mapping in $X^*$, has a maximal element for the inclusion.

If $\Delta$ is a half-plane, consider the set of all open half-planes containing $\Delta$ and on which $f$ has an affine extension to $X^*$. Consider their union: it is either the whole plane $\C$ or a half-plane $\hat H$.
In any case, for any two open half-planes $H_1$, $H_2$ containing the half-plane $\Delta$, one must contain the other. If $f$ has affine extensions to $H_1$, $H_2$, then they must coincide on the smallest of $H_1$ or $H_2$, by the affine identity theorem since they coincide on $\Delta$.
It follows that we can define $\hat f(z)$ for $z\in\C$ or $\hat H$ as the union of all these extensions, i.e.\ as $f_1(z)$ for any extension $f_1$ defined on a half-plane $H_1$ containing $z$.

If $\Delta$ is a disk, consider the supremum $R\in(0,+\infty]$ of the radii of disks $\Delta'$ containing $\Delta$ and on which $f$ can extend.
There exists a sequence $\Delta_n$ whose radii sequence $R_n$ tends to $R$.
If $R$ is finite then by extraction we may assume that their centers converge too, to some $c\in\C$.
As in the case for planes, we can then use the affine identity theorem to define an extension $\hat f$ on $B(c,R)$ using an exhaustion by connected compact subsets containing a closed disk contained in $\Delta$.
If $R=+\infty$ then either the distance from the center $z_1$ of $\Delta$ to the circle $\partial\Delta'$ tends to infinity, in which case we can extract a nested subsequence and $f$ extends to $\C$. 
Or this distance stays bounded and we can extract a subsequence such that the $\Delta_n$ tend to a half-plane $H$ and we can proceed as in the case $R$ is finite to get an extension to $H$ by exhaustion.
It may happen that $H$ is not maximal.
But we know that the extension to $H$ has a maximal extension because we treated the case $\Delta$ is a half-plane.
\end{proof}

\subsection{Maximal elements in the Delaunay category}\label{sub:maximal}

In the following two statements, we prove that maximal immersion have specific boundary behavior.

But first let us recall what we proved concerning the boundary behavior of general immersions $f$ of disk and half-planes in finite type affine surfaces.
By \Cref{lem:diskImmersion} and the end of \Cref{prop:halfPlaneImmersion}, disk and half-plane affine immersions continuously extend to the whole boundary except maybe one point (in which case $f$ factors through a Reeb cylinder embedding).
Denote $\bar f$ the continuous extension.
Denote
\[A=\bar f^{-1}(\cal S)\]
the set of points (on the boundary) that are mapped to a singularity by the extension.
For half-planes, using \Cref{lem:b} and the classification in \Cref{prop:halfPlaneImmersion}, we get:
if there is a discontinuity point $z_0$, then $A$ is a finite union of orbits under a common dilation fixing $z_0$;
if there is no discontinuity point, either $A$ is finite or $f$ factors through a cylinder embedding and $A$ a finite union of orbits under a common translation of the boundary line.
For disks, by \Cref{lem:diskImmersion}, in the presence of a discontinuity then $A$ is empty, while if there is no discontinuity then $A$ is finite.

The results of the present section will in particular prove that, given a finite type affine surface $(X,\nabla,\mathcal{S})$, an element $(\Delta,f)$ of the Delaunay category is maximal if and only if it is one of the following:
\begin{itemize}
    \item[--] the immersion of a disk with at least two points in $A$;
    \item[--] the immersion of a half-plane with at least one point in $A$;
    \item[--] the immersion of a half-plane that does not extend at some point of the boundary;
    \item[--] the immersion of the whole plane.
\end{itemize}

\begin{lem}\label{lem:DiskMaxImmersion}
Consider a finite type affine surface $(X,\nabla,\mathcal{S})$ and an affine immersion $f:\Delta \to X^{\ast}$ of a disk.
If $(\Delta,f)$ is maximal then $f$ has no discontinuity point on $\partial \Delta$ and $\bar f$ maps at least two (boundary) points to singularities.
Conversely, an affine immersion of a disk such that $\bar f$ sends at least two points to singularities is maximal (hence has no discontinutiy point).
\end{lem}

\begin{figure}[htbp]
\begin{center}
\begin{tikzpicture}
\node at (0,0) {\includegraphics[scale=0.8]{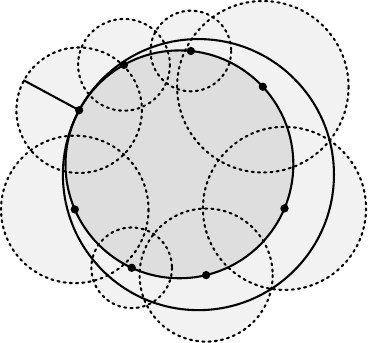}};
\end{tikzpicture}
\end{center}
\caption{A disk immersion that has an affine extension near every point of its boundary except one where it has a slit extension, can be extended to a bigger tangent disk.}
\label{fig:eds}
\end{figure}

\begin{proof}
Assume $(\Delta,f)$ is maximal.
In the presence of a discontinuity, we saw $f$ would factor to a Reeb cylinder embedding, which implies that $f$ extends to a half-plane containing the disk, so $f$ would be non-maximal.
So there is no discontinuity: $f$ extends continuously to $\partial \Delta$.

If the whole boundary maps to regular points, there is a local affine extension near every boundary point, and using a compactness argument and the affine identity theorem, we can extend to a disk of the same center but bigger radius.

If there is only one point $z_0\in\partial \Delta$ mapped to a singularity, then a similar procedure can be applied to extend to a bigger disk tangent to $\Delta$ at $z_0$, using at $z_0$ that we have an extension to a slit neighborhood (\Cref{lem:SlitDisk}). See \Cref{fig:eds}

So there are at least two points on $\partial \Delta$ that are mapped to singularities by the continuous extension of $f$.

Conversely, if two boundary points map to singularities, any extension of $f$ must be defined on a set avoiding these two points. But there is no set among open disks, open half-planes or the whole plane, that strictly contains $\Delta$ but does not contain any of the two points.
\end{proof}

\begin{lem}\label{lem:HPMaxImmersion}
Consider a finite type affine surface $(X,\nabla,\mathcal{S})$ and affine immersion $f:H \to X^{\ast}$ of an open half-plane.
Then $(H,f)$ is maximal if and only if either it has a discontinuity point on $\partial H$, or its extension $\bar f$ to $\partial H$ maps a point\footnote{$\partial H$ denotes the boundary of $H$ in $\C$, i.e.\ the point at infinity is excluded.} to a singularity (these two cases are not exclusive).
\end{lem}
\begin{proof}
Discontinuities and singularities obviously prevent extension.

For the converse, we use \Cref{prop:halfPlaneImmersion}.
By hypothesis, $f$ may satisfy any case of that statement.

In Cases~\eqref{item:embed}  ($f$ is an embedding) or~\eqref{item:antiConeEmbed} ($f$ factors to an embedding of an anti-conical Fuchsian singularity standard neighborhood, of opening angle $<\pi$), \Cref{prop:halfPlaneImmersion} tells us in particular that, as $z$ tends to $\infty$, $f(z)$ tends to an anti-conical singularity or a swath of an irregular singularity.
Applying \Cref{lem:DeformCone} to some half-lines contained in $H$ that are parallel to $\partial H$, we show that $f$ extends on a sector $S$ of angle strictly larger than $\pi$ that contains all the boundary line except possibly a segment $I$ of finite length.
Points mapped to a singularity are necessarily in $I$. If there were no such point, we may combine the extension of $f$ to sector $S$ and a finite union of disks covering $I$ as in the disk case above, and deduce that $f$ extends to an open half-plane $H'$ containing $H$ as a strict subset.

In Case~\eqref{item:exc}, half-plane immersion are very explicit:
for full cylinders and translation tori, they always extend to a bigger half-plane, for affine tori, this is the case too if the lift of the half-plane in $\cet$ has a boundary that does not contain $0$, otherwise continuous extension fails at this point.

In Case~\eqref{item:cylEmbed}, $f$ factors through the embedding of a semi-infinite translation cylinder.
The latter bounded by a chain of finitely many saddle connections (see \Cref{cor:CylinderBoundary}) or by a periodic geodesic.
In the second case, $f$ would extend to a bigger half-plane (by a similar argument as for disks above), so we are in the first case.

Case~\eqref{item:ReebEmbed} is when $f$ factors through the embedding of a Reeb cylinder of angle $\pi$, more precisely there exists a non-trivial dilation $\phi\in\AutC$ whose fixed point does not belong to $H$, and a factorization $f = g \circ \pi|_H$ where $\pi : H' \to H'/\langle \phi \rangle$ is the quotient map to a Reeb cylinder of angle $\pi$, $H'$ is the half-plane containing $H$ and whose boundary passes through the fixed point of $\phi$ and $g$ is an affine embedding from the Reeb cylinder into $X^*$.
Since $g\circ\pi$ defines an extension of $f$ to $H'$, we actually have $H'=H$. So there is a discontinuity on $\partial H$.
\end{proof}

As a consequence:

\begin{prop}\label{prop:maxHalfPlane}
Consider a finite type affine surface $(X,\nabla,\mathcal{S})$ that is not an affine torus.
Let $(H,f)$ be a maximal element of the Delaunay category with $H$ a half-plane.
Then $X$ does not admit an immersion of the whole plane.
Moreover $(H,f)$ satisfies one of these three mutually exclusive statements:
\begin{enumerate}[(i)]
    \item\label{item:ehp:other} $f$ tends at infinity to an anti-conical singularity or a irregular singularity swath, it extends continuously at every point of the boundary line, the number of points mapped to points of $\mathcal{S}$ is finite and at least one, it maps the segment between two consecutive such points to a saddle connection.\item\label{item:ehp:cyl} $f$ is the universal covering of a semi-infinite translation cylinder, it extends continuously at every point of the boundary line, the boundary line is mapped to a chain of finitely many saddle connections;
    \item\label{item:ehp:reeb} $f$ is the universal covering of an embedded Reeb cylinder of angle $\pi$, it extends continuously at every point of the boundary line excepted exactly one point, the image of the boundary is the image of the two boundary circles of the cylinder, each of which being a closed geodesic or a finite union of saddle connections;
\end{enumerate}
\end{prop}

\begin{proof}
We cannot be in Case~\eqref{item:exc} of \Cref{prop:halfPlaneImmersion}: we excluded affine tori and for full cylinders and translation tori, the only maximal immersions are that of the whole plane.

Then, the only case for which a discontinuity can arise on $\partial H$ is, by the end of \Cref{prop:halfPlaneImmersion}, Case~\eqref{item:ReebEmbed} of that statement, when the fixed point of the generator $\phi$ is on $\partial H$. With this and \Cref{cor:CylinderBoundary}, we get Case~\eqref{item:ehp:reeb} of the present statement.

Otherwise there is continuous extension to $\partial H$ and we are not in Case~\eqref{item:ReebEmbed} of \Cref{prop:halfPlaneImmersion} for otherwise we would be able to extend $f$, as already noted in the proof of \Cref{lem:HPMaxImmersion}.

In Cases~\eqref{item:embed} or~\eqref{item:antiConeEmbed} of \Cref{prop:halfPlaneImmersion}, we saw that $f$ extends near $\infty$ to a sector containing the boundary line and tends to an anti-conical singularity or an irregular singularity swath.
By \Cref{lem:b} there is a finite number of point on $\partial H$ mapped to singularities. By \Cref{lem:HPMaxImmersion} there is at least one.
We get Case~\eqref{item:ehp:other} of the present statement.

In Case~\eqref{item:cylEmbed} of \Cref{prop:halfPlaneImmersion}, we can apply \Cref{cor:CylinderBoundary,lem:HPMaxImmersion} and get Case~\eqref{item:ehp:cyl} of the present statement.
\end{proof}

We distinguish between \textit{flexible} and \textit{rigid} immersions among the maximal immersions of the Delaunay category of a finite type affine surface.

\begin{defn}\label{defn:flexible}
A maximal element of the Delaunay category of a finite type affine surface is said to be \textit{flexible} in the following cases:
\begin{itemize}
    \item[Type A:] (disk, flexible) the immersion of an open disk that extends continuously to its boundary circle where exactly two points are mapped to points of $\mathcal{S}$;
    \item[Type B:] (half-plane, non-cylinder, flexible) the immersion of an open half-plane that extends continuously to its boundary line where exactly one point is mapped to a point of $\mathcal{S}$;
    \item[Type C:] (Reeb, flexible) the immersion of an open half-plane that extends continuously to its boundary line excepted in one point and where no boundary point is mapped to a point of $\mathcal{S}$.
\end{itemize}
\end{defn}

We recall that an immersion of a disk or half-plane always extends continuously to all but maybe one point of the boundary and that the set of points mapping to a singularity is discrete (all points are isolated).

\begin{defn}\label{defn:rigid}
A maximal element of the Delaunay category of a finite type affine surface is said to be \textit{rigid} in the following cases:
\begin{itemize}
    \item[Type I:] (disk, rigid) the immersion of an open disk that extends continuously to its boundary circle where at least three of points are mapped to points of $\mathcal{S}$;
    \item[Type II:] (half-plane, non-cylinder, rigid) the immersion of an open half-plane that extends continuously to its boundary line where only a finite number of points, at least two, are mapped to points of $\mathcal{S}$;
    \item[Type III:] (semi-cylinder, rigid) the immersion of an open half-plane that extends continuously to its boundary line where an infinite and translation-periodic subset of points are mapped to points of $\mathcal{S}$;
    \item[Type IV:] (Reeb, rigid, $>\pi$) the immersion of an open half-plane that extends continuously to its boundary line excepted in one point $z_0$ and where the set of boundary points mapped to points of $\mathcal{S}$ is non-empty, invariant by non-trivial dilation fixing $z_0$, and contained on one side of $z_0$;
    \item[Type V:] (Reeb, rigid, $\pi$) same as Type IV but with points mapped to $\cal S$ on both sides of $z_0$;
    \item[Type VI:] (whole plane) the immersion of a whole plane.
\end{itemize}
\end{defn}

Using the previous results, we check immediately that, for finite type affine surface, every maximal element of the Delaunay category is either rigid or flexible.

\subsection{Affine Delaunay spine and pivoting}\label{sub:AffineSpine}

Given a maximal element $(\Delta,f)$ of the Delaunay category where $\Delta$ is not the whole plane, we define a way to deform it which we call \emph{pivoting}.
Consider the boundary $\partial \Delta$ from which we remove the discontinuity point if there is one, and every point mapped to $\cal S$ by the continuous extension (then by \Cref{lem:b} these singularities are conical Fuchsian or irregular).
One gets an open subset of $\partial \Delta$, which has at least two connected components, possibly infinitely many if $\Delta$ is a half-plane.
Choose one component, denote it $J$ and call it the \emph{egress component}.
It is an open arc of circle, an open segment, or an open half-line.

\smallskip

\noindent Case 1: the component $J$ is bounded. Then $J$ is bounded by two continuity points $A$ and $B$ mapped to singularities. There is thus an extension of $f$ (see \Cref{lem:SlitDisk}) to a set $V$ containing $\Delta$, a neighborhood of the open arc/segment $J$ and two slit neighborhoods of $A$ and $B$, and we can choose the slits tangent to $\partial \Delta$ at these points.
As a consequence, consider pencil of circles/line passing through $A$ and $B$, which form a continuous family we can parametrize by the angle $\theta\in \R/\pi\Z$ they make at $A$ with $\partial \Delta$, with $\theta$ small and positive corresponding to the existence of an arc outside $\Delta$ but close to $J$. Call $C_\theta$ the circle/line of parameter $\theta$ in the pencil.
Then there is some $\eps>0$ such that for all $\theta\in[0,\eps]$, $V$ contains the disk $\Delta'$ bounded by the circle of parameter $\theta$.
By \Cref{lem:DiskMaxImmersion}, $(\Delta',\hat f|_{\Delta'})$ is a maximal element of the Delaunay category.

\smallskip

\noindent Case 2: the component $J$ is unbounded, \ie it is a half-line. If it is bounded by a discontinuity point $z_0$ then $f\circ\phi=f$ for some non-trivial dilation $\phi$ fixing $z_0$. Let $[a,b)\in J$ be a fundamental domain of the action of $\phi$ on $J$. Every point of $J$ has a neighborhood on which $f$ extends, and by compactness of $[a,b]$ and $\phi$-invariance we can extend $f$ to an affine immersion defined on the union of $H$ and of an open infinite sector bisected by $J$.
This immersion still satisfies $\hat f\circ\phi=\hat f$.
If $J$ is bounded by a point $z_0$ mapped to a singularity, then by \Cref{lem:SlitDisk} near $z_0$ and regular points, and \Cref{prop:halfPlaneImmersion} near infinity (we are in Cases~\eqref{item:embed} or~\eqref{item:antiConeEmbed} of that statement), we can also affinely extend $f$ to the union of $H$ and of an open infinite sector bisected by $J$ and the extension extends continuously to the boundary of its domain.
In both cases the domain of $\hat f$ contains every half-plane $\Delta'$ containing $J$ and obtained by rotation $H$ around $z_0$ by a small enough amount. the half-plane bounded by the rotated line and containing $J$.
Then $(\Delta',\hat f|_{\Delta'})$ is still a maximal element of the Delaunay category.

\smallskip

\begin{figure}[htbp]
\begin{tikzpicture}
\node at (0,0) {\includegraphics[scale=0.75]{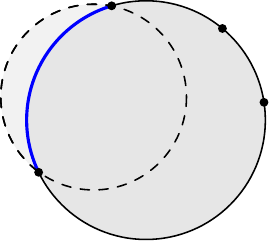}};
\node at (0,-3) {\includegraphics[scale=0.75]{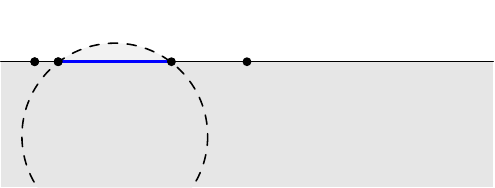}};
\node at (0,-6) {\includegraphics[scale=0.75]{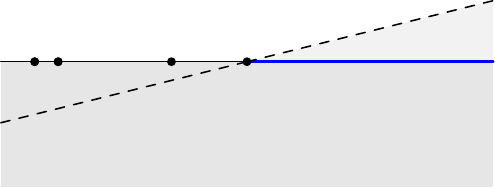}};
\end{tikzpicture}
\caption{Pivoting maximal disk/half-plane immersions. The egress component $J$ is in blue.}
\label{fig:pivot}
\end{figure}

See \Cref{fig:pivot} for an illustration of these deformations.

\medskip

If $J$ is a bounded segment then:
If $\Delta$ is a plane and $J$ is a bounded segment, then for all $\theta \in (0,\pi)$, $\Delta_\theta$ is a disk and we define $\theta_0=\pi$ and $\Delta_{\theta_0}$ as the open half-plane bounded by $\partial \Delta$ but opposite to $\Delta$.
If $\Delta$ is a disk, then there is some $\theta=\theta_0$ for which $C_\theta$ is the line through $A$ and $B$. We define $\Delta_{\theta_0}$ as the half-plane bounded by $C_{\theta_0}$ and which contains $J$.
Consider then the supremum $\theta_1\in(0,\theta_0]$ of $\theta\in(0,\theta_0)$ such that $f$ extends to an affine immersion of $\Delta \cup\Delta_{\theta}$.
Then $f$ extends to an immersion $\hat f:\Delta_{\theta_1}\cup J\cup \Delta \to X^*$ and $(\Delta_{\theta_1},\hat f|_{\Delta_{\theta_1}})$ is a rigid maximal element of the Delaunay category.

If $\Delta$ is a plane and $J$ is an infinite half-line, then for all $\theta \in (0,+\infty)$, $\Delta_\theta$ is a half-plane.
There is a maximal $\theta_1\in(0,+\infty]$ and for all $\theta\in I:=[0,\theta_1]\cap\R$ an immersion $f_\theta:\Delta_\theta\to X^*$ such that for all $\theta,\theta' \in I$ with $|\theta-\theta'|<2\pi$, $f_\theta$ and $f_{\theta'}$ coincide on $\Delta_\theta\cap \Delta_{\theta'}$.

We get a family of maximal elements of the Delaunay category parametrized by $\theta \in [0,\theta_1]\cap \R$. We call this process \emph{pivoting}.

\loctitle{The affine Delaunay spine}

This process naturally endows the set of isomorphism classes of maximal elements of the Delaunay category with a graph structure.
The vertices are the classes of rigid element, and the edge interiors points are the classes of flexible elements: indeed, flexible element have exactly two egress components, while rigid ones have at least three (with the exception of whole plane immersions, in which case the graph is a single point).
It also comes with a natural metric along edges given by the pivoting angle.
In this case there may be infinitely long edges in the presence of semi-infinite Reeb cylinders (\ie of a Reeb type Fuchsian singularity) or irregular singularities.
This graph is an analogue of the Voronoi graph of a finite set of marked points in the Euclidean plane. Actually in that setting, the Delaunay decomposition we defined coincides with the classical one, and it is known to be dual to the Voronoi decomposition: the center of all maximal disks trace the edges and vertices of the Voronoi decomposition associated to the marked point.
In our setting of a finite type affine surface, the center of maximal disks also trace a graph on $X^*$, which can be called the Voronoi decomposition, and is homeomorphic to the part of the Delaunay spine excluding the half-plane immersions.
This is interesting because distance is not well-defined on affine surfaces (this is already the case for dilation surfaces, see \cite{BGT25}).
To embed the full Delaunay spine, one would have to blow up non-conical singularities and irregular singularities in a certain way.

\subsection{Delaunay segments}\label{sub:delSeg}

We recall here the definition of a Delaunay segment introduced in \Cref{defn:DelaunayEdge}. Immersions of type A are immersions of open disks that extend continuously to the boundary circle, with exactly two points on the boundary mapped to singularities of $\nabla$. Delaunay segments are the images of the chords drawn between the two singular points under these immersions. In this section, we show that Delaunay segments are non-self-intersecting saddle connections with mutually disjoint interiors. These claims follow from a lemma of Euclidean geometry, which state without proof.

\begin{figure}[htbp]
  \begin{center}
    \begin{tikzpicture}
      \node at (0,0) {\includegraphics[scale=0.8]{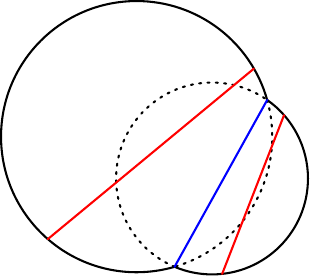}};
      \node at (1.75,0.8) {$I$};
      \node at (0.2,-2.1) {$J$};
      \node at (-1.7,-1.6) {$A_1$};
      \node at (1.6,1.3) {$B_1$};
      \node at (2.25,0.5) {$A_2$};
      \node at (1,-2.2) {$B_2$};
      \node at (-2.1,1.2) {$\alpha_1$};
      \node at (2.4,-1.3) {$\alpha_2$};
    \end{tikzpicture}
  \end{center}
  \caption{Illustration of \Cref{lem:TwoChords}}
  \label{fig:chords}
\end{figure}

\begin{lem}\label{lem:TwoChords}
Let $D_{1},D_{2}$ be two open disks in the plane with a nonempty intersection $D_{1} \cap D_{2}$ and such that neither $D_{1} \subset D_{2}$ nor $D_{2} \subset D_{1}$. We denote by $I,J$ the intersection points of their boundary circles $\partial D_{1}$ and $\partial D_{2}$. Let $\alpha_{1}$ (resp. $\alpha_{2}$) be the circular arc drawn between $I,J$ contained in $\partial D_{1}$ and disjoint from $D_{2}$ (resp. contained in $\partial D_{2}$ and disjoint from $D_{1}$).
Consider two open segments $(A_{1},B_{1})$, $(A_{2},B_{2})$ with endpoints respectively in $\alpha_{1}$ and $\alpha_{2}$, see \Cref{fig:chords}. If $(A_{1},B_{1})$ and $(A_{2},B_{2})$ intersect, then they both coincide with $(I,J)$.
\end{lem}

\begin{cor}\label{cor:DelaunayEmbedded}
In any affine surface $(X,\nabla,\mathcal{S})$, the interior of Delaunay segments are embedded finite straight open segments drawn between points of $\mathcal{S}$ (they have no self-intersection outside the endpoints).
\end{cor}

\begin{proof}
Consider a Delaunay open segment $e$ of $(X,\nabla,\mathcal{S})$ that is the image of a chord $(A,B)$ by the affine immersion of an open disk $\Delta$. We assume by contradiction that $e$ admits a self-intersection in some regular point $z \in X^{\ast}$ and consider two distinct points $x_{1},x_{2}$ in $(A,B)\cap f^{-1}(\{z\})$.
\par
By \Cref{sub:hol}, there exists $g\in\AutC$ such that $g(x_{1})=x_{2}$ and $f=f\circ g^{-1}$ holds on $\Delta\cap g(\Delta)$.
This allows to extend $f$ to an immersion $\hat f:\Delta\cup g(\Delta)\to X^*$.
Both immersions $\hat f\vert_{\Delta}=f$ and $\hat f\vert_{g(\Delta)}=f\circ g^{-1}$ are immersions of type~A (exactly two points of their boundary circles are mapped to singularities of the affine structure; for the second these are $g(A)$ and $g(B)$). It follows that neither $\Delta \subset g(\Delta)$ nor $g(\Delta) \subset \Delta$ can hold: indeed if we had, say, $g(\Delta) \subset \Delta$ then since neither $g(A)$ nor $g(B)$ can belong to $\Delta$, this means $g(\Delta)=\Delta$. In that case, the points mapped to $\cal S$ by $\ov f\circ g^{-1}$ are only $A$ and $B$, so $\{g(A),g(B)\} = \{A,B\}$. So either $g=\on{id}_\C$, which contradicts $g(x_1)=x_2$, or $g$ permutes $A$ and $B$, which means $[A,B]$ is a diameter of $\Delta$ but then $g$ fixes the midpoint of $[A,B]$ in $\Delta$, which is impossible by \Cref{lem:HolFixNeg}.
\par
By definition, the intersection $\Delta \cap g(\Delta)$ is nonempty because it contains a neighborhood of $x_{2}$. Then the open chord $(A,B)$ of $\Delta$ and the open chord $g((A,B))$ of $g(\Delta)$ meet (at $x_2$).
On the other hand, their endpoints cannot belong to $\Delta\cup g(\Delta)$ for then they would map by $\hat f$ to $X^*$ but also to $\cal S$, a contradiction.
They thus form a configuration of segments in $\Delta \cup g(\Delta)$ that is forbidden by \Cref{lem:TwoChords}.
\end{proof}

Similar arguments allow us to show below that Delaunay segments do not intersect each other.

\begin{cor}\label{cor:DelaunayDisjoint}
In any affine surface $(X,\nabla,\mathcal{S})$, Delaunay segments have disjoint interiors.
\end{cor}

\begin{proof}
We assume by contradiction that two Delaunay segments $e_{1},e_{2}$ of  $(X,\nabla,\mathcal{S})$ intersect transversely in some common interior point $z \in X^{\ast}$. Then, there is an affine immersion $f$ of a open disk $D_{0} \subset \mathbb{C}$ into $X^{\ast}$ such that $z \in f(D_{0})$.
\par
Then $f^{-1}(e_{1})$ and $f^{-1}(e_{2})$ form two chords of $D_{0}$ we extend $f$ along two disks $D_{1}$ and $D_{2}$ in such a way that $f\vert_{D_{1}}$ and $f\vert_{D_{2}}$ are immersions of type A. Observe that neither $D_{1} \subset D_{2}$ nor $D_{2} \subset D_{1}$ are possible since the boundary circle of each of them contains points mapped to singularities. The intersection $D_{1} \cap D_{2}$ is nonempty because it contains a neighborhood of $z$. It follows that $f^{-1}(e_{1})$ and $f^{-1}(e_{2})$ form a configuration of segments in $D_{1} \cup D_{2}$ that is forbidden by \Cref{lem:TwoChords}.
\end{proof}

Delaunay segments form a system of disjoint arcs drawn between the singularities of the meromorphic connection. 

\begin{lem}\label{lem:DelaunayFinite}
  There are finitely many Delaunay segments.
\end{lem}
\begin{proof}
  This will follow from \Cref{lem:SystemTopologicalArcs2} once we have checked that any non-empty and finite subset $\cal A$ of the set of Delaunay segments forms an apical arc system (\Cref{def:aas}) and that in the complement in $X$ of the support of $\cal A$, every ordinary face (topological disk not containing an interior singularity or swath neighborhood) has at least three sides.
  If such a face $F$ had less, consider a developing map $\phi:F\to\C$.
  The map $\phi$ has a continuous extension to the closure of $F$.
  This is immediate near interior points of the Delaunay segments bounding $F$, since they are regular points. Near a conical point, this comes from the local model of such point. Near an irregular point $p$, this comes from the local model of $p$ and from the fact that $F$ contains no swath neighborhood of $p$.
  The image $\phi(F)$ is thus a bounded open subset of $\C$, whose boundary is contained in the union of one or two segments.
  Such a set does not exist.
\end{proof}

\begin{lem}\label{lem:bdyhpdel}
  Any saddle connection image of a boundary segment by a continuous extension of a half-plane immersion, is a Delaunay segment.
\end{lem}
\begin{proof}
  Denote $(\Delta,f)$ the half-plane immersion, and $[A,B]$ the boundary segment.
  It is possible to apply the pivot process of \Cref{sub:AffineSpine} to $(\Delta,f)$ with egress component $(A,B)$.
  We get a type~A Delaunay element whose chord $[A,B]$ defines a Delaunay segment according to \Cref{defn:DelaunayEdge}.
\end{proof}

%Chercher dans le texte où on a implicitement utilisé ce résultat

\subsection{Delaunay polygons and the core}\label{sub:core}

\begin{defn}\label{defn:delPol}
  An \emph{open Delaunay polygon} is the image by a rigid maximal immersion $(\Delta,f)$ of type~I (disk, rigid) of the interior of the convex hull $P$ of $\bar f^{-1}(\cal S)$ (the points of $\partial \Delta$ that map to singular points by the extension of $f$ to $\partial \Delta$).
\end{defn}

\begin{lem}
  In \Cref{defn:delPol}, the image by $\bar f$ of the sides of $P$ are Delaunays segments.
  No diagonal of the convex polygon $P$ is mapped  a Delaunay segment.
\end{lem}
\begin{proof}
  For the first claim, note that we can always pivot (see \Cref{sub:AffineSpine}) the immersion around the two endpoints of a side of $P$ so as to get a type~A (disk, flexible) immersion. 

  We proved in \Cref{sub:delSeg} that Delaunay segments do not self-intersect and that no two Delaunay segments can cross either.
  The only possibility for the presence of a point in a Delaunay segment in an open Delaunay polygon would be that $P$ has at least $4$ sides and that a diagonal of $P$ maps to a Delaunay segment.
  There would be a type~A (disk, flexible) immersion $(\Delta',f')$ for which only two boundary points $A',B'$ map to $\cal S$ and such that $f'((A,B))=f(D)$.
  Then by composing with an appropriate automorphism of $\C$ and possibly permuting $A'$ and $B'$ we may assume that $A'=A$, $B'=B$ and $f'=f$ on a neighborhood of $(A,B)$.
  The points $A$ and $B$ are thus on $\partial \Delta'$.
  Then $\Delta'=\Delta$ for otherwise $\Delta'$ would contain a point $C\in f^{-1}(\cal S)$, and by affine continuation and continuity we would get $f'(C)\in\cal S$ which is forbidden for affine immersions.
  From $\Delta'=\Delta$ we get that $f'=f$ is not a type~A immersion, but $f$ is of type~I: this is a contradiction.
\end{proof}

The next result tells that the interior of Delaunay polygons are always embedded.

\begin{cor}\label{cor:DelaunayPolygon}
Given a finite type affine surface $(X,\nabla,\mathcal{S})$, we consider an immersion $f:\Delta \to X^{\ast}$ of an open disk into $X^{\ast}$ that extends continuously to its boundary circle, where a finite number $m \geq 3$ of boundary points $A_{1},\dots,A_{m}$ are mapped to points of $\mathcal{S}$. Then, the restriction of $f$ to the interior $\mathcal{P}$ of the convex hull of $\lbrace{ A_{1},\dots,A_{m}\rbrace}$ is an embedding.
\end{cor}

\begin{proof}
We consider such an immersion and assume by contradiction that $\left. f \right|_{\mathcal{P}}$ is not injective. Then there exists $\phi \in \Aut(\mathbb{C})$ such that:
\begin{itemize}
    \item $\phi\neq\on{id_\C}$;
    \item $\phi^{-1}(\mathcal{P}) \cap \mathcal{P} \neq \emptyset$;
    \item for any $z \in \phi^{-1}(\Delta) \cap \Delta$, $f \circ \phi (z) = f(z)$.
\end{itemize}
The open convex subset $\mathcal{P}$ and $\phi(\mathcal{P})$  of $\mathbb{C}$ have nonempty intersection.
Moreover we saw that no point in $\cal P$ or $\cal \phi(P)$ can be on a Delaunay segment.
The only possibility is that $\phi(\cal P)=\cal P$.
Then $\phi$ has a fixed point contained in the open polygon $\mathcal{P}$. This contradicts \Cref{lem:HolFixNeg}.
\end{proof}

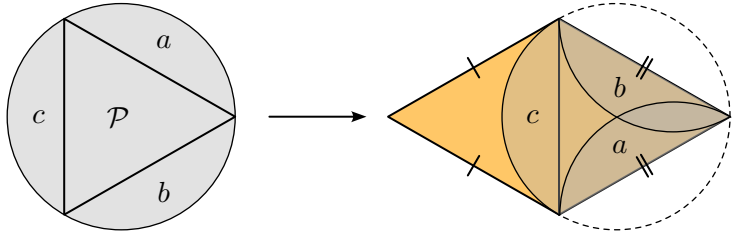
\begin{figure}[htbp]
  \begin{center}
    \begin{tikzpicture}
      \node at (0,0) {\includegraphics[scale=0.75]{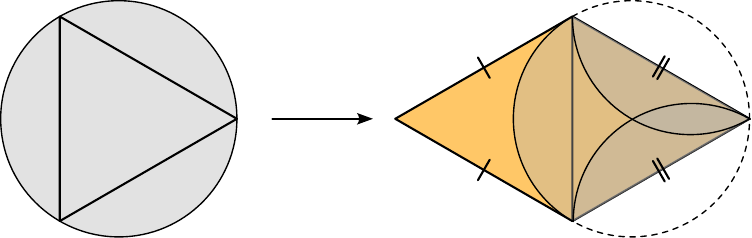}};
      \node at (-3.3,0) {$\cal P$};
      \node at (-2.7,1) {$a$};
      \node at (-2.7,-1) {$b$};
      \node at (-4.35,0) {$c$};
      \node at (3.3,-0.4) {$a$};
      \node at (3.3,0.4) {$b$};
      \node at (2.15,0) {$c$};
    \end{tikzpicture}
  \end{center}
  \caption{The affine surface on the right is formed of two equilateral triangles glued together in a specific way. The disk on the left is immersed in it as indicated. This map injective on the interior of $\cal P$ but not on its boundary.}
  \label{fig:noninjonbdy}
\end{figure}

The map $\bar f$ in not necessarily injective on $\partial\cal P$, see \Cref{fig:noninjonbdy}.

\begin{defn}\label{defn:core}
  Given a finite type affine surface $(X,\nabla,\mathcal{S})$, we call \emph{exterior} the union of the non-conical Fuchsian singularities and of the image of all affinely immersed half-plane.
  We call \emph{core} the complement of the exterior.
\end{defn}

The core thus always contains the conical Fuchsian singularities and the irregular singularities and no other singularity.
It may be reduced to these singularities, in which case the surface is exceptional, see \Cref{sub:Exceptional}.
Actually we have slightly better: in a non-exceptional surface, every singularity in the core is at the end of at least one Delaunay edge (\Cref{lem:de0}).
Reeb cylinders are in the exterior, as are appropriate standard neighborhoods of non-conical singularities and swathes.
Note that for translation surfaces, the core is always connected, see \cite{Ta18}. However, for affine surfaces, the core can be disconnected by Reeb cylinders.

\begin{lem}\label{lem:finSegNum}
  Consider a disk immersion $(\Delta,f)$ that extends continuously to the boundary.
  Let $\cal A$ be a finite set of geodesic segments of $X$, of finite life span (so ending on regular points, conical Fuchsian singularities or foci of irregular singularites).
  Let $A$ be the union of the support of the elements of $\cal A$.
  Then $f^{-1}(A)$ consists in finitely many segments.
\end{lem}
\begin{proof}
  By \Cref{lem:b} there are only finitely many points in $\partial \Delta$ mapped to singularities, where there are slit extensions of $f$, and at all the other points there are affine extensions.
  It follows we can extend $f$ to a neighborhood of $\ov\Delta$, minus a finite number of slits form the points mapped to singularities.
  For the extension, the number of segments is locally finite as it is for any affine immersions.
  This is also the case near the tip of a slit on $\partial\Delta$, because of the local model of the singularity (which is either conical or a focus of an irregular singularity): finitely many segment ends go straight to it, and they have at most one preimage in the slit in the focus case, and finitely many in the conical case; the other segments stay at some distance of the singularity.
\end{proof}

\begin{prop}\label{prop:DelPolCoverCore}
The core is the union of the open Delaunay polygon interiors, the Delaunay segment interiors, the conical Fuchsian singularities and the irregular singularities.
\end{prop}

\begin{proof}
Given a finite type affine surface $(X,\nabla,\cal S)$,
assume $x\in X^*$ belongs to the core.
Consider a local chart $f_0:\Delta_0\to X^*$ with $\Delta_0=B(z_0,r_0)$
and $f_0(z_0)=x$.
\Cref{lem:charMax2} ensures that there exists a maximal immersion $(\Delta_1,f_1)$ of a disk/half-plane/plane containing $\Delta_0$ and extending $f_0$.
Then $\Delta_1$ must be a disk for otherwise $x$ would belong to the exterior.
By \Cref{sub:maximal}, $(\Delta_1,f_1)$ is of type~A (disk, flexible) or~I (disk, rigid).
If it is of type~A, this disk defines a particular Delaunay segment $[A,B]\subset\ov\Delta_1$.
If $z_0\in[A,B]$, then we are done.
Otherwise it is on one side of the segment.
By \Cref{lem:DelaunayFinite} there are only finitely many Delaunay segments. The preimage of these segments in $\Delta_1$ by $f_1$ forms a  set $\cal C$ of open chords, including $(A,B)$.
This set is finite by \Cref{lem:finSegNum}.
Let us pivot the maximal disk towards $z_0$ (see \Cref{sub:AffineSpine}) until we meet a rigid maximal immersion $(\Delta_2,f_2)$.
The domain $\Delta_2$ necessarily contains $z_0$.
So it cannot be a half-plane, for $z_0$ is not in the exterior.
Being rigid, it necessarily defines a Delaunay polygon $P_2$, moreover, $[A,B]$ is one of its edges.
If $z_0$ belongs to the closed polygon $P_2$, we are done.
Otherwise, it belongs to a part of the disk $\Delta_2$ cut by an edge $[A_2,B_2]$ of $P_2$ that is closer to $z_0$ than $[A,B]$.
Moreover $[A_2,B_2]\cap \Delta_1$ is necessarily a chord in the set $\cal C$.
We then repeat the process and pivot the maximal disk $\Delta_2$ to get a new rigid maximal immersion $(\Delta_3,f_3)$, etc.
At each step $n$, if $z_0$ is not in the closed polygon $P_n$, then the side $[A_n,B_n]$ of $P_n$ separating $z_0$ from the interior of $P_n$ in $\Delta_n$, satisfies that $[A_n,B_n]\cap \Delta_1$ is one of the chords in $\cal C$, and moreover these chords get closer and closer to $z_0$ in the sense that the component of $\Delta_1$ minus the chord that contains $z_0$ is strictly decreasing for the inclusion.
This implies by finiteness of $\cal C$ that the process must stop at some point, i.e.\ that $z_0\in P_n$ for some $n$.
\end{proof}

We did not state above that the core is the union of Delaunay polygons. Indeed it may happen that the core has connected components that are graphs, or single points: consider for instance the whole plane from which a square is removed, and the sides of the square are glued together (we can get a sphere or a torus according which pairing is chosen). Its core consist in the two Delaunay segments on the boundary of the removed square.

\subsection{Exceptional affine surfaces}\label{sub:Exceptional}

These surfaces have several characterization, one being that they have no Delaunay segments. We start with a useful lemma. 

\begin{lem}\label{lem:de0}
  Assume a finite type affine surface $(X,\nabla,\cal S)$ has a conical Fuchsian singularity $s\in \cal S$ with no Delaunay segment ending on it, or an irregular singularity with a focus $s$ with no Delaunay segment ending on it.
  Then $X^*$ is either isomorphic to the exponential-affine plane or to a skew cylinder.
\end{lem}

\begin{proof}
Every maximal element of the Delaunay category is a half-plane immersion.
By the local model of $s$, there exists a disk immersion $(\Delta,f)$ which extends continuously at $\partial \Delta$ (in the case of a focus, we use the extended topology defined in \Cref{sub:helix}) and such that one point $z_0\in\partial \Delta$ is mapped to $s$.
Consider a maximal Delaunay element $(\Delta',f')$ extending $(\Delta,f)$.
Then $\Delta'$ is a half-plane and $z_0\in\partial\Delta'$.
There can be no discontinuity point $z_1\in\partial \Delta'$ (for the usual topology on $X$) for otherwise we would be in the situation where $f'\circ\phi=f'$ for some dilation $\phi$ fixing $z_1$, and every point of the form $\phi^n(z_0)$, $n\in \Z$, would be mapped to $s$ so there would be, by \Cref{lem:bdyhpdel}, Delaunay segments obtained by taking the image by the extension of $f'$ of the portion of $\partial \Delta'$ between two consecutive points mapped to $\cal S$.
Then, apart from $z_0$, there can not be any other point of $\partial \Delta'$ mapped to $\cal S$ for otherwise there would be a Delaunay segment for the same reason.
It means we can pivot $\Delta'$ in the two directions, by a rotation around $z_0$ (see \Cref{sub:AffineSpine}).
This rotation goes on forever, in both directions, for an hypothetical half-plane immersion at the end would contradict absence of Delaunay segment for the same reasons as above.
This implies that $\cet$ immerses into $X^*$.
We then conclude with \Cref{prop:EImmersion} that $X^*$ either admits an affine immersion of the plane, a situation which is classified in \Cref{prop:PlaneImmersion}, or is isomorphic the exponential-affine plane, to a full Reeb cylinder, to a skew cylinder or to an affine torus. The translation and affine tori have no singularity so are excluded. A whole plane has only one singularity, which is anti-conical Fuchsian, a full translation cylinder has only two, which are cylindrical Fuchsian, and a full Reeb cylinder only two, which are of Reeb type, so these cases are excluded too.
\end{proof}

\begin{prop}\label{prop:NoDelaunayEdge}
If a finite type affine surface $(X,\nabla,\mathcal{S})$ has no Delaunay segments, then it is isomorphic to one of the following models:
\begin{itemize}
    \item[(a)] \textit{whole plane}: $X$ is a sphere, $X^{\ast}=\mathbb{C}$ and $\nabla$ has a unique Fuchsian singularity of residue equal to $2$;
    \item[(b)] \textit{translation cylinder}: $X$ is a sphere, $X^{\ast} = \mathbb{C}/\mathbb{Z}$ and $\nabla$ has two Fuchsian singularities with residues equal to $1$ (the two ends of the cylinder);
    \item[(c)] \textit{translation torus}: $X=\mathbb{C}/\Lambda$ where $\Lambda$ is a lattice of $(\mathbb{C},+)$ and $\nabla$ has no singularity.
    \item[(d)] \textit{infinite angle cone}: $X$ is a sphere, $X^{\ast}$ coincides with the universal covering $\mathcal{E}$ of the punctured plane $\mathbb{C}^{\ast}$ (also called the \textit{exponential-affine plane}) and the only singularity of $\nabla$ is a double pole, of residue equal to $2$;
    \item[(e)] \textit{affine cylinder} (skew cylinder or full Reeb cylinder): $X$ is a sphere and $X^{\ast}$ is a quotient of $\mathcal{E}$ by a linear transformation; 
    \item[(f)] \textit{affine torus}: $X$ is a quotient $\mathcal{E}$ by a lattice and $\nabla$ has no singularity.
\end{itemize}
\end{prop}

\begin{proof}
If $X$ has a singularity, it is covered by \Cref{lem:de0} and we are in situation (a), (b), (d) or (e).
If $X$ has no singularity then consider a maximal Delaunay category element $(\Delta,f)$.
The set $\Delta$ cannot be a disk for there would be a Delaunay edge. If $\Delta$ is a half-plane we are in the type~B case (Reeb, flexible) of \Cref{defn:flexible}.
In this case we can proceed exactly as in the proof of \Cref{lem:de0} and deduce that there is an immersion of the exponential-affine plane $\cal E$ in $X^*$ and we get the same conclusion, except that case (f) is possible.
Finally, \Cref{prop:PlaneImmersion} classifies the case $\Delta$ is the whole plane, which leads to cases (a), (b) or (c).
\end{proof}

Actually the converse holds. It can be directly seen on each models. Alternatively, it follows from \Cref{thm:ExceptionalSurfaces} (for instance all the examples satisfy criterion (iii): $2g+n\leq 2$).

\begin{proof}[Proof of \Cref{thm:ExceptionalSurfaces}]
We have to prove that for a finite type affine surface $(X,\nabla,\mathcal{S})$ of genus $g$ with $n$ poles (counted with multiplicity), the following statements are equivalent:
\begin{enumerate}[(i)]
    \item[(o)] $X$ is in the list of \Cref{prop:NoDelaunayEdge};
    \item every maximal geodesic on $X^{\ast}$ is infinite in the future or in the past;
    \item $X$ has no Delaunay segment;
    \item $2g+n \leq 2$;
    \item the universal covering of $X^{\ast}$ is either the flat plane $\mathbb{C}$ or the exponential-affine plane $\mathcal{E}$;
    \item there exists an affine immersion $\cal E\to X^*$;
    \item the core  is contained in $\cal S$.
\end{enumerate}
Let us justify
(o) $\Rightarrow$ (iv) $\Rightarrow$ (i) $\Rightarrow$ (ii) $\Rightarrow$ (o), (o) $\Leftrightarrow$ (iii), (iv) $\Leftrightarrow$ (v) and (ii) $\Leftrightarrow$ (vi).

\medskip

(o) $\Rightarrow$ (iv) follows from inspection of the list: (a), (b), (c) are quotients of $\C$ and (d), (e), (f) of $\cal E$.

(iv) $\Rightarrow$ (i): every geodesic lifts to the universal covering $\C$ or $\cal E$. In $\C$, every geodesic have infinite lifespan in the past and future. In $\cal E\cong \cet$ every geodesic too, except those aimed at $0$, which hit $0$ in finite time and are infinite lived in the other direction.

(i) $\Rightarrow$ (ii): every Delaunay segment is a goedoesic whose lifespan is finite in both directions.

(ii) $\Rightarrow$ (o): this is \Cref{prop:NoDelaunayEdge}.

\medskip

(o) $\Rightarrow$ (iii) follows by inspection of the list: (a) and (d) $(g,n) = (0,1)$, (b) and (e) $(g,n) = (0,2)$, (c) and (f) $(g,n) = (1,0)$.

(iii) $\Rightarrow$ (o): this is \Cref{prop:LowComplexity}.

\medskip

(iv) $\Rightarrow$ (v): if there is a universal covering $\pi: \cal E\to X^* $ then this is trivial: $\pi$ is an affine immersion from $\cal E\to X^*$. If there is a universal covering $\pi: \C \to X^*$ then $\pi\circ\exp:\C\to X^*$ is an affine immersion (where $\exp:\cal E\to \C^*\subset\C$).

(v) $\Rightarrow$ (iv): see \Cref{prop:EImmersion} and \Cref{prop:PlaneImmersion}.

\medskip

(ii) $\Leftrightarrow$ (vi): immediate from the definitions and the fact that Delaunay polygons are bounded by Delaunay segments.
\end{proof}

\subsection{Classification of components of the Delaunay decomposition}\label{sub:DelaunayComponents}

In this section we prove \Cref{thm:ClassificationComponents}.

Consider a non-exceptional finite type affine surface $(X,\nabla,\cal S)$.
By (ii) of \Cref{thm:ExceptionalSurfaces} there is at least one Delaunay segment.
By \Cref{cor:DelaunayEmbedded,cor:DelaunayDisjoint} the Delaunay segment interiors are embedded and disjoint.
Consider the union of all Delaunay segments.
By \Cref{lem:de0}, the endpoints contains all conical Fuchsian singularities and all irregular singularities.
The converse inclusion holds by definition of Delaunay segments.
By \Cref{lem:DelaunayFinite} there are finitely many Delaunay segments.
Their union is thus a closed subset of $X$ and the complement of this union is an open subset of $X^*$.
Its connected components are called the \emph{Delaunay components}.
\Cref{thm:ClassificationComponents} states that these components are of the following types:
\begin{itemize}
    \item[(i)] \textbf{open Delaunay polygon}: defined in \Cref{defn:delPol}, they are in particular embedded convex polygons with $p \geq 3$ sides drawn between points of $\mathcal{S}$;
    \item[(iia)] \textbf{finite angle Reeb cylinder}: dilation cylinders of angle at least $\pi$;
    \item[(iib)] \textbf{Semi-infinite Reeb cylinder} and its end, a Reeb-type Fuchsian singularity;
    \item[(iii)] \textbf{Semi-infinite translation cylinder} and its end, a cylindrical Fuchsian singularity;
    \item[(iv)] \textbf{anti-conical domain} associated to an anti-conical Fuchsian singularity;
    \item[(v)] \textbf{swath domain} associated to an irregular singularity swath.
\end{itemize}
We refer to \Cref{sub:cylindersAndSkewCones} for the description of the cylinders, to \Cref{subsub:UnboundedDomains} for anti-conical and swath domains and to \Cref{sub:core} for Delaunay polygons.

\begin{lem}\label{lem:uddadc}
An anti-conical or swath domain is a Delaunay component.
\end{lem}
\begin{proof}
They are bounded by simple saddle connections $\gamma$ that do not cross each other. It is enough to prove that these connections are Delaunay segments.
For each $\gamma$, \Cref{lem:Unbounded} allows to define a half-plane immersion, that will be rigid of type~II, and whose boundary extension has a segment mapping to $\gamma$.
By \Cref{lem:bdyhpdel}, $\gamma$ is a Delaunay segment.
%Applying the pivot process with this segment as egress component, we get a flexible maximal disk immersion which defines $\gamma$ as a Delaunay segment.
\end{proof}

Similarly, any maximal embedded and non-full Reeb cylinder and its eventual end singularity form a Delaunay component, and this is the case too for a maximal embedded non-full translation cylinder and its eventual end singularity.

\medskip

Let us prove \Cref{thm:ClassificationComponents}.
Consider a Delaunay component $\mathcal{C}$. Let $x\in \cal C$.
If $x$ is in the core of $X$ (defined in \Cref{defn:core}), then \Cref{prop:DelPolCoverCore} proves that $x$ is in an open Delaunay polygon.
Otherwise, $x$ is in the exterior of $X$ (the complement of the core).
Let $(\Delta,f)$ be a half-plane immersion such that $x\in f(\Delta)$.
Without loss of generality, we can assume that $f$ is maximal in the sense of \Cref{defn:DelaunayMaximal}. These maximal immersions are classified in \Cref{sub:maximal}, see \Cref{defn:flexible,defn:rigid}.
Note that $f(\Delta)$ is contained in the exterior, hence in $\cal C$.
\par
If $f$ is a rigid immersion of type~III (translation cylinder), then its image is a semi-infinite translation cylinder open neighborhood of a cylindrical Fuchsian singularity bounded by a chain of saddle connections that are Delaunay segments by \Cref{lem:bdyhpdel} (see point \Cref{item:cylEmbed} in \Cref{prop:halfPlaneImmersion}). It follows that $\mathcal{C}$ coincides with the translation cylinder and its end.
\par
If $f$ is  a flexible immersion of type~C (Reeb) or a rigid immersion of type~IV (Reeb, opening angle $>\pi$) or~$V$ (Reeb, opening angle $\pi$), then \Cref{prop:halfPlaneImmersion} shows that its image is contained in a dilation cylinder.
Consider the maximal dilation cylinder $\cal C'$ containing $f(\Delta)$.
It has two ends, which may be finite or infinite (in terms of the angle span).
An infinite end is a neighborhood of a Reeb type Fuchsian singularity by \Cref{lem:embedFuchsian}.
A finite end is necessarily bounded by Delaunay segments by \Cref{prop:halfPlaneImmersion,lem:bdyhpdel}.
Again, $\mathcal{C}$ coincides with the Reeb cylinder.
\par
In the remaining cases, for any $x \in \mathcal{C}$, $f$ is either a flexible immersion of type~B (non-cylinder) or a rigid immersion of type~II (non-cylinder).
In these cases, the map $f$ extends continuously to $\partial \Delta$ and $f^{-1}(\cal S)$ consists in finitely many points.
They cut $\partial \Delta$ into two half-lines and finitely many segments ($0$ for type~B).
We can apply the pivot process described in \Cref{sub:AffineSpine} with egress component either of the two half-lines.
Initially we get type~B immersions. In the process we may meet rigid immersions, but note that they have to be of type~III because in the other cases, we already classified the situation and they exhibit no type~B immersions.
We get a continuous family of immersions $(\Delta_\theta,f_\theta)$ parametrized by the pivot angle $\theta\in I$, with $I$ an open interval contained in $\R$ (we will soon prove that $I$ can be taken equal to $\R$).
These immersions all map to the exterior, so, by connectedness, to $\cal C$.

Applying \Cref{prop:halfPlaneImmersion} to each $(\Delta_\theta,f_\theta)$, we get that we are in situation~\eqref{item:embed} or~\eqref{item:antiConeEmbed} of that statement.
In particular, at infinity $f_\theta$ tends to a singularity $y\in X$, which is either anti-conical Fuchsian or irregular.
Since the $f_\theta$ overlap, we gradually get that all $f_\theta$ tend to the same $y$ at infinity.

Rigid immersions occur for isolated parameters $\theta$, and each defines one or more, but finitely many, consecutive Delaunay segments sides on $\partial \cal C$.
Moreover, two consecutive rigid parameters define consecutive Delaunay segments sides on $\partial \cal C$.
Since there are only finitely many Delaunay segment, $\cal C$ is bounded by finitely many Delaunay segment sides (one segment can bound $\cal C$ on at most two of its sides).

\begin{rem}
Though not strictly necessary, it is interesting to interpret this in terms of the the affine Delaunay spine introduced in \Cref{sub:AffineSpine}, a graph whose points represent all the maximal immersion equivalence classes.
The real number $\theta$ parametrizes a point moving on a succession of edges on the part of the Delaunay spine corresponding to half-plane immersions.
We recall that interior points of the edge are the flexible immersions, while vertices are the rigid ones.
\end{rem}

There is at least one rigid parameter, for or otherwise by the proof of \Cref{lem:de0}, $X$ would be a skew cylinder, hence exceptional.

Two situations may occur:
\begin{itemize}
\item[--] Case 1: there are two different parameters $\theta$, $\theta'$ such that the associated immersions follow the same side of a given Delaunay segment.
Then these immersions are equivalent ($\phi(\Delta_\theta) = \Delta_{\theta'}$ and $f_{\theta'}\circ\phi=f_\theta$ or some $\phi\in\AutC$).
\item[--] Case 2: if we are not in Case~1, this implies that the number of rigid parameters $\theta$ is at most the number of sides of $\cal C$. This implies we can take $I=\R$: indeed after the last rigid immersion, the pivoting only happens around a single point, so at a finite upper/lower bound of $I$, there is an immersion, which can then be still pivoted since it cannot be rigid, leading to a contradiction.
\end{itemize}

If the function $\theta\in I\mapsto (\Delta_\theta,f_\theta)$ is not injective then using affine continuation we can periodize: we can take $I=\R$ and up to equivalence, the function $\theta\in\R\mapsto (\Delta_\theta,f_\theta)$ is necessarily $\theta_p$-periodic. Hence in both Cases~1 and~2, we can take $I=\R$.
And Case~1 is equivalent to the fact that $\theta\in\R\mapsto (\Delta_\theta,f_\theta)$ is periodic up to equivalence.

\loctitle{Let us assume we are in Case~1}
Take $\theta_p$ to be the minimal period.
%There exits $\theta_p>0$ such that the immersion $(\Delta_{\theta_p},f_{\theta_p})$ equivalent to the initial one, $(\Delta_0,f_0)$.
%By equivalent we mean $f_{\theta_p} = f_0\circ\phi_p$ and $\Delta_{\theta_p}=\phi_p^{-1}(\Delta_0)$ for some $\phi_p\in\AutC$.
%Using affine continuation we get that, up to equivalence, the function $\theta\in\R\mapsto (\Delta_\theta,f_\theta)$ is $\theta_p$-periodic.
Call
\[\cal U = \bigcup_{\theta\in[0,\theta_p)} f_\theta(\Delta_\theta) \subset \cal C\setminus\{y\}\]
note that if $y$ is an irregular singularity then $y\notin\cal C$; however we will prove in a moment that $y$ is anti-conical, that $\cal C=\cal U\cup\{y\}$ and that $\cal C$ is an anti-conical domain, as per \Cref{defn:Unbounded}. We will also prove as a bonus that $y$ has opening angle $\theta_p$.

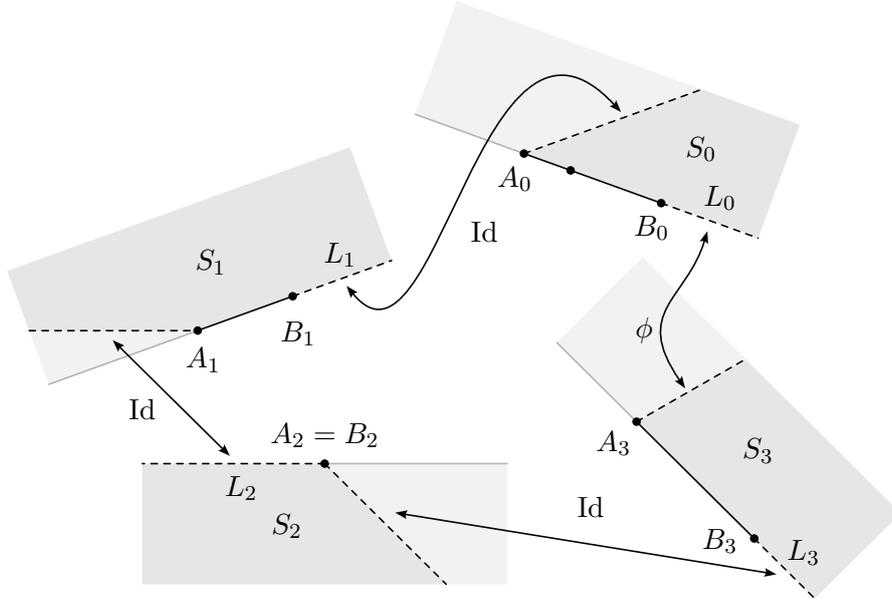
\begin{figure}
\begin{center}
\begin{tikzpicture}
\node at (0,0) {\includegraphics[scale=0.8]{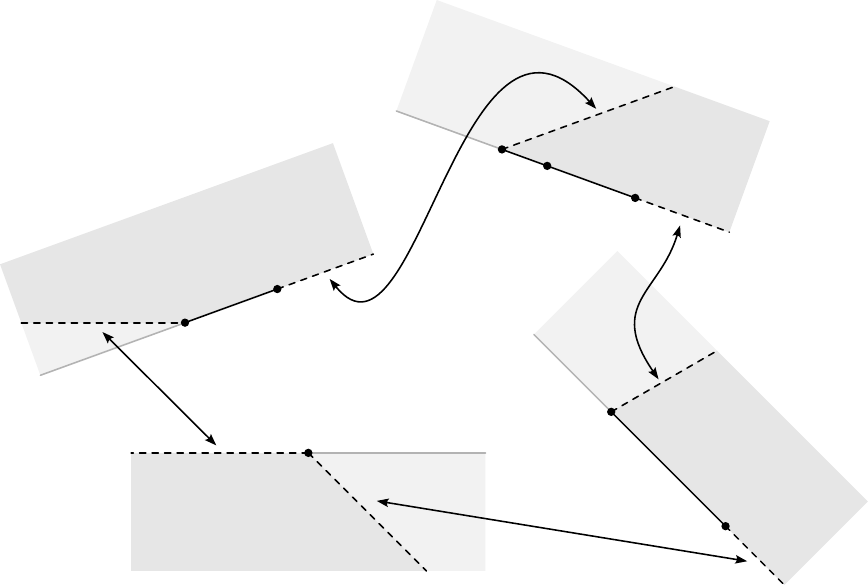}};
\node at (-4.1,-1.45) {$\Id$};
\node at (0.4,0.9) {$\Id$};
\node at (1.8,-2.75) {$\Id$};
\node at (2.5,-0.4) {$\phi$};
\node at (3.25,2) {$S_0$};
\node at (-3.2,0.5) {$S_1$};
\node at (-2.2,-3) {$S_2$};
\node at (4,-2) {$S_3$};
\node at (0.8,1.6) {$A_0$};
\node at (2.6,0.95) {$B_0$};
\node at (3.5,1.35) {$L_0$};
\node at (-3.3,-0.8) {$A_1$};
\node at (-2.05,-0.45) {$B_1$};
\node at (-1.5,0.6) {$L_1$};
\node at (-1.7,-1.8) {$A_2=B_2$};
\node at (-2.8,-2.5) {$L_2$};
\node at (2.1,-1.9) {$A_3$};
\node at (3.5,-3.2) {$B_3$};
\node at (4.6,-3.35) {$L_3$};
\end{tikzpicture}
\end{center}
\caption{Finite angle unbounded domain reconstruction.}
\label{fig:unbddModel}
\end{figure}

There are only finitely many $\theta\in[0,\theta_p)$ such that $(\Delta_\theta,f_\theta)$ is rigid, \ie such that the extension of $f$ to $\partial\Delta_\theta$ passes more that once through a singularity, and there is at least one.
Without loss of generality we assume that $\theta=0$ is one of them.
Consider the circular list of rigid parameters in $[0,\theta_p]$.
If necessary, insert finitely many (flexible) $\theta$ in the list so that every consecutive elements of the list differ by less than $\pi$.
Call $\theta_0=0<\theta_1<\ldots<\theta_{k-1}<\theta_p$ the obtained finite list.
For each $\theta_i\neq\theta_p$, orient the line $\partial \Delta_{\theta_i}$ so that $\Delta_{\theta_i}$ is on the left of the line, and denote $A_i$ the first point on $\partial \Delta_{\theta_i}$ that is mapped to a singularity and $B_i$ the last.
Then $\cal U$ decomposes as the union of $f(S_i)$ where $S_i\subset\Delta_i$ is a sector based on $A_i$, with one boundary line going through $B_i$ and the other being the rotation of this line by $\theta_{i+1}-\theta_i\in(0,\pi]$, and of $f(L_i)$ where $L_i\subset\partial\Delta_i$ is the half-line from $B_i$ and onward.
See \Cref{fig:unbddModel}.
The disjoint union of the $S_i$ and $L_i$ defines an affine surface $\cal A$ that surjectively immerses in $\cal U$ via the $f_{\theta}$:
\[F: \cal A\to \cal U.\]

Let us show that $F$ is proper from $\cal A$ to $\cal C\setminus\{y\}$.
Consider a sequence $a_n\in\cal A$ that leaves every compact subset of $\cal A$.
We must prove that $F(a_n)$ tends to $\partial \cal C\cup\{y\}$, i.e.\ that a subsequence tends to a point of that set.
If $a_n$ tends to infinity in the $S_i$ or the seams between them then we already know that $F(a_n)$ tends to $y$.
Otherwise, a subsequence of $a_n$ tends to a point in one of the segments $[A_i,B_i]$ and by continuous extension of the $f_{\theta_i}$ (and of possibly the adjacent $f_{\theta_j}$) we get the claim.

As we saw, $F$ extends continuously by mapping $\infty$ to $y$, hence holomorphically.
In particular, $y$ is separated from $\partial C$ by an annulus, so $y$ cannot be one of the irregular singularities (which all belong to the union of Delaunay segments), hence $y$ is anti-conical.

As a holomorphic proper map without critical points over a connected set is a finite degree covering, so $F$ is a covering from $\cal A$ to $\cal C\setminus\{y\}$.
It is hence surjective and $\cal C\setminus\{y\}$ must be homeomorphic to a covering quotient of the pointed disk, so is to a pointed disk too, so $\cal C$ is simply connected.
It is bounded by a finite chain of Delaunay segments (hence saddle connections), with possible repetitions, and the angle between them is $\geq \pi$ according to the model $\cal A$.

\begin{rem}
We do not stricto sensu need it but we find it interesting to prove that $F$ has degree $\leq 1$, so is injective, hence an isomorphism from $\cal A$ to $\cal C\setminus\{y\}$.
% For this, it is enough to find a point that has at most one preimage. We will seek such a point near $y$.
By the same argument as above, if $x\in\cal C\setminus\{y\}$ is close to $y$ then all the elements in $F^{-1}(\{x\})$ are close to $\infty$.
% Note that infinity in $\cal A$ has a neighborhood that is isomorphic to a standard neighborhood of anti-conical singularity.
A neighborhood of $\infty$ in $\cal A$ is isomorphic to an anti-conical standard neighborhood with angle $\theta_p$ and ratio some $s>0$, with possibly a shift (see shifted anti-conical singularities in \Cref{sub:Fuchsian}).
If $\deg F>1$ then the deck transformation generator of the covering $F$ necessarily is a rotation of angle $\theta_p/\deg F$ (a priori it may be $>2\pi$ and it may include a translation in shifted case).
But then this would imply that the period of the family $(\Delta_\theta,f_\theta)$ up to equivalence, is strictly smaller than $\theta_p$, in contradiction with the minimality of $\theta_p$.
\end{rem}

\loctitle{Let us now assume we are in Case 2} No two $(\Delta_\theta,f_\theta)$, $\theta\in\R$, are equivalent.
We can still define an affine surface $\cal A$ with an immersion $F: \cal A\to \cal C$.
This times it needs infinitely many flexible $\theta_i$, more precisely: finitely many between rigid parameters; after the last and before the first, an infinite sequence that we can take  spaced by $\pi/2$ for instance. Let $\theta_i$, $i\in \Z$ the increasing bi-infinite sequence thus obtained.

\begin{lem}\label{lem:SiDisjoint}
Let $i<j$ be such that $\theta_j<\theta_i+\pi$. Then $S_i\cup L_{i+1}$ is disjoint from $S_j\cup L_{j+1}$.
\end{lem}
\begin{proof}
If $j=i+1$ this is easy. Otherwise, we proceed as follows: the segment $[B_k,A_k]$ is supported on $\partial \Delta_{\theta_k}$.
We have $B_k=A_{k-1}$ for all $k\in\Z$.
The vector $A_k-B_k = A_k-A_{k-1}$ has argument $\theta_k$
The vectors $A_k-B_k$ for $k$ between $i$ and $j-1$ included, have an argument in $(\theta_j-\pi,\theta_j]$.
This implies that the distance from $A_k$ to $\Delta_{\theta_j}$ is strictly increasing as $k$ decreases from $j-1$ to $i-1$.
In particular, the closure of the set $S_i$, which is a closed sector based on $A_{i}$ and whose elements are linear combination with non-negative coefficients of the vectors $A_{i-1}-A_i$ and $A_i-A_{i+1}$ lies in the complement of $\ov{\Delta}_{\theta_j}$.
\end{proof}

\begin{figure}
\begin{center}
\begin{tikzpicture}
\node at (0,0) {\includegraphics[scale=0.6]{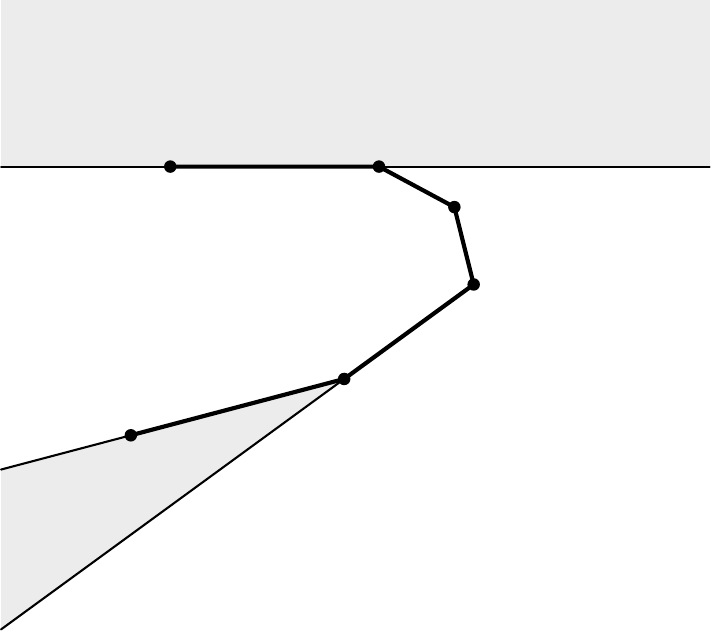}};
\node at (-2.5,-1.8) {$S_1$};
\node at (0,2.3) {$\Delta_5$};
\node at (-2.3,-0.9) {$A_0$};
\node at (-0.3,-0.3) {$A_1$};
\node at (1.55,0.1) {$A_2$};
\node at (1.4,1.1) {$A_3$};
\node at (0.25,1.1) {$A_4$};
\node at (-1.9,1.1) {$A_5$};
%  \draw[step=.1,overlay,lightgray,very thin] (-3,-3) grid (3,3);
%  \draw[step=.5,overlay,gray,thin] (-3,-3) grid (3,3);
%  \draw[step=1,overlay] (-3,-3) grid (3,3);
\end{tikzpicture}
\end{center}
\caption{Proof of \Cref{lem:SiDisjoint}}
\label{fig:SiDisjoint}
\end{figure}

Let us prove that $F$ is injective.
By contradiction assume there is some $i,i'\in\Z$, $z\in S_i\cup L_{i+1}$ and $z'\in S_{i'}\cup L_{i'+1}$ such that $f_{\theta_i}(z)=f_{\theta_{i'}}(z')$ but $z\neq z'$ or $i\neq i'$.
There exists $\phi\in \AutC$ such that $\phi(z')=z$ and such that
$\forall w$ close to $z$, $f_{\theta_{i}}\circ\phi(w) = f_{\theta_{i'}}(w)$.
The map $f_{\theta_i}\circ\phi$ is defined on the half-plane $H=\phi^{-1}(\Delta_{\theta_i})$ and is a rigid maximal Delaunay element, since its continuous extension maps at least two boundary point to a singularity.
There is a unique $\theta\in [\theta_{i'}-\pi,\theta_{i'}+\pi)$ such that $H$ is oriented like $\Delta_\theta$, i.e.\ is a translate thereof.
We cannot have $\theta=\theta_{i'}-\pi$ for otherwise, it would allow to extend $f_{\theta_{i'}}$ to the whole plane, and $(X,\nabla,\cal S)$ would be exceptional.
The two sets $H$ and $\Delta_\theta$ they are either equal or one contains strictly the other.
By affine continuation, the map $f_{\theta_i}\circ\phi$ coincides $f_{\theta_{i'}}$ with on the infinite sector $H\cap\Delta_{\theta_{i'}}$.
By definition, the map $f_\theta$ coincides $f_{\theta_{i'}}$ on the infinite sector $\Delta_{\theta}\cap\Delta_{\theta_{i'}}$.
It follows that $f_{\theta_i}\circ\phi$ and $f_\theta$ coincide on a non-empty open set, hence on the smallest of $H$ and $\Delta_\theta$.
By maximality of the immersions $f_{\theta_i}\circ\phi$ and $f_\theta$, this implies that $H=\Delta_{\theta}$ and
$f_{\theta_i}\circ\phi = f_\theta$, \ie $\phi$ is an equivalence between $(\Delta_\theta,f_\theta)$ and $(\Delta_{\theta_i},f_{\theta_i})$.
By hypothesis, this means $\theta=\theta_i$.
So $\phi(\Delta_{\theta_i}) = \Delta_{\theta_i}$ and $f_{\theta_i}\circ\phi=f_{\theta_i}$.
This implies that $\phi$ fixes $A_i$ and $B_i$ hence $\phi$ is the identity, so $z=z'$, hence $i\neq i'$ by hypothesis.
Moreover by one of the above claims, $|\theta_i-\theta'_i|<\pi$.
So the sets $S_i\cup L_{i+1}$ and $S_{i'}\cup L_{i'+1}$ are disjoint by \Cref{lem:SiDisjoint}, which contradicts $z=z'$. This proves the claim.

We then prove that $y$ is an irregular singularity.
Otherwise, it would be anti-conical Fuchsian.
If $y$ is pure anti-conical,
consider a pointed neighborhood $W$ of $y$ isomorphic to a standard neighborhood $V$ of a skew cylinder $\cet/\langle\lambda\rangle$ via an isomorphism $\psi:V\to W$.
The immersion $f_0$ maps $z\in\Delta_0$ to $W$ as soon as $|z|$ is big enough.
Fix a real $M>0$.
It follows by affine continuation that 
we can factor $F$ as $F = \psi\circ \pi\circ\Phi$ on a subset $U\subset \cal A$ with:
$U$ is the union of $(S_i\cup L_i)\setminus \ov B(0,R)$ for the indices $i$ satisfying $\theta_i \in [-M,M]$,
$R>0$,
$\Phi:U\to\cet$ is affine,
$\pi:\cet\to\cet/\langle\lambda\rangle$ is the canonical projection.
% there is $R>0$ and an affine map $\Phi:U\to\cet$ from the portion $U$ of $\cal A$ defined by the $S_i\cup L_i$ whose index $i$ satisfies $\theta_i \in [-M,M]$,
% such that $\forall\theta_i\in(-M,M)$ and $\forall z\in (S_i\cup L_i)\setminus \ov B(0,R)$, $F(\hat z) = f_{\theta_i}(z) = \psi\circ\pi\circ\Phi(\hat z)$ where $\hat z$ is the point of $U$ defined by $z$ and $i$ and $\pi:\cet\to\cet/\langle\lambda\rangle$.
If $M$ is chosen big enough, this contradicts injectivity of $F$.
In the shifted anti-conical case, the proof is similar.

We then claim that $F$ is proper from $\cal A$ to $\cal C$.
As in Case~1, we already know that $\hat z\in\cal A$ tends to infinity within one $S_i\cup L_i$ then $F(\hat z)$ tends to $y$, and that if $\hat z$ tends to a boundary segment of $\cal A$ for one of the finitely rigid $\theta_i$, then $F(\hat z)$ tends to a Delaunay segment.
Compared to Case~1 there is a supplementary possibility up to extraction for a sequence of points $\hat z\in\cal A$ leaving every compact subset: that $\theta_i$ tends to $\pm\infty$ where $i$ is the index such that $\hat z \in S_i\cup L_i$.
To treat that case we use that a neighborhood of $\infty$ in $\Delta_0$ maps to a neighborhood of type $B_{n+1/2}$ of the swath, in the nomenclature of \Cref{ss:model}.
This allows by analytic continuation to prove that for some $M>0$ big enough, for all $i$ such that $\theta_i>M$, $F(S_i\cup L_i)$ maps to the $A_n$ piece and that as $\theta_i\to+\infty$, $F(\hat z)$ tends to $y$.
A similar statement holds for $\theta_i$ tending to $-\infty$.
We also get that the first and last points on $\partial A$ map to adjacent foci of $y$ (the same focus if $y$ is a degree $2$ pole of $\nabla$), and that $F(\cal A)$ contains a neighborhood of a swath of $y$, and no neighborhood of any other swath of any irregular singularity.
Moreover we get that the closure of $F(\cal A)$ consists in $F(\cal A)$ union the image of all boundary $[A_i,B_i]$ of $\cal A$ by the extensions of the finitely many rigid $f_{\theta_i}$.

A proper injective holomorphic map over a connected set is a homeomorphism, so $\cal C$ is isomorphic to $\cal A$, so it is simply connected.
From the model of $\cal A$, it follows that the angles at the singularities are $\geq \pi$.
We have proved that $\cal C$ a swath domain, as per \Cref{defn:Unbounded}.

\medskip

In the classification we just proved, we have in every case a precise enough description so that we know an exterior Delaunay component $D$ that is not a finite angle Reeb cylinder contains the neighborhood of exactly one swath or non-conical singularity and is disjoint from small enough neighborhoods of all the others,
while a finite angle Reeb cylinder is disjoint from  small enough neighborhoods of all swathes and non-conical singularities.
(If the boundary of $D$ contains an irregular singularity $p$, disjointness of $D$ from small enough swath neighborhoods of $p$ (except for the swath associated to $D$ if $D$ is a swath domain) follows from the local model given in \Cref{ss:model} and the fact that the angle of the domain boundary at $p$ is finite).

\medskip

Conversely, near any anti-conical, Reeb or cylindrical Fuchsian singularity, or irregular singularity swath, we can immerse a half-plane.
The Delaunay component $D$ that contains the image of this immersion is necessarily an exterior domain.
Every neighborhood of the singularity or swath we started from intersects $D$.
By the previous paragraph, $D$ is not a finite angle Reeb cylinder, and thus is associated to a singularity or swath which, again by the previous paragraph, has to be the one we started from.

\subsection{A complexity bound}\label{sub:ComplexityBound}

The number of components in the Delaunay decomposition of a finite type affine surface $(X,\nabla,\mathcal{S})$ is controlled by topological data: the number $n$ of poles of $\mathcal{S}$ (counted with multiplicity) and the genus $g$ of the underlying Riemann surface $X$.
We assume that $2g+n \geq 3$, so that there are Delaunay segments.
\Cref{thm:ComplexityBound} is actually finer than just a bound on the number of components.
It gives an equality, which to be expressed requires to decompose Delaunay polygons without adding vertices, which is always possible.
A Delaunay polygon with $s$ sides will give $s-2$ triangles.
Some of the triangles may be self-folded but we allow this.
We call $t$ the total number of triangles.
Concerning exterior domains, we count the total number of sides on all these domains and call it $\beta$.
Then \Cref{thm:ComplexityBound} states that
\[ t+\beta = 4g-4+2n
.\]

\begin{proof}[Proof of \Cref{thm:ComplexityBound}]

For a non-exceptional finite type affine surface $(X,\nabla,\mathcal{S})$, we draw the union $\mathcal{A}$ of the Delaunay segments.
The set $X \setminus \mathcal{A}$ decomposes according to \Cref{thm:ClassificationComponents}.
We draw additional arcs to decompose further Delaunay polygons into topological triangles $t$.
We obtain a system of arcs $\mathcal{A}'$ (containing $\mathcal{A}$) satisfying the hypotheses of \Cref{lem:SystemTopologicalArcs}.
The special faces of the decomposition are exactly the semi-infinite Reeb cylinders, semi-infinite translation cylinders, anti-conical domains and swath domains, together with their singularity (except in the swath domain).
The ordinary faces of the decomposition form a topological triangulation of the union of the Delaunay polygons and Reeb cylinders of finite angle.
From \Cref{lem:SystemTopologicalArcs} we obtain:
\begin{align*}
|\mathcal{A}| &= 2g - 2 + n + o \\
\text{and} \quad |\mathcal{A}| &= 6g - 6 + 3n - \sigma
\end{align*}
where $o$ is the number of ordinary faces and $\sigma$ the total side count of the special faces.
A finite angle embedded Reeb cylinder with side-count $s$ is triangulated by exactly $s$ ordinary faces, whence $o+\sigma = t+\beta$.
Combining these identities, we obtain $t+\beta=4g-4+2n$.
\end{proof}

\section{Existence of meromorphic connections with prescribed irregular singularities}\label{sec:ExistenceProblem}

\subsection{A Riemann-Hilbert problem}

Meromorphic abelian differentials,\footnote{Which take the form $f(z)dz$ in charts, $f$ meromorphic.} are equivalent to translation surfaces with polar domains.
A complete system of local invariants are given by the orders of the zeroes and the order and residues at the poles.
The complete characterization of the profiles of local invariants that can be realized by such a differential on a Riemann surface of genus $g$ was given in \cite{getaab}.
Analogous results were proven for quadratic differentials\footnote{Which take the form $f(z)dz^2$ in charts.} and higher order differentials\footnote{Which take the form $f(z)dz^q$ in charts for some integer $q\geq 3$.} in \cite{getaquad} and \cite{getakdiff}, respectively.
\par
All these differentials define affine structures corresponding to special subclasses of Fuchsian meromorphic connections\footnote{With Christoffel symbol $\Gamma(z)=f'(z)/qf(z)$ in charts.} on the tangent bundle of a Riemann surface.
In each case, the list of obstructions is remarkably rich, and the algorithms for constructing flat surfaces with prescribed local invariants are highly intricate.
\par
The global realization problem for a profile of local invariants on affine surfaces is a Riemann–Hilbert-type problem, one that has never been solved at the above of generality, let-alone in the presence of irregular singularities:

\begin{pbm}
Given a profile of local analytic invariants for poles, does there exist a Riemann surface $X$ of genus $g$ and a meromorphic connection $\nabla$ on the tangent bundle $TX$ of $X$ whose poles realize this profile of local invariants?
\end{pbm}

In this paper, we address a significantly simpler problem.

\begin{pbm}
Given a double pole with prescribed local invariants, what is the minimal number $m$ such that there exists a meromorphic connection $\nabla$ on the tangent bundle of $\mathbb{CP}^{1}$ such that $\nabla$ has a double pole realizing this profile of local invariants, $m$ Fuchsian singularities and no other singularity?
\end{pbm}

Let us introduce the following terminology. We recall that the asymptotic value family of an irregular singularity is a $\Z$-indexed sequence $u_n$ of complex numbers defined in \Cref{prop:di}.

\begin{defn}\label{def:centered}
An irregular singularity is \emph{centered} when the asymptotic value family of a developing map (hence of every developing map) is constant.
\end{defn}
According to Section~\ref{subsub:doublepoles}, for a double pole of and a given residue $\res \in \mathbb{C}$, there are exactly two possible local affine isomorphism classes, and they can be distinguished by the pole being centered or not.
The choice of the invariant $\iota\in\cal I_{2,\res}$, which characterizes the class for a fixed residue, amounts to the choice between centered or not.

\begin{rem}\label{rem:ResInZCen}
When $\res \in \mathbb{Z}$, a double pole is centered if and only if $b=0$, where the affine holonomy of a developing map of the pole is $L(z)=z+b$. See \Cref{def:is}.
\end{rem}

\begin{thm}\label{thm:CharacterizationDoublePole}
Given a residue $\res \in \mathbb{C}$ and an asymptotic values invariant $\iota \in \mathcal{I}_{2,\res}$, the minimal number $s \in \mathbb{N}$ such that there exists a meromorphic connection on the tangent bundle of $\mathbb{CP}^{1}$ with $s$ simple poles and a double pole of residue $\res$ and invariant $\iota$ is:
\begin{itemize}
    \item $0$ if $\res=2$ and $\iota$ is centered;
    \item $1$ if $\res \in \mathbb{Z}\cap[3,+\infty)$ and $\iota$ is centered;
    \item $1$ if $\res \in \C\setminus\{2,3,\ldots\}$ and $\iota$ is not centered;
    \item $2$ if $\res \in \mathbb{Z}\cap[2,+\infty)$ and $\iota$ is not centered; 
    \item $2$ if $\res \in \C\setminus\{2,3,\ldots\}$ and $\iota$ is  centered;
\end{itemize}
\end{thm}

We can present this list differently:

%\SetTblrInner{hspan=even}
\begin{center}
\begin{tblr}{
  vlines,hlines
 ,hline{3}={1}{lr}
 ,hline{4}={1}{lr}
 ,vline{2,3,4}={1}{abovepos=-1,belowpos=-1}
 ,colspec={cX[c]X[c]X[c]}
 ,width=0.8\textwidth
}
$\res$: & $\C\setminus\{2,3,\ldots\}$ & $2$ & $3,4,\ldots$
\\ 
centered & \SetCell[c=3]{}
\\ 
yes & 2 & 0 & 1
\\ 
no & 1 & 2 & 2
\\ 
\end{tblr}
\end{center}

\medskip

In the remainder of \Cref{sec:ExistenceProblem}, we prove \Cref{thm:CharacterizationDoublePole} by providing the obstructions in \Cref{sub:Obstruction} and explicitly constructing the connections in \Cref{sub:DoublePoleConstruction}. Note that the case $s=0$ is already completely characterized as the exponential-affine plane $\mathcal{E}$ by \Cref{lem:ExpoAffineUniqueness}.

The analogous problem for poles of order $d \geq 3$ is significantly more complicated due to the wide variety of possible asymptotic sequences.

\subsection{Obstructions}\label{sub:Obstruction}

Here we identify and prove obstructions to the existence of a meromorphic connection on $\RS$ with one double pole and a small number of simple poles (Fuchsian singularities).
Recall (see \Cref{prop:sumRes}) that the sum of residues of a meromorphic connection on $\RS$ is equal to $2$.

\begin{prop}\label{prop:DoublePoleObstruction}
We assume that there exists a meromorphic connection $\nabla$ on $\mathbb{CP}^{1}$ with exactly two poles, one being a simple pole while the other is a double pole of residue $\res \in \mathbb{C} \setminus \lbrace{ 2 \rbrace}$ and invariant $\iota \in \cal I_{\res,2}$. Then the invariant $\iota$ is centered if and only if $\res \in \mathbb{Z}\cap[3,+\infty)$.
\end{prop}

\begin{figure}[htbp]
  \begin{center}
    \begin{tikzpicture}
% \draw[step=.1,overlay,lightgray,very thin] (-3,-2) grid (3,2);
% \draw[step=.5,overlay,gray,thin] (-3,-2) grid (3,2);
% \draw[step=1,overlay] (-3,-2) grid (3,2);
      \node at (0,0) {\includegraphics[scale=1.33]{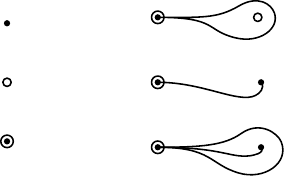}};
      \node at (-1.6,1.5) {\begin{minipage}{2cm}\small conical\\ Fuchsian\end{minipage}};
      \node at (-1.6,0.1) {\begin{minipage}{2cm}\small non-conical\\ Fuchsian\end{minipage}};
      \node at (-1.6,-1.2) {\begin{minipage}{2cm}\small irregular\\ degree $2$\end{minipage}};
      \node at (1.25,2.1) {A};
      \node at (2.2,1.5) {B};
      \node at (1.25,0.5) {A};
      \node at (1.25,-0.9) {A};
      \node at (2.35,-1.25) {C};
    \end{tikzpicture}
  \end{center}
  \caption{The three kind of Delaunay decompositions we can meet for a meromorphic connection as in \Cref{prop:DoublePoleObstruction}.
  The Delaunay components (connected components of the complement of the graph) are labeled A, B, C, and this is related to their classification (i), (ii), \ldots (v) as in \Cref{thm:ClassificationComponents} as follows: A for infinite angle unbounded type (v); B for either semi-infinite Reeb (ii) or translation (iii) cylinder, or finite angle unbounded type (iv); C for open Delaunay polygon (i).}
  \label{fig:dd1}
\end{figure}
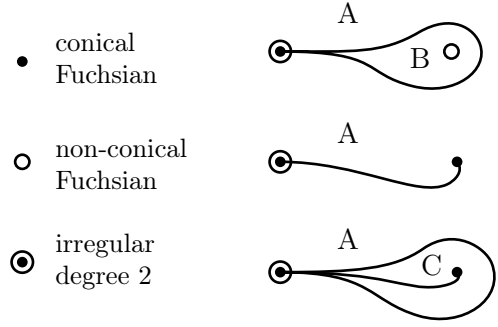

\begin{proof}[Geometric proof of \Cref{prop:DoublePoleObstruction}]
We assume that such a connection exists. 
Call $p$ the double pole and $q$ the simple pole.
The surface is non-exceptional since $(g,n)=(0,3)$ in the notations of \Cref{thm:ExceptionalSurfaces}.
Call
\[\rho = 2- \res\]
the residue of the simple pole $p$.
The complement $U$ of the two poles in $\RS$ has fundamental group $\pi_1(U)$ isomorphic to $\Z$.
Let $\gamma$ be a loop generating it.
Let $\pi:\tilde U \to U$ be a universal covering.
Its deck transformation group is isomorphic to $\pi_1(U)$, hence to $\Z$.
Let $\tau$ be the generator associated to $\gamma$.
Consider a developing map $\phi:\tilde U\to\C$.
Then according to \Cref{sub:Turning}, $\phi\circ\tau = L\circ\phi$ for some $L\in\AutC$ (called the holonomy) with $L(z)=az+b$ and $a=\exp(-2\pi i\res)=\exp(2\pi i\rho)$ or $a= \exp(2\pi i\res)=\exp(-2\pi i\rho)$ depending on the choice of generator $\gamma$.
\par
If the simple pole $q$ is conical, then it has a punctured neighborhood $V$ affine isomorphic to the following local model: a connected neighborhood of $0$ in a closed sector in $\C$ (or in $\cet$ if the angle of the sector needs to be $>2\pi$), based on $0$, with the two boundary lines identified by an affine map fixing $0$.
The developing map sends any lift of $V$ in $\tilde U$ to an affine image of the sector. Choose one lift and call $z_0$ the tip of the sector.
Then $L$ sends one of the half-lines bounding the sector to the other one, in particular it fixes $z_0$.
\par
If $\res \in \mathbb{Z}\cap[3,+\infty)$, then $\rho \in\Z\cap (-\infty,-1]$.
The simple pole is therefore a conical singularity of angle $(1-\rho)2\pi$.
The coefficient $a$ of the affine holonomy $L$ is equal to $1$.
Moreover, we get from the local model of a conical singularity that the affine holonomy of a small loop around it always has a fixed point.
So $L$ is the identity.
But $L$ is also the affine holonomy of a small loop winding once around the double pole $q$ (in the other direction).
So $q$ is centered by \Cref{rem:ResInZCen}.

\medskip

Below we use the Delaunay decomposition of $RS$ associated to the affine structure. The possible cases that we will meet are illustrated on \Cref{fig:dd1}.

\medskip

If $\Re(\res)>1$ but $\res \notin \mathbb{Z}$, then we have $\Re(\rho)<1$ so the simple pole is a conical singularity. Since $\Re(\rho) \notin \mathbb{Z}$, the coefficient $a$ of the affine holonomy $L$ is not equal to $1$.
\par
In the Delaunay decomposition induced by the connection (see Theorem~\ref{thm:ClassificationComponents}), there is exactly one swath domain $D$ (containing the swath of the double pole). Following Theorem~\ref{thm:ComplexityBound}, we have $t+\beta=2$ where $t$ is the number of triangles in any  triangulation of the rest of the surface with vertices in $\cal S$ while $\beta$ is the number of boundary sides of the exterior domains (including $D$).
There are two cases.
\par If $\beta=1$, then $D$ is the only exterior domain and the double pole is non-centered because, following a developing map along the unique Delaunay segment bounding $D$ gives a segment, whose two endpoints are different, which shows that the double pole $q$ has at least two different asymptotic values (the path defined by the Delaunay segment can be shrunk within the closure of $D$ in $X^*$ to an arbitrarily small neighborhood of the double pole, while still preserving a initial and final portion of the path). So $q$ is not centered.
We are in the bottom right situation of \Cref{fig:dd1}.
\par
If $\beta=2$, then $D$ is still the only exterior domain because the only Fuchsian singularity is conical and a Reeb cylinder of finite angle would have at least two boundary sides. It follows that $t=0$ and the swath domain has two boundary sides and we obtain the affine surface by taking a model $\cal A$ as in Case~2 of the proof of \Cref{thm:ClassificationComponents} in \Cref{sub:DelaunayComponents}, gluing on each other the two boundary segments of $\cal A$, corresponding to the two sides of the Delaunays segment, with one end being the conical singularity, while the other end is the focus of the the irregular singularity.
Following the developing map $\phi$ from the simple to the double pole along this segment, we trace a straight segment in $\C$ from the fixed point $z_0$ of $L$ to a term $u_n$ of the asymptotic value family of the double pole $q$.
Hence $L(u_n)\neq u_n$.
But $L(u_n)=u_{n+1}$ (or $u_{n-1}$ according to the choices of generators).
It follows that $q$ is non-centered.
We are in the middle right situation of \Cref{fig:dd1}.

\medskip

If $\Re(\res) \leq 1$, then $\Re(\rho) \geq 1$, so the Fuchsian singularity $p$ is cylindrical, of Reeb-type or anti-conical, so no Delaunay segment can reach $p$, and $p$ belong contained to a component $D_1$ of the Delaunay decomposition that is either a semi-infinite translation cylinder, a Reeb cylinder an anti-conical domain or a swath domain (together with their singularities in the non-swath case).
Following Theorem~\ref{thm:ClassificationComponents}, the Delaunay decomposition of the affine structure defined by the connection must also contain a swath domain $D_2$ containing the swath of the double pole and different from $D_1$.
So the total number $\beta$ of boundary sides of these domains satisfies $\beta \geq 2$. Theorem~\ref{thm:ComplexityBound} tells that $t+\beta = 2$ where $t$ is some non-negative integer related to the other domains $D_i$ of the decomposition, $i\geq 3$, if there are.
It follows that $\beta=2$ and $t=0$, so there are in fact no domains beyond $D_1$ and $D_2$.
The surface is formed by these two domains glued along the two sides of a single closed saddle connection from the double pole to itself.
We are in the top right situation of \Cref{fig:dd1}.
%By a homographic change of variable, we may assume that the two poles in $\RS$ are $0$ and $\infty$.
%Consider now a developing map $\phi$, defined on the universal covering of $\C^*$.
Then the value of the developing map $\phi$ along any lift of the Delaunay segment to $\tilde U$ (the universal covering of the complement in $\RS$ of the two poles) follow a straight segment whose ends are two consecutive terms of the asymptotic value family associated to the double pole $q$.
Consequently, $q$ is not centered.

Coming back to the case $\res \in \mathbb{Z}\cap[3,+\infty)$, which we saw is centered, the same argument as for the case $\Re(\res)>1$ implies that the domain $D$ is 2-sided so the Delaunay graph is a single simple curve between the two poles: we are in the middle right situation of \Cref{fig:dd1}.
\end{proof}

\subsection{Sector grafting}\label{ss:grafting}

\begin{figure}[htbp]
  \begin{center}
    \begin{tikzpicture}
      \node at (0,0) {\includegraphics[scale=1.0]{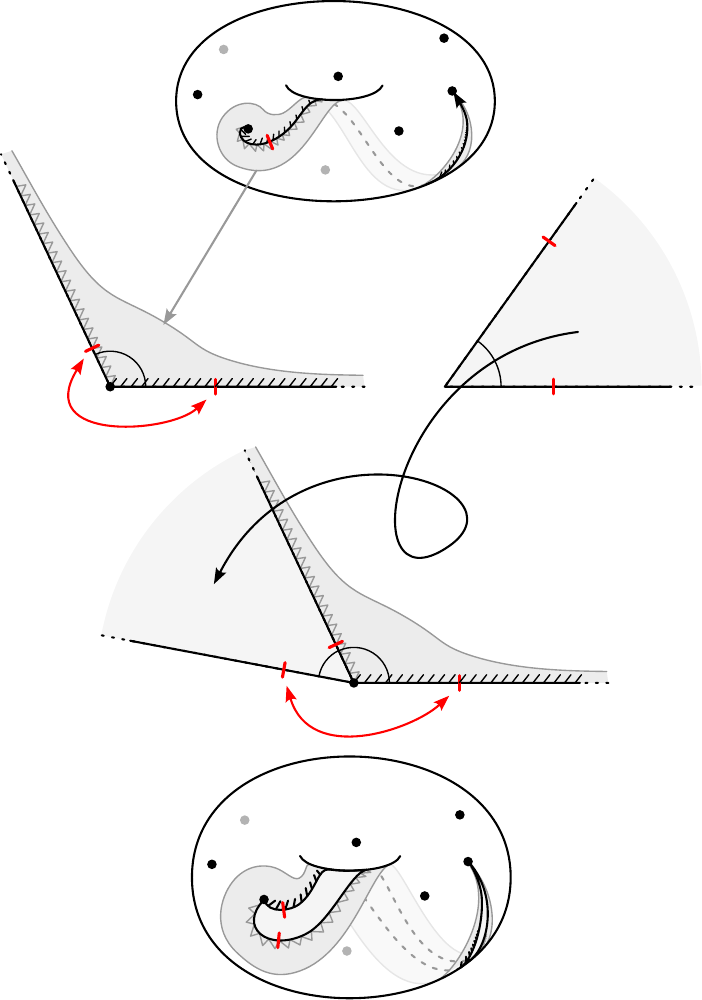}};
      \node at (2.9,6.8) {$X$};
      \node at (2.05,7.15) {$S_2$};
      \node(B) at (-1.4,7.8) {$S_1$};
      \draw[-] (B) -- (-1.71,6.3);
      \node(A) at (-3.5,5.9) {$\gamma$};
      \draw[-] (A) -- (-1.7,6.1);
      \node at (4.5,3.5) {$C_\theta$};
      \node at (3.2,-6.4) {$X'$};
    \end{tikzpicture}    
  \end{center}
  \caption{Grafting on a torus (\Cref{ss:grafting})}
  \label{fig:grafting}
\end{figure}

We introduce a surgery that allows us to modify the residues of the singularities.
\par
Assume we are given, on a finite type affine surface $(X,\nabla,\mathcal{S})$, a geodesic $\gamma: (0,+\infty) \to X^{\ast}$ such that:
\begin{itemize}
    \item $\gamma(t)$ converges to a conical singularity $S_{1}$ of residue $\res_{1}$ as $t \tends 0$;
    \item $\gamma(t)$ converges to a singularity $S_{2}$ of residue $\res_{2}$ as $t \tends +\infty$;
    \item $\gamma$ is injective (not self-crossing).
\end{itemize}
We make no asumption on the order of the pole at $S_2$ but note that $S_2$ cannot be a conical Fuchsian singularity, for it only would only take finite time for a geodesic to reach it, nor a Reeb type singularity, which is unreachable.
\par
We cut $X$ along the line traced by $\gamma$: this gives a surface with boundary, the slit having a left side and a right side with respect to the orientation of $\gamma$.
We then glue an infinite sector $C_{\theta}$ of angle $\theta\in(0,+\infty)$, $C_\theta=\{z\in\wt\C^*\,;\,\arg z\in[0,\theta]\}$, as follows: a point in $C_{\theta}$ of polar coordinate $(r,0)$ is glued to the right boundary point of the slit at $\gamma(r)$; the point of $C_{\theta}$ of polar coordinate $(s r,\theta)$ is glued to the left boundary point of the slit at $\gamma(r)$, where $s  \in \mathbb{R}_{>0}$ is some dilation factor. See \Cref{fig:grafting}.
\par
The obtained topological space is again a surface, and actually a finite type affine surface $X'$.
The conical point remains a conical point and $S_2$ remains a singularity, of the same polar order.
In $X'$, the residues at the singularities $S_{1}$ and $S_{2}$ become respectively
\[
\begin{split}
\res'_1 & = \res_{1}-\frac{\log(s )+i\theta}{2\pi i} = \res_{1}+i\frac{\log(s )}{2\pi}-\frac{\theta}{2\pi} \text{ and }
\\
\res'_2 & = \res_{2}+\frac{\log(s )+i\theta}{2\pi i} = \res_{2}-i\frac{\log(s )}{2\pi}+\frac{\theta}{2\pi}.
\end{split}\]
This operation does not affect the other singularities.

\begin{figure}[htbp]
  \begin{center}
    \begin{tikzpicture}
      \node at (0,0) {\includegraphics[scale=1.0]{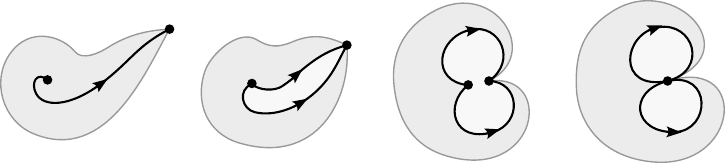}};
    \end{tikzpicture}    
  \end{center}
  \caption{Grafting an infinite angle seen as a limit process.}
  \label{fig:grafting:infini}
\end{figure}

We can also graft an infinite angle, by grafting two semi exponential half-planes to each side of the slit $\gamma$.
In other words, cut $\cet$ into two pieces along the radial line of argument $0$.
Glue the top piece to the left of the slit and the bottom piece to the right of the slit.
This time there is no dilation parameter, and the two singularities merge into a single singularity $S'$, whose order is $d+1$ where $d$ is the order of $S_2$.

\begin{lem}
The residue of $S'$ is $\res_1+\res_2$.
\end{lem}
\begin{proof}
Take a simple closed loop near $\gamma$ winding anticlockwise around it: its holonomy is $\exp(-2\pi i(\res_1+\res_2))$ and it does not change after the grafting, then using that the linear holonomy only depends on free homotopy class, we get
$\exp(-2\pi \res') = \exp(-2\pi i(\res_1+\res_2))$
so $\res' \equiv \res_1+\res_2) \bmod \Z$.
To get rid of the modulo, one looks at the turning number of the loop, as defined in \Cref{sub:Turning} and the deformation invariance stated in \Cref{lem:tninv} to reduce to a (non-injective) $C^1$ curve that loops around $S_1$, then rides along $\gamma$ then loops around $S_2$, then rides back; the small loops have classical turning number $0$ in a chart, not $1$, so \cref{eq:tn1} gives that their turning number for the affine structure is $-\Re \res$ (not $1-\Re\res$ as in \cref{eq:tn2}).

Alternatively one could try and prove convergence of the finite angle grafting (with no dilation) to the infinite angle grafting and use the computation of residues we made in the finite angle case.
See \Cref{fig:grafting:infini}.
\end{proof}

Informally, the conical singularity becomes a focus, while for the other singularity: if it is Fuchsian, it becomes a swath; if it is irregular, the swath to which $\gamma$ converges before the infinite angle grafting, is cut into two swathes.

\subsection{Construction of connections with prescribed irregular singularity class}\label{sub:DoublePoleConstruction}

\subsubsection{Constructions with one additional Fuchsian singularity}\label{ss:ConstructionAdditionalOne}

\begin{figure}[htbp]
  \begin{center}
    \begin{tikzpicture}
      \node at (0,0) {\includegraphics[scale=0.5]{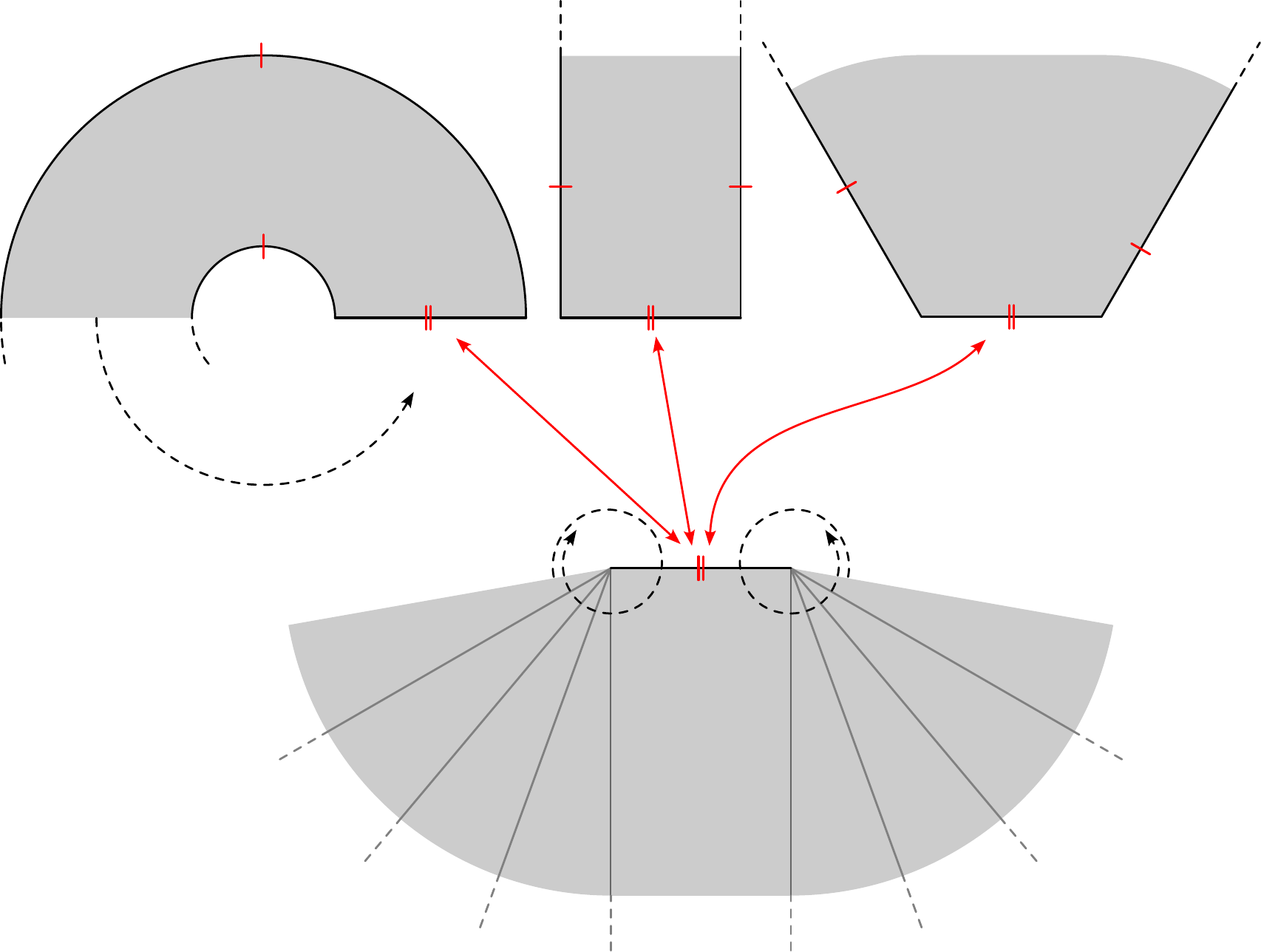}};
      \node at (-.25,0.7) {or};
      \node at (1,0.7) {or};
%      \node at (2.5,0.2) {or};
    \end{tikzpicture}    
  \end{center}
  \caption{First and second construction of \Cref{ss:ConstructionAdditionalOne}}
  \label{fig:constr1}
\end{figure}

% \begin{figure}[htbp]
%   \begin{center}
%     \begin{tikzpicture}
%       \node at (0,0) {\includegraphics[scale=0.66]{23.pdf}};
%     \end{tikzpicture}    
%   \end{center}
%   \caption{The third construction of \Cref{ss:ConstructionAdditionalOne} consists in applying the grafting of \Cref{ss:grafting} to the second construction. The sector is inserted along the line indicated by scissors.}
%   \label{fig:constr2}
% \end{figure}

\begin{figure}[htbp]
  \begin{center}
    \begin{tikzpicture}
      \node at (0,0) {\includegraphics[scale=0.66]{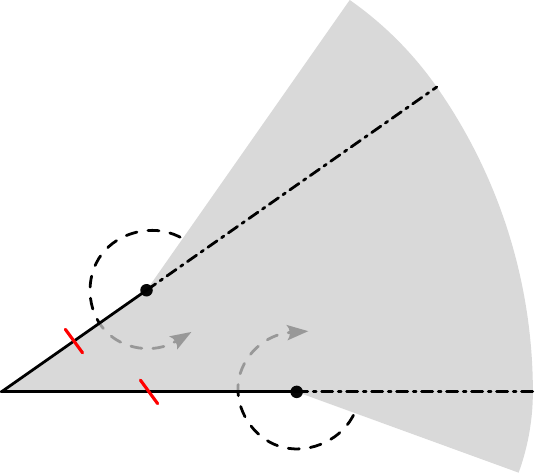}};
    \end{tikzpicture}    
  \end{center}
  \caption{One way to present the second  construction of \Cref{ss:ConstructionAdditionalOne}.}
  \label{fig:constr2}
\end{figure}

The first construction consists in gluing together, along a closed segment that will end up being a saddle connection, two domains:
\begin{itemize}
    \item a swath domain\footnote{Using the terminology of \Cref{sub:DelaunayComponents}. It is not necessary to have read that section to understand the construction here.} that will contain a swath of the double pole. It is constructed by taking the universal covering $\cet$ of $\C^*$, slitting it along a half-line from $0$, and inserting a semi-infinite strip. See the bottom of \Cref{fig:constr1}.
    \item An anti-conical domain \textit{or} a semi-infinite translation cylinder \textit{or} a semi-infinite Reeb cylinder (an exterior domain containing the simple pole). 
    See the top part of \Cref{fig:constr1}.
\end{itemize}
Both ends of the segment tend to the focus the double pole. In this construction, the residue of the simple pole can be freely chosen in $\lbrace \rho \in \mathbb{C} \,;\, \Re(\rho) \geq 1 \rbrace$. Then the residue $\res$ at the double pole can be any complex number satisfying $\Re(\res) \leq 1$ since $\res+\rho=2$ by \Cref{prop:sumRes}. Note that in this construction, the double pole is never centered.

\medskip

In the second construction, for any $\rho$ with $\Re\rho<1$, start from a skew cylinder where the conical singularity has residue. Take a half-line embedded in the regular part, for instance a radial line minus an initial portion from the conical singuarity. One end of this geodesic is a regular point, the other is the anti-conical point at infinity.
Perform an infinite angle sector grafting as in \Cref{ss:grafting}.
This merges the regular and anti-conical points into a double pole, without changing the conical singularity, and the resulting compact Riemann surface is still homeomorphic to a sphere.

% The second construction is the same where for the second domain we use a self-folded triangle, where two sides are identified by an affine map.
% The residue $\rho$ of the simple pole, corresponding to a conical singularity whose angle lies in $(0,\pi)$, can be freely chosen in $\lbrace{ z \in \mathbb{C} \,;\, \frac{1}{2}<\Re(z) < 1 \rbrace}$. Then the residue $\res$ at the double pole can be any complex number satisfying $1<\Re(\res) < \frac{3}{2}$ since $\res+\rho=2$.
% Again, since the segment on which the gluing is made is a saddle self-connection of the double pole, the latter is non-centered.

% \vspace{0.5cm}

% In the third and final construction, we start with a surface produced by the second construction, choose an infinite geodesic $\gamma$ starting from the simple pole (which is a conical singularity) and passes through an interior point of the saddle connection: it converges in infinite time toward the unbounded sector of the double pole, and has no self-intersections. Grafting along $\gamma$ a cone of arbitrary angle and arbitrary dilation ratio, we realize every residue $\res$ satisfying $\Re(\res) \geq \frac{3}{2}$.

\medskip

For each of these two constructions, we can either use \Cref{prop:DoublePoleObstruction} to determine whether the double pole is centered or not, depending on the value of $\res$, or see it directly on the geometric model.

\subsubsection{Constructions with a non-centered double pole and two additional Fuchsian singularities}\label{ss:constr3}

\begin{figure}[htbp]
  \begin{center}
    \begin{tikzpicture}
      \node at (0,0) {\includegraphics[scale=0.66]{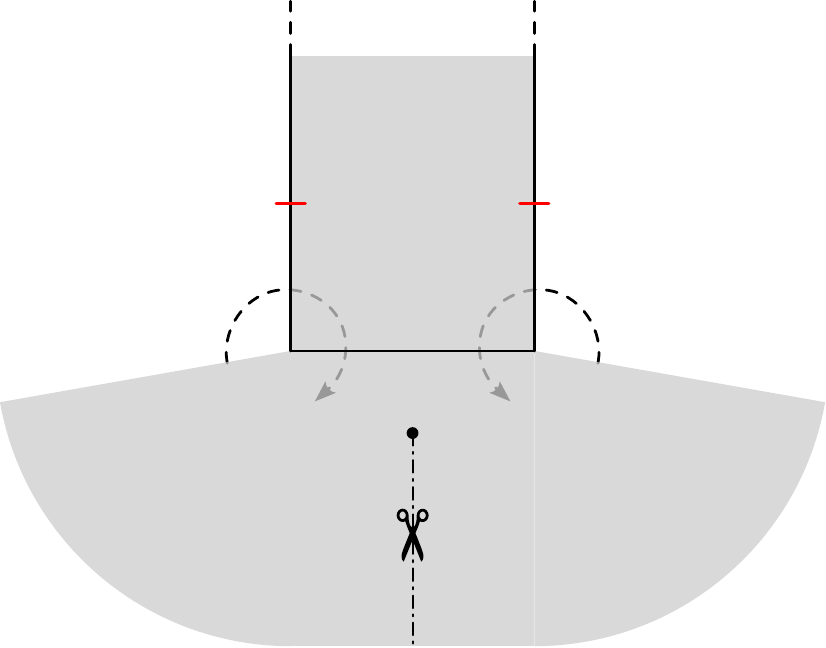}};
    \end{tikzpicture}    
  \end{center}
  \caption{The construction of \Cref{ss:constr3} consists in applying the sector grafting of \Cref{ss:grafting} to the first construction of \Cref{ss:ConstructionAdditionalOne} (with a translation cylinder for the second domain). The sector is inserted along the line indicated by scissors and has an angle that is a positive integer number of turns.}
  \label{fig:constr3}
\end{figure}

We have to construct a meromorphic connection on $\mathbb{CP}^{1}$ with a non-centered double pole of residue $\res \in \mathbb{Z}\cap[2,+\infty)$.
We start from the first construction of \Cref{ss:ConstructionAdditionalOne} where for the second domain we chose a semi-infinite translation cylinder.
We then take any point in the first domain draw a half-line to infinity within that domain (see an example on \Cref{fig:constr3}), on which we perform a grafting (see \Cref{ss:grafting}) of angle $2\pi (\res-1)$ and no dilation.
We obtain a meromorphic connection on $\mathbb{CP}^{1}$ with a Fuchsian singularity of residue $1$ at the end of the translation cylinder and a Fuchsian singularity of residue $1-\res$, being a conical singularity of angle $2\pi \res$. Then, the double pole has a saddle self-connection so is non-centered and its residue is equal to $\res$ since the sum of the three residues has to be $2$.

\subsubsection{Constructions with a centered double pole and two additional Fuchsian singularities}\label{ss:constr4}

\begin{figure}[htbp]
  \begin{center}
    \begin{tikzpicture}
      \node at (0,0) {\includegraphics[scale=0.66]{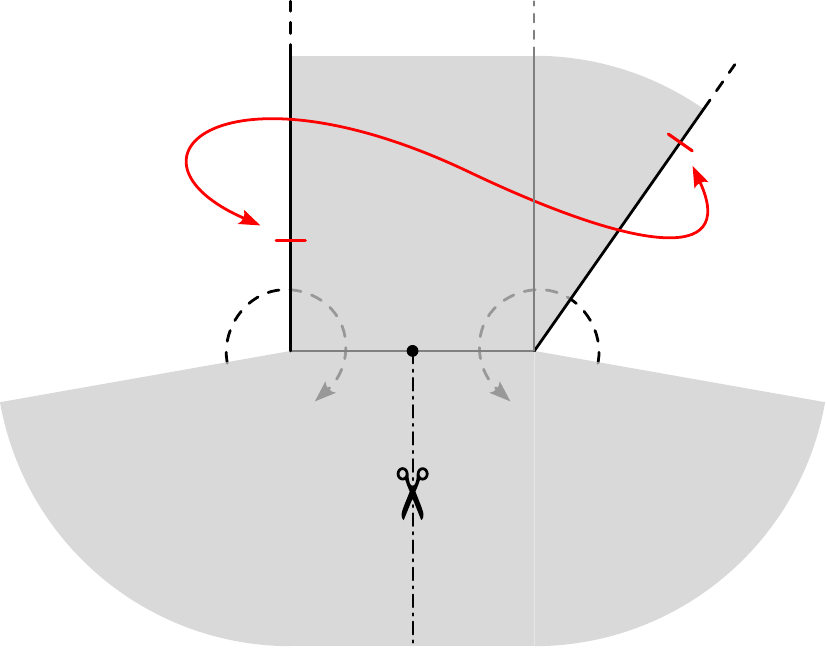}};
      \node at (2,2) {$2\pi a$};
    \end{tikzpicture}    
  \end{center}
  \caption{The construction of \Cref{ss:constr4} consists in applying two sector graftings to the first construction of \Cref{ss:ConstructionAdditionalOne} (with a translation cylinder for the second domain). The sectors are inserted along: the line indicated by scissors; the vertical line along which a strip was initially glued to give a translation cylinder.}
  \label{fig:constr4}
\end{figure}

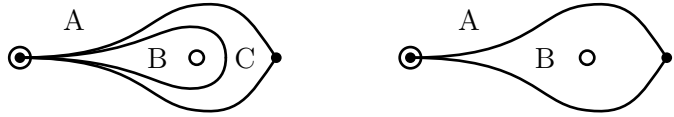
\begin{figure}[htbp]
  \begin{center}
    \begin{tikzpicture}
      \node at (0,0) {\includegraphics[scale=1.33]{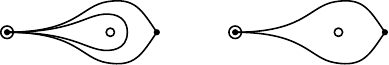}};
      \node at (-3.5,0.5) {A};
      \node at (-2.4,0) {B};
      \node at (-1.25,0) {C};
      \node at (1.7,0.5) {A};
      \node at (2.7,0) {B};
    \end{tikzpicture}    
  \end{center}
  \caption{Delaunay decompositions corresponding to the constructions of \Cref{ss:constr3,ss:constr4}. Same conventions as on \Cref{fig:dd1}.}
  \label{fig:dd2}
\end{figure}

We have to construct a meromorphic connection on $\mathbb{CP}^{1}$ with a centered double pole of residue $\res \notin \mathbb{Z}\cap [2,+\infty)$.
We also start here from the first construction of \Cref{ss:ConstructionAdditionalOne} where for the second domain we chose a semi-infinite translation cylinder.
We then take the mid-point of the saddle connection between the two domains and draw a vertical half-line to infinity within the first domain (see on \Cref{fig:constr4}), on which we perform a grafting of angle $2\pi k-\pi$ and no dilation, $k\in\N^*$.
We also perform a grafting on the vertical line bounding the upper half-strip (initially glued to the other vertical side to give a translation cylinder), inserting an angle of $2\pi a$, $a\in [0,+\infty)$ and this time allowing for a dilation of factor $s >0$.
We obtain a meromorphic connection on $\mathbb{CP}^{1}$ with:
\begin{itemize}
    \item at infinity in the former cylinder, a Fuchsian singularity of residue $\rho_1$ with $\rho_1 = 1+a +i \frac{ \log s }{2\pi}$;
    \item at the midpoint, a Fuchsian singularity of residue $1/2-k$, since it is conical with no dilation and angle  $2\pi k+\pi$.
\end{itemize}
By construction, the double pole is centered.
Its residue is equal to 
\[\res := k+\frac{1}{2}-a-i \frac{ \log s }{2\pi}\]
since the three residues add up to $2$.
Since $k$ and $a$ can be chosen arbitrarily in $\mathbb{N}^{\ast}$ and $[0,+\infty)$ respectively, $\Re(\res)$ can take any prescribed value.
Finally $\Im(\rho)$ can take any prescribed value since $s $ can be freely chosen in $(0,+\infty)$

\subsection{Analytic approach and Euler's Gamma function}\label{subsub:PrescribedAnalyticApproach}

Here we reprove part of \Cref{thm:CharacterizationDoublePole} with another approach. More precisely we get the minimal value $s$ on the number of supplementary Fuchsian singularities when this value is $0$ or $1$ but instead of $s=2$ we only get $s\geq 2$. Interestingly, part of the argument is based on Euler's Gamma function.

\begin{defn}
Let $\cal C$ be the class of finite type meromorphic connections on $\RS$ with two singularities: a double pole and either a simple pole or a marked point. 
\end{defn}
Up to a homographic change of variable we can assume that the double pole is at $0$ and that the other singularity is $\infty$.

We denote here $\zeta$, instead of $\Gamma$, the Christoffel symbol of the connection in the canonical chart $\C$ of $\RS$.
It is a holomorphic function on $\C^*$.
It thus has a Laurent power series expansion
\[\zeta(z) = \sum_{n\in\Z} a_n z^n\]
with domain of convergence containing $\C^*$.
By hypothesis, it has at $0$ a double pole, so the sum may be restricted to $n\geq -2$.
Under the change of variable $w=1/z$, the Christoffel symbol transforms into
\[\zeta_w = - \frac{1}{w^2}\,\zeta(1/w)-\frac{2}{w} = -\frac{2}{w} -\sum_{n\geq -2} a_n w^{-n-2}\]
and is supposed to have at $0$ a regular point or a simple pole, so actually $a_n=0$ for all $n\geq 0$ and:
\[\zeta(z) = \frac{a_{-2}}{z^2}+\frac{a_{-1}}{z}\]
with $a_{-2}\neq 0$ and $a_{-1}$ is the opposite of the double pole residue $\res$ (see \Cref{sub:minicours}).
By a linear change of variable on $z$ we may assume that $a_{-2}=1$:
\[\zeta(z) = \frac{1}{z^2}-\frac{\res}{z}.\]

We immediately deduce that:
\begin{lem*}
  For each value of the residue, there is one and only one local isomorphism class of double poles with this residue realized in the class $\cal C$.
\end{lem*}
We saw in \Cref{subsub:doublepoles} that there are exactly two possible classes of double poles for each value of the residue, one centered and one non-centered.
This already shows that \emph{not all double pole local isomorphism classes can be realized in class $\cal C$}.

\medskip

We will figure out with purely analytical methods which of the two classes is realized. We will see that it is not easy, outside some special cases.

The developing map of any germ of affine chart is a map $\phi:\wt\C^*\to \C$ that satisfies $\phi''/\phi' = \zeta = \frac{1}{z^2}-\frac{\res}{z}$.
So $\log\phi' = -\frac{1}{z} - \res\times\log z$ up to an additive constant and
\[\phi' = z^{-\res}\exp(-1/z)\]
up to a multiplicative constant.

When $\res = m\in\Z$, we are in a favorable situation:
$\phi'$ is well-defined on $\C^*$, has a classical Laurent power series expansion, and $\phi$ is a logarithm plus a Laurent power series:
\[\phi'=\sum_{n\geq 0} \frac{(-1)^n}{n!} z^{-m-n}\]
up to a multiplicative constant so, up to post-composition with an affine bijection of $\C$:\\
A. If $m\leq 1$, write $m=1-n_0$ and we have
\[\phi = \frac{(-1)^{n_0}}{n_0!}\log z+ \sum_{n\geq 0, n\neq n_0} \frac{(-1)^n}{n!(-m-n+1)} z^{-m-n+1}.\]
B. If $m\geq 2$:
\[\phi = \sum_{n\geq 0} \frac{(-1)^n}{n!(-m-n+1)} z^{-m-n+1}.\]
In case~A the affine holonomy of $\phi$ is a non-zero translation and the double-pole is non-centered.
In case~B is is the identity and the double pole is centered.

When $-\res \notin\Z$, things are more complicated.
The map $\phi$ can still be computed as a generalized Laurent power series (without log term):
\[\phi = \sum_{n\in\N} \frac{(-1)^n}{n!(-\res-n+1)} z^{-\res-n+1} = z^{-\res+1} \sum_{n\in\N} \frac{(-1)^n}{n!(-\res-n+1)}z^{-n},\]
in which $z^\alpha$ for $\alpha\in\C\setminus\Z$ is to be understood as $\exp(\alpha\log z)$, where $\log z$ is well-defined if $z\in\wt\C^*$.
However, the absence of the log term does not imply anything on the isomorphism class of the double pole.
With the choice of developing map $\phi$ above, the affine holonomy of the expression is $L:z\mapsto \lambda z+0$ with $\lambda = \exp(-2\pi i \res)$, in particular it fixes the origin.
The question amounts to determining whether or not the limit of $\phi(x)$ as the real $x>0$ tends to $0$, is $0$ or a non-zero complex number.
With the expansion above, it seems hard.

Instead of using series expansions, we will try and determine the class using a path integral of $\phi'=z^{-\res}\exp(-1/z)$ as follows.
Let $r>0$ and consider the path $\gamma_r$ that is composed of three pieces: the first is the segment from $0$ to $r$, the second is the circle of center $0$ and radius $r$ followed counterclockwise, and the last one is initial segment followed backwards.
Then the double pole local isomorphism class is completely determined by whether or not the following integral vanishes:
\[I = \int_{\gamma_r} \phi'(z) dz\]
where $\phi'(z)$ is analytically continued along $\gamma_r$.
By classical complex analysis, this integral is independent of $r$.
Since as $x\to+\infty$, $\phi'(x)\sim x^{-\res}$,
we get that if $\Re(\res)>1$, then the integral along the circle tends to $0$ when its radius $r$ tends to infinity, hence $I$ is equal to the difference of the integral of two determinations of $\phi'$ along the real line:
\[\Re(\res)>1 \implies I = (1-e^{-2\pi i \res})\int_0^{+\infty} x^{-\res} e^{-1/x} dx\]
and the change of variable $y=1/x$ yields
$I = (1-e^{-2\pi i \res})\int_0^{+\infty} y^{\res-2} e^{-y}\, dy$ and we recognize \emph{Euler's Gamma function}.
So
$\Re(\res)>1 \implies I = (1-e^{-2\pi i \res})\Gamma_{\text{Euler}}(\res-1)$.
The quantity $\phi'(z)$ depends analytically on $\res$ and for a fixed $r$, the integral along $\gamma_r$ is analytic in $\res$, hence $I$ is an analytic function of $\res$. By the analytic identity theorem:
\[\forall\res\in\C \setminus \Z,\ I = (1-e^{-2\pi i \res})\,\Gamma_{\text{Euler}}(\res-1).\]
Thanks to the non trivial fact that the Gamma function never vanishes, it follows that
\[\forall\res\in\C \setminus \Z,\ I\neq 0.\]
In these cases, the double pole is non-centered.

As a conclusion:
\begin{prop*}
  For all $\res\in\C$, for the unique (up to isomorphism) affine surface in the class $\cal C$ where the double pole $p$ has residue $\res$, then $p$ is centered if and only if $\res\in\Z\cap[2,+\infty)$.
\end{prop*}

Again, \Cref{thm:CharacterizationDoublePole} contains this result and gives more information.

\smallskip

Last, we can reverse the line of arguments: the analysis of this section plus the geometric argument of \Cref{prop:DoublePoleObstruction} give \emph{a new proof that Euler's Gamma function does not vanish on $\C \setminus \Z$}.\footnote{Classical proofs include one via Euler's reflection formula and several via infinte products.
One may also use the recursion formula $\Gamma(z+1)=z\Gamma(z)$ together with an integral estimate when $\Re(z)$ is much bigger than $\Im(z)$.
The reflection formula gives the alternate expression $I = - 2i/(e^{\pi i\res}\,\Gamma_{\text{Euler}}(2-\res))$.}

\bibliographystyle{alpha}
\bibliography{bib}

\stopcontents[all]

\end{document}